\documentclass[a4paper,12pt,pagesize,twoside,BCOR1mm,numbers=noenddot,toc=bibliography]{scrreprt}
\pdfoutput=1
\usepackage[utf8]{inputenc}
\usepackage[english,ngerman]{babel}
\usepackage[T1]{fontenc}
\usepackage{lmodern}
\usepackage{dsfont}
\usepackage{microtype}
\usepackage{amsmath}
\usepackage{amssymb}
\usepackage{hyperref}
\usepackage[all]{hypcap}
\usepackage[hyperref,amsmath,thmmarks]{ntheorem}
\usepackage{tikz}
\usepackage[all]{xy}
\usetikzlibrary{matrix,arrows}

\newcommand{\R}{\mathbb R}
\newcommand{\Q}{\mathbb Q}
\newcommand{\Z}{\mathbb Z}

\newcommand{\K}{\mathbb K}

\newcommand{\C}{\mathbb C}
\newcommand{\F}{\mathbb F}

\newcommand{\bP}{\mathbb P}
\newcommand{\one}{\mathds 1}
\newcommand{\Sph}{\mathbb S}
\newcommand{\bH}{\mathbb H}
\newcommand{\cA}{\mathcal A}

\newcommand{\cC}{\mathcal C}
\newcommand{\caD}{\mathcal D}
\newcommand{\cF}{\mathcal F}
\newcommand{\cW}{\mathcal W}
\newcommand{\cP}{\mathcal P}
\newcommand{\cQ}{\mathcal Q}
\newcommand{\cX}{\mathcal X}
\newcommand{\cY}{\mathcal Y}
\newcommand{\cZ}{\mathcal Z}
\newcommand{\cU}{\mathcal U}
\newcommand{\caL}{\mathcal L}
\newcommand{\caH}{\mathcal H}
\newcommand{\cI}{\mathcal I}

\newcommand{\fs}{\mathfrak s}
\newcommand{\ft}{\mathfrak t}

\newcommand{\stern}{*}
\DeclareMathOperator{\id}{id}
\DeclareMathOperator{\im}{im}

\DeclareMathOperator{\Ad}{Ad}
\DeclareMathOperator{\Aut}{Aut}

\DeclareMathOperator{\Gr}{Gr}

\DeclareMathOperator{\Gl}{GL}

\DeclareMathOperator{\PGl}{PGL}
\DeclareMathOperator{\BGl}{BGL}
\DeclareMathOperator{\Sl}{SL}
\DeclareMathOperator{\SO}{SO}

\DeclareMathOperator{\St}{St}
\DeclareMathOperator{\Sp}{Sp}
\DeclareMathOperator{\U}{U}
\DeclareMathOperator{\BU}{BU}
\DeclareMathOperator{\E}{E}
\DeclareMathOperator{\Orth}{O}
\DeclareMathOperator{\Sym}{Sym}

\DeclareMathOperator{\Flag}{Flag}

\DeclareMathOperator{\inv}{inv}

\DeclareMathOperator{\rank}{rank}
\DeclareMathOperator{\type}{type}
\DeclareMathOperator{\supp}{supp}
\DeclareMathOperator{\ind}{ind}
\DeclareMathOperator{\lk}{lk}
\DeclareMathOperator{\Lk}{Lk}
\DeclareMathOperator{\st}{st}
\DeclareMathOperator{\Cone}{Cone}
\DeclareMathOperator{\Tot}{Tot}
\DeclareMathOperator{\Ch}{Ch}
\DeclareMathOperator{\Shapes}{Shapes}
\DeclareMathOperator{\diam}{diam}
\DeclareMathOperator{\pr}{pr}

\newcommand{\coloneq}{\mathrel{\mathop :}=}
\newcommand{\alter}[2]{\left\{\genfrac{}{}{0pt}{}{#1}{#2} \right\}}

\theoremstyle{plain}
\newtheorem{Definition}{Definition}[section]
\newtheorem{Theorem}[Definition]{Theorem}
\newtheorem{Lemma}[Definition]{Lemma}
\newtheorem{Proposition}[Definition]{Proposition}
\newtheorem{Corollary}[Definition]{Corollary}

\theorembodyfont{\rm}
\newtheorem{Examples}[Definition]{Examples}
\newtheorem{Example}[Definition]{Example}
\newtheorem{ConstructionN}[Definition]{Construction}

\theoremclass{Theorem}
\theoremstyle{break}
\newtheorem{BreakTheorem}[Definition]{Theorem}
\newtheorem{BreakProposition}[Definition]{Proposition}

\theoremstyle{nonumberbreak}
\newtheorem{SlnTheorem}{Theorem \protect\ref{th:sln_result}}
\newtheorem{UnTheorem}{Theorem \protect\ref{th:un_result}}
\newtheorem{SOnTheorem}{Theorem \protect\ref{th:son_result}}
\theoremstyle{nonumberplain}
\newtheorem{IntroTheorem}{Theorem}
\newtheorem{SpecSTheorem}{Theorem \protect\ref{th:stability_pair_spectral_sequence}}
\newtheorem{A2generalTheorem}{Theorems \protect\ref{th:a2_is_building} and \protect\ref{th:a2_classification}}
\newtheorem{A2specialTheorem}{Theorem \protect\ref{th:cyclic_lattices}}
\newtheorem{A2homologyTheorem}{Theorem \protect\ref{th:a2_homology}}
\newtheorem{C2Theorem}{Theorems \protect\ref{th:c2_building} and \protect\ref{th:c2_result2}}
\newtheorem{C2homologyTheorem}{Theorems \protect\ref{th:c2_rational_homology} and \protect\ref{th:c2_homology}}
\theorembodyfont{\rmfamily}
\newtheorem{IntroExamples}{Examples}
\newtheorem{Remark}{Remark}
\newtheorem{Remarks}{Remarks}
\newtheorem{Construction}{Construction}
\newtheorem{Terminology}{Terminology}
\theoremsymbol{\ensuremath\Box}
\newtheorem{Proof}{Proof}

\addtokomafont{paragraph}{\rmfamily}
\addtokomafont{subparagraph}{\rm\it}

\addtokomafont{subject}{\sf}

\subject{Mathematik}
\title{Buildings, Group Homology and Lattices}
\author{Jan Essert}
\date{}%\today}
\publishers{Inaugural-Dissertation\\\normalsize
zur Erlangung des akademischen Grades\\ eines Doktors der Naturwissenschaften\\ im Fachbereich Mathematik und Informatik\\ der Westfälischen Wilhelms-Universität Münster\\[2em]vorgelegt von\\Jan Essert\\aus Seeheim-Jugenheim\\2010}

\lowertitleback{\begin{tabular}{ll}
Dekan & Prof.~Dr.~Dr.~h.c.~Joachim Cuntz\\
Erster Gutachter & Prof.~Dr.~Linus Kramer\\
Zweiter Gutachter & Prof.~Dr.~Benson Farb\\
Tag der mündlichen Prüfung& 18.05.2010\\
Tag der Promotion& 18.05.2010
\end{tabular}}

\xdefinecolor{myblue}{rgb}{0.1,0.1,0.6}
\xdefinecolor{mygreen}{rgb}{0.1,0.4,0.3}

\hypersetup{
pdftitle={Buildings, Group Homology and Lattices},
pdfauthor={Jan Essert},
pdfkeywords={Buildings, Group homology, Lattices, Classical groups, Group theory},
colorlinks=true,
linkcolor=myblue,
citecolor=mygreen,
urlcolor=myblue
}

\begin{document}
\selectlanguage{english}
\maketitle

\phantom{a}\vspace{2cm}
\begin{center}
{\Large\sffamily\bfseries Abstract}
\end{center}
\vspace{0.25em}

\noindent This thesis discusses questions concerning the homology of groups related to buildings. We also give a new construction of lattices in such groups and investigate their group homology.

Specifically, we construct a simplicial complex called the \emph{Wagoner complex} associated to any group of Kac-Moody type. For a 2-spherical group $G$ of Kac-Moody type, we show that the fundamental group of the Wagoner complex is almost always isomorphic to the Schur multiplier of the little projective group of $G$.

Furthermore, we present a general method to prove \emph{homological stability for groups with weak spherical Tits systems}, that is, groups acting strongly transitively on weak spherical buildings. We use this method to prove strong homological stability results for special linear groups over infinite fields and for unitary groups over division rings, improving the previously best known results in many cases.

Finally, we give a new construction method for buildings of types $\tilde A_2$ and $\tilde C_2$ with \emph{cocompact lattices} in their full automorphism groups. Almost all of these buildings are exotic in the case $\tilde C_2$. In the $\tilde A_2$-case, it is not known whether these buildings are classical. The advantages of our construction are very explicit presentations of the lattices as well as very explicit descriptions of the buildings. Using these, we give partial results about the structure of the $\tilde A_2$-buildings and calculate group homology of all the lattices. To the author's knowledge, these are the first known presentations of lattices in buildings of type $\tilde C_2$.

A detailed overview of the results can be found in Chapter \ref{ch:results}.
\cleardoublepage

\tableofcontents
\clearpage
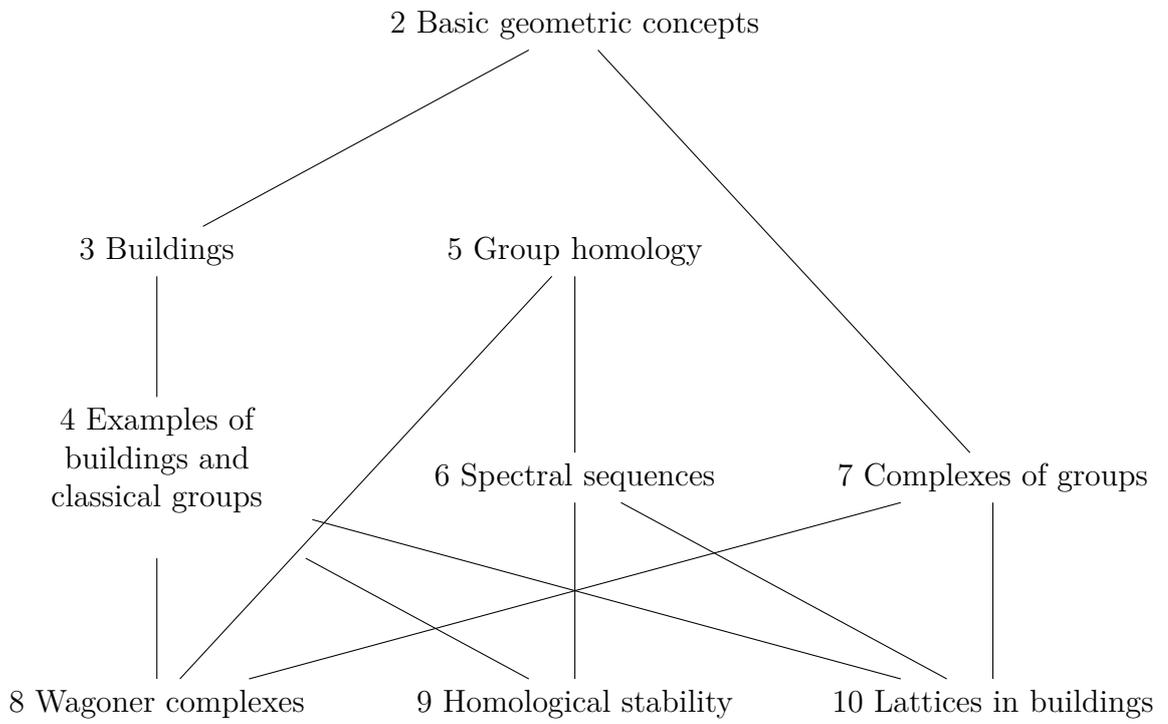
\begin{figure}[hb]
    \centering
        \phantom{a}\vspace{2cm}
        \begin{tikzpicture}
            \node (Basic) {2 Basic geometric concepts};
            \node[above of=Basic,node distance=2cm] (None) {};
            \node[below of=Basic,node distance=3cm] (Homol) {5 Group homology};
            \node[below of=Homol,node distance=3cm] (Spec) {6 Spectral sequences};

            \node[left of=Homol,node distance=5.5cm] (Build) {3 Buildings};
            \node[below of=Build,node distance=3cm] (Examp) {\parbox{3.8cm}{\begin{center}4 Examples of buildings and classical groups\end{center}}};

            \node[right of=Spec,node distance=5.5cm] (Cplx) {7 Complexes of groups};

            \node[below of=Examp,node distance=3cm] (Wagon) {8 Wagoner complexes};
            \node[right of=Wagon,node distance=5.5cm] (Stab) {9 Homological stability};
            \node[right of=Stab,node distance=5.5cm] (Latt) {10 Lattices in buildings};

            \draw (Basic) -- (Build) -- (Examp) -- (Stab);
            \draw (Examp) -- (Wagon); \draw (Examp) -- (Latt);
            \draw (Homol) -- (Spec) -- (Stab);
            \draw (Spec) -- (Latt);
            \draw (Homol) -- (Wagon);
            \draw (Cplx) -- (Wagon);
            \draw (Basic) -- (Cplx) -- (Latt);
        \end{tikzpicture}
    \caption{The architecture of this thesis.}
    \label{fig:leitfaden}
\end{figure}
\pagebreak
\phantom{a}\vfill
\section*{Acknowledgements}

I would like to thank Sarah for supporting me during these past years despite the long distance between us. Many thanks go to my family for their constant support.

It is a great pleasure to thank Linus Kramer for the amount of time he had for me, for all the good advice and for simply being a great advisor. I wish to thank Benson Farb for his hospitality during my stay at the University of Chicago and for his ongoing support afterwards.

I was supported financially by the \selectlanguage{ngerman}\emph{Graduiertenkolleg: \glqq Analytische Topologie und Me\-ta\-geo\-me\-trie\grqq}\selectlanguage{english} and the \selectlanguage{ngerman}\emph{SFB 478: \glqq Geometrische Strukturen in der Mathematik\grqq}\selectlanguage{english}.

Many people have contributed expert knowledge and valuable advice to this project. I wish to thank Uri Bader, Pierre-Emmanuel Caprace, Ruth Charney, Walter Freyn, Theo Grundhöfer, Guntram Hainke, Clara Löh, Frederick Magata, Hendrik van Maldeghem, Lars Scheele, Petra Schwer, Markus Stroppel, Koen Thas, Anne Thomas, Alain Valette, Karen Vogtmann, Richard Weiss and Stefan Witzel for their help.

In addition to mathematical discussions with and advice from Frederick, Guntram, Lars and Petra, I would like to thank my friends for the good times we had in the last three years.

\clearpage
\addtocounter{chapter}{-1}

\chapter{Introduction}

In the mathematical universe, this thesis is located in the area of geometric group theory, a topic at the intersection of algebra, geometry and topology. Geometric group theory studies groups by investigating their actions on `nice' geometric objects. In this thesis, these geometric objects will usually be buildings or geometries related to buildings. The groups will be classical groups over division rings or, later, lattices in the automorphism groups of buildings. The questions we shall investigate are related to group homology, focusing on homological stability, and to the construction of lattices.

\paragraph{} Buildings were discovered by Tits in the fifties as he studied semisimple complex Lie groups in a uniform way using geometry. Later, he generalised the concept of a building to semisimple algebraic groups over arbitrary fields. The first type of buildings to be studied were so-called spherical buildings, which can be thought of as certain simplicial complexes covered by subcomplexes which are tessellated spheres.

Spherical buildings were later classified by Tits in \cite{Tit:BsT:74}. In higher rank, there is a one-to-one correspondence between spherical buildings and a class of groups consisting (more or less) of semisimple isotropic linear algebraic groups and of isotropic classical groups. The term `classical group' is not rigorously defined, but we take it to mean groups of matrices over division rings preserving certain forms. Typical examples are linear and unitary groups.

Algebraic groups over finite fields, and hence finite spherical buildings play an important role in the classification of finite simple groups. In general, the theory of spherical buildings is very rich and a powerful tool to prove strong results for the associated groups.

If the algebraic group associated to a spherical building is defined over a field with a discrete valuation, then Bruhat and Tits constructed in \cite{BT:GRC:72} a second building associated to this group, a so-called affine building. These buildings are simplicial complexes which are covered by tessellated Euclidean spaces.

Affine buildings were classified by Bruhat and Tits in several following papers (see \cite{Wei:SAB:09} for a complete account). Higher rank affine buildings (more or less) correspond bijectively to semisimple isotropic algebraic groups over fields with discrete valuations. The classification heavily uses the fact that such buildings have a spherical building `at infinity' such that the classification of spherical buildings can be used as an ingredient.

It is common to endow spherical and affine buildings with metrics which come from the canonical metrics on the spheres and Euclidean spaces by which they are covered. With these metrics, spherical buildings are naturally CAT(1)-spaces whereas Euclidean buildings are naturally CAT(0)-spaces. It turns out that buildings can also be characterised by these metric properties. All buildings have a geometric realisation as CAT(0)-spaces, the Davis realisation.

The notion of a building has been generalised in various directions: to twin buildings, to non-discrete buildings and to $\Lambda$-buildings. On the metric side, buildings which are covered by tessellated hyperbolic spaces are also considered. In this thesis, we will mostly restrict our attention to the `classical' notions of spherical and affine buildings.

\paragraph{} Group homology is a concept arising at the intersection of group theory and topology. It is a homology theory associated to groups, which can be defined using two different approaches. On the one hand, group homology can be defined using algebraic constructions. This leads to interpretations of low-dimensional homology and cohomology groups as important algebraic invariants, for instance related to central extensions.

On the other hand, group homology can be considered as the (singular) homology of classifying spaces. This leads to topological interpretations and applications of group homology. An equivalent, completely topological point of view is the study of homology groups of aspherical spaces. Combining these two approaches, a rich interplay between topology and group theory has been achieved, leading to important results in both areas.

We will mostly be concerned with classical groups, which usually come in infinite series with inclusions such as
\[
    \Gl_1(R) \rightarrow \Gl_2(R) \rightarrow \cdots \rightarrow \Gl_k(R) \rightarrow \cdots.
\]
Homological stability is the question whether the induced maps on homology groups become isomorphisms for large $k$. This happens for many series of groups, for instance for symmetric groups, for mapping class groups and for a large class of linear groups. Homological stability results for general linear groups have been proved by various authors. The direct limit of the above sequence is called the stable general linear group $\Gl(R)$, whose homology is isomorphic to almost all of the homology groups of the unstable groups in the above sequence by the aforementioned homological stability results.

Homological stability is useful since stable homology groups are usually easier to compute. For general linear groups, stable group homology is interesting because of its connection to algebraic $K$-Theory via the Hurewicz map:
\[
    K_i(R) \rightarrow H_i(\Gl(R)).
\]
A related theory called hermitian $K$-theory is related to the homology of stable unitary groups.

Homological stability results for classical groups have been proved by different authors over the years, usually over rings with additional conditions such as finite stable rank, see \cite[Chapter 2]{Knu:HLG:01} for an overview. In this thesis, we restrict ourselves to classical groups over division rings. These groups have associated spherical buildings whose structure allows us to prove stronger stability results.

\paragraph{} A lattice in a locally compact topological group is a discrete subgroup of finite covolume. Traditionally, lattices arise as arithmetic subgroups in algebraic groups over local fields. There are many strong results in this area, for instance Margulis' super-rigidity \cite{Mar:DSS:91} stating roughly that any group homomorphism from a lattice into an algebraic group with non-compact image can be extended to the ambient algebraic group. These strong properties of lattices make their investigation a highly interesting topic.

Algebraic groups over local fields act on affine buildings and it can be shown that the lattice theory of the full automorphism group of the building does not differ much from the lattice theory of the algebraic group.

It is hence natural to extend the theory of lattices in algebraic groups to lattices in the automorphism groups of affine buildings, or even of other polyhedral complexes of non-positive curvature. A survey on results and problems in this direction by Farb, Hruska and Thomas can be found in \cite{FHT:PAG:09}.

Examples of this extension of lattice theory are lattices in trees studied by Bass and Lubotzky in \cite{BL:TL:01} and lattices in products of trees investigated by Burger and Mozes in \cite{BM:LPT:00}, but also lattices in hyperbolic or right-angled buildings studied by various authors. Any Kac-Moody group over a finite field is a lattice in the automorphism group of its twin building, see \cite{Rem:CRK:99}.

Here, we focus on lattices in two-dimensional affine buildings. More specifically, we give a construction of two-dimensional affine buildings with uniform lattices. Such a construction of lattices provides either new geometric descriptions of known arithmetic lattices or interesting examples of exotic buildings with cocompact lattices.

\minisec{The structure of this thesis} This thesis is organised as follows. We start with a detailed discussion of the main results in the \hyperref[ch:results]{first} chapter. In the \hyperref[ch:basic]{second} chapter, we introduce some basic geometric concepts which are required throughout the thesis such as simplicial complexes, metric polyhedral complexes and bounded curvature for metric spaces. Buildings and groups acting on buildings are introduced in the \hyperref[ch:buildings]{third} chapter along with all the required concepts from the theory of buildings. In Chapter \ref{ch:building_examples}, we define the classical groups and the buildings with which we are working.

The \hyperref[ch:group]{fifth} chapter introduces group homology and the concept of homological stability. It also shows the connection to algebraic $K$-theory. In the \hyperref[ch:spectral]{sixth} chapter, we define spectral sequences in general and the spectral sequences with which we will be working. Complexes of groups, which will be our basic tool for the construction of lattices, are introduced in Chapter \ref{ch:complexes_of_groups}.

The last three chapters form the part of the thesis where we discuss our results. In Chapter \ref{ch:wagoner}, we construct Wagoner complexes associated to groups of Kac-Moody type. In Chapter \ref{ch:stability}, we prove homological stability for large classes of classical groups. Finally, in Chapter \ref{ch:lattices}, we construct new lattices in two-dimensional affine buildings and calculate their homology.

Chapter \ref{ch:problems} lists open problems related to the results of this thesis. It may serve as the starting point for further research in this direction.

\pagestyle{headings}
\chapter{Main results}\label{ch:results}

We begin this thesis with a high-level overview of the main results, which are contained in Chapters \ref{ch:wagoner}, \ref{ch:stability} and \ref{ch:lattices}. Readers that are unfamiliar with the terminology should consult Figure \ref{fig:leitfaden} on page \pageref{fig:leitfaden} to find the location of the prerequisites or start directly with Chapter \ref{ch:basic} and come back for this overview at a later point.

\section{Wagoner complexes}\label{sec:wagoner_results}

Wagoner complexes are certain simplicial complexes first constructed by Wagoner in \cite{Wag:BSK:73}. In his original construction, they were associated to the group $\Gl_n(R)$ over an arbitrary ring $R$. Wagoner used the homotopy groups of these complexes to give a tentative definition of higher algebraic $K$-theory. Later Anderson, Karoubi and Wagoner proved that this definition is equivalent to Quillen's definition of $K$-theory in \cite{AKW:HAK:73,AKW:HAK:77,Wag:EAK:77}.

The complexes Wagoner defined are closely related to the building associated to the general linear group. Unlike the building, the Wagoner complex is not homotopy equivalent to a bouquet of spheres. Instead it has interesting homotopy and homology groups. Wagoner complexes have already been generalised to reductive $F$-groups by Petzold in \cite{Pet:HDW:04}. There, rational homology groups of Wagoner complexes are used to construct rational representations of reductive groups.

We define Wagoner complexes in the context of groups of Kac-Moody type (or groups with root data), which we believe to be the appropriate level of generality. All definitions are very natural and the connection to the theory of buildings is clearly visible. Important examples of such groups are semi-simple Lie groups and algebraic groups as well as Kac-Moody groups.

We exhibit interesting substructures in Wagoner complexes analogous to apartments in buildings and present some general properties of Wagoner complexes. The main result is the following characterisation of low-dimensional homotopy groups of Wagoner complexes in the case where the associated root datum is of $2$-spherical type.

\begin{IntroTheorem} Let $G$ be a group of Kac-Moody type and assume that its root datum is of $2$-spherical type. Let $G^\dagger$ be the little projective group and let $\hat G$ be the associated Steinberg group in the sense of Tits. Let $\cW(G)$ be the associated Wagoner complex.
    
    Then, in almost all cases
        \begin{align*}
                \pi_0(|\cW(G)|) &\cong G/G^\dagger. \\
                \intertext{If $G$ is of non-spherical type or $\rank(G)\geq 3$, we have}
                \pi_1(|\cW(G)|) &\cong \ker(\hat G\twoheadrightarrow G^\dagger).
        \end{align*}
        By a recent result by Caprace in \cite{Cap:O2S:07}, this means $\pi_1(|\cW(G)|)\cong H_2(G^\dagger;\Z)$ except over very small fields.
\end{IntroTheorem}
For a precise statement of the above theorem see Theorem \ref{th:wag_main_result} and Corollary \ref{cor:wag_main_result}.

In addition to this result, we also make the connection to algebraic $K$-theory explicit in Section \ref{sec:ktheory}. Finally, we give a definition of affine Wagoner complexes associated to groups admitting root data with valuations in Section \ref{sec:affine}.

\section{Homological stability}\label{sec:stability_results}

There are many homological stability results for various series of classical groups over rings with finite stable rank. Usually, these proofs are tailored to the specific series of groups. In this thesis, we restrict ourselves to division rings in order to have the theory of buildings as a tool at hand. This enables us to construct a spectral sequence which can be used for homological stability results for various series of groups.

More concretely, using buildings, we develop a method to prove homological stability for a large class of groups: those which act strongly transitively on a weak spherical building. Among these groups are general and special linear groups as well as unitary, symplectic and orthogonal groups. We apply our method to concrete series of groups to obtain homological stability results which improve the previously known stability ranges in many cases.

Homological stability proofs are usually based on the following idea: Consider the action of a group on a highly connected simplicial complex, find smaller groups of the same series as stabiliser subgroups and inspect the right spectral sequence.

The starting point for our investigation are the homological stability proofs for linear, symplectic and orthogonal groups by Charney in \cite{Cha:HSD:80} and \cite{Cha:gtV:87} and by Vogtmann in \cite{Vog:HSO:79} and \cite{Vog:SPH:81}. Our method is based on the observation that the simplicial complexes used there are closely related to the theory of buildings --- they are the \emph{opposition complexes} studied by von Heydebreck in \cite{vH:HPC:03}.

The opposition complexes admit Levi subgroups as vertex stabilisers. Using these, we construct a spectral sequence involving relative group homology of Levi subgroups. These Levi subgroups often split as direct of semidirect products of smaller groups which can be used for stability proofs.

We will use this spectral sequence for homological stability proofs for special linear and for unitary groups. For these groups, the Levi subgroups split as direct or semidirect products, both of special linear or unitary groups and, interestingly, of general linear groups. Using strong stability results for general linear groups, we can prove stability results for special linear or unitary groups.

Our method can probably also be used to show low-dimensional homological stability for groups of types $E_6$, $E_7$ and $E_8$. Additionally, one could try to compare group homology of groups of different types using this method. Finally, homological stability results for all reductive algebraic groups should be possible, albeit with a rather weak stability range.

\subsection{Homological stability results}

The method outlined above has originally been used by Charney in \cite{Cha:HSD:80} to prove homological stability of general and special linear groups, but yielding a comparatively weak stability range. For special linear groups, however, it is an interesting observation that terms involving general linear groups appear in the spectral sequence. This allows us to apply a strong theorem by Sah in \cite{Sah:HcL:86} on homological stability for general linear groups to prove homological stability for special linear groups.

\begin{SlnTheorem}[Homological stability of special linear groups]
  If $D$ is an infinite field, then $n\geq 2k-1$ implies
  \[
	  H_k(\Sl_{n+1}(D),\Sl_{n}(D);\Z)=0.
  \]
\end{SlnTheorem}

\noindent For fields of characteristic zero, there is a far better result by Hutchinson and Tao in \cite{HT:HSS:08} with stability range $n\geq k$. Up to now, the best result known to the author applicable to other infinite fields has been a result by van der Kallen in \cite{vdK:HSL:80} for rings with stable rank one. It guarantees a stability range of $n\geq 2k$.

Vogtmann originally used a version of this construction to prove homological stability for orthogonal and symplectic groups in \cite{Vog:HSO:79} and \cite{Vog:SPH:81}. In this thesis, we investigate the general situation of unitary groups associated to a hermitian form of Witt index $n+1$ on a vector space $V$. This vector space then splits non-canonically as an orthogonal sum of a hyperbolic module $\caH_{n+1}$ and an anisotropic complement $E$. We consider the unitary group induced on the subspace $\caH_n\perp E$, where $\caH_n$ is a hyperbolic module of rank $n$ and ask for homological stability. Again, the spectral sequence we consider has terms involving the relative homology of general linear groups. We can hence apply Sah's theorem again to obtain

\begin{UnTheorem}[Homological stability of unitary groups]
  For a division ring $D$ with infinite centre, the relative homology modules
  \[
    H_k\bigl(\U(\caH_{n+1}\perp E), \U(\caH_n\perp E);\Z\bigr)
  \]
  vanish for $n\geq 2$ if $k=1$ and for $n\geq k\geq 2$. If the centre of $D$ is finite, relative homology vanishes for $n\geq 2k$.
\end{UnTheorem}

\noindent This is an improvement over the results by Mirzaii and van der Kallen in \cite{MaB:HSU:02} and \cite{Mir:HSU:05}, where homological stability for unitary groups with stability range $n\geq k+1$ has already been proved. Their result is valid for a much larger class of rings, namely local rings with infinite residue fields, but only for the case of maximal Witt index, that is for $E=\{0\}$.

Rather surprisingly, the following strong result can also be proved using this method.

\begin{SOnTheorem}[Homological stability of special orthogonal groups]
  For an infinite field $D$, we have
  \[
    H_k\bigl(\SO_{n+1,n+1}(D),\SO_{n,n}(D);\Z\bigr) =0
  \]
  for $n\geq 2$ if $k=1$ and for $n\geq k\geq 2$. If $D$ is a finite field, then the relative homology groups vanish for $n\geq 2k$.

\end{SOnTheorem}

\subsection{The construction of a relative spectral sequence}

In the following, we give an outline of the method used to prove these results and we state the main theorem which is the basis for the homological stability proofs. Consider a group $G$ with a weak spherical Tits system of rank \mbox{$n+1$}. We enumerate the type set $I=\{i_1,\ldots,i_{n+1}\}$ arbitrarily. For $1\leq p\leq n+1$, denote by $L_p$ certain Levi subgroups of $G$ of type $I\backslash\{i_p\}$.

For the applications discussed in the previous section, the group $G$ is of type $A_{n+1}$ or $C_{n+1}$ with a linear ordering of the type set. The resulting Coxeter diagrams of $G$ and $L_p$ are illustrated in the following picture.
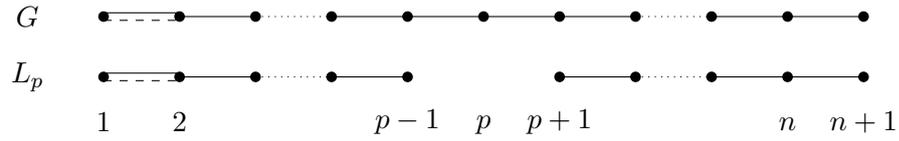
\begin{center}
\begin{tikzpicture}[font=\small]
  \node (G) at (-1,.8) {$G$};
  \node (L) at (-1,0) {$L_p$};
  \foreach \y in {0,.8} {
  \foreach \x in {0,1,2,3,4,6,7,8,9,10} { \fill (\x,\y) circle (.7mm);}
  \draw (1,\y) -- (2,\y);
  \draw[dotted] (2,\y) -- (3,\y);
  \draw (3,\y) -- (4,\y);
  \draw (6,\y) -- (7,\y);
  \draw[dotted] (7,\y) -- (8,\y);
  \draw (8,\y) -- (10,\y);
  \draw (0,\y + .05) -- (1,\y + .05);
  \draw[dashed]  (0,\y - .05) -- (1,\y -.05);
  }
  \draw (4,.8) -- (6,.8);
  \fill (5,.8) circle (.7mm);
  \node (0) at (0,-.6) {$1$}; \node (1) at (1,-.6) {$2$}; \node (3) at (4,-.6) {$p-1$}; \node (4) at (5,-.63) {$p$}; \node (5) at (6,-.6) {$p+1$}; \node (7) at (9,-.62) {$n$}; \node (8) at (10,-.6) {$n+1$};
\end{tikzpicture}
\end{center}
We choose a subgroup $G'\leq L_{n+1}$ of type $I\backslash \{i_{n+1}\}$ and write $L'_p=L_p\cap G'$. Again, in the concrete applications, we have the following situation.
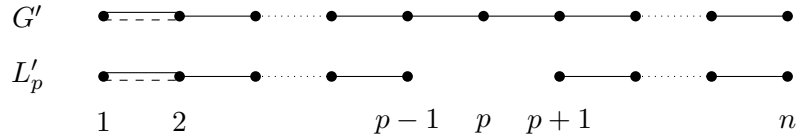
\begin{center}
\begin{tikzpicture}[font=\small]
  \node (G) at (-1,.8) {$G'$};
  \node (L) at (-1,0) {$L'_p$};
  \foreach \y in {0,.8} {
  \foreach \x in {0,1,2,3,4,6,7,8,9} { \fill (\x,\y) circle (.7mm);}
  \draw (1,\y) -- (2,\y);
  \draw[dotted] (2,\y) -- (3,\y);
  \draw (3,\y) -- (4,\y);
  \draw (6,\y) -- (7,\y);
  \draw[dotted] (7,\y) -- (8,\y);
  \draw (8,\y) -- (9,\y);
  \draw (0,\y + .05) -- (1,\y + .05);
  \draw[dashed]  (0,\y - .05) -- (1,\y -.05);
  }
  \draw (4,.8) -- (6,.8);
  \fill (5,.8) circle (.7mm);
  \node (0) at (0,-.6) {$1$}; \node (1) at (1,-.6) {$2$}; \node (3) at (4,-.6) {$p-1$}; \node (4) at (5,-.63) {$p$}; \node (5) at (6,-.6) {$p+1$}; \node (7) at (9,-.6) {$n$};
  \node (empty) at (10,0) {\phantom{$n+1$}};
\end{tikzpicture}
\end{center}
Using a filtration by types of vertices of the opposition complex, we construct two exact chain complexes of $G$- and $G'$-modules. From these chain complexes, we obtain a spectral sequence involving relative homology of Levi subgroups with coefficient modules $M_p$ which are top-dimensional homology modules of opposition complexes of type $\{i_1,\ldots,i_{p-1}\}$, except for $M_1=\Z$.

\begin{SpecSTheorem}[Relative spectral sequence]
  There is a spectral sequence with first page
  \[
  E^1_{p,q}=\begin{cases}
    H_q(G,G';\Z) & p=0 \\
    H_q(L_p,L'_p;M_p) & 1\leq p \leq n\\
    H_q(L_{n+1},G';M_{n+1}) & p=n+1
  \end{cases}
  \]
  which converges to zero.
\end{SpecSTheorem}

This can be used to prove homological stability for groups of type $A_{n+1}$ and $C_{n+1}$ in the following way: We want to prove that $H_q(G,G';\Z)$ vanishes for all $q\leq k$ for a given $k$. Hence we must show that $H_q(L_p,L'_p;M_p)$ vanishes for $p+q\leq k+1$. For $2\leq p\leq n-1$ the Levi subgroups, having disconnected diagrams, usually split as direct or semidirect products of two groups whose types belong to the connected components of the Coxeter diagrams.
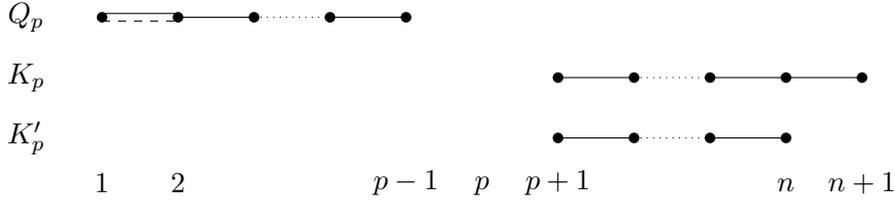
\begin{center}
\begin{tikzpicture}[font=\small]
  \node (Qp) at (-1,1.6) {$Q_p$};
  \node (K) at (-1,0.8) {$K_p$};
  \node (Kp) at (-1,0) {$K'_p$};
  \foreach \x in {0,1,2,3,4} { \fill (\x,1.6) circle (.7mm);}
  \draw (1,1.6) -- (2,1.6);
  \draw[dotted] (2,1.6) -- (3,1.6);
  \draw (3,1.6) -- (4,1.6);
  \draw (0,1.6 + .05) -- (1,1.6 + .05);
  \draw[dashed]  (0,1.6 - .05) -- (1,1.6 -.05);
  \foreach \y in {0,.8} {
  \foreach \x in {6,7,8,9} { \fill (\x,\y) circle (.7mm);}
  \draw (6,\y) -- (7,\y);
  \draw[dotted] (7,\y) -- (8,\y);
  \draw (8,\y) -- (9,\y);
  }
  \draw (9,.8) -- (10,.8);
  \fill (10,.8) circle (.7mm);
  \node (0) at (0,-.6) {$1$}; \node (1) at (1,-.6) {$2$}; \node (3) at (4,-.6) {$p-1$}; \node (4) at (5,-.63) {$p$}; \node (5) at (6,-.6) {$p+1$}; \node (7) at (9,-.62) {$n$}; \node (8) at (10,-.6) {$n+1$};
\end{tikzpicture}
\end{center}
This means that there are groups $Q_p$, $K_p$ and $K'_p$ of types $\{i_1,\ldots,i_{p-1}\}$, $\{i_{p+1},\ldots,i_{n+1}\}$ and $\{i_{p+1},\ldots,i_n\}$, respectively, such that $L_p=Q_p\ltimes K_p$ and $L'_p=Q_p\ltimes K'_p$. The modules $M_p$ are constructed in such a fashion that the groups $K_p$ and $K'_p$ act trivially on $M_p$. If we know that relative \emph{integral} homology of the subgroups $K_p$ and $K'_p$  vanishes, we can produce zeroes in the first page $E^1_{p,q}$ by using a relative Lyndon\slash Hochschild-Serre spectral sequence.

The structure of the Levi subgroups and the corresponding semidirect product decompositions depend on the specific series of groups. But note that we always require relative integral homology of groups of type $A_*$. For special linear groups over fields and for unitary groups over division rings, these subgroups of type $A_*$ are general linear groups. Hence, as mentioned above, we can use strong results on the homological stability of general linear groups to obtain homological stability of these series of groups.

\paragraph{}Variations of this method using an appropriate type filtration can probably be used to prove homological stability results for different series of reductive groups. As mentioned above, this could be used to study groups of type $E_6$, $E_7$ or $E_8$, or to compare group homology of groups of different types. In particular, relations between algebraic K-theory and hermitian K-theory could also be studied by choosing a different type enumeration of a group of type $C_{n+1}$, forcing $G'$ to be of type $A_n$ instead of type $C_n$.

\section{Lattices in two-dimensional affine buildings}\label{sec:lattices_results}

A \emph{uniform lattice} is a discrete and cocompact subgroup in a locally compact topological group. We describe a new geometric construction method producing buildings of type $\tilde A_2$ and $\tilde C_2$ with a uniform lattice in their full automorphism groups. The lattices act regularly, that is transitively and freely, on panels of the building having the same type. The advantages of this construction are very simple presentations for the lattices, very explicit descriptions of the buildings and a very good understanding how these lattices act on the building. Using the latter fact, we can calculate the group homology of the lattices. To the author's knowledge, these are the first known presentations of lattices in buildings of type $\tilde C_2$.

It is known that there are countably many buildings of types $\tilde A_2$ and $\tilde C_2$ coming from algebraic groups over local fields, as well as uncountably many other so-called \emph{exotic} buildings of these types.  For the buildings of type $\tilde C_2$, all but possibly one of the buildings we construct are exotic. In the case of $\tilde A_2$, this is not known yet, except for a single building whose associated lattice we can realise explicitly in $\Sl_3(\F_2(\!(t)\!))$.

If the building $X$ we construct is classical, then there is an associated algebraic group $G\leq \Aut(X)$. By a result of Tits in \cite{Tit:BsT:74}, the algebraic group $G$ is always cocompact in $\Aut(X)$ and it is hence conceivable that our lattice (or a finite-index subgroup) is already contained in $G$. Then Margulis' Arithmeticity \cite{Mar:DSS:91} asserts that the lattice is \emph{arithmetic}, that is, comes from an algebraic construction. Our construction would then provide simple presentations for arithmetic lattices, along with a very explicit description of the action on the building.

\minisec{Details} For the construction, we use \emph{Singer polygons}, generalised polygons with a point-regular automorphism group, to construct small complexes of groups. The local developments of these complexes of groups are cones over the generalised polygons with which we have started. These cones are automatically non-positively curved in a natural metric, so the complexes of groups are developable by a theorem of Bridson and Haefliger in \cite{BH:NPC:99}. Then, by a recognition theorem by Charney and Lytchak in \cite{CL:MC:01}, we know that their universal covers are buildings. The fundamental groups of these complexes of groups are then uniform lattices on two-dimensional locally finite affine buildings.

In principle, our construction might even extend to lattices in buildings of type $\tilde G_2$, but unfortunately there are no known Singer hexagons.

\minisec{Previous constructions} There are three geometric constructions of chamber-transitive lattices in buildings of type $\tilde A_2$ by Köhler, Meixner and Wester in \cite{KMW:A2l:84, KMW:2ab:85} and by Kantor in \cite{Ka:SLF:87}, all of these in characteristic two. All chamber-transitive lattices in classical buildings have been classified by Kantor, Liebler and Tits in \cite{KLT:CTL:87}. In particular, the examples of Köhler, Meixner, Wester and Kantor are the only $\tilde A_2$-buildings with chamber-transitive lattices.

Cartwright, Mantero, Steger and Zappa construct vertex-transitive lattices on buildings of type $\tilde A_2$ in \cite{CMSZ1, CMSZ2}. They prove that some of the lattices they construct are contained in $\Sl_3(\F_q(\!(t)\!))$, while others correspond to exotic buildings.

In \cite{Ron:TG:84}, Ronan constructs (possibly exotic) buildings of type $\tilde A_2$ admitting a lattice acting regularly on vertices of the same type.

Except for the single lattice which we can realise in $\Sl_3(\F_2(\!(t)\!))$, it is not clear whether any of our lattices are commensurable to the lattices constructed by Cartwright-Mantero-Steger-Zappa or by Ronan.

\paragraph{} For buildings of type $\tilde C_2$, there is a free construction by Ronan in \cite{Ron:CBR:86} which produces very unstructured examples of exotic buildings and does not give any control over the automorphism groups. In \cite{Kan:GPS:86}, Kantor gives a construction of exotic buildings of type $\tilde C_2$ admitting uniform lattices acting freely on vertices of the building. However, no presentations of these lattices are given.

Again, we do not know whether the lattices we construct here are commensurable to the examples given by Kantor.

\subsection{Buildings of type \texorpdfstring{$\tilde A_2$}{\textasciitilde A2}}

For buildings of type $\tilde A_2$, we obtain the following result. Fix three Singer projective planes (Singer generalised triangles) of order $q$ along with three \emph{Singer groups} $S_i$, $i\in\{1,2,3\}$, that is, groups acting regularly on points (and hence on lines) of these planes. For each of these planes, fix a point $p_i$ and a line $l_i$ and write
\[
    D_i = \{ d\in S_i: p_i\text{ is incident to } d(l_i) \}.
\]
These sets are called \emph{difference sets}. Finally, write $J=\{0,1,\ldots,q\}$ and fix three bijections $d_i:J\rightarrow D_i$.

\begin{A2generalTheorem}
    There is an $\tilde A_2$-building $X$ such that the group
    \[
    \Gamma_1 = \Bigl\langle S_1,S_2,S_3 \,\Big| \begin{array}{c}\text{ all relations in the groups }S_1,S_2,S_3,\\ d_1(j)d_2(j)d_3(j) = d_1(j')d_2(j')d_3(j') \quad\forall j,j'\in J\end{array}\Bigr\rangle
    \]
is a uniform lattice in the full automorphism group of the building. The union of all chambers containing a fixed panel is a fundamental domain for the action.

Conversely, every panel-regular lattice on a building of type $\tilde A_2$ arises in this way and admits a presentation of the above form with $d_1(0)d_2(0)d_3(0)=1$.
\end{A2generalTheorem}

\noindent A special case arises as follows: Classical projective planes admit cyclic Singer groups. Denote the generators of the three Singer groups $S_i$ by $\sigma_i$, and write
\[
\Delta_i=\{ \delta \in \Z / |S_i| : \sigma_i^\delta \in D_i\}.
\]
These are difference sets in the classical sense.

\begin{A2specialTheorem}
    For any prime power $q$, for any three classical difference sets $\Delta_1$, $\Delta_2$, $\Delta_3$ containing 0, and for any bijections $\delta_i:J\rightarrow \Delta_i$ satisfying $\delta_i(0)=0$, the group $\Gamma_2$ with presentation
	\[
		\Gamma_2 = \langle \sigma_1, \sigma_2, \sigma_3 \,|\, \sigma_1^{q^2+q+1}=\sigma_2^{q^2+q+1}=\sigma_3^{q^2+q+1}=1,\, \sigma_1^{\delta_1(j)}\sigma_2^{\delta_2(j)}\sigma_3^{\delta_3(j)} = 1 \quad\forall j \in J\rangle
\]
	is a uniform lattice in a building of type $\tilde A_2$.
\end{A2specialTheorem}

\begin{IntroExamples}
 Two examples are:
\begin{align*}
    \Lambda &= \langle \sigma_1,\sigma_2,\sigma_3 \,|\, \sigma_1^7=\sigma_2^7=\sigma_3^7 = \sigma_1\sigma_2\sigma_3 = \sigma_1^3\sigma_2^3\sigma_3^3 = 1\rangle,\\
    \Lambda' &= \langle \sigma_1,\sigma_2,\sigma_3 \,|\, \sigma_1^{13}=\sigma_2^{13}=\sigma_3^{13} = \sigma_1\sigma_2^3\sigma_3^9 = \sigma_1^3\sigma_2^9\sigma_3 =\sigma_1^9\sigma_2\sigma_3^3 = 1\rangle.
\end{align*}
Further examples can easily be constructed using the list of difference sets found in the La Jolla Difference Set Repository \cite{LaJolla}.
\end{IntroExamples}

If we start from cyclic Singer groups, we obtain very explicit descriptions of the associated buildings. Using these, we can give explicit descriptions of the spheres of radius two in these buildings, which lead to an interesting incidence structure called a \emph{Hjelmslev plane of level two}. Using a result of Cartwright-Mantero-Steger-Zappa in \cite{CMSZ2}, we show that the lattices from the latter theorem cannot belong to the building associated to $\Sl_3(\Q_p)$ if $\Delta_1=\Delta_2=\Delta_3$.

For the homology groups, we obtain the following result.

\begin{A2homologyTheorem} If $\Gamma_2$ is a lattice as in the latter theorem, then
    \[
        H_j(\Gamma_2;\Z) \cong\begin{cases}
            \Z & j=0\\
            \ker(\caD) & j=1 \\
            \Z^q & j=2 \\
            (\Z/(q^2+q+1))^3 & j\geq 3 \text{ odd} \\
            0 & \text{ else,}
        \end{cases}
    \]
    where $\caD:(\Z/(q^2+q+1))^3\rightarrow (\Z/(q^2+q+1))^q$ is given by the matrix $(\delta_j(i))_{i,j}$.

    In addition
    \[
        H_2(\Gamma_i;\Q)=\Q^q,\qquad H_j(\Gamma_i;\Q)=0 \text{ for $j\not \in \{0,2\}$}
    \]
    for any lattice $\Gamma_i$ as in one of the theorems above.
\end{A2homologyTheorem}

\subsection{Buildings of type \texorpdfstring{$\tilde C_2$}{\textasciitilde C2}}

For buildings of type $\tilde C_2$, we give two different constructions of panel-regular lattices, depending on the types of panels the lattice acts regularly on. We use \emph{slanted symplectic quadrangles} as defined in \cite{GJS:SSQ:94} for the vertex links. Since almost all of these quadrangles are exotic, we necessarily construct exotic buildings admitting panel-regular lattices.

Starting from a slanted symplectic quadrangle of order $(q-1,q+1)$, where $q>2$ is a prime power, we write $J=\{0,1,\ldots,q+1\}$ and consider a Singer group $S$ acting regularly on points of the quadrangle. Fix the set of lines $L$ through a point and a bijective enumeration function $\lambda:J\rightarrow L$. All line stabilisers $S_l$ are isomorphic to $\Z/q$, and we fix isomorphisms $\psi_j:\Z/q\rightarrow S_{\lambda(j)}$. We repeat this construction for a second quadrangle of order $(q-1,q+1)$ and obtain a Singer group $S'$ and functions $\lambda'$ and $\{\psi'_j\}_{j\in J}$.

\begin{C2Theorem}
    The finitely presented groups
    \begin{align*}
    \Gamma_1 &= (S * S')/ \langle [S_{\lambda(j)}, S'_{\lambda'(j)}] : j\in J\rangle \\
    \intertext{as well as}
    \Gamma_2 &= (S * S' * \langle c\rangle) / \langle c^{q+2},\, c^j \psi_j(x) c^{-j} \psi_j'(x)^{-1} : j\in J, x\in \Z/q\rangle
    \end{align*}
    are uniform panel-regular lattices in buildings of type $\tilde C_2$.

    \begin{itemize}
        \item If $q$ is an odd prime power, then $S$ and $S'$ are three-dimensional Heisenberg groups over $\F_q$.

            If $q$ is an odd prime, the above presentations can be made more explicit by writing out presentations for $S$ and $S'$.
        \item If $q$ is even, then $S$ and $S'$ are isomorphic to the additive groups of $\F_q^3$. The set of lines can be identified with $\bP\F_q^2\sqcup \{0\}$ and the stabilisers have the following structure:
            \begin{itemize}
                \item For a point $[a:b]$ in the projective plane $\bP\F_q^2$, the stabiliser $S_{[a:b]}$ is the $\F_q$-subspace of $S=\F_q^3$ spanned by $(a,b,0)^T$.
                \item The stabiliser $S_0$ is the $\F_q$-subspace spanned by $(0,0,1)^T$.
            \end{itemize}
    \end{itemize}
\end{C2Theorem}

\noindent If $q>3$, then the associated buildings are necessarily exotic. In any case, the above presentations imply very simple descriptions of these buildings. Finally, we calculate group homology.

\begin{C2homologyTheorem}
    If $\Gamma_i$ for $i\in\{1,2\}$ is any of the two lattices constructed in the above theorem, we have $H_j(\Gamma_i;\Q)=0$ for $j\neq 0$. In addition, for the first type of lattices we have
    \[
        H_1(\Gamma_1;\Z)\cong (\Z/q)^6,\qquad H_2(\Gamma_1;\Z) \cong H_2(S) \oplus H_2(S').
    \]
\end{C2homologyTheorem}

\chapter{Basic geometric concepts}\label{ch:basic}

\section{Simplicial complexes}

The basic geometric object with which we will be dealing throughout this thesis is a simplicial complex.

\begin{Definition}
	An \emph{abstract simplicial complex} consists of a set of \emph{vertices} $V$ together with a collection $\Sigma$ of finite subsets of $V$ such that
    \begin{itemize}
        \item $\{v\} \in \Sigma$ for every $v\in V$.
        \item If $\sigma\in\Sigma$ and $\tau\subseteq \sigma$ then $\tau\in\Sigma$.
    \end{itemize}
    The elements of $\Sigma$ are called \emph{simplices}. The cardinality of a simplex $\sigma$ is called its \emph{rank}, denoted by $\rank(\sigma)$. The \emph{dimension of $\sigma$} is defined by $\dim(\sigma) = \rank(\sigma) - 1$. The \emph{dimension} and \emph{rank} of $\Sigma$ is defined to be the supremum of the dimensions and ranks of simplices, respectively. Usually, we write $\Sigma^0$ for the vertex set $V$ and we will identify vertices $v$ with the corresponding singletons $\{v\}$ in $\Sigma$.
\end{Definition}

\noindent Note that a simplicial complex always contains the empty set, which is called the \emph{empty simplex}.

\begin{Example}
    The power set of the finite set $\{0,1,\ldots,n\}$ is a simplicial complex of dimension $n$, called the \emph{abstract standard $n$-simplex}, sometimes denoted by $\Delta^n$.

    This example should also explain the terminology: Compare this abstract complex to the \emph{$n$-dimensional standard simplex}, the convex hull of the $n+1$ standard basis vectors in $\R^{n+1}$. Observe that each face of the standard simplex is uniquely determined by its set of vertices.
\end{Example}

\noindent Since we are interested in group theory, we require the concept of automorphisms of simplicial complexes.

\begin{Definition}
    A \emph{simplicial map} or morphism of simplicial complexes is a map between the vertex sets which maps simplices to simplices. An action of a group by simplicial automorphisms will be called a \emph{simplicial action}. The group of all simplicial automorphisms of $\Sigma$ is denoted by $\Aut(\Sigma)$.
\end{Definition}

\noindent We say that two simplices $\sigma,\tau\in\Sigma$ are \emph{joinable} if the set $\sigma\cup \tau$ is a simplex, that is $\sigma\cup \tau\in\Sigma$.

\begin{Definition}
    Let $\Sigma$ be a simplicial complex and let $\sigma\in\Sigma$ be a simplex. We define the following subsets of $\Sigma$:
    \begin{itemize}
        \item The \emph{link} of $\sigma$ is given by $\lk_\Sigma(\sigma) = \{\tau\in\Sigma: \tau,\sigma \text{ are joinable}, \sigma\cap \tau=\emptyset \}$.
        \item The \emph{residue} (or \emph{open star}) is given by $R_\Sigma(\sigma) = \{ \tau\in\Sigma : \sigma\subseteq \tau\}$.
        \item The \emph{(closed) star} is given by $\st_\Sigma(\sigma) = \{\tau\in\Sigma: \tau,\sigma \text{ are joinable}\}$.
    \end{itemize}
    For a vertex $v$, we have $\st_\Sigma(v) = \lk_\Sigma(v) \sqcup R_\Sigma(v)$.
\end{Definition}

\noindent Links and closed stars are always subcomplexes of $\Sigma$. There is a simple way to construct a simplicial complex out of every partially ordered set.

\begin{Definition}
    Let $(P,\leq)$ be a partially ordered set. The \emph{flag complex $\Flag(P)$} has $P$ as its vertex set. A simplex of $\Flag(P)$ is a finite, totally ordered subset of $P$, that is, a set whose elements can be arranged in an ascending chain. The flag complex is always a simplicial complex.
\end{Definition}

\noindent Since simplicial complexes always include the empty simplex $\emptyset$, we make the following convention.

\begin{Definition}
    If $\Sigma$ is a simplicial complex partially ordered by inclusion, then $\Flag(\Sigma\setminus\{\emptyset\})$ is called the \emph{barycentric subdivision of $\Sigma$}. It is denoted by $\Sigma'$.
\end{Definition}

\noindent We are interested in simplicial complexes where the maximal dimension is attained everywhere.

\begin{Definition}
    A simplicial complex is \emph{pure} if all maximal simplices have the same dimension. In this case, maximal simplices are called \emph{chambers} and simplices of codimension one are called \emph{panels}.

    Associated to every pure simplicial complex $\Sigma$, there is the \emph{chamber graph} whose vertices are the chambers of $\Sigma$ and where two vertices are joined by an edge if the corresponding chambers share a common panel. A path in the chamber graph is called a \emph{gallery}.

    A pure simplicial complex is called a \emph{chamber complex} if its chamber graph is connected.
\end{Definition}

\noindent All the chamber complexes we will encounter are colourable in the following sense.

\begin{Definition}\label{def:type_function}
    Let $\Sigma$ be a chamber complex of rank $n$ with vertex set $\Sigma^0$. Let $I$ be a set with $n$ elements. A \emph{type function} on $\Sigma$ with values in $I$ is a function
    \[
        \type: \Sigma^0 \rightarrow I
    \]
    such that the vertices of each chamber are mapped bijectively onto $I$. A chamber complex that admits a type function is called \emph{colourable}.
\end{Definition}

\begin{Remark}
    Obviously, a type function is unique up to a bijection of the type set $I$. Moreover, given a type function, we define the type of a simplex to be the set of types of its vertices. In this way, we obtain an extension of the type function to $\Sigma$ taking its values in the power set of $I$.
\end{Remark}

\begin{Definition}
    A simplicial automorphism $\varphi\in\Aut(\Sigma)$ of a colourable chamber complex $\Sigma$ is called \emph{type-preserving} if $\type= \type\circ\varphi$. The set of type-preserving automorphisms forms a subgroup of $\Aut(\Sigma)$ which is usually denoted by $\Aut^0(\Sigma)$.
\end{Definition}

\noindent There is a way to construct new simplicial complexes out of old ones.

\begin{Definition}
    Let $\Sigma_1$ and $\Sigma_2$ be two abstract simplicial complexes. The \emph{simplicial join} $\Sigma_1 * \Sigma_2$ is the abstract simplicial complex with vertex set $(\Sigma_1*\Sigma_2)^0 = \Sigma_1^0 \sqcup \Sigma_2^0$ and with simplices
    \[
        \Sigma_1 * \Sigma_2 = \{ \sigma \subseteq \Sigma_1^0 \sqcup \Sigma_2^0 \,:\, \sigma \cap \Sigma_1^0 \in \Sigma_1,\, \sigma \cap \Sigma_2^0 \in \Sigma_2\}.
    \]
\end{Definition}

\begin{Examples}
    The join of a simplicial complex with a single vertex is the \emph{cone} over the simplicial complex. The join of the abstract standard simplices $\Delta^n$ and $\Delta^m$ is the abstract standard simplex $\Delta^{n+m}$.
\end{Examples}

\noindent Associated to every simplicial complex, there is a topological space, the geometric realisation.

\begin{Definition}
    Let $\Sigma$ be a simplicial complex with vertex set $\Sigma^0$. Consider an $\R$-vector space $E$ with basis $\Sigma^0$. For every simplex $\sigma\in\Sigma$, the \emph{geometric realisation} $|\sigma|$ is given by the convex hull of its vertices
    \[
    |\sigma| = \bigl\{ \sum_{v\in \sigma} \lambda_v v \,:\, \lambda_v\geq 0,\, \sum_{v\in \sigma} \lambda_v = 1 \bigr\} \subset E.
    \]
    Note that $|\sigma|$ inherits a topology from the finite-dimensional subspace of $E$ spanned by the vertices of $\sigma$. The \emph{geometric realisation} of the full complex $|\Sigma|$ is given as the union
    \[
    |\Sigma| = \bigcup_{\sigma\in\Sigma} |\sigma|.
    \]
    We endow $|\Sigma|$ with the weak topology, that is, a subset of $|\Sigma|$ is closed if and only if all its intersections with geometric realisations of simplices are closed.
\end{Definition}

\begin{Remark}
	We could equally well consider the standard scalar product on $E$ and the induced metric on $|\Sigma|$. If $\Sigma^0$ is infinite, then the corresponding topology called the \emph{metric topology} is not the same as the weak topology. However, the identity map $\id: |\Sigma|_{\text{weak}}\rightarrow |\Sigma|_{\text{metric}}$ induces a homotopy equivalence by \cite[Proposition IV.4.6]{LW:CW:69}. Since we are usually only interested in the homotopy type of the simplicial complex, the choice of the metric topology would not change the results.

    In this thesis, the prominent examples of simplicial complexes are Coxeter complexes and buildings, where we will not only consider the weak topology. On spherical and affine buildings, there is a natural metric, which is discussed in Section \ref{subsec:metrics_on_buildings}.
\end{Remark}

\begin{Lemma}
    The geometric realisation of a simplicial complex is homeomorphic to the geometric realisation of its barycentric subdivision.
\end{Lemma}

\begin{Proof}
    Let $\Sigma$ be a simplicial complex. The vertices in the barycentric subdivision $\Sigma'$ correspond to non-empty simplices of $\Sigma$. Let $E$ be the $\R$-vector space with basis $\Sigma^0$ and let $E'$ be the $\R$-vector space with basis $(\Sigma')^0=\Sigma\setminus\{\emptyset\}$. The linear map $\varphi:E'\rightarrow E$ given on the basis vectors by
    \[
        \varphi(\sigma) = \sum_{v\in \sigma} \tfrac{1}{\rank \sigma} v,
    \]
    where $v$ runs through all vertices of $\sigma$, obviously restricts to a homeomorphism on the realisation $|S|$ of each simplex of $S\in\Sigma'$. Since both $\Sigma$ and $\Sigma'$ are endowed with the weak topology, the map $\varphi$ restricts to a homeomorphism $|\Sigma'|\rightarrow |\Sigma|$.
\end{Proof}

\section{Non-positive curvature}

In addition to simplicial complexes as combinatorial geometric objects and their realisations as topological spaces, we will also consider metrics on simplicial complexes. We will require a few concepts from metric geometry which are discussed in this section. The following material can be found in much greater detail in \cite{BH:NPC:99}.

\begin{Definition}
    Let $(X,d)$ be a metric space. A \emph{geodesic} joining $x\in X$ to $y \in X$ is a map $c:[0,l]\rightarrow X$ from a closed interval in $\R$ to $X$ such that
    \[
        c(0)=x,\quad c(l)=y,\quad d(c(t),c(t'))=|t-t'|\quad\text{ for all } t,t'\in [0,l].
    \]
    The image of a geodesic is called a \emph{geodesic segment}.

    A subset $C\subseteq X$ is said to be \emph{convex} if every pair of points of $C$ can be joined by a geodesic such that the corresponding geodesic segment is contained in $C$.

    The space $X$ is said to be an \emph{$r$-geodesic space} if every two points in $X$ of distance less than $r$ can be joined by a geodesic. The space $X$ is \emph{geodesic} if it is $r$-geodesic for all $r>0$.

    The \emph{diameter $\diam(X)$} is defined as
    \[
    \diam(X) \coloneq \sup_{x,y\in X}d(x,y).
    \]
\end{Definition}

\noindent We will introduce two basic spaces which will serve as comparison spaces throughout this thesis. The first one is the standard Euclidean space.

\begin{Definition}
    Consider the standard $n$-dimensional Euclidean space $\R^n$ endowed with the usual scalar product $\langle -,-\rangle$. We will write $M^n_0=(\R^n,d)$ for the associated metric space. \emph{Hyperplanes} in $M^n_0$ are defined to be affine subspaces of codimension one. The geodesic segment between two points $x,y\in M^n_0$ is given by
    \[
        [x,y] \coloneq \{ ty + (1-t)x : 0\leq t\leq 1\}.
    \]
    Given a point $p\in M^n_0$ and two geodesic segments $c=[p,x]$ and $c'=[p,y]$ the \emph{angle at $p$ between $c$ and $c'$} is given by the usual Euclidean angle: the unique number $\alpha$ in $[0,2\pi)$ such that
    \[
    \cos \alpha = \frac{\langle x-p, y-p \rangle}{\sqrt{\langle x-p,x-p\rangle}\sqrt{\langle y-p,y-p\rangle}}.
    \]
    We write $\angle_p^{(0)}(x,y)=\alpha$.
\end{Definition}

\noindent The second basic spaces are spheres with the angular metric.

\begin{Definition}
    We write $M^n_1=\Sph^n$ for the $n$-sphere in $\R^{n+1}$, endowed with the \emph{angular metric}, which assigns to each pair $(x,y)\in\Sph^n\times\Sph^n$ the unique number $d(x,y)\in [0,\pi]$ such that
    \[
        \cos d(x,y) = \langle x,y \rangle,
    \]
    which is the standard scalar product in $\R^{n+1}$.
\end{Definition}

\noindent A proof that this is a metric can be found in \cite[Proposition I.2.1]{BH:NPC:99}.

\begin{Definition}
    A \emph{hyperplane} in $M^n_1$ is the intersection of an $n$-dimensional linear subspace of $\R^{n+1}$ with $\Sph^n$.
\end{Definition}

\noindent Any geodesic segment is given by a starting point $p$ and an initial unit vector $u\in \R^{n+1}$ such that $\langle p,u\rangle=0$, see \cite[I.2]{BH:NPC:99}. The geodesic segment then has the form
    \[
        \{ \cos(t)p + \sin(t)u : t\in [0,l] \}
    \]
    for $l\in[0,\pi]$. For any two points $v,w\in \Sph^n$ of distance strictly less than $\pi$, there is a unique geodesic segment joining them denoted by $[v,w]$.

\begin{Definition}
The \emph{spherical angle} between two geodesic segments starting at $p$ with initial vectors $u$ and $v$ joining $p$ to $x$ and $y$ is given by the unique number $\alpha$ in $[0,\pi]$ such that
    \[
        \cos \alpha = \langle u, v\rangle.
    \]
    We write $\angle_p^{(1)}(x,y)=\alpha$.
\end{Definition}

\begin{Proposition}[Proposition I.2.11 in \cite{BH:NPC:99}]
    The spaces $M^n_{\kappa}$ are geodesic for $\kappa\in\{0,1\}$.
\end{Proposition}

\noindent We will also require the appropriate notion of subspaces.

\begin{Definition}
    Let $0<m\leq n$ be integers. A subset of $M^n_\kappa$ is called an \emph{$m$-plane} if it is isometric to $M^m_\kappa$.
\end{Definition}

\begin{Lemma}[Corollary I.2.22 in \cite{BH:NPC:99}]
    If $m<n$, then every $m$-plane in $M^n_\kappa$ is the intersection of $(n-m)$ hyperplanes. Every subset of $M^n_\kappa$ is contained in a unique $m$-plane of minimal dimension.
\end{Lemma}

\noindent In classical Riemannian geometry, the notion of curvature is a central aspect. Spaces with curvature bounded from above or below are studied extensively in this area. For metric spaces, there is a notion of non-positive curvature due to Alexandrov, which is very important in metric geometry.

We will first introduce the concept of a comparison triangle.

\begin{Lemma}[Lemma I.2.14 in \cite{BH:NPC:99}]
    Choose $\kappa\in\{0,1\}$. Let $x$, $y$ and $z$ be three points in a metric space $(X,d)$, where we assume the perimeter $d(x,y)+d(y,z)+d(z,x)$ of the triangle $(x,y,z)$ to be less than $2\pi$ if $\kappa=1$. Then there is a triple of points $(\bar x, \bar y, \bar z)$ in $M^2_\kappa$ such that
    \[
        d(x,y)= d(\bar x, \bar y),\quad d(y,z) = d(\bar y, \bar z),\quad d(z,x) = d(\bar z, \bar x).
    \]
    Such a triple is called an \emph{$M^2_\kappa$-comparison triangle}.
\end{Lemma}

\noindent In this situation, for a point $p$ on a geodesic segment from $x$ to $y$, the \emph{comparison point} is the point $\bar p \in [\bar x, \bar y]$ such that
\[
    d(p,x) = d(\bar p, \bar x) \quad\text{and}\quad d(p,y) = d(\bar p, \bar y).
\]
Comparison points on $[y,z]$ and on $[z,x]$ are defined analogously.

\begin{Definition}
    Let $(X,d)$ be a metric space and let $\kappa\in\{0,1\}$. Let $\Delta$ be a geodesic triangle in $X$, of perimeter less than $2\pi$ if $\kappa=1$, and let $\bar \Delta$ be an $M^2_\kappa$-comparison triangle. Then $\Delta$ is said to \emph{satisfy the CAT($\kappa$)-inequality} if for all $p,q\in\Delta$ and for all comparison points $\bar p,\bar q\in \bar \Delta$, we have
    \[
        d(p,q) \leq d (\bar p, \bar q).
    \]
    The space $X$ is called a \emph{CAT(0)-space} if it is geodesic and all geodesic triangles satisfy the CAT(0)-inequality. It is called a \emph{CAT(1)-space} if it is $\pi$-geodesic and all geodesic triangles of perimeter less than $2\pi$ satisfy the CAT(1)-inequality.
\end{Definition}

\noindent We will also require a local version of these conditions.

\begin{Definition}
    A metric space $X$ is said to be \emph{locally CAT($\kappa$)} or of \emph{curvature at most $\kappa$}, if for every $x\in X$ there exists $r_x>0$ such that the ball $B(x,r_x)$ with the induced metric is CAT($\kappa$).

    A metric space which is locally CAT(0) is said to be \emph{non-positively curved}.
\end{Definition}

\noindent It turns out that in simply connected spaces the local condition implies the global condition.

\begin{BreakTheorem}[Cartan-Hadamard-Theorem, Theorem II.4.1 in \cite{BH:NPC:99}]\label{th:cartan_hadamard}
    A complete, connected and simply connected metric space which is non-positively curved is already CAT(0).
\end{BreakTheorem}

\noindent CAT($\kappa$)-spaces have many interesting and strong properties. In this thesis, we will only use the following result about the homotopy type.

\begin{Theorem}[Corollary II.1.5 in \cite{BH:NPC:99}]
    A CAT(0)-space is contractible.
\end{Theorem}

\noindent In the following, we will sometimes consider quotients of metric spaces. The following construction provides a pseudometric on quotients.

\begin{Definition}
    Let $(X,d)$ be a metric space, let $\sim$ be an equivalence relation on $X$, set $\bar X = X/\sim$ and denote the projection by $p:X\rightarrow \bar X$. For $\bar x$, $\bar y\in \bar X$, the \emph{quotient pseudometric $\bar d$} is given by
    \[
    \bar d(\bar x, \bar y) \coloneq \inf \sum_{i=1}^n d(x_i,y_i),
    \]
    where the infimum is taken over all sequences $(x_1,y_1,x_2,y_2,\ldots,x_n,y_n)$ of points of $X$ such that $x_1\in \bar x$, $y_n\in \bar y$ and $y_i\sim x_{i+1}$ for $i\in\{1,\ldots,n-1\}$.
\end{Definition}

\section{Metric polyhedral complexes}

The quotient of the geometric realisation of a simplicial complex modulo a group action is, in general, no longer simplicial. Even though the cells are still simplices, it may happen that two simplices intersect in a union of faces, which is forbidden for a simplicial complex. On the other hand, it is interesting to consider complexes whose cells are not simplices but other polytopes. Cube complexes are a particularly interesting class. So, even though we will only use complexes whose cells are simplices, we introduce the larger class of \emph{polyhedral complexes}, which is closed under quotients.

Roughly speaking, a polyhedral complex is a cell complex endowed with a metric such that the cells are isometric to convex polyhedra in a sphere, a Euclidean space or a hyperbolic space. Since we will only require piecewise Euclidean and piecewise spherical complexes, we will restrict the definitions to these cases.

We shall first define convex $M_\kappa$-polyhedral cells, which are the models for cells in polyhedral complexes.

\begin{Definition}
    Fix $\kappa\in\{0,1\}$. A \emph{convex $M_\kappa$-polyhedral cell $C$} is the convex hull of a finite set of points in $M^n_\kappa$. If $\kappa=1$, these points are required to lie in an open ball of radius $\pi/2$. The \emph{dimension} of $C$ is the dimension of the smallest $m$-plane containing it. The \emph{interior of $C$} is the interior of $C$ as a subset of this $m$-plane.

    Let $H$ be a hyperplane in $M^n_\kappa$. If $C$ lies in one of the closed half-spaces bounded by $H$ and if $H\cap C \neq\emptyset$, then $F=H\cap C$ is called a \emph{face} of $C$. Again, the dimension of a face is the dimension of the smallest $m$-plane containing it and its interior is the interior of $F$ in this plane. By convention, the whole cell $C$ is also a face of itself.

    The zero-dimensional faces of $C$ are called its \emph{vertices}. The \emph{support} of $x\in C$, denoted by $\supp(x)$ is the unique face containing $x$ in its interior.
\end{Definition}

\noindent We obtain a polyhedral complex by glueing polyhedral cells along their faces.

\begin{Definition}
    Let $(C_\lambda : \lambda \in \Lambda)$ be a family of convex $M_\kappa$-polyhedral cells and let $X=\coprod_{\lambda\in\Lambda}(C_\lambda \times \{\lambda\})$ denote their disjoint union. Let $\sim$ be an equivalence relation on $X$ and let $K=X/\sim$. Let $p:X\rightarrow K$ be the natural projection and define $p_\lambda: C_\lambda \rightarrow K$ by $p_\lambda(x)\coloneq p(x,\lambda)$.

    \noindent Then $K$ is called an \emph{$M_\kappa$-polyhedral complex} if
    \begin{itemize}
        \item For every $\lambda\in\Lambda$, the restriction of $p_\lambda$ to the interior of each face of $C_\lambda$ is injective.
        \item For all $\lambda_1,\lambda_2$ and all $x_1\in C_{\lambda_1}$ and $x_2\in C_{\lambda_2}$ satisfying $p_{\lambda_1}(x_1) = p_{\lambda_2}(x_2)$, there is an isometry $h:\supp(x_1)\rightarrow \supp(x_2)$ such that $p_{\lambda_1}(y)=p_{\lambda_2}(h(y))$ for all $y\in\supp(x_1)$.
    \end{itemize}
    The set of isometry classes of the faces of the cells of $C_\lambda$ is denoted by $\Shapes(K)$.
\end{Definition}

\noindent We endow $K$ with the \emph{intrinsic pseudometric}, which is the quotient pseudometric induced by the projection $p:\coprod_{\lambda\in\Lambda}C_\lambda\rightarrow K$. By the following theorem, this is really a metric if $\Shapes(K)$ is finite.

\begin{Theorem}[Theorem I.7.50 in \cite{BH:NPC:99}]\label{th:finite_shapes_poly_complexes_are_geodesic}
    If $\Shapes(K)$ is finite and $K$ is connected, then $K$ is a complete geodesic space with respect to the intrinsic (pseudo-)metric.
\end{Theorem}

\noindent We will also require the notion of geometric links in $M_\kappa$-polyhedral complexes, which is the analogue of combinatorial links in abstract simplicial complexes in this metric context. Links are defined for all points of a cell, but not for faces which are not vertices. In the following definition, note that geodesic segments are uniquely defined by the properties of the ambient spaces $M^n_\kappa$.

\begin{Definition}\label{def:geom_links_in_simplices}
    Let $C$ be a convex polyhedral cell of dimension $m$ and fix $x\in C$. The \emph{geometric link of $x$} is the set of $\sim$-equivalence classes of geodesic segments joining $x$ to points of $C$, that means
    \[
        \Lk(x,C) \coloneq \bigl\{ [x,y] : y\in C\setminus\{x\} \bigr\} / \sim,
    \]
    where $[x,y_1]\sim [x,y_2]$ if $[x,y_1]\subseteq [x,y_2]$ or vice versa.

    The geometric link is metrised in the following way: The distance between the equivalence classes $[x,y_1]$ and $[x,y_2]$ is the angle $\angle^{(\kappa)}_x(y_1,y_2)$ in the model space $M^n_\kappa$. If $x$ is a vertex of $C$, then with this metric $\Lk(x,C)$ is a convex polyhedral cell in $\Sph^{m-1}$.
\end{Definition}

\noindent Links in polyhedral complexes are defined by glueing links in simplices.

\begin{Definition}\label{def:geometric_link}
    For any point $x$ in an $M_\kappa$-polyhedral complex $K$, the \emph{open star} $\st(x)$ is defined to be the union of the interiors of cells containing $x$.

    Analogously to Definition \ref{def:geom_links_in_simplices}, given two geodesic segments $[x,y_1]$ and $[x,y_2]$ where $y_1,y_2\in\st(x)\setminus\{x\}$, we define an equivalence relation $[x,y_1]\sim[x,y_2]$ if $[x,y_1]\subseteq [x,y_2]$ or vice versa. We define the \emph{geometric link} to be
    \[
    \Lk(x,K) \coloneq \bigl\{ [x,y] : y\in\st(x)\setminus\{x\}\bigr\} / \sim.
    \]
    We use the projection $\coprod_{\lambda\in\Lambda} \Lk(x_\lambda, C_\lambda) \rightarrow \Lk(x,K)$ induced by the projection
    \[
        p:\coprod_{\lambda\in\Lambda} C_\lambda \rightarrow K
    \]
    to endow $\Lk(x,K)$ with the quotient pseudometric. Here, we set $\Lk(x_\lambda,C_\lambda)=\emptyset$ if $x\not\in p(C_\lambda)$ and otherwise $x_\lambda$ is chosen to satisfy $x=p_\lambda(x_\lambda)$.
\end{Definition}

\noindent Usually, we will only consider simplices instead of more general polyhedra.

\begin{Definition}
    A set of $(m+1)$ points in $M^n_\kappa$ is \emph{in general position} if it is not contained in any $(m-1)$-plane.

    A \emph{geodesic $m$-simplex $S\subset M^n_\kappa$} is the convex hull of $m+1$ points in general position. If $\kappa=1$, then these points are required to lie in an open ball of radius $\pi/2$.
\end{Definition}

\noindent An $M_\kappa$-simplicial complex is not only constructed solely of simplices, but it also has the usual property of an abstract simplicial complex: that the intersection of any two simplices is empty or a simplex.

\begin{Definition}
    The polyhedral complex $K$ is an \emph{$M_\kappa$-simplicial complex} if each of the cells $C_\lambda$ is a geodesic simplex, each of the maps $p_\lambda$ is injective and the intersection of any two cells in $K$ is empty or a single face.

    An $M_0$-simplicial complex is often called a \emph{piecewise Euclidean complex}, an $M_1$-simplicial complex is called a \emph{piecewise spherical complex}.
\end{Definition}

\noindent Note that geometric links of vertices in piecewise Euclidean complexes are piecewise spherical complexes.

\begin{Remark}
    If $\Shapes(K)$ is finite, then $\Shapes(\Lk(x,K))$ is automatically finite for any vertex $x\in K$. By Theorem \ref{th:finite_shapes_poly_complexes_are_geodesic}, geometric links in piecewise Euclidean complexes with a finite set of shapes are also complete geodesic spaces.
\end{Remark}

\noindent For $M_\kappa$-polyhedral complexes, there are many strong results implying CAT($\kappa$). We will only state the following special result.

\begin{Theorem}[Theorem II.5.4 in \cite{BH:NPC:99}]
    A piecewise Euclidean complex with a finite set of shapes is CAT(0) if and only if it is simply connected and the geometric link of every vertex is CAT(1).
\end{Theorem}

\chapter{Buildings}\label{ch:buildings}

The goal of geometric group theory is to study groups by investigating their actions on `nice' geometric objects. Buildings are such `nice' geometric objects with a rich structure. They are associated to many groups appearing naturally in group theory, in particular to classical groups and algebraic groups.

Not every building corresponds to a `nice' group, but those which do are usually related to a group of matrices, an algebraic group or a classical group over a division ring. For the buildings in which we are interested most, so-called \emph{spherical} and \emph{affine} buildings, it can be shown that all higher-dimensional buildings are \emph{classical}, that is, associated to such a matrix group.

On the other hand, it is known that there are uncountably many so-called \emph{exotic} buildings of small dimensions which do not fit into this classical picture. It will be a key problem in Chapter \ref{ch:lattices} whether we can find out if a given building is classical or exotic.

\begin{figure}[hbt]
    \centering
        \begin{tikzpicture}
            \node (1) at (0,0) {$\langle e_1\rangle$};
            \node (2) at (4,0) {$\langle e_2\rangle$};
            \node (3) at (2,3.464) {$\langle e_3\rangle$};
            \node (12) at (2,0) {$\langle e_1,e_2\rangle$};
            \node (23) at (3,1.732) {$\langle e_2,e_3\rangle$};
            \node (13) at (1,1.732) {$\langle e_1,e_3\rangle$};
            \draw (1) -- (12) -- (2) -- (23) -- (3) -- (13) -- (1);
        \end{tikzpicture}\caption{An apartment in 2-dimensional projective space. This is the barycentric subdivision of the boundary of a 2-simplex, hence a triangulated 1-sphere.}
    \label{fig:a2_apartment}
\end{figure}
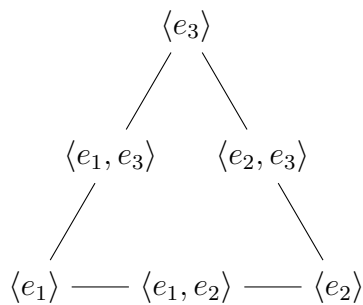

\paragraph{Example} For a spherical building, the basic example of which to think is the $n$-dimensional projective space over a field $k$. This is the flag complex over the set of non-trivial proper subspaces of $k^{n+1}$ which is partially ordered by inclusion. It is a simplicial complex which has very simple subcomplexes, called \emph{apartments}. For every basis of $k^{n+1}$, the set of all proper subspaces which can be expressed as spans of subsets of this basis form an apartment. It is not hard to see that this subcomplex is in fact the barycentric subdivision of the boundary of an $n$-simplex, in particular a triangulated $(n-1)$-sphere. See Figure \ref{fig:a2_apartment} for an illustration for $n=2$.

Buildings contain many apartments. In this example, for every two subspaces, there is a basis such that these subspaces can be expressed as spans of subsets of this basis. This implies that every two vertices are contained in a common apartment.

Note also that the group $\Gl_{n+1}(k)$ acts naturally on this projective space. Since it acts transitively on bases, it acts transitively on apartments (even on pairs consisting of an apartment and a maximal simplex contained in the apartment).

Now, in general, a spherical building can be thought of as a highly symmetric simplicial complex covered by apartments which are not necessarily boundaries of simplices, but some triangulated spheres. These special complexes are called \emph{Coxeter complexes}, which we will have to define first in order to introduce buildings.

\paragraph{} Buildings can be introduced in many different ways, among these as simplicial complexes, as chamber systems and as Coxeter-group-valued metric spaces. In this thesis, we will almost exclusively focus on the point of view of simplicial complexes, since this is the most useful one for our applications.

We will first construct Coxeter complexes, which are the basic `building blocks' for buildings. For most of the cases we are interested in, these can be thought of as triangulated spheres or triangulated Euclidean spaces. Using these complexes, we will construct buildings and investigate natural metrics. Afterwards, we will investigate the `right' kind of group actions on buildings and discuss implications on the group level.

\paragraph{} For the reader not familiar with buildings, the book by Brown \cite{Bro:Bdg:89} provides a very good introduction. This book has later been extended by Abramenko and Brown to \cite{AB:B:08} which is very comprehensive and which we will use for almost all standard references. The books by Ronan \cite{Ron:LoB:89}, by Weiss \cite{Wei:SSB:03, Wei:SAB:09} and by Garrett \cite{Gar:BCG:97} provide additional information. The book by Tits \cite{Tit:BsT:74} is of course the classical reference.

\section{Coxeter complexes}

In order to define Coxeter complexes, we first introduce abstract Coxeter groups. It will take a bit of work to see that these abstract Coxeter groups have faithful linear representations. We will mostly be interested in those Coxeter groups that can be realised as reflection groups on spheres or as affine reflection groups on Euclidean spaces.

\subsection{Coxeter groups}\label{subsec:coxeter_groups}

Coxeter groups in this sense were first introduced by Tits to generalise the notion of finite groups with similar presentations studied earlier by Coxeter.

\begin{Definition}
    Let $I$ be a finite set. A \emph{Coxeter matrix} $M=(m_{i,j})_{i,j\in I}$ over $I$ is a symmetric matrix with entries 1 on the diagonal and with entries in $\{2,3,\ldots,\infty\}$ elsewhere.

    We associate to $M$ the \emph{Coxeter diagram}, an edge-labelled graph with vertex set $I$ and with edges between $i$ and $j$ if $m_{i,j}>2$. We label all edges by the corresponding matrix elements $m_{i,j}$. It is customary to omit the label if $m_{i,j}\in\{3,4\}$ and to draw a double edge for $m_{i,j}=4$.
\end{Definition}

\noindent A list of diagrams and their names can be found in \ref{fig:spherical_coxeter_diagrams} and \ref{fig:affine_coxeter_diagrams} in the appendix. Associated to every Coxeter matrix there is a finitely presented group, the Coxeter group.

\begin{Definition}
    Let $I$ be a finite set and let $M$ be a Coxeter matrix over $I$. The associated \emph{Coxeter group $W$} is a finitely presented group with generating set
    \[
        S=\{s_i:i\in I\}
    \]
    and with presentation
    \[
    W = \langle S \;|\; (s_is_j)^{m_{i,j}}=1\quad\forall i,j\in I \rangle,
    \]
    where $m_{i,j}=\infty$ will mean that we omit the corresponding relation.
    The pair $(W,S)$ is called a \emph{Coxeter system}, the cardinality $|I|$ is called its \emph{rank}.
\end{Definition}

\begin{Remark}
    We will try to be precise by always specifying the Coxeter system $(W,S)$ instead of just the Coxeter group $W$. In fact, there are different Coxeter diagrams leading to isomorphic Coxeter groups. Determining the right condition on the diagram to ensure that it is the only diagram associated to a Coxeter group is a difficult problem, the so-called \emph{isomorphism problem} for Coxeter groups.
\end{Remark}

\begin{Example}\label{ex:dihedral_group}
    For the matrix $M=\begin{pmatrix}
        1 & m \\ m & 1
    \end{pmatrix}$ with corresponding Coxeter diagram
	
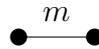
\begin{center}
		\begin{tikzpicture}
			\node (1) at (.5,.3) {$m$};
			\tikzstyle{every node}=[draw,circle,fill,inner sep=0pt,minimum width=6pt];
			\draw (0,0) node {} -- (1,0) node {};
		\end{tikzpicture}
	\end{center}
the associated Coxeter group $W$ admits the presentation
    \[
        W = \langle s_1,s_2 \;|\; s_1^2 = s_2^2 = (s_1s_2)^m = 1 \rangle.
    \]
    If $m$ is finite, it is easy to see that $W$ is isomorphic to the dihedral group of order $2m$ and $W$ is said to be of type $I_2(m)$. The types $I_2(3)$, $I_2(4)$ and $I_2(6)$ are also usually called $A_2$, $C_2$ and $G_2$, respectively. If $m=\infty$, then $W$ is isomorphic to the infinite dihedral group and $W$ is said to be of type $\tilde A_1$.
\end{Example}

\begin{Example}\label{ex:an_coxeter_group}
    Let $I=\{1,2,\ldots,n\}$. For the Coxeter diagram $A_n$ on $n$ vertices as follows
	\begin{center}
		\begin{tikzpicture}
			\tikzstyle{every node}=[draw,circle,fill,inner sep=0pt,minimum width=6pt];
			\draw (0,0) node {} -- (1,0) node {};
			\draw[dotted] (1,0) -- (2,0);
			\draw (2,0) node {} -- (3,0) node {};
		\end{tikzpicture}
    \end{center}
    the associated Coxeter group has the presentation
    \[
    W = \langle s_1,s_2,\ldots, s_n \;|\; s_i^2=1, (s_is_{i+1})^3=1\quad\forall 1\leq i<n, s_is_j = s_js_i \quad\forall |i-j|>2 \rangle
    \]
    It is not difficult to verify that this is a presentation for the symmetric group on $n+1$ letters $\Sym(n+1)$, where the elements $s_i$ correspond to the transpositions $(i, i+1)$.
\end{Example}

\begin{Example}\label{ex:cn_coxeter_group}
    Let $I=\{1,2,\ldots,n\}$ and consider the Coxeter diagram $C_n$ on $n$ vertices as shown below
    \begin{center}
		\begin{tikzpicture}
			\tikzstyle{every node}=[draw,circle,fill,inner sep=0pt,minimum width=6pt];
            \node (0) at (0,0) {};
            \node (1) at (1,0) {};
			\draw (0,0.07) -- ++(1,0);
			\draw (0,-0.07) -- ++(1,0);
            \draw (1,0) -- ++(1,0) node {};
			\draw[dotted] (2,0) -- ++(1,0);
            \draw (3,0) node {} -- ++(1,0) node {};
		\end{tikzpicture}
    \end{center}
    The associated Coxeter group has the presentation
    \[
    W = \Biggl\langle s_1,s_2,\ldots,s_n \,\bigg|\, \begin{array}{cccc}s_i^2=1&\forall i\in I,& (s_is_{i+1})^3 =1&\forall 1 < i < n,\\  s_is_j=s_js_i&\forall |i-j|>2,& (s_1s_2)^4=1\end{array}\Biggr\rangle.
    \]
    This group can be interpreted as the group of \emph{signed permutations} of the set
    \[
        [\pm n]\coloneq\{-n,\ldots,-1,1,\ldots,n\}.
    \]
    These are those permutations $w$ satisfying $-w(k)=w(-k)$ for all $k\in[\pm n]$. The elements $(s_i)_{2\leq i\leq n}$ interchange $\pm (i-1)$ with $\pm i$ and leave all other numbers fixed, while $s_1$ flips $1$ and $-1$ and leaves all other elements fixed.
\end{Example}

\noindent We also consider subgroups generated by subsets of the generating set $S$.

\begin{Definition}
    Let $(W,S)$ be a Coxeter system with associated type set $I$ and let $J\subseteq I$ be a subset. We consider the corresponding generator set $S_J = \{ s_i : i\in J\}$ and the subgroup $W_J=\langle S_J \rangle$. Subgroups of this form are called \emph{standard subgroups} or \emph{standard parabolic subgroups}. Any coset $wW_J$ will be called a \emph{standard coset}.
\end{Definition}

\noindent It is easy to see that subgroups of this form are also Coxeter groups. For any $J\subseteq I$, the subgroup $W_J$ is a Coxeter group with generator set $S_J$ with respect to the restricted Coxeter matrix $M_J \coloneq (m_{i,j})_{i,j\in J}$.

\begin{Definition}
    Let $(W,S)$ be a Coxeter system with type set $I$. Then $(W,S)$ is said to be \emph{irreducible} if there is no decomposition $I=J_1\sqcup J_2$ into two non-empty sets such that $W= W_{J_1}\times W_{J_2}$.
\end{Definition}

\noindent Obviously, a Coxeter system is irreducible if and only if the associated Coxeter diagram is connected.

\subsection{Coxeter complexes}

Coxeter complexes are the basic geometric objects associated to Coxeter groups. For the Coxeter groups in which we are interested most, we will later see that the associated complexes are naturally triangulated spheres or triangulated Euclidean spaces.

\begin{Definition}
    We denote the set of all standard cosets by
	\[
		\Sigma(W,S)=\{wW_J: w\in W, J\subseteq I\},
	\]
	and order it partially by reverse inclusion.
\end{Definition}

\begin{Theorem}[Theorem 3.5 in \cite{AB:B:08}]
     For any Coxeter system $(W,S)$, the partially ordered set $\Sigma(W,S)$ is a colourable chamber complex, the \emph{Coxeter complex} associated to $(W,S)$. The Coxeter group $W$ acts by simplicial automorphisms on $\Sigma(W,S)$ via left multiplication. The type function is given by
     \[
        \type( wW_J ) = I\setminus J.
     \]
\end{Theorem}

\begin{Terminology}
    As usual for chamber complexes, top-dimensional simplices of $\Sigma(W,S)$ are called \emph{chambers}, codimension-1-simplices are called \emph{panels}. The set of chambers of $\Sigma(W,S)$ will be denoted by $\Ch(\Sigma(W,S))$. Two chambers $c$ and $c'$ are called \emph{adjacent}, $c\sim c'$, if they share a common panel.
\end{Terminology}

\begin{Remark}
    Note that the chambers of $\Sigma(W,S)$ are singleton subsets of $W$. We will hence usually identify $W$ with the set of chambers of $\Sigma(W,S)$.

    Note also that we gave the type function directly for simplices, which is equivalent to the procedure in Definition \ref{def:type_function}, where we first construct a type function for vertices and then extend it to all simplices.
\end{Remark}

\begin{Example}\label{ex:i2_complex}
    For the Coxeter group of type $I_2(m)$, given by
    \[
    W=\langle s_1,s_2 \;|\; s_1^2=s_2^2=(s_1s_2)^m=1\rangle
    \]
    as in Example \ref{ex:dihedral_group}, the only non-trivial proper standard subgroups are $\langle s_1\rangle=\{1,s_1\}$ and $\langle s_2\rangle=\{1,s_2\}$. In this case, the Coxeter complex is obviously one-dimensional. As remarked above, chambers, that is 1-simplices, correspond to elements of $W$, whereas vertices correspond to cosets of the two standard subgroups.

    Figure \ref{fig:coxeter_complex_i2_4} shows the Coxeter complex of type $I_2(4)=C_2$. Here, vertices which are cosets of $\langle s_1 \rangle$ are drawn white, vertices which are cosets of $\langle s_2\rangle$ are drawn black. Black vertices are of type $\{1\}$ and white vertices are of type $\{2\}$. The elements of each coset correspond directly to the chambers which contain the associated vertex.
    \begin{figure}[hbt]
        \centering\begin{tikzpicture}
            \tikzstyle{every node}=[circle, draw, inner sep=0pt, minimum width=6pt]
                \draw \foreach \x in {22.5,112.5,202.5,292.5}
                {
                (\x-45:2) node[circle,draw,fill=white] {} -- (\x:2)
                (\x:2) node[circle,draw,fill=black] {} -- (\x+45:2)
                };
                \tikzstyle{every node}=[]
                \node (s1) at (135:2.5) {$s_1s_2s_1$};
                \node (s1) at (180:2.5) {$s_1s_2$};
                \node (s1) at (225:2.5) {$s_1$};
                \node (1) at (270:2.5) {$1$};
                \node (s2) at (315:2.5) {$s_2$};
                \node (s2) at (0:2.5) {$s_2s_1$};
                \node (s2) at (45:2.5) {$s_2s_1s_2$};
                \node (s2) at (90:2.5) {$w_0$};
                \node (w0) at (7, 9 pt) {where};
                \node (w0) at (7, -9 pt) {$w_0=s_1s_2s_1s_2=s_2s_1s_2s_1$};
            \end{tikzpicture}\caption{The Coxeter complex of type $I_2(4)$.}\label{fig:coxeter_complex_i2_4}
    \end{figure}
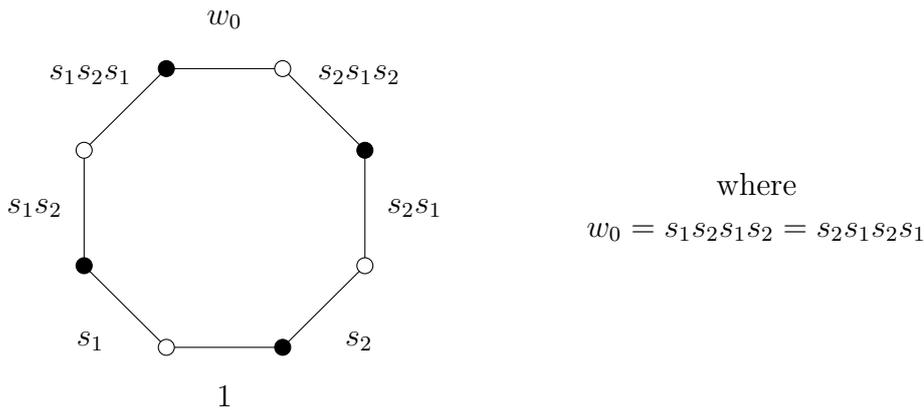
\end{Example}

\begin{Example}\label{ex:an_coxeter_complex}
    The Coxeter complex of type $A_n$, which is associated to the symmetric group on $n+1$ letters as in Example \ref{ex:an_coxeter_group}, can be thought of as the barycentric subdivision of the boundary of an $n$-simplex in the following way:

    Note first of all that the simplicial complex associated to a single $n$-simplex is in bijection to the power set of the set $\{1,\ldots,n+1\}$, since every simplex is uniquely determined by its set of vertices.

    Let $\Sigma$ be the simplicial complex of proper subsets of $\{1,\ldots,n+1\}$, which corresponds to the boundary of the $n$-simplex. This is obviously a chamber complex of dimension $(n-1)$. The symmetric group $W=\Sym(n+1)$, which is a Coxeter group of type $A_n$ with the standard generating set, acts on $\Sigma$ by permutation of the vertices.

    Consider now the action of $W$ induced on the barycentric subdivision $\Sigma'$. The subcomplex of all subsimplices of the maximal flag
    \[
        c_0 = (\{1\} \subset \{1,2\} \subset \cdots \subset \{1,\ldots,n\})
    \]
    is obviously a fundamental domain of the action and $W$ acts regularly on maximal flags. By calculating the stabilisers of subsimplices of $c_0$, it is not hard to see that $\Sigma'$ is isomorphic to the Coxeter complex $\Sigma(W,S)$. The type of a vertex is simply given by its cardinality.
\end{Example}

\begin{Example}\label{ex:cn_coxeter_complex}
    We have seen in Example \ref{ex:cn_coxeter_group} that the Coxeter group of type $C_n$ can be identified with the group of signed permutations on $[\pm n]$. Similarly to the example of type $A_n$, the Coxeter complex of type $C_n$ can be interpreted as the barycentric subdivision of the boundary of the \emph{$n$-hyperoctahedron} or \emph{$n$-cross polytope}, which is the convex hull of the $2n$ vectors $\{\pm e_i\}$ in $\R^n$, where $e_i$ is the standard basis. The $W$-action is given by permutation of the `signed' basis vectors. Equivalently, one can also consider the barycentric subdivision of the boundary of the dual polytope, which is the $n$-dimensional hypercube. Again, types of vertices are given by their cardinalities.
\end{Example}

\noindent Links in Coxeter complexes are again Coxeter complexes.

\begin{Lemma}[Proposition 3.16 in \cite{AB:B:08}]
    Let $\sigma\in\Sigma(W,S)$ be a simplex of type $I\setminus J$. Then $\lk_\Sigma(\sigma)$ is isomorphic to the Coxeter complex $\Sigma(W_J,S_J)$, where we set $S_J=\{ s_j: j\in J\}$ and $W_J=\langle S_J\rangle$.
\end{Lemma}

\begin{Remark}
    In particular, the Coxeter diagram of this link is obtained very easily: Just take the subdiagram induced by all the vertices in $J$.
\end{Remark}

\noindent Once we know that a simplicial complex $\Sigma$ with a given type function $\type: \Sigma^0 \rightarrow I$ is a Coxeter complex, we can reconstruct the Coxeter matrix.

\begin{Proposition}[Corollary 3.20 in \cite{AB:B:08}]\label{prop:coxeter_complex_determines_matrix}
    If $\Sigma=\Sigma(W,S)$ is a Coxeter complex with type function $\type: \Sigma^0\rightarrow I$, the Coxeter matrix is given by
    \[
    m_{i,j} = \diam( \lk_{\Sigma(W,S)}(\sigma)),
    \]
    where $\sigma\in\Sigma$ is a simplex satisfying $\type(\sigma)=I\setminus\{i,j\}$. In particular, a Coxeter complex $\Sigma$ determines the associated Coxeter system.
\end{Proposition}

\noindent Remember that $W$ acts on $\Sigma(W,S)$ by left multiplication. One can see that all type-preserving automorphisms arise in this way.

\begin{Lemma}[Proposition 3.32 in \cite{AB:B:08}]\label{lem:simplicial_automorphisms_of_coxeter_complexes}
	The monomorphism
    \[
        W \hookrightarrow \Aut(\Sigma(W,S))
    \]
    surjects onto the group of type-preserving automorphisms $\Aut^0(\Sigma(W,S))$.
\end{Lemma}

\noindent Since $\Sigma(W,S)$ is a chamber complex, we can construct the \emph{chamber graph} as in the second chapter: Its vertex set is the set of chambers $\Ch(\Sigma(W,S))$ where two chambers are connected by an edge if they are adjacent. We endow the chamber graph with the usual distance function $d$, where each edge has length $1$.  Remember that paths in the chamber graph are called \emph{galleries}. 

\begin{Definition}
    A gallery between two chambers $c$ and $c'$ is \emph{minimal} if its length is the distance between $c$ and $c'$ in the chamber graph.
\end{Definition}

\begin{Remark}
    The chamber graph of $\Sigma(W,S)$ is isomorphic to the Cayley graph of $W$ with respect to the generating set $S$.
\end{Remark}

\noindent It can be shown that the simplicial complex $\Sigma(W,S)$ can be recovered from the set of chambers with the adjacency relation if we keep track of the types of panels. This structure is usually called a \emph{chamber system}. Some facts about Coxeter complexes are easier to state in the language of chamber systems. When translated to the language of simplicial complexes, the formulation requires the notion of induced subcomplexes.

\begin{Definition}
    The \emph{subcomplex $\bar C$ induced by a set of chambers $C\subseteq \Ch(\Sigma(W,S))$} is given by
    \[
	\bar C = \{ \sigma\in\Sigma : \sigma\subseteq c \in C \}.
    \]
    By abuse of notation, we will often denote the induced subcomplex also by $C$.
\end{Definition}

\noindent Now we can define convex subcomplexes.

\begin{Definition}
    A set of chambers is called \emph{convex} if it contains all minimal galleries between its chambers. A subcomplex induced by a convex set of chambers is also called convex.
\end{Definition}

\begin{Remark}
	We will later also define a metric on the geometric realisations of Coxeter complexes. The notion of convexity of sets of chamber is different from metric convexity, however. A metrically convex subcomplex need not be convex in the sense above.
\end{Remark}

\noindent We will require another important feature of Coxeter complexes --- the existence of projections.

\begin{Proposition}[Proposition 3.103 and 3.105 in \cite{AB:B:08}]\label{prop:projections}
    Let $\sigma\in\Sigma(W,S)$ be a simplex and let $c\in\Sigma(W,S)$ be a chamber. Then there is a unique chamber $e$ in the residue $R_\Sigma(\sigma)$ called the \emph{projection of $c$ onto $R_\Sigma(\sigma)$} such that for any chamber $d\in R_\Sigma(\sigma)$ there is a minimal gallery from $c$ to $d$ through $e$.
\end{Proposition}

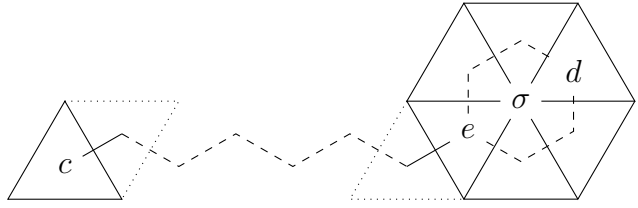
\begin{figure}
    \centering
        \begin{tikzpicture}
            \node (s) at (0,0) {$\sigma$};
            \foreach \x in {0,60,120,180,240,300} {
                \draw (\x+60:1.5) -- (\x:1.5) -- (s);
            };
            \node (e) at (210:0.8) {$e$};
            \node (c) at (-6,-0.866) {$c$};
            \draw (-6,0) -- ++(-0.75,-1.3) -- ++(1.5,0) -- cycle;
            \draw[dotted] (-6,0) -- ++(1.5,0) -- ++(-0.75,-1.3);
            \draw[dotted] (-1.5,0) -- ++(-0.75,-1.3) -- ++(1.5,0);
            \draw (c) -- (-5.25,-0.433);
            \draw[dashed] (-5.25,-0.433) -- ++(0.75,-0.433) -- ++(0.75,0.433) -- ++(0.75,-0.433) -- ++(0.75,0.433) -- ++(0.75,-0.433);
            \draw (-1.5,-0.866) -- (e);

            \node (d) at (30:0.8) {$d$};
            \draw[dashed] (e) -- (150:0.8) -- (90:0.8) -- (d);
            \draw[dashed] (e) -- (270:0.8) -- (330:0.8) -- (d);
        \end{tikzpicture}\caption{Projections in a Coxeter complex} \label{fig:projections}
\end{figure}

\noindent This property is also called the \emph{gate property}. See Figure \ref{fig:projections} for an illustration.

\subsection{Roots}\label{subsec:roots}

Roots are very important subcomplexes of Coxeter complexes. In the case of spherical or affine Coxeter complexes, roots can be thought of as simplicial half-spheres or half-spaces, respectively. We will define roots via their sets of chambers and take the induced subcomplexes.

\begin{Definition}
	Let $c\sim c'$ be adjacent chambers in a Coxeter complex $\Sigma(W,S)$. The set of chambers
	\[
    \Ch {\mathbf\alpha}(c,c') = \{ e\in\Ch(\Sigma(W,S)) : d(e,c)<d(e,c') \}
	\]
	is called a \emph{root} or \emph{half-apartment}. The associated induced subcomplex is denoted by $\alpha(c,c')$ and also called root or half-apartment. The set of roots of $\Sigma(W,S)$ is denoted by $\Phi(W,S)$. We will usually denote roots by $\alpha$ or $\beta$ without specifying the chambers $c$ and $c'$ if they are not important.
\end{Definition}

\begin{Terminology}
    If $\alpha = \alpha(c,c')$ and $\beta = \alpha(c',c)$, then we say that $\beta$ is the root \emph{opposite of $\alpha$} and we write $\beta=-\alpha$.

	The intersection of the two subcomplexes associated to roots $\alpha$ and $-\alpha$ is called a \emph{wall}, denoted by $\partial\alpha$.
\end{Terminology}

\begin{Lemma}[Proposition 3.79 in \cite{AB:B:08}]
    Let $\Sigma=\Sigma(W,S)$ be a Coxeter complex and let $\sigma\in\Sigma$ be a simplex. Then the roots of $\lk_\Sigma(\sigma)$ are given by
    \[
    \{ \alpha \cap \lk_\Sigma(\sigma) : \alpha\in\Phi(W,S), \sigma\in\partial\alpha\}.
    \]
    The walls of $\lk_\Sigma(\sigma)$ are given by
    \[
    \{ \partial\alpha \cap \lk_\Sigma(\sigma) : \alpha\in\Phi(W,S), \sigma\in\partial\alpha\}.
    \]
\end{Lemma}

\begin{Lemma}[Lemmas 3.44 and 3.45 in \cite{AB:B:08}]\label{l:roots_are_convex}
	Roots are convex. Additionally, we have
    \[
        \Ch(\alpha) \sqcup \Ch(-\alpha) = \Ch(\Sigma(W,S)).
    \]
\end{Lemma}

\noindent In fact, roots can be used to define convex hulls and convexity can be expressed using roots as well as minimal galleries.

\begin{Definition}
    For any set of chambers $C\subseteq \Ch(\Sigma(W,S))$, its \emph{convex hull} is the intersection of all roots containing $C$.
\end{Definition}

\begin{Lemma}[Theorem 3.131 in \cite{AB:B:08}]
    A set of chambers is convex if and only if its induced subcomplex is equal to its convex hull.
\end{Lemma}

\noindent The following Lemma is a very useful tool.

\begin{Lemma}[Proposition 2.7 in \cite{Ron:LoB:89}]\label{lem:roots_separating_chambers}
    Let $c$ and $d$ be chambers in $\Sigma(W,S)$ and fix a minimal gallery $c=d_0,d_1,\ldots,d_n=d$ from $c$ to $d$. Write $\alpha_i=\alpha(d_i,d_{i+1})$ for the roots containing one chamber of the gallery but not the next one. Then these roots are mutually distinct and the set $\{\alpha_i:i=0,1,\ldots,n-1\}$ is precisely the set of roots containing $c$ but not $d$.
\end{Lemma}

\noindent The notion of a prenilpotent pair of roots plays a large role later on.

\begin{Definition}
    A pair of roots $\{\alpha,\beta\}$ is called \emph{prenilpotent}, if both intersections $\alpha\cap\beta$ and $-\alpha\cap-\beta$ contain a chamber.
\end{Definition}

\begin{Definition}\label{def:root_intervals}
    Let $\{\alpha,\beta\}$ be a prenilpotent pair of roots of a Coxeter complex $\Sigma(W,S)$. We define the \emph{closed interval} by
    \[
		[\alpha,\beta] \coloneq \{ \gamma \in \Phi(W,S) : \alpha\cap\beta \subseteq \gamma, (-\alpha\cap-\beta)\subseteq -\gamma \}
    \]
    and the \emph{half-open} and \emph{open intervals} by
    \[
		[\alpha,\beta) = [\alpha,\beta] \setminus \{\beta\},\qquad (\alpha,\beta] = [\alpha,\beta] \setminus \{\alpha\}, \qquad (\alpha,\beta) = [\alpha,\beta] \setminus \{\alpha,\beta\}.
    \]
\end{Definition}

\begin{Remark}
    The notion of a prenilpotent set of roots is particularly simple for finite Coxeter systems. A pair of roots $\{\alpha,\beta\}$ in a finite Coxeter system is prenilpotent if and only if $\alpha\neq-\beta$. In this case, we have
    \[
    [\alpha,\beta] = \{\gamma\in\Phi(W,S) : \alpha\cap\beta\subseteq\gamma\}.
    \]
    See \cite[Example 8.43]{AB:B:08} for details.
\end{Remark}

\noindent For finite rank two Coxeter complexes, root intervals are not difficult to describe.

\begin{Example}\label{ex:root_intervals_rank_two}
    Let $\Sigma$ be a finite Coxeter complex of rank two, hence of dimension one. Then $\Sigma$ is an ordinary $2m$-gon as in Example \ref{ex:i2_complex}. Its roots are paths of length $m$. We denote the vertices of $\Sigma$ by $(x_i)_{i\in \Z/2m}$. This leads to an enumeration of the roots $(\alpha_i)_{i\in \Z/2m}$ such that every root $\alpha_i$ contains the vertices $\{x_i,x_{i+1},\ldots,x_{i+m}\}$.
    Then it is easy to see that
    \[
        [\alpha_0,\alpha_i]=\{\alpha_0,\alpha_1,\ldots,\alpha_i\},
    \]
    for $1\leq i\leq m-1$, which justifies the name `root interval'.
\end{Example}

\begin{Remark}
    For any root $\alpha\in\Phi(W,S)$, there is a reflection $s_\alpha$ which fixes the wall $\partial\alpha$ and interchanges $\alpha$ and $-\alpha$. These reflections are conjugate to elements of $S$. We will not construct the reflections explicitly in this thesis. They only appear in the formulation of the next Lemma and briefly in Chapter \ref{ch:wagoner}. In any case, the elements $s_\alpha$ are not required for the parts of the statements used in this thesis, but instead to understand the references given in the text.
\end{Remark}

\noindent In infinite Coxeter complexes, there are two possibilities for prenilpotent pairs of roots.

\begin{Lemma}[Proposition 3.165 and Lemma 8.45 in \cite{AB:B:08}]\label{lem:nested_or_finite_order}
    Let $\{\alpha,\beta\}$ be a prenilpotent pair of roots in $\Phi(W,S)$. Then there are two possibilities
\begin{itemize}
    \item The product $s_\alpha s_\beta$ has infinite order. In this case, the roots $\alpha$ and $\beta$ are \emph{nested}, that is, we have $\alpha\subseteq\beta$ or $\beta\subseteq\alpha$.
    \item The product $s_\alpha s_\beta$ has finite order. Then every maximal simplex $\sigma$ of the intersection of the walls $\partial\alpha\cap\partial\beta$ has codimension two and the link $\lk_\Sigma(\sigma)$ is finite. We write $\bar\alpha=\alpha\cap\lk_\Sigma(\sigma)$, $\bar \beta=\beta\cap\lk_\Sigma(\sigma)$ for the induced roots in the link. The map given by
    \begin{align*}
        [\alpha,\beta] &\rightarrow [\bar\alpha,\bar\beta]\\
        \gamma & \mapsto \gamma\cap\lk_\Sigma(\sigma)
    \end{align*}
    is a bijection.
\end{itemize}
\end{Lemma}

\noindent This can be used to prove the following characterisation for half-open intervals which we will require later on.

\begin{Lemma}\label{l:root_intervals}
    For any prenilpotent, non-nested pair of roots $\{\alpha,\beta\}$, there are either roots $\alpha'$ and $\beta'$ such that
    \[
        (\alpha,\beta] = [\alpha',\beta],\qquad [\alpha,\beta) = [\alpha,\beta']
    \]
    or $(\alpha,\beta] = \{\beta\}$ and $[\alpha,\beta)=\{\alpha\}$.
\end{Lemma}

\begin{Proof}
    If the rank of $\Sigma$ is two and $\Sigma$ is finite, we fix an enumeration of the vertices and of the roots of $\Sigma$ as in Example \ref{ex:root_intervals_rank_two}. For root intervals of the form $[\alpha_0,\alpha_i]$, we obtain by the description: $[\alpha_0,\alpha_1]=\{\alpha_0,\alpha_1\}$ and $(\alpha_0,\alpha_i] = [\alpha_1,\alpha_i]$ and $[\alpha_0,\alpha_i) = [\alpha_0,\alpha_{i-1}]$ for $i>1$. Since the numbering of vertices was arbitrary, the Lemma follows for finite rank two Coxeter complexes.

    If the rank of the Coxeter complex $\Sigma$ is greater than two, and if $\{\alpha,\beta\}$ is a prenilpotent, non-nested pair of roots, then by Lemma \ref{lem:nested_or_finite_order}, we can reduce to the rank two case as follows. Fix a maximal simplex $\sigma$ in the intersection of the walls $\partial\alpha\cap\partial\beta$ and write $\bar\alpha=\alpha\cap\lk_\Sigma(\sigma)$, $\bar \beta=\beta\cap\lk_\Sigma(\sigma)$ as above. Then consider the bijection
    \begin{align*}
        [\alpha,\beta] &\rightarrow [\bar\alpha,\bar\beta]\\
        \gamma & \mapsto \gamma\cap\lk_\Sigma(\sigma).
    \end{align*}
    Now the lemma follows by applying the finite rank two case to $\lk_\Sigma(\sigma)$ and using the bijection to transport the result to the ambient building.
\end{Proof}

\noindent Finally, root intervals can always be ordered in the following fashion.

\begin{Lemma}[2.2.6 in \cite{Rem:GKM:02}]\label{lem:root_intervals_order}
    Let $\{\alpha,\beta\}$ be a prenilpotent pair of roots. Fix a chamber $c\in\alpha\cap\beta$ and a chamber $d\in (-\alpha\cap-\beta)$. Fix a minimal gallery from $c$ to $d$ and assume that it crosses the wall $\partial\alpha$ before the wall $\partial\beta$. Then $\partial\alpha$ is the first wall crossed among the walls $\partial\gamma$ where $\gamma\in[\alpha,\beta]$ and $\partial\beta$ is the last.

    The order of the walls crossed by the gallery induces an order on the root interval, which depends on the choice of the gallery. 
\end{Lemma}

\subsection{The canonical linear representation}

In this section, we will see that every abstract Coxeter group acts as a group of reflections on a finite-dimensional vector space. This will be used to define spherical and affine Coxeter groups, which act as groups of reflections on spheres and as groups of affine reflections on Euclidean spaces, respectively.

Let $M$ be a Coxeter matrix with type set $I$ and let $(W,S)$ denote the associated Coxeter system. We construct a real vector space $E$ with a bilinear form $b$ on which the Coxeter group $W$ acts faithfully, preserving $b$.

\begin{Construction}
    Let $E=\R^I$ be the real vector space with basis $(e_i)_{i\in I}$. Let $b$ be the symmetric bilinear form on $E$ such that
    \[
    b(e_i,e_j) = - \cos( \pi / m_{i,j} ).
    \]
    Consider the map
    \begin{align*}
        \rho: W &\rightarrow \Gl(E)
        \intertext{defined on $S$ by}
        s_i&\mapsto \Bigl( x\mapsto \frac{x - 2b(e_i,x)e_i}{b(x,x)}\Bigr).
    \end{align*}
    This is well defined, since the images of generators are obviously involutions and it is a simple calculation to verify that the other relations hold. See \cite[Section 2.5]{AB:B:08} for additional details.

    Then $\rho$ is a faithful representation of $W$ in $\Gl(E)$ which preserves the bilinear form $b$.
\end{Construction}

\subsection{Spherical Coxeter complexes}\label{ss:spherical_coxeter_complexes}

The simplest Coxeter groups are finite Coxeter groups. The classification of finite groups generated by reflections on a real vector space is a classical, famous result by Coxeter in \cite{Cox:DGR:34}. It turns out that such finite reflection groups are precisely the finite Coxeter groups.

Throughout this section, let $(W,S)$ be a Coxeter system and let $\rho:W\rightarrow \Gl(E)$ be the canonical representation, where $b$ is the associated bilinear form on $E$.

\begin{Proposition}[Corollary 2.68 in \cite{AB:B:08}]
    The following two conditions are equivalent:
    \begin{itemize}
        \item The bilinear form $b$ is positive definite.
        \item The Coxeter group $W$ is finite.
    \end{itemize}
    In this case, $W$ is called \emph{spherical}.
\end{Proposition}

\noindent The spherical Coxeter groups are classified.

\begin{BreakTheorem}[VI.\S 4, Theorem 1 and 2 in \cite{Bou:LaL:02}]
    An irreducible Coxeter group is finite if and only if its diagram appears on the list \ref{fig:spherical_coxeter_diagrams} in the appendix.
\end{BreakTheorem}

\noindent A spherical Coxeter group $W$ preserves the unit sphere $\Sph(E)$ of $E$ with respect to the bilinear form $b$. Denote by $\caH$ the set of all hyperplanes in $E$ which are left invariant by conjugates of $S$ in $W$.

For all hyperplanes $H\in \caH$ we fix a linear functional $h$ on $E$ such that $H=H^0\coloneq h^{-1}(0)$ and write $H^+\coloneq h^{-1}( [0,\infty))$ and $H^-\coloneq h^{-1}( (-\infty,0])$ for the corresponding closed half-spaces of $E$.

\begin{Construction}
    A \emph{cell in $E$} with respect to $\caH$ is an intersection
    \[
    C = \bigcap_{H\in \caH} H^{\sigma(H)}
    \]
    such that $\sigma : \caH \rightarrow \{+,0,-\}$ is any function, called the corresponding \emph{sign function}.

    Now consider the set of cells
    \[
    \cC(W,S) = \{  C\cap \Sph(E) : C \text{ is a cell in $E$}\}.
    \]
\end{Construction}

\begin{Proposition}[Theorem 1.111 in \cite{AB:B:08}]\label{prop:coxeter_complex_triangulated_sphere}
    The set $\cC(W,S)$ partially ordered by inclusion is isomorphic to the Coxeter complex $\Sigma(W,S)$.

    In particular, vertices of the simplicial complex $\cC(W,S)$ are points on the sphere $\Sph(E)$. We denote the unit vector corresponding to the vertex $v\in \Sigma(W,S)^0$ by $e_v$.
\end{Proposition}

\noindent The simplicial complex $\cC(W,S)$ can hence be seen as a triangulation of the sphere $\Sph(E)$.

\begin{Construction}
Now consider the following map
\begin{align*}
    \lvert\cC(W,S)\rvert & \rightarrow \Sph(E) \\
    \intertext{given on the cell $C$ by}
    \sum_{v\in C} \lambda_v v & \mapsto \frac{ \sum_{v\in C} \lambda_v e_v } {\sqrt{b(\sum_{v\in C} \lambda_v e_v,\sum_{v\in C} \lambda_v e_v)}}
\end{align*}
which is obviously a bijection. We use this bijection and the isomorphism of Proposition \ref{prop:coxeter_complex_triangulated_sphere} to pull back the angular metric of the sphere and obtain a metric on $|\Sigma(W,S)|$.
\end{Construction}

\begin{Remark}
    An equivalent method is given by endowing all simplices of $|\Sigma(W,S)|$ with the spherical metric pulled back from a single cell $|C\cap\Sph(E)|$ in this fashion and considering the intrinsic metric, see \cite[Example 12.39]{AB:B:08}.
\end{Remark}

\begin{Definition}
    From the construction of the complex $\cC(W,S)$, it is obvious that the involution $E\rightarrow E$ given by $v\mapsto -v$ induces a simplicial map on $\cC(W,S)$. We consider the induced involutory automorphism on $\Sigma(W,S)$. A simplex and its image under this automorphism are said to be \emph{opposite}.

    Note that this opposition involution need not be type-preserving. It induces a bijection on the type set $I$. For any $J\subseteq I$, we denote the image under this type bijection by $J^0$.
\end{Definition}

\begin{Remark}
    The concept of opposition is a very strong tool for spherical buildings, but it does not make sense in non-spherical buildings. However, there is a way to combine two buildings by means of a so-called \emph{co-distance}, such that simplices in one of the building can be said to have opposite simplices in the other building. This construction is known as a \emph{twin building}. Twin buildings are intimately related to groups with twin root data or groups of Kac-Moody type, which we will define later. However, since we do not make use of twin buildings, we will omit their definition.
\end{Remark}

\noindent The concept of opposition can also be expressed in a purely combinatorial fashion. For chambers, this is particularly simple.

\begin{Lemma}[Proposition 1.57 in \cite{AB:B:08}]
    Two chambers in $\Sigma(W,S)$ are opposite if and only if they are at maximal distance in the chamber graph.
\end{Lemma}

\noindent For arbitrary simplices, opposition can be expressed best in terms of the \emph{longest word}.

\begin{Lemma}[Proposition 1.77 in \cite{AB:B:08}]
    Let $(W,S)$ be a spherical Coxeter system. There is a unique element $w_0\in W$ whose word length is maximal with respect to the generating set $S$, called the \emph{longest element of $W$}. The longest element $w_0$ is of order two.
\end{Lemma}

\begin{Proposition}[Lemma 5.111 in \cite{AB:B:08}]\label{prop:opposition_via_w0}
    Let $(W,S)$ be a spherical Coxeter system. Two chambers $w$ and $w'$ in $\Sigma(W,S)$ are opposite if and only if $w=w'w_0$. The simplex opposite of the simplex $wW_J$ is given by $ww_0W_{J^0}$.
\end{Proposition}

\begin{Example}\label{ex:an_opposition}
    Let $(W,S)$ be the Coxeter system of type $A_n$ with type set $I$. Consider the associated Coxeter complex, which can be considered as the flag complex over the partially ordered set of all non-trivial proper subsets of $\{1,\ldots,n+1\}$ as in Example \ref{ex:an_coxeter_complex}. The type of a vertex is given by its cardinality as a subset of $\{1,\ldots,n+1\}$.

    By \cite[5.7.4]{AB:B:08}, in this case the opposite type $J^0$ is given by applying the unique non-trivial diagram automorphism to the diagram $A_n$. For a vertex of type $j$, the opposite type is hence given by $n+1-j$. By Proposition 1.5.2 in \cite{BB:CoG:05}, the word length in the symmetric group is given by the number of inversions. Clearly, the element $w_0$ of the symmetric group $\Sym(n+1)$ with the maximal number of inversions is given by $w_0(k)=n+2-k$.

    Hence the chamber opposite of the standard chamber
    \[
        c_0 = (\{1\} \subset \{1,2\} \subset \cdots \subset \{1,\ldots,n\})
    \]
    is given by
    \[
        w_0(c_0) = ( \{n+1\} \subset \{n, n+1\} \subset \cdots \subset \{ 2,\ldots, n+1\}).
    \]
    Since the type of a vertex is given by its cardinality, we obtain in particular that the vertex opposite of $\{1,\ldots,k\}$ is given by $\{k+1,\ldots,n+1\}$.
\end{Example}

\begin{Example}\label{ex:cn_opposition}
    Now let $(W,S)$ be the Coxeter system of type $C_n$ with type set $I$. The associated Coxeter complex can be considered as the flag complex over the partially ordered set of all `signed' subsets of $[\pm n]$ as in Example \ref{ex:cn_coxeter_complex}. By 5.7.4 in \cite{AB:B:08}, types and their opposites coincide for type $C_n$, so $J^0=J$ for all $J\subseteq I$.

    By Proposition 8.1.1 in \cite{BB:CoG:05}, the word length in $W$ is given by:
    \[
    l(w) = \inv(w) + \sum_{ \{k\in\{1,\ldots,n\} \,:\, w(k)<0\}} (-w(k)),
    \]
    where $\inv(w)$ is the number of inversions of the permutation $w$ on $[\pm n]$. Clearly, the second sum is maximal if $w(k)<0$ for all $k>0$. In this case, the number of inversions is maximal for the longest word given by $w_0(k)=-k$ for all $k\in[\pm n]$. Applying this to the Coxeter complex, we obtain that the vertex $\{1,\ldots,k\}$ is opposite of the vertex $\{-1,\ldots,-k\}$.
\end{Example}

\subsection{Affine Coxeter complexes}\label{ss:affine_coxeter_complexes}

Affine Coxeter groups are precisely those which can be interpreted as a group of affine isometries of a Euclidean space, generated by affine reflections.

\begin{Definition}
    An irreducible Coxeter system $(W,S)$ and its associated Coxeter complex $\Sigma(W,S)$ are said to be \emph{affine} or \emph{Euclidean} if the bilinear form $b$ from the canonical linear representation is positive but degenerate.
\end{Definition}

\noindent These affine Coxeter groups can be classified as well.

\begin{Theorem}[VI.\S 4, Theorem 4 in \cite{Bou:LaL:02}]
    A Coxeter system is affine if and only if its Coxeter diagram is in list \ref{fig:affine_coxeter_diagrams} in the appendix.
\end{Theorem}

\noindent The construction of a Euclidean space on which the group $W$ acts by affine reflections is more involved than in the spherical case. We will hence omit the explicit construction here.

\begin{Definition}
    Assume that $V$ is a Euclidean vector space of dimension $n\geq 1$. Let $W$ be a group of affine isometries of $V$. Then $W$ is an \emph{affine reflection group} if there is a set of affine hyperplanes $\caH$ of $V$ such that
    \begin{itemize}
        \item The group $W$ is generated by the orthogonal reflections $s_H$ at the hyperplanes $H\in \caH$.
        \item The set of hyperplanes $\caH$ is $W$-invariant.
        \item The set $\caH$ is locally finite, meaning that every point of $V$ has a neighbourhood that meets only finitely many hyperplanes $H\in\caH$.
    \end{itemize}
\end{Definition}

\begin{Theorem}[Theorem 10.42 and Proposition 10.44 in \cite{AB:B:08}]
    Any affine Coxeter system $(W,S)$ has a canonical realisation as an affine reflection group acting on an affine hyperplane $A$ in the dual space $E^*$ of the image of the canonical linear representation $\rho: W \rightarrow \Gl(E)$.
\end{Theorem}

\noindent This realisation as an affine reflection groups allows a concrete realisation of the Coxeter complex.

\begin{Construction}
    Let $(W,S)$ be an affine Coxeter system and consider the action of $W$ on the associated affine space $A$. For every hyperplane $H\in \caH$, we write $H^0=H$ and $H^+$ and $H^-$ for the two closed half-spaces bounded by $H$, making arbitrary choices.

    If we consider the set of cells
    \[
    \tilde \cC(W,S) = \Bigl\{ \bigcap_{H\in \caH} H^{\sigma(H)} \,\Big|\, \sigma: \caH \rightarrow \{-,0,+\}\text{ sign function}\Bigr\},
    \]
    we obtain a simplicial complex as in the spherical case above. Vertices of cells correspond to elements of $A$, which are denoted by $a_v$ for $v\in \Sigma(W,S)^0$.
\end{Construction}

\noindent We visualise this simplicial complex as a triangulation of the affine hyperplane $A$.

\begin{Proposition}[Proposition 10.13 in \cite{AB:B:08}]\label{prop:coxeter_complex_triangulated_euclidean_space}
    Let $(W,S)$ be an affine Coxeter system. The simplicial complex $\tilde \cC(W,S)$ is isomorphic to $\Sigma(W,S)$.
\end{Proposition}

\begin{Construction}
 Now consider the following map
\begin{align*}
    \lvert \tilde\cC(W,S)\rvert & \rightarrow A \\
    \intertext{given on the cell $C$ by}
    \sum_{v\in C} \lambda_v v & \mapsto \sum_{v\in C} \lambda_v a_v
\end{align*}
which is obviously a bijection. We use this bijection and the isomorphism of Proposition \ref{prop:coxeter_complex_triangulated_euclidean_space} to pull back the Euclidean metric of the hyperplane $A$ and obtain a metric on the geometric realisation $\lvert\Sigma(W,S)\rvert$.
\end{Construction}

\noindent We will show that an affine Coxeter complex has a simplicial sphere `at infinity'. For this, we first require the notions of special vertices and gems.

\begin{Definition}
    A vertex of type $\{j\}$ of the affine Coxeter complex $\Sigma(W,S)$ is called \emph{special} if the removal of the vertex of the Coxeter diagram corresponding to $j$ creates a Coxeter diagram of the same name, but removing the tilde $\sim$.
    A residue $R_\Sigma(v)$ of a special vertex $v$ is called a \emph{gem}.
\end{Definition}

\noindent For example, any vertex in a Coxeter system of type $\tilde A_n$ is special and only those vertices corresponding to the two rightmost vertices in the Coxeter diagram of type $\tilde B_n$ are special.

\begin{Remark}
    We will later see that affine buildings always have associated spherical buildings. In order to minimize confusion, roots in affine Coxeter complexes will usually be called \emph{half-apartments}, while the term `roots' will be reserved for roots in the associated spherical Coxeter complex.
\end{Remark}

\noindent Using gems, we can define sectors.

\begin{figure}[hbt]
    \centering\begin{tikzpicture}
            \fill[lightgray!45] (0:4) -- (0,0) -- (60:4) -- (0:4);
            \fill[lightgray!17] (0:4) -- (0,0)  -- (240:4) -- (300:4) -- (0:4);
            \fill[lightgray!17] (60:4) -- (120:4) -- (180:4) -- (0,0) -- (60:4);
            \node (s) at (0,0) {$v$};
            \foreach \x in {0,60,120,180,240,300} {
                \draw (\x+60:1.5) -- (\x:1.5) -- (s);
                \draw[dotted] (\x:1.5) -- (\x+30:2.598);
                \draw[dotted] (\x:1.5) -- (\x-30:2.598);
                \draw[dotted] (\x:1.5) -- (\x:3) -- (\x+60:3);
            };
            \node (c) at (30:0.9) {$c$};
            \node (a) at (170:3.2) {$\alpha$};
            \node (aa) at (250:3.2) {$\alpha'$};
            \fill[lightgray!45] (20:3) rectangle (45:2.1);
            \node (si) at (30:2.5) {$\fs(R,c)$};
            \node (R) at (210:1.8) {$R$};
            \draw (0:4) -- (s) -- (60:4);
            \draw (180:4) -- (s) -- (240:4);
        \end{tikzpicture}\caption{Sectors in affine Coxeter complexes}
    \label{fig:def_sector}
\end{figure}
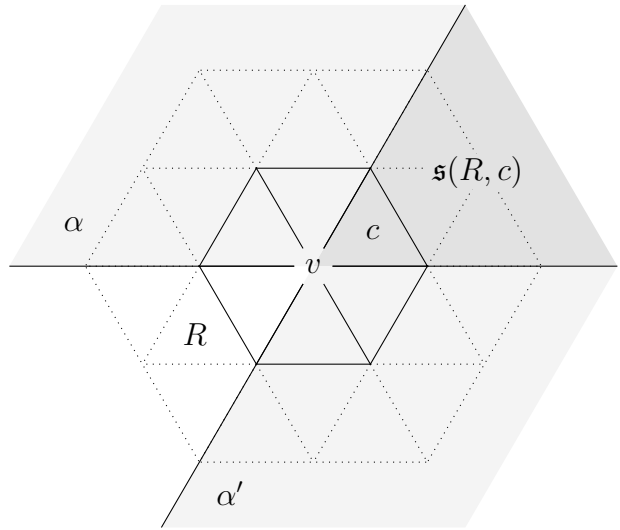

\begin{Definition}
    Let $R$ be a gem and let $c$ be a chamber in $R$. Let $[R,c]$ be the set of half-apartments that contain $c$ and whose opposite half-apartment contains one of the chambers of $R$ adjacent to $c$. A \emph{sector} $\fs$ of $\Sigma(W,S)$ is a subcomplex of the form
    \[
    \fs(R,c) \coloneq \bigcap_{\alpha\in[R,c]} \alpha.
    \]
    A sector is a simplicial cone and its faces are called \emph{sector faces}.

    Figure \ref{fig:def_sector} illustrates the definition. There, $R=R_\Sigma(v)$ is a gem and $\alpha$ and $\alpha'$ are the two roots separating $c$ from adjacent chambers in $R$.
\end{Definition}

\noindent Later, we will consider an equivalence relation on sectors and faces given by `parallelism'. It can be shown that the equivalence classes of sectors and sector faces of an affine Coxeter complex form a spherical Coxeter complex called the Coxeter complex \emph{at infinity}. We will discuss this directly for affine buildings in the next section.

\section{Buildings}

Using Coxeter complexes as subcomplexes, we can now give a definition of buildings. As mentioned in the introduction of this chapter, buildings can be defined in many ways. The following definition uses the language of simplicial complexes we have introduced.

It turns out that we will only require very few results on buildings themselves for this thesis. It will be enough to know about their homotopy type and the canonical metrics.

\subsection{The definition of buildings}

Roughly speaking, a building is a highly symmetric simplicial complex covered by the union of many so-called apartments, which are copies of a fixed Coxeter complex.

\begin{Definition}
	A simplicial complex $X$ together with a collection of subcomplexes $\cA$ called \emph{apartments} is a \emph{building} if
	\begin{itemize}
		\item Every apartment is a Coxeter complex.
		\item Every two simplices of $X$ are contained in a common apartment.
		\item For any two apartments $A_1,A_2\in \cA$ containing a common chamber, there is an isomorphism $A_1\rightarrow A_2$ fixing $A_1\cap A_2$ pointwise.
	\end{itemize}
    The building $X$ is called \emph{thick} if every panel is contained in at least three chambers. We will call a building \emph{weak} if it is not necessarily thick.
\end{Definition}

\noindent Most of the language we have introduced for Coxeter complexes translates directly to buildings: \emph{Chambers}, \emph{panels}, \emph{galleries} and \emph{convexity} are defined as for Coxeter complexes. Roots of apartments are called \emph{roots of the building}. This notion is independent of the apartment: If a set is a root contained in an apartment, it is a root as a subset of any other apartment in which it is contained.

An immediate consequence of the axioms and of Proposition \ref{prop:coxeter_complex_determines_matrix} is the following:

\begin{Proposition}[Propositions 4.6 and 4.7 in \cite{AB:B:08}]
    All apartments have the same Coxeter matrix $M$ with index set $I$ and the building is colourable via a type function $\type:X^0\rightarrow I$. There is in particular a unique Coxeter system $(W,S)$ associated to $X$ and we say that $X$ is \emph{of type $(W,S)$}. If $W$ is spherical or affine, we say that $X$ is spherical or affine.
\end{Proposition}

\noindent The following terminology will become clearer in Section \ref{sec:generalised_polygons}.

\begin{Definition}
    A building of type $I_2(m)$ is called a \emph{generalised polygon} or {generalised $m$-gon}.
\end{Definition}

\noindent Like for Coxeter complexes, links in buildings are again buildings of smaller rank.

\begin{Proposition}[Propositions 4.9 and 3.79 in \cite{AB:B:08}]\label{prop:links_are_buildings}
    Let $\sigma\in X$ be a simplex of type $\type(\sigma)= I \setminus J$. Then $\lk_X(\sigma)$ is a building of type $(W_J,S_J)$. Its apartments, roots and walls are intersections of apartments, roots and walls of $X$ with the link.
\end{Proposition}

\begin{Remark}[Theorems 4.54 and 4.70 in \cite{AB:B:08}]
    It can be shown that the union of any family of apartment systems for $X$ is again an apartment system. Hence there is a largest system of apartments, called the \emph{complete system of apartments} for $X$. For spherical buildings, there is a unique system of apartments which is precisely the set of convex hulls of opposite chambers.
\end{Remark}

\noindent This classical result on the geometric realisations of buildings is central to many of our arguments.

\begin{Theorem}[Solomon-Tits, Theorem 4.127 in \cite{AB:B:08}]\label{th:solomon_tits}
    Let $X$ be a building. If $X$ is spherical of rank $n$, then the geometric realisation $|X|$ is homotopy equivalent to a bouquet of $(n-1)$-spheres. Furthermore, in this case, the top-dimensional homology group is generated by the apartments. If $X$ is not spherical, then $|X|$ is contractible.
\end{Theorem}

\subsection{Moufang spherical buildings}\label{subsec:moufang_spherical_buildings}

The Moufang condition ensures a large automorphism group of a building. For spherical buildings, this condition is both strong enough to allow for a complete classification and weak enough to be satisfied automatically in higher rank. The Moufang condition can also be defined for non-spherical buildings, but we restrict ourselves to spherical buildings here.

Let $X$ be a thick spherical building of type $(W,S)$ and of rank at least two. Fix an apartment $\Sigma_0$ in $X$ and denote the set of roots of the fixed apartment $\Sigma_0$ by $\Phi_0$.

\begin{Definition}
	For any root $\alpha\in\Phi_0$, the \emph{root group} $U_\alpha$ is the pointwise stabiliser of the union of all stars $\st_X(\sigma)$ over all panels $\sigma\in(\alpha\setminus\partial\alpha)$ in the full group of type-preserving automorphisms of $X$. The subgroup
	\[
		G^\dagger \coloneq \langle U_\alpha : \alpha\in\Phi_0\rangle
	\]
	is called the \emph{little projective group} of $X$.
\end{Definition}

\begin{Definition}
	The building $X$ is said to be \emph{Moufang} if every root group $U_\alpha$ for $\alpha\in\Phi_0$ acts transitively on the set of apartments containing $\alpha$. If this action is regular, the building is said to be \emph{strictly Moufang}.
\end{Definition}

\begin{Remarks}
	If the Coxeter diagram of $(W,S)$ has no isolated nodes, then these conditions are equivalent by \cite[Lemma 7.25]{AB:B:08}. In various other texts, the authors require Coxeter diagrams to have no isolated nodes or even irreducibility of $(W,S)$ when defining the Moufang condition so that the Moufang condition and the strict Moufang condition agree.

    It is not hard to see that the little projective group is transitive on the apartments of a Moufang building. This can be used to show that the definition above does not depend on the choice of $\Sigma_0$, see \cite[Remarks 7.29]{AB:B:08}.
\end{Remarks}

\noindent The Moufang condition is inherited by links.

\begin{Proposition}[Proposition 7.32 in \cite{AB:B:08}]\label{prop:links_are_moufang}
    If $\sigma$ is a simplex in a Moufang building $X$, then $\lk_X(\sigma)$ is also Moufang. If $X$ is strictly Moufang and the Coxeter diagram of $\lk_X(\sigma)$ has no isolated nodes, then $\lk_X(\sigma)$ is also strictly Moufang.
\end{Proposition}

\noindent It turns out that the Moufang condition is no restriction in higher rank.

\begin{Theorem}[Tits, Satz 1 in \cite{Tit:ESW:77}]\label{th:higher_rank_sph_buildings_are_moufang}
    If $X$ is a thick, irreducible and spherical building of rank at least $3$ then $X$ is strictly Moufang.
\end{Theorem}

\noindent In rank two, the Moufang condition imposes a severe restriction on the Coxeter diagram of the building. Much earlier, an analogous result for thick finite generalised polygons was proved by Feit and Higman, see Theorem \ref{th:feit_higman}.

\begin{Theorem}[Tits-Weiss, Theorem 17.1 in \cite{TW:MP:02}]
    Moufang generalised $m$-gons, $m\geq 3$, can only exist for $m\in\{3,4,6,8\}$.
\end{Theorem}

\noindent This theorem is the starting point for the full classification of Moufang polygons by Tits and Weiss. Another consequence is that, since links in Moufang buildings are again Moufang by Proposition \ref{prop:links_are_moufang}, thick spherical buildings of type $H_3$ or $H_4$ do not exist. This implies

\begin{Corollary}[Corollary 7.61 in \cite{AB:B:08}]
    If $X$ is a thick, irreducible and spherical building of rank $n\geq 3$, then $X$ is of type $A_n$, $C_n$, $D_n$, $E_6$, $E_7$, $E_8$ or $F_4$.
\end{Corollary}

\begin{Remark}
	In fact, irreducible spherical buildings of rank at least three have been classified by Tits in \cite{Tit:BsT:74}. At the heart of this classification is a rigidity result: A certain map between small parts of two buildings can always be extended to an isomorphism. With this result and Coxeter's classification of finite reflection groups, Tits uses a case-by-case analysis for each possible Coxeter diagram.

    Later, Moufang polygons were classified by Tits and Weiss in \cite{TW:MP:02}. Using this classification and Theorem \ref{th:higher_rank_sph_buildings_are_moufang}, a simpler and more structured classification of all Moufang spherical buildings of rank at least two is given by Tits and Weiss in \cite[Chapter 40]{TW:MP:02}. A detailed overview and the complete classification can be found in \cite{Wei:SSB:03}.

	The result of the classification is usually given by stating which groups correspond to the building. Hence we will give the result of the classification in Section \ref{sec:group_actions}, once we have explained the notion of strong transitivity.
\end{Remark}

\begin{Remark}
	The Moufang condition is necessary for the classification. Rank two buildings include all projective planes, for instance, which cannot be classified because of the existence of free constructions.
\end{Remark}

\subsection{Metrics on buildings}\label{subsec:metrics_on_buildings}

In Sections \ref{ss:spherical_coxeter_complexes} and \ref{ss:affine_coxeter_complexes}, we have seen that there are canonical metrics on spherical and affine Coxeter complexes. In this section, we will extend these metrics to spherical and affine buildings.

\begin{Theorem}[Theorem II.10A.4 in \cite{BH:NPC:99}]\label{th:metrics_on_buildings}
    The intrinsic metric on a spherical or affine building is the unique metric inducing the canonical metrics on the apartments. With this metric, the building is a CAT(1) or CAT(0) space, respectively. If the building is thick, then it is complete as a metric space.
\end{Theorem}

\noindent Since every two points in a building are contained in a common apartment, and the natural metric on a spherical Coxeter complex comes from a unit sphere, the diameter of a spherical building is $\pi$. The metric structure is naturally inherited by links.

\begin{Proposition}[Proposition 2.3 in \cite{CL:MC:01}]
    If $X$ is a spherical or affine building of dimension at least two with the intrinsic metric as above, then the geometric link of every vertex is naturally a spherical building. The induced metric on the geometric link is exactly the intrinsic metric of the spherical building, in particular every link is connected and has diameter $\pi$.
\end{Proposition}

\noindent In \cite{CL:MC:01}, Charney and Lytchak show that spherical and affine buildings can actually be characterised by these properties. We will give the metric recognition theorem for Euclidean buildings here, since we will use it in Chapter \ref{ch:lattices}. Since we will only consider two-dimensional affine buildings there, we quote the theorem only in the two-dimensional case.

\begin{Theorem}[Theorem 7.3 in \cite{CL:MC:01}]\label{th:recognition}
    Let $X$ be a connected two-dimensional piecewise Euclidean complex satisfying
    \begin{itemize}
        \item The complex $X$ is CAT(0).
        \item Every 1-cell is contained in at least three 2-cells ($X$ is thick).
        \item The link of every vertex is connected and has diameter $\pi$.
    \end{itemize}
    Then $X$ is a two-dimensional affine building or a product of metric trees.
\end{Theorem}

\subsection{Affine buildings}

Let $X$ be a thick building of type $(W,S)$, where the Coxeter system $(W,S)$ is of affine type. Then, as we have seen, apartments of the building are triangulated Euclidean spaces. Affine buildings have a canonically associated spherical building `at infinity', which we will describe now. The building at infinity itself will not play a large role in this thesis, but it is related to groups with \emph{root data with valuations}, which we will define later. It is one of the most important structures associated to an affine building and plays a large role in the classification of affine buildings.

\begin{Definition}
    A \emph{sector (sector face)} of the building is a sector (sector face) of any of its apartments. Two sectors (sector faces) $\fs$ and $\ft$ are said to be \emph{parallel} if the sets
    \[
    \{d(x,\ft) : x\in\fs\},\qquad \{d(\fs,y) : y\in \ft\}
    \]
    are bounded as subsets of $\R$.

    For sectors, this can be reformulated: Two sectors are parallel if and only if they contain a common subsector.
\end{Definition}

\noindent Obviously, the set of parallel classes of sectors and sector faces, denoted by $X^\infty$, inherits a simplicial structure. Let $I$ be the type set of $X$ and let $j\in I$ be the type of a special vertex.

\begin{Theorem}[Theorem 9.6 in \cite{Ron:LoB:89}]
    The simplicial complex $X^\infty$ is a thick building of type $I\setminus\{j\}$, called the \emph{building at infinity} of $X$.
\end{Theorem}

\noindent A type function on the building at infinity can be constructed by realising each parallel class of sector faces at a fixed gem $R_X(v)$ and defining the type of this sector face to be the type of the intersection with the link $\lk_X(v)$.

\begin{Remark}
    The building at infinity of a higher rank affine building is always Moufang by Theorem \ref{th:higher_rank_sph_buildings_are_moufang}. Moufang spherical buildings are classified, as we have seen above. Using some additional information constructed from the affine building, affine buildings of higher rank can also be classified. A more detailed explanation of the classification can be found in Section \ref{subsec:root_datum_valuation}.
\end{Remark}

\section{Groups acting on buildings}\label{sec:group_actions}

As we have seen in the previous section, every higher rank spherical building has a large automorphism group, its little projective group. Buildings with large automorphism groups are easier to understand and can often be classified.

In this section, we define the concept of a `good' group action on a building. We study the implications for the group and see that there is an algebraic datum, a so-called Tits system or BN pair, that allows us to reconstruct the building. In the same spirit, we discuss group theoretic counterparts of Moufang spherical buildings and of affine buildings with a Moufang building at infinity.

\subsection{Strongly transitive actions and Tits systems}\label{subsec:str_tr_and_bn_pairs}

We will now study groups acting on buildings in an interesting fashion. It turns out that the right condition is \emph{strong transitivity}.

\begin{Definition}
    Let $X$ be a building and let $G$ be a group that acts simplicially on $X$. The $G$-action on $X$ is called \emph{strongly transitive} if the action is transitive on pairs of chambers and apartments containing them.
\end{Definition}

\begin{ConstructionN}\label{con:tits_system}
    We fix a standard apartment $\Sigma_0$ and a standard chamber $c_0$ in $\Sigma_0$. Associated to $G$ acting strongly transitively on $X$, there are the following subgroups: The \emph{Borel subgroup $B$} is the stabiliser of the chamber $c_0$. The subgroup $N$ is the setwise stabiliser of the apartment $\Sigma_0$.

    Their intersection $H=B\cap N$ is the pointwise stabiliser of $\Sigma_0$ and the group $H$ is normal in $N$. The quotient $N/H$ acts chamber-regularly on the apartment $\Sigma_0$ and is hence isomorphic to the associated Coxeter group $W$, where we have fixed the generating set $S$.
    \begin{center}
    \begin{tikzpicture}
        \node (G) at (0,2) {$G$};
        \node[anchor=mid] (B) at (-1,1) {$B$};
        \node[anchor=mid] (N) at (1,1) {$N$};
        \node (H) at (0,0) {$H$};
        \node[anchor=mid] (W) at (3,1) {$W=\langle S\rangle$};
        \draw (G) -- (B) -- (H);
        \draw (G) -- (N) -- (H);
        \draw[->] (N) -- ++(1,0);
    \end{tikzpicture}
    \end{center}
    To sum up, from a strongly transitive action of a group $G$ on a building $X$, we obtain a quadruple $(G,B,N,S)$ with the above properties by choosing a standard apartment $\Sigma_0$ and a standard chamber $c_0\in\Sigma_0$. This associated quadruple will play a large role from now on.
\end{ConstructionN}

\begin{Definition}\label{def:parabolic_levi}
    The stabilisers of the subsimplices of the standard chamber are called \emph{standard parabolic subgroups}. Stabilisers of arbitrary simplices are called \emph{parabolic subgroups}.

    If the building $X$ is spherical, then the intersections of parabolic subgroups associated to opposite simplices are called \emph{Levi subgroups}.
\end{Definition}

\begin{Remark}
    If $X$ is thick, parabolic subgroups can also be characterised as all subgroups of $G$ containing a conjugate of $B$.
\end{Remark}

\noindent One of the consequences for a group $G$ that acts strongly transitively on a building is the \emph{Bruhat decomposition} which, after having been discovered by Bruhat long before the definition of buildings, is the starting point for Chevalley's investigation of simple algebraic groups which led to Tits' definition of `Tits systems' or BN pairs.

\begin{Theorem}[Bruhat decomposition, 6.1.5 in \cite{AB:B:08}]\label{th:bruhat}
	Let $G$ be a group acting strongly transitively on a weak building, where we fix an apartment and a chamber. Construct $B$, $N$, $H$ and $W$ as above. Then
	\[
	G = BWB = \coprod_{w\in W} BwB
	\]
	as a disjoint union of double cosets.
\end{Theorem}

\begin{Remark}
    Note that, actually, one has to pick preimages of $W\cong N/H$ in $N$ such that this description makes sense. But since $H\subseteq B$, the decomposition is independent of these choices and hence commonly written in this fashion.
\end{Remark}

\noindent Using the Bruhat decomposition, we obtain the following purely group-theoretic properties of the quadruple $(G,B,N,S)$ we have constructed above.

\begin{Definition}
    Let $G$ be a group. Let $B$ and $N$ be subgroups of $G$ such that their intersection $H=B\cap N$ is normal in $N$. Assume also that the quotient group $W\coloneq N/H$ is generated by a set $S\subset W$. The quadruple $(G,B,N,S)$ is called a \emph{weak Tits system} if
    \begin{itemize}
        \item $G=\langle B \cup N\rangle$,
        \item $(W,S)$ is a Coxeter system and
        \item for $s\in S$ and $w\in W$, we have $BsBwB\subseteq BwB \sqcup BswB$.
    \end{itemize}
    If in addition $sBs\neq B$ for all $s\in S$, we call $(G,B,N,S)$ a \emph{Tits system} for $G$.
\end{Definition}

\noindent It turns out that the concept of a strongly transitive action and of a weak Tits system are actually equivalent.

\begin{Theorem}\label{th:bn_pair_equiv_strongly_transitive}
    If a group $G$ acts strongly transitively on a weak building, the quadruple $(G,B,N,S)$ as constructed in Construction \ref{con:tits_system} is a weak Tits system. If the building is thick, then $(G,B,N,S)$ is a Tits system.

    The converse is also true: If $(G,B,N,S)$ is a (weak) Tits system for $G$, we call a subgroup of the form $\bigsqcup_{w\in W_J} BwB$ a \emph{standard parabolic subgroup}. Then the set of cosets
    \[
        X \coloneq \{ G/P : P \text{ standard parabolic} \}
    \]
    ordered by reverse inclusion is a (weak) building on which $G$ acts strongly transitively. In this building, chambers are cosets of the Borel group $B$. The apartment
    \[
        \Sigma_0\coloneq \{ wB : w\in W\}
    \]
    is stabilised by $N$. The chamber $c_0\coloneq B$ is obviously fixed by $B$.
\end{Theorem}

\begin{Proof}
    Assuming that a group acts strongly transitively on a weak building, we obtain the first property of a Tits system by the Bruhat decomposition in Theorem \ref{th:bruhat}. The second property is fulfilled by construction. The third property is the content of Corollary 6.12 and Theorem 6.21 in \cite{AB:B:08}.

    The second statement of the theorem is Proposition 6.52 in \cite{AB:B:08}. Finally, thickness of the building $X$ is equivalent to $sBs\neq B$ for all $s\in S$ by \cite[Remark 6.26]{AB:B:08}.
\end{Proof}

\begin{Remark}
    If the building is thick, the situation is easier. Then we do not have to assume that $W$ is a Coxeter group and the generating set $S$ can be recovered from the other data. This is why, in this situation, one speaks of the \emph{BN pair $(B,N)$}, which is enough to construct the Tits system $(G,B,N,S)$. Since we often deal with buildings which are not necessarily thick, we will use Tits systems in this thesis.
\end{Remark}

\noindent Theorem \ref{th:bn_pair_equiv_strongly_transitive} is the heart of the strong interplay between buildings and groups. This can be used to state the result of the classification of Moufang spherical buildings.

\paragraph{The classification result} The classification roughly states that irreducible spherical buildings of rank at least three are associated to simple isotropic algebraic groups or to simple isotropic classical groups (with a large overlap between those classes) plus some additional series of buildings which are called \emph{mixed} buildings.

Several examples of classical groups with their associated spherical buildings will be described in Chapter \ref{ch:building_examples}.

\subsection{Groups of Kac-Moody type}\label{subsec:root_data}

The group-theoretic analogue of a Moufang spherical building is a group with a spherical root datum. Instead of defining groups with spherical root data, we will introduce the more general notion of groups with twin root data, also called \emph{groups of Kac-Moody type}. Their geometric counterparts are so-called twin buildings, which will not be defined in this thesis, since this would enlarge it unnecessarily.

For further details on spherical root data see \cite[Section 7]{AB:B:08}. For the more general notion of twin root data, see \cite[Section 8]{AB:B:08} or \cite{Rem:GKM:02}.

Let $\Sigma(W,S)$ be a Coxeter complex and let $\Phi=\Phi(W,S)$ be the set of its roots. As usual, we identify the set of chambers $\Ch(\Sigma(W,S))$ with the elements of $W$ via the regular action. For any $s\in S$, denote by $\alpha_s\in \Phi$ the unique root containing the chamber $1$ but not the chamber $s$. A root containing $1$ is called \emph{positive}.

\begin{Definition}\label{def:root_datum}
	Let $G$ be a group and assume that there is a system of subgroups $(U_\alpha)_{\alpha\in \Phi}$ of $G$ indexed by the set of roots $\Phi$. For all roots $\alpha$, we write $U^*_\alpha\coloneq U_\alpha\setminus\{1\}$. We define
    \[
        U_+ \coloneq \langle U_\alpha : \alpha\text{ positive}\rangle.
    \]
    Let $H\leq G$ be a subgroup that normalises all groups $U_\alpha$. The triple $(G,(U_\alpha)_{\alpha\in\Phi},H)$ is called a \emph{root datum} of type $(W,S)$ if the following hold
	\begin{enumerate}
        \item For all roots $\alpha$, we have $U_\alpha\neq\{1\}$.
        \item For every prenilpotent pair of roots $\{\alpha,\beta\}$ with $\alpha\neq\beta$, we have
            \[
            [U_\alpha,U_\beta] \subseteq U_{(\alpha,\beta)} \coloneq \langle U_\gamma : \gamma\in(\alpha,\beta)\rangle.
            \]
        \item For all $s\in S$ and all $u\in U_{\alpha_s}^*$, there are $u'$ and $u''$ in $U_{-\alpha_s}$ such that $m(u)\coloneq u' u u''$ conjugates $U_\beta$ onto $U_{s\beta}$ for all roots $\beta\in\Phi$.

            In addition, for all $u,v\in U_{\alpha_s}$, we have $m(u)H=m(v)H$.
        \item For all $s\in S$, we have $U_{-\alpha_s}\not\subset U_+$.
        \item $G = H\langle U_{\alpha} : \alpha\in\Phi\rangle$.
	\end{enumerate}
    The subgroups $U_\alpha$ are then called \emph{root groups}. The subgroup $G^\dagger\coloneq\langle U_{\alpha} : \alpha\in\Phi\rangle$ is called the \emph{little projective group}. The root datum is called \emph{spherical} if the Coxeter system $(W,S)$ is spherical.
\end{Definition}

\begin{Remark}
    It can be shown that, under these assumptions, we have $H = \bigcap_{\alpha\in\Phi}N_G(U_\alpha)$, see \cite[1.5.3]{Rem:GKM:02}.
\end{Remark}

\noindent From a root datum, we obtain a Tits system and hence a thick building.

\begin{Proposition}[Proposition 8.54 in \cite{AB:B:08}]\label{prop:root_datum_yields_bn_pair}
    Let $G$ be a group with a root datum $(G,(U_\alpha)_{\alpha\in\Phi},H)$ of type $(W,S)$. Consider the subgroups
    \begin{align*}
        B_+ &= HU_+, \\
        N &= \langle H \cup \{ m(u) : u\in U_{\alpha_s}^*, s\in S \}\rangle.
    \end{align*}
    Then $(B_+ \cap N) = H$, $H$ is normal in $N$ and there is an isomorphism $N/H \rightarrow W$. The quadruple $(G,B_+,N,S)$ forms a Tits system of type $(W,S)$ for $G$.

    In particular, we obtain a thick building $X$ on which the group $G$ acts strongly transitively. The set of chambers of $X$ can be identified with the cosets $G/B_+$. We denote the chamber identified with $B_+$ by $c_0$ and the apartment with chambers $\{wB_+ : w\in W\}$ by $\Sigma_0$.
\end{Proposition}

\noindent If the root datum is spherical, the building is Moufang in the sense of Section \ref{subsec:moufang_spherical_buildings}.

\begin{Proposition}[Theorem 7.116 in \cite{AB:B:08}]\label{prop:root_datum_building_is_moufang}
    Assume that $(W,S)$ is of spherical type. Then the above building $X$ is Moufang. The homomorphism $G\rightarrow \Aut^0(X)$ injects the root groups $U_\alpha$ into the corresponding root groups $U_{\alpha_0}$ coming from the apartment $\Sigma_0$. If the Coxeter diagram of $(W,S)$ has no isolated nodes, the building $X$ is strictly Moufang and the root groups $U_\alpha$ are mapped isomorphically onto the groups $U_{\alpha_0}$.
\end{Proposition}

\noindent It can be shown that Levi subgroups inherit a root datum. Since we have defined Levi subgroups only for spherical buildings, we state the following result only in this case.

\begin{Proposition}[6.2.3 in \cite{Rem:GKM:02}]\label{prop:levi_factors_have_root_datum}
    Assume that $(W,S)$ is of spherical type. Fix a pair of opposite simplices $(s,t)$ in $X$. Remember that the stabiliser of these simplices in $G$ is called a \emph{Levi subgroup}, denote it by $L_{s,t}$.

    The group $L_{s,t}$ admits a spherical root datum, where the root groups are precisely those root groups of $X$ where the associated wall contains the simplices $s$ and $t$. It acts strongly transitively on the links $\lk_X(s)$ and $\lk_X(t)$.
\end{Proposition}

\noindent As stated above, the natural geometric object corresponding to a twin root datum is a twin building, which is also Moufang in an appropriate sense. We only include the one specific property of the positive half of the twin building we require:

\begin{Proposition}[Proposition 8.56 in \cite{AB:B:08}]\label{prop:general_moufang}
    The action of the system of root groups $(U_\alpha)_{\alpha\in\Phi}$ on the building $X$ as in Proposition \ref{prop:root_datum_yields_bn_pair} satisfies the following properties:
    \begin{itemize}
        \item Each root group $U_\alpha$ stabilises the corresponding root in $\Sigma_0$ pointwise.
        \item For each root $\alpha\in\Phi$ and for each panel $p\in\partial\alpha$, the root group $U_\alpha$ acts regularly on the set of chambers which contain $p$ but which are not contained in $\alpha$.
    \end{itemize}
\end{Proposition}

\begin{Definition}
    For any of the root intervals $\Psi \in \{ [\alpha,\beta], (\alpha,\beta], [\alpha,\beta), (\alpha,\beta)\}$ for a prenilpotent pair of roots $\alpha\neq\beta$, we define the \emph{root group interval}
\[
U_{\Psi}\coloneq \langle U_\gamma : \gamma \in \Psi \rangle.
\]
\end{Definition}

\noindent We will need the following explicit description of root group intervals.

\begin{Lemma}\label{lem:root_group_interval_as_cyclic_product}
    Let $\{\alpha,\beta\}$ be a prenilpotent pair of roots. Then there is an ordering of the interval $[\alpha,\beta]=\{\alpha=\alpha_1,\ldots,\alpha_k=\beta\}$ such that
    \[
    U_{[\alpha,\beta]} = U_{\alpha_1}\cdots U_{\alpha_k}.
    \]
\end{Lemma}

\begin{Proof}
    For the proof, we will use the following piece of terminology: A set of roots $\Psi$ is \emph{convex}, if both $\bigcap_{\alpha\in\Psi}\alpha$ as well as $\bigcap_{\alpha\in\Psi}(-\alpha)$ contain a chamber. Closed root intervals are by definition convex.

    By Lemma 8.14 in \cite{AB:B:08}, every convex set of roots admits a so-called \emph{admissible ordering}, that is, an ordering $\{\alpha_1,\ldots,\alpha_k\}$, where each subset $\{\alpha_i,\ldots,\alpha_m\}$ is convex. As one can see from the proof of this lemma, the admissible ordering is constructed as in Lemma \ref{lem:root_intervals_order}. In particular, the first root in the order is $\alpha$, the last one is $\beta$.

    Finally, by Proposition 8.33 in \cite{AB:B:08}, we obtain the description of $U_{[\alpha,\beta]}$.
\end{Proof}

\subsection{Root data with valuations}\label{subsec:root_datum_valuation}

Just as Tits systems encode general buildings and as groups with spherical root data encode Moufang buildings, \emph{root data with valuations} encode affine buildings with Moufang buildings at infinity. Since we mention root data with valuations only at the very end of Chapter \ref{ch:wagoner} and since we will make very little use of the precise structure of these root data, we will omit some technical details of the definition.

\begin{Definition}[3.21 in \cite{Wei:SAB:09}]\label{def:root_datum_valuation}
    Let $(W,S)$ be a spherical Coxeter system whose diagram has no isolated nodes. Let $\Phi=\Phi(W,S)$ be its set of roots and let $G$ be a group with a spherical root datum $(G,(U_\alpha)_{\alpha\in\Phi},H)$. A \emph{valuation} of the root datum is a family of surjective maps $\varphi_\alpha : U_\alpha^* \rightarrow \Z$ for all roots $\alpha\in\Phi$ satisfying
    \begin{enumerate}
        \item For each $\alpha\in\Phi$,
            \[
            U_{\alpha,k} \coloneq \{u\in U_\alpha: \varphi_\alpha(u)\geq k \}
            \]
            is a subgroup of $U_\alpha$ for each $k\in\Z$, where we set $\varphi_\alpha(1)=\infty$, so that $1\in U_{\alpha,k}$ for all $k$.
        \item For all $\alpha,\beta\in\Phi$ such that $\alpha\neq\pm\beta$, we have
            \[
            [U_{\alpha,k},U_{\beta,l}] \subset \prod_{\gamma\in (\alpha,\beta)} U_{\gamma, p_\gamma k + q_\gamma l}
            \]
            for all $k,l\in\R$. The coefficients $p_\gamma $ and $q_\gamma$ depend on a \emph{root map} which realises the set of roots $\Phi$ as a \emph{root system}, a set of vectors in a real vector space. The product is taken in the natural ordering of the roots. For details, see \cite[Chapter 3]{Wei:SAB:09}.
        \item For all $\alpha,\beta\in\Phi$ and for all $u\in U_{\alpha}^*$, there exists $t\in \Z$ such that
            \[
            \varphi_{s_{\alpha}(\beta)}(m(u)^{-1}xm(u)) = \varphi_{\beta}(x) + t
            \]
            for all $x\in U_\beta^*$, where $s_\alpha$ is the unique reflection interchanging $\alpha$ and $-\alpha$.
        \item If $\alpha=\beta$ in 3., we have $t=-2\varphi_\alpha(u)$.
    \end{enumerate}
\end{Definition}

\noindent A root datum with valuation leads to an affine building with a Moufang building at infinity.

\begin{Theorem}[Theorems 13.30 and 14.47 in \cite{Wei:SAB:09}]\label{th:root_data_with_valuations_have_affine_buildings}
    Let $G$ be a group with a root datum with valuation $(G,(U_\alpha)_{\alpha\in\Phi},(\varphi_\alpha)_{\alpha\in\Phi},H)$. Then there is an affine building $X$ on which the group $G$ acts strongly transitively. In addition, the spherical building at infinity of $X$ is isomorphic to the Moufang spherical building associated to the root datum.

    There is an affine apartment $\Sigma\subset X$ and a special vertex $x\in \Sigma$ such that the half-apartments of $\Sigma$ can be parametrised as $H_{\alpha,k}$, where $\alpha\in\Phi$ is a spherical root and $k\in \Z$, such that $x\in \partial H_{\alpha,0}$ and $H_{\alpha,k} \subset H_{\alpha,l}$ for $k<l$. The subgroups $U_{\alpha,k}$ then stabilise the half-apartments $H_{\alpha,k}$.

    Conversely, every root datum with valuation arises in this way.
\end{Theorem}

\paragraph{Classification} Affine buildings of rank at least three with a Moufang building at infinity have been classified by Bruhat and Tits, see Tits' article \cite{Tit:RLF:79} or Weiss' book \cite{Wei:SAB:09}. The additional condition of the building at infinity being Moufang is void for the case of rank at least four by Theorem \ref{th:higher_rank_sph_buildings_are_moufang}. Using this spherical building, a group with a spherical root datum as in Section \ref{subsec:root_data} is constructed.

The spherical root groups $U_\alpha$ act on the affine building, each element is contained in a root group corresponding to a half-apartment in a parallel class associated to the root $\alpha$. By ordering the elements of the parallel classes by inclusion and by fixing an origin, one obtains an enumeration of each parallel class of roots. The valuation evaluated at a root group element is simply the number of the half-apartment the element stabilises.

It turns out that spherical root data with valuations can be classified algebraically which leads to a complete classification of affine buildings. Roughly, affine buildings come from algebraic groups or classical groups, which correspond to the Moufang building at infinity, over a field which admits a \emph{discrete valuation}. This valuation induces a valuation of the root datum, which consequently leads to an affine building.

\begin{Remark}
    The Moufang condition at infinity is necessary in rank three. There is a free construction by Ronan in \cite{Ron:CBR:86} constructing rank three buildings whose rank two residues can be arbitrary generalised polygons. Since these cannot be classified, there is no hope of classifying affine buildings of rank three whose building at infinity is not Moufang.
\end{Remark}

\chapter{Examples of buildings and classical groups}\label{ch:building_examples}

In this section, we give a variety of examples of spherical and affine buildings and their associated groups which we will use throughout this thesis.

The smallest interesting examples of spherical buildings are rank two buildings: generalised polygons, which will be described in the first section.

The next two sections introduce those buildings and classical groups which are the main examples in this thesis. We will apply our results in Chapter \ref{ch:stability} to these examples to obtain concrete homological stability results.

In the last section, we give a quick overview of the classical two-dimensional affine buildings and known constructions of exotic two-dimensional affine buildings.

\section{Generalised polygons}\label{sec:generalised_polygons}

Buildings were originally defined in terms of incidence geometries. For rank two buildings, this point of view leads to classical point-line incidence geometries such as projective planes. This description is more common and natural for rank two buildings.

\begin{Definition}
    An \emph{incidence structure} $\cI$ is a triple of sets $(\cP,\caL,\cF)$, where $\cP\cap\caL=\emptyset$ and $\cF\subseteq \cP \times \caL$. The elements of $\cP$ are called \emph{points}, the elements of $\caL$ are called \emph{lines}. Elements of $\cF$ are called \emph{flags}. If $(p,l)\in\cF$ for a point $p$ and a line $l$, we say that $p$ and $l$ are \emph{incident}.

    Associated to every incidence structure is its \emph{incidence graph} with vertex set $\cP\sqcup\caL$ and with edge set $\cF$.
\end{Definition}

\noindent In the previous section, we have defined generalised polygons to be buildings of type $I_2(m)$. There is a canonical incidence structure associated to such buildings, also called a generalised polygon.

\begin{Definition}
	A \emph{generalised $m$-gon} is an incidence structure $\cI=(\cP,\caL,\cF)$ whose associated incidence graph has diameter $m$ and girth $2m$. Remember that the \emph{girth} of a graph is the length of the shortest non-trivial cycle.
\end{Definition}

\noindent It is not hard to see that the Coxeter complex of type $I_2(m)$ is a $2m$-cycle, compare Example \ref{ex:i2_complex}. Using this fact, the correspondence between the generalised polygon as an incidence geometry and as a building can be established fairly easily.

\begin{Proposition}[7.14 and 7.15 in \cite{Wei:SSB:03}]
    Consider a generalised $m$-gon in the sense of the last definition. Its associated incidence graph, considered as a one-dimensional simplicial complex, is a spherical building of type $I_2(m)$. Conversely, every incidence graph of a generalised $m$-gon arises in this way.
\end{Proposition}

\noindent The simplest examples for generalised polygons are complete bipartite graphs.

\begin{Example}
	It is easy to verify that every generalised $2$-gon is a complete bipartite graph and vice versa.
\end{Example}

\begin{Definition}
    A thick generalised $3$-gon or triangle is called a \emph{projective plane}. A projective plane can equivalently be characterised as an incidence structure satisfying the following axioms.
    \begin{itemize}
        \item Every two distinct lines intersect in a unique point.
        \item There is a unique line through every two distinct points.
        \item There are four points such that no line is incident with more than two of them. (Existence of a quadrangle)
    \end{itemize}
\end{Definition}

\noindent For every division ring there is a so-called \emph{classical} projective plane. Its construction is the special case $n=2$ of the construction of classical spherical buildings of type $A_n$, which can be found in the next section.

We will only require a few results concerning projective planes.

\begin{Lemma}[3.2.1 in \cite{De:FG:68}]
    For a finite projective plane, the number of points on a line is constant throughout the plane. In addition, it is equal to the number of lines through a point, which is also constant.

    If there are $q+1$ points on a line, we say that the projective plane has \emph{order $q$}. In buildings terminology, this means that every panel is contained in $q+1$ chambers.

    A projective plane of order $q$ has $q^2+q+1$ points and $q^2+q+1$ lines.
\end{Lemma}

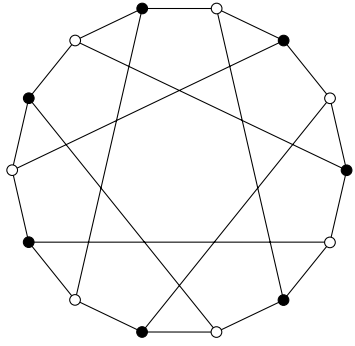
\begin{figure}
    \centering
            \begin{tikzpicture}[scale=1.1]
                \tikzstyle{every node}=[circle, draw, fill=black!50,
                        inner sep=0pt, minimum width=4pt]
                \draw \foreach \x in {0,51.42,...,309.68}
                {
                (\x:2) node[circle,draw,fill=black] {} -- (\x+25.71:2)
                (\x+25.71:2) node[circle,draw,fill=white] {} -- (\x+51.42:2)
                (\x:2) -- (\x+128.65:2)
                };
            \end{tikzpicture}\caption{The Fano plane}
    \label{fig:fano_plane}
\end{figure}

\begin{Example}
    The incidence graph of the smallest projective plane, the \emph{Fano plane} can be found in Figure \ref{fig:fano_plane}. The Fano plane has order $2$. Points and lines are indicated with black and white vertices, respectively.
\end{Example}

\begin{Remark}
    A generalised $4$-gon is usually called a \emph{generalised quadrangle}. For quadrangles, the number of points per line and the number of lines through a point must each be constant, but these are not necessarily equal. A generalised quadrangle is said to be \emph{of order $(s,t)$} if there are $s+1$ points on a line and $t+1$ lines through a point. Explicit examples will be constructed in Section \ref{sec:ex_buildings_of_type_bn_cn} and also in Section \ref{subsec:GQs}.
\end{Remark}

\noindent Similar to the corresponding result by Tits and Weiss for Moufang polygons, the following famous theorem by Feit and Higman shows that thick finite generalised polygons can only exist for very few parameters. This shows in particular that there cannot be finite thick buildings of types $H_3$ and $H_4$.

\begin{Theorem}[Feit-Higman, \cite{FH:NEP:64}]\label{th:feit_higman}
    Thick finite generalised $m$-gons only exist for $m\in\{2,3,4,6,8\}$. In addition, there are explicit bounds on the thickness of panels in the various special cases.
\end{Theorem}

\noindent There is a free construction by Tits that produces infinite generalised $m$-gons for any $m$.

\section{The spherical buildings of type \texorpdfstring{$A_n$}{An} and linear groups}\label{sec:an_example}

Let $D$ be a division ring, and let $V$ be a right vector space of dimension $n\geq 2$ over $D$.

\begin{Definition}
    The \emph{projective space $P$} over $V$ is the set of all non-trivial proper subspaces of $V$, partially ordered by inclusion. A \emph{frame in $V$} is a set of one-dimensional subspaces $\cF=\{L_1,\ldots,L_n\}$ such that $V=L_1\oplus\cdots\oplus L_n$.
\end{Definition}

\noindent We now construct a simplicial complex associated to the projective space.

\begin{Definition}
    Denote the flag complex over $P$ by $X \coloneq \Flag(P)$. For any given frame $\cF$, the subcomplex consisting of flags of subspaces spanned by elements of $\cF$ is denoted by $\Sigma(\cF)$.
\end{Definition}

\begin{Theorem}
    The complex $X$ is a thick building of type $A_{n-1}$. Its apartments are the subcomplexes $\Sigma(\cF)$ for all frames $\cF$. A type function on vertices is given by the dimension function:
    \[
        \type( (V')) = \dim(V'),
    \]
    where $0\lneq V' \lneq V$ and the type set $I$ is given by $I=\{1,\ldots, n-1\}$.
\end{Theorem}

\begin{Proof}
    A sketch of a proof where $D$ is a field can be found in \cite[4.3]{AB:B:08}. Instead of giving a complete proof, we consider the action of $\Gl_n(D)$ on $X$. We shall see in Theorem \ref{th:an_tits_system} that the corresponding stabilisers form a Tits system for $G$, which implies that $X$ is isomorphic to the building corresponding to the Tits system. This shows in particular that $X$ is a building.
\end{Proof}

\begin{Remark}
    In the special case of $n=3$, we obtain a building of type $A_2$ which is a generalised triangle in the sense of Section \ref{sec:generalised_polygons}. This is of course just the definition of a classical projective plane over $D$.
\end{Remark}

\noindent Later on, it will be important for us to characterise opposite vertices.

\begin{Lemma}\label{lem:an_opposition}
    Two vertices $(V')$ and $(V'')$ in $X$ are opposite if and only if
    \[
        V'\bigoplus V''= V.
    \]
\end{Lemma}

\begin{Proof}
    Let $(V')$ and $(V'')$ be opposite vertices. From the above description of apartments, we see that there is a frame $\cF$ such that the vertices are both contained in the corresponding apartment: $(V'),(V'')\in\Sigma(\cF)$. Obviously, we can choose an ordered basis $(e_i)_{i\in\{1,\ldots,n\}}$ such that $\cF =\{\langle e_1\rangle,\ldots,\langle e_n\rangle\}$ and $V'=\langle e_1,\ldots,e_k\rangle$ for some $k\in I$.

    The Coxeter group of type $A_{n-1}$ acts on the apartment $\Sigma(\cF)$, and we have already seen that it is isomorphic to the symmetric group on $n$ letters acting on the basis $(e_i)_i$. The apartment $\Sigma(\cF)$ is canonically isomorphic to the flag complex over non-trivial proper subsets of $\{1,\ldots,n\}$ which we have discussed in Example \ref{ex:an_coxeter_complex}.

    In Example \ref{ex:an_opposition}, we have seen that the vertex opposite of $\{1,\ldots,k\}$ is $\{k+1,\ldots,n\}$, which shows that $V''=\langle e_{k+1},\ldots,e_n\rangle$. So we have $V'\oplus V''=V$.

    For the other direction, find a suitable basis and use the corresponding apartment.
\end{Proof}

\noindent Now fix a basis $e_1,\ldots,e_n$ of $V$ and consider the corresponding action of $G\coloneq\Gl_n(D)$ on $V$, which induces an order-preserving action on $P$ and hence a simplicial action on $X$. We fix the chamber
\[
    c_0 = ( \langle e_1\rangle \subset \langle e_1,e_2\rangle\subset\cdots\subset\langle e_1,\ldots,e_{n-1}\rangle)
\]
and the frame
\[
    \cF_0 = \{\langle e_1\rangle,\ldots,\langle e_n\rangle\},
\]
which leads to an apartment $\Sigma_0=\Sigma(\cF_0)$.

\begin{Definition}
    Denote by $B\leq G$ the group of upper triangular matrices and by $N\leq G$ the group of monomial matrices, i.e.~matrices with exactly one non-zero entry in every row and column.

    The group $N$ acts as a group of permutations on the frame $\cF_0$ and we obtain a surjection of $N$ onto the symmetric group on $n$ letters with kernel $H=N\cap B$, which is the subgroup of diagonal matrices. Hence $W\coloneq N/H$ can be identified with the symmetric group on $n$ letters and we fix $S$ to be the set of standard generators $S=\{s_i=(i,i+1): i\in I\}$.
\end{Definition}

\noindent Obviously, $B$ is the stabiliser of $c_0$ and $N$ is the setwise stabiliser of $\Sigma_0$.

\begin{Theorem}[Section 6.5 in \cite{AB:B:08}]\label{th:an_tits_system}
    The quadruple $(G,B,N,S)$ is a Tits system of type $A_{n-1}$ for $G=\Gl_n(D)$. The associated building is isomorphic to $X$.
\end{Theorem}

\begin{Remark}
    In fact, the group $\Gl_n(D)$ admits a root datum of type $A_{n-1}$. It can be shown that the roots in the Coxeter complex of type $A_{n-1}$ can be indexed by  $\{(i,j) : 1\leq i\neq j \leq n\}$, such that the positive roots are exactly indexed by those pairs $(i,j)$ where $i<j$. The corresponding root groups $U_{i,j}$ have the form
    \[
    U_{i,j} = \{ \one_n + \varepsilon_{i,j}(d) : d\in D \},
    \]
    where $(\varepsilon_{i,j}(d))_{i,j}=d$ and $(\varepsilon_{i,j}(d))_{k,l}=0$ for $i\neq k$ or $j\neq l$.

    In particular, the corresponding little projective group is the group generated by elementary matrices, denoted by $\E_n(D)$. It is a perfect group by \cite[Proposition 2.1.4]{Ros:AKa:94}. The group $U_+$ is the subgroup of all upper triangular matrices with $1$ on the diagonal.
\end{Remark}

\noindent We make the standard parabolic and Levi subgroups explicit, since we will use this later.

\begin{Remark}
    For the vertex $(\langle e_1,\ldots, e_k\rangle)$ of the chamber $c_0$, the corresponding parabolic subgroup is given by
    \[
    P_k \coloneq\Biggl\{ \begin{pmatrix}
        A & B \\ 0 & C
    \end{pmatrix} : A \in \Gl_{k}(D), C\in \Gl_{n-k}(D), B\text{ arbitrary} \Biggr\}.
    \]
    By Lemma \ref{lem:an_opposition}, the opposite vertex in the apartment $\Sigma_0$ is $(\langle e_{k+1},\ldots,e_n\rangle)$. The corresponding Levi subgroup is given by
    \[
    L_k \coloneq\Biggl\{ \begin{pmatrix}
        A & 0 \\ 0 & C
    \end{pmatrix} : A \in \Gl_{k}(D), C\in \Gl_{n-k}(D) \Biggr\}.
    \]
\end{Remark}

\noindent If $D$ is a field, there is a determinant function $\Gl_n(D)\rightarrow D^\times$. Its kernel is the special linear group $\Sl_n(D)$. Special linear groups can also be defined for division rings using the Dieudonné determinant, but the case of fields is enough for our purposes. We also have $\Sl_n(D)=\E_n(D)$ by \cite[2.2.6]{HoM:CGK:89}.

\begin{Proposition}
    The special linear group $\Sl_n(D)$ also acts strongly transitively on the building $X$. Consequently, the quadruple $(\Sl_n(D), B \cap \Sl_n(D), N\cap\Sl_n(D),S)$ is a Tits system for $\Sl_n(D)$.
\end{Proposition}

\begin{Proof}
    Note that strong transitivity on the building is equivalent to saying that every pair of an apartment $\Sigma$ and a chamber $c\in\Sigma$ can be mapped to the apartment $\Sigma_0$ and to the chamber $c_0$. If there is an element $g\in\Gl_n(D)$ with this property, then the diagonal matrix $h$ with diagonal $(\det(g)^{-1},1,\ldots,1)$ is contained in $H$ and fixes hence both $\Sigma_0$ and $c_0$. In particular, $hg$ still maps the pair $(\Sigma,c)$ to $(\Sigma_0,c_0)$ and $hg\in\Sl_n(D)$, which proves strong transitivity. The structure of the Tits system follows from Theorem \ref{th:bn_pair_equiv_strongly_transitive}.
\end{Proof}

\noindent Again, the structure of the standard parabolic subgroups and of the Levi subgroups can easily be obtained.

\begin{Remark}
    If $(\langle e_1,\ldots, e_k\rangle)$ is a vertex of the chamber $c_0$, the corresponding parabolic subgroup in $\Sl_n(D)$ is given by
    \[
    P_k \coloneq\Biggl\{ \begin{pmatrix}
        A & B \\ 0 & C
    \end{pmatrix} : \begin{array}{l} A \in \Gl_{k}(D), C\in \Gl_{n-k}(D),\\ \det(A)\det(C)=1, B\text{ arbitrary}\end{array} \Biggr\}.
    \]
    Again, the opposite vertex in the apartment $\Sigma_0$ is $(\langle e_{k+1},\ldots,e_n\rangle)$. The corresponding Levi subgroup is given by
    \[
    L_k \coloneq\Biggl\{ \begin{pmatrix}
        A & 0 \\ 0 & C
    \end{pmatrix} : A \in \Gl_{k}(D), C\in \Gl_{n-k}(D), \det(A)\det(C)=1 \Biggr\}.
    \]
\end{Remark}

\noindent It turns out that all buildings of type $A_n$ for $n\geq 3$ are isomorphic to one of the buildings we have constructed in this section. We will not make use of this classification, but we state the following result, which is the first step in Tits' famous classification of higher-rank spherical buildings.

\begin{Theorem}
    Every building of type $A_n$ for $n\geq 3$ arises as the flag complex over a projective space over a division ring. Buildings of type $A_2$ are projective planes of which there are many exotic ones. If we additionally assume that the $A_2$-building is Moufang, then it arises as the flag complex over a projective plane over a division ring or a Cayley division algebra.
\end{Theorem}

\begin{Proof}
    For $n\geq 3$, the main result is Theorem 6.3 in \cite{Tit:BsT:74}, which shows that every building of type $A_n$ is a projective space in the sense of incidence geometries. The Moufang condition translates into the Desarguesian condition for the projective space. The fundamental theorem of projective geometry, see for instance \cite{Lun:BHG:65}, then shows that every such projective space is a projective space over a division ring.

    The case of Moufang triangles, that is $n=2$, is discussed in Theorems 17.2 and 17.3 in \cite{TW:MP:02}.
\end{Proof}

\noindent For a survey on the general linear group and the associated building from the point of view of incidence geometry, see \cite[Part 1]{Kra:BCG:00}.

\section{The spherical buildings of type \texorpdfstring{$C_n$}{Cn} and unitary groups}\label{sec:ex_buildings_of_type_bn_cn}

Let $D$ be a division ring. In this section, we will define the building associated to an $(\varepsilon,J)$-hermitian form on a right vector space over $D$. A very detailed overview of such forms and their theory can be found in \cite[Chapter 6]{HoM:CGK:89}. We refer to this book for further details.

\begin{Definition}
    An \emph{involution $J$ of $D$} is an additive map $J: D\rightarrow D$, written as $a\mapsto a^J$ satisfying $(ab)^J = b^Ja^J$ and $J^2=\id_D$.
\end{Definition}

\noindent Fix an involution $J:D\rightarrow D$ and fix $\varepsilon\in\{\pm 1\}$. Let $V$ be a $d$-dimensional right vector space over $D$.

\begin{Definition}
	An $(\varepsilon,J)$-hermitian form is a biadditive map $h:V\times V\rightarrow D$ that satisfies
\begin{align*}
    h(va,w) &= a^Jh(v,w),& h(v,wa)&=h(v,w)a,& h(w,v)&= h(v,w)^J\varepsilon.
\end{align*}
    The form is said to be \emph{non-degenerate} if there is no $v\in V\setminus\{0\}$ such that $h(v,w)=0$ for all $w\in V$.
\end{Definition}

\noindent For the rest of this section, let $h$ be a non-degenerate $(\varepsilon,J)$-hermitian form. If the characteristic of $D$ is two, we assume additionally that $h$ is \emph{trace-valued}, that is:
\[
    h(v,v) \in \{ a + a^J\varepsilon : a\in D\}.
\]
If the characteristic of $D$ is not two, then every hermitian form is automatically trace-valued by \cite[5.1.10]{HoM:CGK:89}.

\begin{Definition}
	A subspace $E\leq V$ is \emph{totally isotropic} if $h(v,w)=0$ for all $v,w\in E$. At the other extreme, a subspace $F\leq V$ is said to be \emph{anisotropic} if $h(v,v)\neq 0$ for all $v\in F\setminus\{0\}$. Denote by $\ind(h)$ the \emph{Witt index of $h$}, the maximal dimension of totally isotropic subspaces of $V$.
\end{Definition}

\noindent We will collect some basic results on these forms in the following proposition.

\begin{BreakProposition}[6.1.12 and 6.2.13 in \cite{HoM:CGK:89}]\label{prop:bases_for_spaces_with_hermitian_form}
    All maximal totally isotropic subspaces have dimension $n=\ind(h)$ and we have
    \[
        2\ind(h)\leq d=\dim(V).
    \]
    There is a basis of $V$ of the form
\[
  e_{-n},\ldots,e_{-1},e_1,\ldots,e_n,f_1,\ldots,f_r
\]
  with $2n+r=d$, such that
\[
  h(e_i,e_j) =\begin{cases}
    \varepsilon & \text{if } i+j=0 \text{ and } i>0 \\
    1 & \text{if } i+j=0\text{ and } i<0 \\
    0 & \text{otherwise,}
  \end{cases}
\]
and that $V=\caH\oplus E$, where
\[
\caH=\langle e_{-n},\ldots,e_{-1},e_1,\ldots,e_n\rangle,\qquad E=\langle f_1,\ldots,f_r\rangle
\]
with $\caH\perp E$ and where $E$ is anisotropic. The subspace $\caH$ is called a \emph{hyperbolic module}.
\end{BreakProposition}

\noindent Again, we define an associated partially ordered set and appropriate subsets.

\begin{Definition}
    The set $P$ of all totally isotropic subspaces of $V$, partially ordered by inclusion is called the \emph{polar space associated to $h$}. For a basis as in Proposition \ref{prop:bases_for_spaces_with_hermitian_form}, the set of totally isotropic subspaces $\cF=\{\langle e_{-n}\rangle,\ldots \langle e_{-1}\rangle, \langle e_1 \rangle,\ldots \langle e_n \rangle\}$ is called a \emph{frame}.
\end{Definition}

\noindent As before, there is a simplicial complex associated to the polar space.

\begin{Definition}
    Denote the flag complex over the polar space $P$ by $X\coloneq \Flag(P)$. For a frame $\cF$ as above, denote by $\Sigma(\cF)$ the subcomplex consisting of all flags of all totally isotropic subspaces spanned by subsets of $\cF$.
\end{Definition}

\begin{Theorem}
    The complex $X$ is a spherical building of type $C_n$. The apartments are the subcomplexes $\Sigma(\cF)$ for all frames $\cF$. The building is thick except if $J=\id$, $\varepsilon\neq -1$ and $\dim(V)=2\ind(h)$.
\end{Theorem}

\begin{Proof}
    This is the content of \cite[7.4]{Tit:BsT:74}, using that the polar space $P$ gives rise to an incidence geometry also called a polar space. Instead of defining such incidence geometries, it is easier to find a group acting on $X$ and showing that it has a Tits system which comes from the right stabilisers, as above.
\end{Proof}

\begin{Remark}
Again, for $n=2$, this is the definition of a classical generalised quadrangle as in Section \ref{sec:generalised_polygons}.
\end{Remark}

\noindent It will be helpful to know what opposition means for vertices.

\begin{Lemma}\label{lem:cn_opposition}
    Two simplices $(V')$ and $(V'')$ in $X$ are opposite if and only if
    \[
        V'\oplus (V'')^\perp = V.
    \]
\end{Lemma}

\begin{Proof}
    This proof is analogous to the proof of Lemma \ref{lem:an_opposition}. Assume that $(V')$ and $(V'')$ are opposite vertices in an apartment $\Sigma$. After re-enumeration, we find a basis $\{e_{-n},\ldots,e_{-1},e_1,\ldots,e_n\}$ as above such that $V'=\langle e_1,\ldots,e_k\rangle$ for some $k$ and such that $\Sigma$ is the apartment associated to this basis. We have seen in Example \ref{ex:cn_coxeter_complex} that the Coxeter complex can be identified with the boundary of a hyperoctahedron or, respectively, the flag complex over `signed' subsets of $\{1,\ldots,n\}$. By the description of opposition in Example \ref{ex:cn_opposition}, we see then that $V''=\langle e_{-1},\ldots, e_{-k}\rangle$.

    For the reverse direction, assume that $V'\oplus (V'')^\perp = V$. Then $V'$ is contained in a totally isotropic subspace $Z'$ and obviously $Z' \cap V''=\{0\}$. Now we can iteratively choose a basis $e_1,\ldots,e_n$ of $Z'$ and associated vectors $e_{-1},\ldots,e_{-n}$ such that $V''=\langle e_{-1},\ldots,e_{-k}\rangle$ and such that $\{e_{-n},\ldots,e_{-1},e_1,\ldots,e_n\}$ is a frame. By the description in \ref{ex:cn_opposition}, the vertices are then opposite.
\end{Proof}

\begin{Definition}
  The \emph{unitary group associated to $h$} is
  \[
  \U(V) = \{ A\in\Gl(V) : h(Av,Aw) = h(v,w) \quad\forall v,w\in V\},
  \]
  the subgroup of $h$-preserving linear automorphisms of $V$.
\end{Definition}

\begin{Remark}
    For $J=\id$, which forces $D$ to be a field, we obtain two important special cases. If $\varepsilon=1$ and $E=\{0\}$, the group $\U(V)$ is usually denoted by $O_{n,n}(D)$ and is called the associated \emph{orthogonal group}.

    If $\varepsilon=-1$ and $E=\{0\}$, the group $\U(V)$ is a \emph{symplectic group} denoted by $\Sp_{2n}(D)$.
\end{Remark}

\noindent Fix a basis $\{e_{-n},\ldots,e_{-1},e_1,\ldots,e_n,f_1,\ldots,f_r\}$ as above. Consider the chamber
\[
c_0 = ( \langle e_{-n} \rangle \subset \langle e_{-n},e_{-(n-1)} \rangle \subset \cdots \subset \langle e_{-n},\ldots,e_{-1} \rangle)
\]
and the associated frame
\[
    \cF_0=\{\langle e_{-n}\rangle,\ldots, \langle e_{-1}\rangle, \langle e_1 \rangle,\ldots, \langle e_n \rangle\}
\]
which leads to the apartment $\Sigma_0=\Sigma(\cF_0)$.

\begin{Construction}
    Consider the following two subgroups of $G=U(V)$: The group
    \[
        B = \Biggl\{\begin{pmatrix}
            A & T & 0 \\
            0 & A^{-J} & 0 \\
            0 & 0 & C
        \end{pmatrix} : A \text{ upper triangular}, AT^J+TA^J\varepsilon = 0, C \in \U(E) \Biggr\}
    \]
    is precisely the stabiliser of $c_0$, which requires some computation. The subgroup $N$ is the stabiliser of the frame $\cF_0$, it has the form
    \[
        N = \Biggl\{\begin{pmatrix}
            M & 0 \\
            0 & C
        \end{pmatrix} : M \text{ a monomial matrix respecting the form $h$}, C\in \U(E) \Biggr\}.
    \]
    Their intersection $H=B\cap N$ has the form
    \[
        H = \Biggl\{\begin{pmatrix}
            \Delta & 0 \\
            0 & C
        \end{pmatrix} : \Delta \text{ diagonal}, C\in \U(E) \Biggr\},
    \]
    which is obviously normal in $N$.

    The quotient subgroup $W\coloneq N/H$ is isomorphic to the group of permutations of the frame $\cF$, respecting the form. This essentially means that if $w(\langle e_i \rangle) = \langle e_j \rangle$, then automatically $w(\langle e_{-i}\rangle) = \langle e_{-j} \rangle$. The group $W$ can hence be identified with the group of signed permutations on $n$ letters, which we have seen to be a Coxeter group of type $C_n$ in Example \ref{ex:cn_coxeter_group}.
\end{Construction}

\begin{Theorem}
    The quadruple $(G,B,N,S)$ is a weak Tits system of type $C_n$ for the unitary group $U(V)$. It is a Tits system except for $J=\id$, $\varepsilon\neq -1$ and $\dim(V)=2\ind(h)$.
\end{Theorem}

\begin{Proof}
    This can be done by explicit calculations. Proof sketches for special cases can be found in \cite[6.6-6.8]{AB:B:08} and \cite[chapter 10]{Gar:BCG:97}.
\end{Proof}

\noindent Again, we make the structure of parabolic and Levi subgroups explicit.

\begin{Remark}
    For each $1\leq k\leq n$, we write $\caH_k = \langle e_{-k},\ldots,e_{-1},e_1,\ldots,e_k \rangle$. For the vertex $(\langle e_{-n},\ldots,e_{-k}\rangle)$ of the chamber $c_0$, the corresponding parabolic subgroup is given by
 \[
P_k =\biggl\{\begin{pmatrix}
    S & * & * & 0 \\
    0 & A & * & B \\
    0 & 0 & S^{-J} & 0 \\
    0 & C & 0 & D
  \end{pmatrix} : S\in \Gl_{n+1-k}(D), \begin{pmatrix}
    A & B \\
    C & D
  \end{pmatrix}\in \U(\caH_{k-1}\perp E)\biggr\}
\]
with some additional conditions on the entries marked $*$.

By Lemma \ref{lem:cn_opposition}, the opposite vertex in $\Sigma_0$ is given by $(\langle e_k,\ldots,e_n\rangle)$. The corresponding Levi subgroup is given by
\[
L_k =\biggl\{\begin{pmatrix}
    S & 0 & 0 & 0 \\
    0 & A & 0 & B \\
    0 & 0 & S^{-J} & 0 \\
    0 & C & 0 & D
  \end{pmatrix} : S\in \Gl_{n+1-k}(D), \begin{pmatrix}
    A & B \\
    C & D
  \end{pmatrix}\in \U(\caH_{k-1}\perp E)\biggr\}.
\]
\end{Remark}

\noindent This class of buildings of type $C_n$ already covers a large part of the class of all buildings of type $C_n$.

\begin{Theorem}
    Every building of type $C_n$ for $n\geq 3$ is a polar space that either comes from an $(\varepsilon,J)$-hermitian form as above or from a non-degenerate pseudo-quadratic form, which yields additional examples only in characteristic two, or belongs to an exceptional series of buildings of type $C_3$ associated to Cayley division algebras.
\end{Theorem}

\begin{Proof}
    By Theorem 7.4 in \cite{Tit:BsT:74}, every weak building of $C_n$ corresponds to an incidence geometry also called a polar space. By Theorem 8.22 in \cite{Tit:BsT:74} due to Veldkamp and Tits, we obtain the stated description.
\end{Proof}

\noindent See \cite[Part 2]{Kra:BCG:00} for a discussion of polar spaces and of the associated classical groups from the point of view of incidence geometry.

\section{Two-dimensional affine buildings}

The classification of affine buildings (due to Bruhat and Tits in a series of papers, see the book by Weiss \cite{Wei:SAB:09}) asserts that all affine buildings in dimensions three and higher come from algebraic data such as an algebraic group over a local field.
In dimension two, there is a free construction of other, exotic buildings by Ronan in \cite{Ron:CBR:86} which is very general but provides no control over the automorphism group. This construction implies in particular that there cannot be any hope to classify these buildings, because any classification would have to include the classification of all projective planes and generalised quadrangles.

\paragraph{Type $\tilde A_2$:} In a locally finite affine building of type $\tilde A_2$, the number $q+1$ of chambers containing a common panel is constant throughout the building. For every prime power $q=p^e$, there are two classical buildings associated to the special linear groups over either a finite extension of the $p$-adic numbers $\Q_p$ or over the Laurent series field $\F_q(\!(t)\!)$. We will not construct these buildings here, since we do not really use them. A construction can be found in \cite[6.9]{AB:B:08}.

There are many more so-called $\emph{exotic buildings}$ of type $\tilde A_2$, some with very large automorphism groups. In \cite{HvM:NTB:90}, for example, van Maldeghem constructs exotic buildings of type $\tilde A_2$ with a vertex-transitive automorphism group.

\paragraph{Type $\tilde C_2$:} The situation is similar in the case of buildings of type $\tilde C_2$. Here, there are more different types of classical buildings, but again there are uncountably many exotic buildings of type $\tilde C_2$. In \cite{Kan:GPS:86}, Kantor constructs exotic buildings of type $\tilde C_2$ with a large automorphism group.

\chapter{Group homology}\label{ch:group}

Group homology is, as the name indicates, a homology theory on the category of groups. It can be defined algebraically using projective resolutions, but also topologically using classifying spaces. We will focus on the first approach and we will only sketch the second. Group homology has many applications, both in group theory and in topology. For us, the concept of homological stability and its connection to $K$-theory will be particularly important.

In this chapter, we first give a quick reminder of all the homology theory required for this thesis. After that, we define group homology along with interpretations of low-dimensional homology groups, with topological interpretations and all the results we require. Then, we introduce the concept of homological stability in which we shall be interested in Chapter \ref{ch:stability}. Finally, we give a quick overview of the connection between homological stability of classical groups and $K$-theory.

The standard reference for group homology is the book by Brown, see \cite{Bro:CoG:82}. The reader might also want to look at the book by Evens, see \cite{Eve:CoG:91}. The book by Knudson \cite{Knu:HLG:01} is dealing with the homology of linear groups and homological stability. Finally, the book by Weibel \cite{Wei:IHA:94} will be our main reference for homological algebra. It also contains a chapter on group homology.

\section{Homological algebra}

We fix a ring $R$. Throughout this section, all modules are left modules over $R$.

\begin{Definition}
    A \emph{chain complex} $(C_*,\partial_*)$ is a sequence of modules $(C_k)_{k\in\Z}$ with a sequence of homomorphisms $\partial_k : C_k \rightarrow C_{k-1}$ for all $k\in \Z$, which are called \emph{boundary operators} or \emph{differentials}, such that $\partial_{k-1}\circ \partial_k=0$ for all $k\in \Z$.

	Associated to every chain complex are \emph{homology modules}
	\[
	H_k(C) = \ker(\partial_k) / \im(\partial_{k+1}).
	\]
	A chain complex is called \emph{exact} if $H_k(C)=0$ for all $k\in \Z$.
\end{Definition}

\noindent We shall often omit the differentials from the notation by simply denoting a chain complex by $C$ and by referring to its differential later as $\partial^C$. In this case, we will often also omit the indices, as in the following definition.

\begin{Definition}
	A \emph{chain map} $f:C' \rightarrow C$ is a sequence of maps $f_k: C'_k\rightarrow C_k$ such that $\partial^C\circ f = f\circ \partial^{C'}$.
\end{Definition}

\noindent Homology groups are functorial: Every chain map induces homomorphisms on homology, denoted by
\[
    H_k(f) : H_k(C') \rightarrow H_k(C).
\]

\begin{Lemma}[5-lemma, page 129 in \cite{Hat:ATp:02}]\label{l:5lemma}
	Consider the following commutative diagram with exact rows
	\[
	\xymatrix{
	M_1 \ar[r]\ar[d]_{\varphi_1} & M_2 \ar[r] \ar[d]_{\varphi_2} &M_3 \ar[r] \ar[d]_{\varphi_3} &M_4 \ar[r]\ar[d]_{\varphi_4} & M_5 \ar[d]_{\varphi_5}\\
	N_1 \ar[r] & N_2 \ar[r] &N_3 \ar[r] &N_4 \ar[r] & N_5\\
	}
	\]
	where $\varphi_1$, $\varphi_2$, $\varphi_4$ and $\varphi_5$ are isomorphisms. Then $\varphi_3$ must be an isomorphism.
\end{Lemma}

\noindent The following proposition is sometimes called the \emph{zig-zag lemma}.

\begin{Proposition}[Theorem 1.3.1 in \cite{Wei:IHA:94}]\label{prop:zigzag}
	Associated to every short exact sequence of chain complexes
	\[
	0 \rightarrow C' \stackrel{i}{\rightarrow} C \stackrel{p}{\rightarrow} C'' \rightarrow 0,
	\]
	there is a \emph{long exact sequence in homology}
	\[
	\cdots \rightarrow H_k(C') \stackrel{H_k(i)}{\rightarrow} H_k(C) \stackrel{H_k(p)}{\rightarrow} H_k(C'') \stackrel{\delta}{\rightarrow} H_{k-1}(C') \rightarrow \cdots.
	\]
	The homomorphisms $\delta$ are called the \emph{connecting homomorphisms}.
\end{Proposition}

\noindent There is a more general construction, where $i$ need not be injective.

\begin{Definition}
	Let $(C_*,\partial^C)$, $(C'_*,\partial^{C'})$ be two chain complexes and let $\varphi_*:C'_*\rightarrow C_*$ be a chain map. The \emph{mapping cone chain complex $\Cone_*(\varphi)$} has the modules
	\[
	\Cone_k(\varphi) = C'_{k-1} \oplus C_k
	\]
	with differential
\[
\partial^{\Cone}_k (c' + c) \coloneq - \partial^{C'}_{k-1}(c') + \varphi_{k-1}(c') + \partial^C_k(c).
\]
\end{Definition}

\noindent As in the situation before, there is an associated long exact sequence. Note that the reference uses the opposite sign convention, but this does not change the proof.

\begin{Proposition}[1.5.2 in \cite{Wei:IHA:94}]\label{prop:les_mapping_cone}
	Let $\varphi : C' \rightarrow C$ be a chain map. Then there is a long exact sequence
	\[
    \cdots \rightarrow H_k(C') \stackrel{H(\varphi)}{\rightarrow} H_k(C) \stackrel{H(0\oplus\id)}{\rightarrow} H_k(\Cone(\varphi)) \stackrel{H(\pr_1)}{\rightarrow} H_{k-1}(C') \rightarrow \cdots.
	\]
\end{Proposition}

\noindent For injective chain maps, these concepts coincide.

\begin{Lemma}[1.5.8 in \cite{Wei:IHA:94}]\label{l:mapping_cone_relative_homology}
  If $\varphi: C'\rightarrow C$ is the inclusion of subcomplexes $C'\subseteq C$, then
  \[
	  H_*(\Cone(\varphi)) \cong H_*(C/C').
  \]
\end{Lemma}

\noindent We remind the reader briefly of the concept of flat modules and exactness.

\begin{Definition}
    A right module $M$ over $R$ is called \emph{$R$-flat} if the functor $(M \otimes_R -)$ is exact, that is, for every short exact sequence
    \[
        0 \rightarrow A \rightarrow B \rightarrow C \rightarrow 0
    \]
    of left $R$-modules, the sequence
    \[
        0 \rightarrow M \otimes_R A \rightarrow M \otimes_R B \rightarrow M \otimes_R C \rightarrow 0
    \]
    is also exact.
\end{Definition}

\begin{Remark}
    This implies in particular that
    \[
        M \otimes_R A / M\otimes_R B \cong M\otimes_R (A/B)
    \]
    for any left $R$-modules $B\leq A$ and that 
    \[
        M\otimes_R H_*(C) \cong H_*(M\otimes_R C_*)
    \]
    for any chain complex $C$ of left $R$-modules.
\end{Remark}

\begin{Lemma}[Section 3.2 in \cite{Wei:IHA:94}]
    Projective modules are flat. In particular, free modules are flat.
\end{Lemma}

\section{Group homology}

Group homology assigns to any group a sequence of abelian groups as invariants. These homology groups can be defined in purely algebraic terms or via topology. We have chosen to use the algebraic definition, since it simplifies the definition of homology with non-trivial coefficient modules. The relation to topology will be indicated in Section \ref{subsec:interpretations}.

\subsection{Definition}

Throughout this section, let $G$ be any group. We fix some basic concepts.

\begin{Definition}
    The \emph{(integral) group ring $\Z G$} is the free abelian group generated by the elements of $G$. The multiplication on $G$ extends uniquely to a $\Z$-bilinear multiplication on $\Z G$.

    We define a \emph{$G$-module} to be a left module over $\Z G$. A $G$-module $M$ is said to be \emph{trivial} if $gm=m$ for all $g\in G$ and all $m\in M$.
\end{Definition}

\begin{Construction} Tensor products $M\otimes N$ of modules are only defined if $M$ is a right module and $N$ is a left module over a common ring. We want to form tensor products of left modules over the group ring $\Z G$. Note that we can canonically make any left $\Z G$-module $M$ into a right $\Z G$-module by setting
\[
	mg := (g^{-1})m \quad\text{for any}\quad m\in M, g\in G.
\]
As in \cite[II.2]{Bro:CoG:82}, using this construction, we can define tensor products of left $G$-modules $M$ and $N$, denoted by $M\otimes_G N$.
\end{Construction}

\begin{Definition}
    Associated to a group $G$, consider the modules $F_k(G)$ which are the free abelian groups over $(k+1)$-tuples of elements of $G$. We define the boundary map on the basis by
    \[
    \partial(g_0,g_1,\ldots,g_k) = \sum_{i=0}^k (-1)^i (g_0,\ldots,\hat g_i,\ldots,g_k),
    \]
    where, as usual, the hat $\hat g_i$ indicates that the entry $g_i$ is omitted from the tuple. We call the associated chain complex of the form
	\[
	\cdots \rightarrow F_2(G)\rightarrow F_1(G) \rightarrow F_0(G) \rightarrow \Z \rightarrow 0,
	\]
    the \emph{standard resolution of $\Z$ over $\Z G$}.
\end{Definition}

\begin{Lemma}[I.5 in \cite{Bro:CoG:82}]
    This is an exact chain complex of $G$-modules.
\end{Lemma}

\begin{Definition}\label{def:group_homology}
	For any $G$-module $M$, the \emph{group homology $H_*(G;M)$} is defined to be the homology of the chain complex $F_*(G) \otimes_G M$.
    
    As usual, we abbreviate $H_*(G)\coloneq H_*(G;\Z)$, where $\Z$ is considered to be a trivial $G$-module.
\end{Definition}

\begin{Remarks}
    Any group monomorphism $G'\rightarrow G$ induces a chain map $F_*(G')\rightarrow F_*(G)$ which is injective on each chain module.

    For any other free (even projective) resolution $F_*$ of $\Z$ over $\Z G$, we also have
    \[
        H_*(G;M)\cong H_*(F_* \otimes_G M)
    \]
    by \cite[Theorem I.7.5]{Bro:CoG:82}.

    We will frequently use the following observation: If $G'\leq G$ and $F_*$ is a free resolution of $\Z$ over $\Z G$, then it is also a free resolution of $\Z$ over $\Z G'$.
\end{Remarks}

\subsection{Interpretations of group homology}\label{subsec:interpretations}

Low-dimensional homology groups have direct algebraic interpretations.

\begin{Remark} A simple calculation shows that
\[
    H_0(G,M)\cong \Z \otimes_G M \cong M/\langle gm - m : g\in G, m\in M\rangle,\label{e:calch0}
\]
which especially means $H_0(G)\cong \Z$. It is not much more difficult to show that
\[
	H_1(G)\cong (G)_{ab} = G/[G,G].
\]
\end{Remark}

\noindent There is also an algebraic interpretation of $H_2(G)$ as follows.

\begin{Definition}
	A \emph{central extension} of a group $G$ by a group $Z$ is a short exact sequence of groups
	\[
		1 \rightarrow Z \rightarrow E \rightarrow G \rightarrow 1,
	\]
	such that the image of $Z$ is central in $E$. A \emph{universal central extension} of $G$ is a central extension
	\[
		1 \rightarrow U \rightarrow X \rightarrow G \rightarrow 1
	\]
	such that, given any other extension as above, there is unique map $X\rightarrow E$ such that the diagram
	\begin{center}\begin{tikzpicture}
		\node (F) at (0,1) {$X$};
		\node (E) at (0,0) {$E$};
		\node (G) at (1,0.5) {$G$};
		\draw[->] (F) -- (E);
		\draw[->] (E) -- (G);
		\draw[->] (F) -- (G);
	\end{tikzpicture}\end{center}
	commutes. As usual for objects satisfying universal properties, the universal central extension is unique up to unique isomorphism.
\end{Definition}

\noindent It turns out that, if a universal central extension exists, its kernel is isomorphic to the second homology group. This kernel is classically called the \emph{Schur multiplier}.

\begin{Theorem}[Theorem 6.9.5 in \cite{Wei:IHA:94}]\label{th:schur_multiplier}
	A group $G$ admits a universal central extension if and only if it is \emph{perfect}, which means $G=[G,G]$ or equivalently $H_1(G)=0$. In this case, the universal central extension is given by
	\[
		1 \rightarrow H_2(G) \rightarrow X \rightarrow G \rightarrow 1
	\]
	with an appropriate group $X$.
\end{Theorem}

\noindent In fact, this description can often be used to calculate $H_2$ by finding the universal central extension.

\paragraph{} The description we gave for group homology is a purely algebraic one. Equivalently, one can give a topological description of group homology.

A connected CW complex $X$ whose higher homotopy groups are trivial except for $\pi_1$ is called an \emph{aspherical space}. Hurewicz proved that the homotopy type of an aspherical space depends only on the fundamental group. In particular, the homology groups of the space $H_*(X)$ can be interpreted as an invariant of the fundamental group $\pi_1(X)$. For any given group $G$, there is an aspherical CW complex $BG$ with this group as fundamental group, called the \emph{classifying space} associated to $G$, see \cite[8.1.7]{Wei:IHA:94}.

\begin{Theorem}[Proposition 4.2 in \cite{Bro:CoG:82}]
	If $X$ is an aspherical CW complex with fundamental group $G$, then the augmented cellular chain complex
	\[
		\cdots\rightarrow C_2(\tilde X) \rightarrow C_1(\tilde X) \rightarrow C_0(\tilde X) \rightarrow \Z \rightarrow 0
	\]
	is a free resolution of $\Z$ over $\Z G$. Here, $\tilde X$ is the universal cover of $X$. Its cellular chains are naturally $G$-modules. In particular, we have
    \[
        H_k(G;\Z) \cong H_k(X;\Z)
    \]
    by observing that $G\backslash C_*(\tilde X)\cong C_*(X)$.
\end{Theorem}

\noindent This shows that the notion of using the homology groups of the space $X$ as an invariant of its fundamental group leads precisely to the concept of group homology. Using this and the Hurewicz theorem, the above descriptions of $H_0$ and $H_1$ are obvious.

The interaction between topology and algebra via group homology has been very fruitful. One of the first applications of group homology in topology is the following result by Hopf.

\begin{BreakTheorem}[Hopf, see Theorem II.5.2 in \cite{Bro:CoG:82}]
	For a connected CW complex $X$ with fundamental group $G$ there is a map $H_*(X)\rightarrow H_*(G)$ and the sequence
	\[
		\pi_2(X) \rightarrow H_2(X) \rightarrow H_2(G) \rightarrow 0
	\]
	is exact.
\end{BreakTheorem}

\subsection{Basic results}

In this section, we will collect some basic results in group homology, which we will use later on. We begin with the explicit calculation of group homology for cyclic groups. The explicit determination of group homology is very hard in general, but for cyclic groups this is not difficult.

\begin{Proposition}[II.3.1 in \cite{Bro:CoG:82}]\label{prop:hom_of_cycl_groups}
	If $G$ is a cyclic group of order $n$, then
	\[
		H_k(G) \cong\begin{cases}
			\Z & k = 0, \\
			\Z / n & k \text{ odd}, \\
			0 & \text{otherwise.}
		\end{cases}
	\]
\end{Proposition}

\noindent This following lemma is a special case of \cite[Theorem 6.1.12]{Wei:IHA:94}.

\begin{Lemma}\label{l:trivial_coefficient_module}
    Let $G$ be a group and let $M$ be a trivial and $\Z$-free $G$-module. Then we have
    \[
        H_*(G;M) \cong H_*(G;\Z) \otimes_\Z M.
    \]
\end{Lemma}

\begin{Proof}
	Since $M$ is trivial, we have $F_k(G) \otimes_G M \cong F_k(G) \otimes_G \Z \otimes_\Z M$. Since $M$ is $\Z$-free, it is in particular $\Z$-flat and the functor $(-\otimes_\Z M)$ is exact, which proves the result.
\end{Proof}

\noindent The lemma of Shapiro allows the calculation of group homology using subgroups. This will be particularly helpful if the $G$-action is transitive on a $\Z$-basis of a module $N$ with a stabiliser $G'$. In this case, we obviously have $N \cong \Z G \otimes_{G'} \Z$.

\begin{Lemma}[Shapiro's lemma, III.6.2 in \cite{Bro:CoG:82}]\label{l:shapiro}
    If $G'\leq G$ is a subgroup and $M$ is a $G'$-module, we have
    \[
    H_*(G';M) \cong H_*(G; \Z G\otimes_{G'} M).
    \]
\end{Lemma}

\noindent In particular, we obtain that higher dimensional homology of induced modules vanishes.

\begin{Corollary}[III.6.6 in \cite{Bro:CoG:82}]\label{cor:induced_modules_are_acyclic}
    Let $G$ be any group and let $M$ be an abelian group. Then $H_k(G;\Z G \otimes_\Z M) = 0$ for $k>0$.
\end{Corollary}

\noindent Finally, using the so-called \emph{transfer map}, one can see that group homology of finite groups with rational coefficients vanishes.

\begin{Proposition}[Corollary III.10.2 in \cite{Bro:CoG:82}]\label{prop:rational_homology_of_finite_groups}
    If $G$ is a finite group, then we have $H_k(G;\Q) = 0$ for all $k > 0$.
\end{Proposition}

\subsection{Relative group homology}

The introduction of relative homology simplifies the formulation of homological stability considerably. For a group $G$ and a subgroup $G'$, note that there is canonically a map
\begin{align*}
    F_*(G') \otimes_{G'} M &\rightarrow F_*(G) \otimes_G M \\
    f \otimes m &\mapsto f\otimes m,
\end{align*}
induced by the inclusion $G'\rightarrow G$. This chain map is easily seen to be injective.

\begin{Definition}
  Consider a group $G$ and a subgroup $G'\leq G$. Let $M$ be a $G$-module. The \emph{relative group homology} $H_*(G,G';M)$ is defined to be the homology of the quotient complex $F_*(G)\otimes_G M / F_*(G') \otimes_{G'} M$.
\end{Definition}

\noindent Note that there is canonically an associated long exact sequence of the form
\[
\cdots \rightarrow H_k(G';M) \rightarrow H_k(G;M) \rightarrow H_k(G,G';M) \rightarrow H_{k-1}(G';M) \rightarrow \cdots.
\]

\noindent The following description will make the calculation of relative homology groups of subgroups much simpler.

\begin{Proposition}\label{prop:relative_homology_of_subgroups}
    Let $G'$ and $H$ be two subgroups of a given group $G$, write $H'=G'\cap H$. Then we have
    \[
    H_*\Bigl( \frac{F_*(G)\otimes_H M}{F_*(G')\otimes_{H'}M}\Bigr) \cong H_*\Bigl( \frac{F_*(H)\otimes_H M}{F_*(H')\otimes_{H'}M}\Bigr) = H_*(H,H';M).
    \]
\end{Proposition}

\begin{Proof}
    Consider the following diagram:
    \[
        \xymatrix{
        0\ar[r] & F_*(H') \otimes_{H'} M\ar[r]\ar[d] & F_*(H) \otimes_H M\ar[r]\ar[d] & \frac{F_*(H)\otimes_H M}{F_*(H')\otimes_{H'}M}\ar[r]\ar[d] & 0 \\
        0\ar[r] & F_*(G') \otimes_{H'} M\ar[r] & F_*(G) \otimes_H M\ar[r] & \frac{F_*(G)\otimes_H M}{F_*(G')\otimes_{H'}M}\ar[r] & 0
        }
    \]
    Here, the top row is exact by construction. The bottom row is exact since $H'=G'\cap H$. The vertical maps are all induced by inclusions. If we now consider the associated long exact sequences on homology, the first and second vertical arrow induce isomorphisms by the remark after Definition \ref{def:group_homology}. By the 5-lemma \ref{l:5lemma}, we get the desired isomorphisms.
\end{Proof}

\noindent It is easy to see that relative $H_0$ vanishes always.

\begin{Lemma}\label{l:relative_h0}
  For any group $G$, any subgroup $G'$ and any $G$-module $M$, we have
  \[
    H_0(G,G';M) = 0.
  \]
\end{Lemma}

\begin{Proof}
	By the remark above, we know that $H_0(G;M) = \Z \otimes_G M$ and $H_0(G';M)= \Z \otimes_{G'} M$. Now by the long exact sequence we have
  \[
    H_0(G,G';M) \cong (\Z \otimes_G M) / i_*(\Z \otimes_{G'} M),
  \]
  where $i_*$ is the map $z \otimes_{G'} m \mapsto z \otimes_G m$. This map is clearly surjective.
\end{Proof}

\section{Homological stability}\label{sec:stability}

Group homology is often difficult to compute. But given a sequence of groups
\[
    G_1\subset G_2 \subset \cdots \subset G_n \subset \cdots,
\]
we can ask whether there is an integer $n(k)$ such that
\[
	H_k(G_n) \rightarrow H_k(G_{n+1})
\]
is an isomorphism for $n\geq n(k)$ or, equivalently, whether relative homology
\[
H_k(G_{n+1},G_n)
\]
vanishes for $n\geq n(k)$. This phenomenon is called \emph{homological stability}. This happens frequently, for instance for symmetric groups with $\Z_p$-coefficients (see \cite{Nak:Dec:60}), for mapping class groups (see \cite{Har:SMC:85,Wah:HSM:08}), for outer automorphism groups of free groups (see \cite{HV:HSO:04,HVW:HSO:06}) and for classical groups. A very good overview of homological stability for linear groups can be found in the book by Knudson, see \cite[Chapter 2]{Knu:HLG:01}.

Usually, homological stability for classical groups is discussed in the generality of linear groups over rings with some additional conditions, for example finite stable rank. We shall not go into details here and we will simply cite several homological stability results which are related to this thesis, always specialising to the case of division rings.

\minisec{General and special linear groups}
For general and special linear groups, the strongest result for arbitrary division rings known to the author is a special case of a strong theorem by van der Kallen:

\begin{Theorem}[van der Kallen, \cite{vdK:HSL:80}]\label{th:vdk}
  For any division ring $D$ we have
	\[
	H_k(\Gl_{n+1}(D),\Gl_n(D);\Z) = 0
	\]
  for $n\geq 2k$. If $D$ is commutative, then also
	\[
	H_k(\Sl_{n+1}(D),\Sl_n(D);\Z) = 0
	\]
    for $n\geq 2k$.
\end{Theorem}

\noindent For division rings with infinite centre, there is a much stronger result due to Sah. A detailed exposition of this result can be found in \cite{Ess:HSG:06}. A different proof of this result can be found in \cite[2.3]{Knu:HLG:01}.

\begin{Theorem}[Sah, Appendix B in \cite{Sah:HcL:86}]\label{th:sah}
	If $D$ has infinite centre, then $n\geq k$ implies
	\[
	H_k(\Gl_{n+1}(D),\Gl_n(D);\Z) = 0.
	\]
\end{Theorem}

\noindent For a long time, there have not been considerable improvements for the case of special linear groups. There is a new result by Hutchinson and Tao improving the known stability range considerably for fields of characteristic zero.

\begin{Theorem}[Hutchinson-Tao, \cite{HT:HSS:08}]\label{th:hutchinson_tao}
  If $D$ is a field of characteristic zero, then $n\geq k$ implies
  \[
	  H_k(\Sl_{n+1}(D),\Sl_{n}(D);\Z)=0.
  \]
\end{Theorem}

\minisec{Unitary groups}

For the unitary groups $\U_n(\K)$ over $\K\in\{\R,\C,\bH\}$ associated to the standard positive definite hermitian form, there is the following stability result due to Sah:

\begin{Theorem}[Sah, Theorem 1.1 in \cite{Sah:HcL:86}]
  If $n\geq k$, then
  \[
  H_k(\U_{n+1}(\K),\U_n(\K);\Z)=0.
  \]
\end{Theorem}

\noindent Again, a more detailed exposition of this result can be found in \cite{Ess:HSG:06}. In this thesis, we are more interested in unitary groups associated to $(\varepsilon,J)$-hermitian forms as in Section \ref{sec:ex_buildings_of_type_bn_cn}, since there are associated buildings. For these unitary groups, there is a homological stability result by Mirzaii and van der Kallen for local rings with infinite residue fields.

\begin{BreakTheorem}[Mirzaii-van der Kallen and Mirzaii, \cite{MaB:HSU:02, Mir:HSU:05}]\label{th:mirzaii}
  If $D$ is a division ring with infinite centre, $\caH_{n+1}$ is the $2(n+1)$-dimensional hyperbolic module over $D$ and $\caH_n$ is a $2n$-dimensional submodule, then
  \[
  H_k(\U(\caH_{n+1}),\U(\caH_n);\Z)=0
  \]
  for $n\geq k+1$.
\end{BreakTheorem}

\noindent Previously, weaker versions of this theorem for the special cases of orthogonal and symplectic groups have been shown by Vogtmann in \cite{Vog:HSO:79} and \cite{Vog:SPH:81}, which were later generalised to Dedekind rings by Charney in \cite{Cha:gtV:87}. Although these results have been superseded by the theorems by van der Kallen and Mirzaii, the proof method Vogtmann and Charney used is the one we will adapt in Chapter \ref{ch:stability}.

\section{Algebraic K-theory}

Algebraic $K$-theory is a sequence of functors from the category of rings to the category of abelian groups, first considered by Grothendieck in the fifties. In the history of $K$-theory, the low-dimensional $K$-groups $K_0$, $K_1$ and $K_2$ were defined first by Grothendieck, Bass and Milnor, respectively. The following definitions are classical.

\begin{Definition}[Grothendieck]
    Let $R$ be any ring. The set of isomorphism classes of projective $R$-modules $P(R)$ forms a commutative monoid with respect to direct sums. It can be completed to an abelian group by the \emph{Grothendieck construction} to obtain $K_0(R)$ as follows:
    \[
    K_0(R) \coloneq F(P(R)) / \langle [P] + [Q] - [P \oplus Q] \,:\, P,Q\text{ projective}\rangle,
    \]
    where $F(P(R))$ denotes the free abelian group over $P(R)$.
\end{Definition}

\begin{Definition}[Bass]
    For any ring $R$, consider the group $\Gl(R)$ which is the direct limit of the groups $\Gl_n(R)$ with respect to the natural inclusions. Define
    \[
        K_1(R) \coloneq \Gl(R) / [\Gl(R),\Gl(R)] \cong H_1(\Gl(R);\Z).
    \]
\end{Definition}

\noindent Denote the commutator subgroup of $\Gl(R)$ by $E(R)$.

\begin{Definition}[Milnor]
    Let $\St(R)$ be the universal central extension of $E(R)$. We define
    \[
        K_2(R) \coloneq \ker( \St(R)\twoheadrightarrow E(R)) \cong H_2(E(R);\Z).
    \]
\end{Definition}

\noindent After that, there was a period of uncertainty on how higher $K$-groups should be defined. Quillen solved this problem by giving two definitions for higher $K$-theory based on the $+$-construction and the $\cQ$-construction. We will only discuss the first possibility by sketching the $+$-construction. The $+$-construction and all of the above definitions are contained in the textbook by Rosenberg \cite{Ros:AKa:94}.

Let $X$ be a CW complex and let $N$ be a normal perfect subgroup of $\pi_1(X)$. The \emph{$+$-construction} then produces a space $X^+$ with the same homology groups as $X$ and with fundamental group $\pi_1(X^+)\cong \pi_1(X)/N$. This can be used to define higher algebraic $K$-theory as follows:

\begin{Definition}[Quillen]
    Let $R$ be any ring, then we define \emph{algebraic $K$-theory} by
    \[
        K_i(R) \coloneq \pi_i( \BGl(R)^+ \times K_0(R)),
    \]
    where $\BGl(R)$ is the classifying space of $\Gl(R)$, and where the $+$-construction is applied to $N=E(R)$ which is normal and perfect.
\end{Definition}

\noindent In particular, we see that the Hurewicz maps are homomorphisms
\[
    K_i(R) = \pi_i( \BGl(R)^+ \times K_0(R)) \rightarrow H_i( \BGl(R)^+ \times K_0(R)) \cong H_i(\Gl(R);\Z)
\]
for $i>0$. Here, we see the direct connection from $K$-theory to the group homology of the stable linear group and hence to homological stability.

\begin{Remark}
    It has been shown by Gersten that actually $K_3(R)\cong H_3(\St(R),\Z)$ which continues the series of definitions above.

    Several other tentative definitions of $K$-theory were given at that time, among these one by Volodin which can be shown to be isomorphic to $K$-theory as defined by Quillen. Another possible definition was given by Wagoner in \cite{Wag:BSK:73} using what we call \emph{Wagoner complexes} here. This definition was the starting point for our investigations in Chapter \ref{ch:wagoner}. Wagoner's definition of higher algebraic $K$-theory can be found in Section \ref{sec:ktheory}. It is also isomorphic to Quillen's definition of $K$-theory.
\end{Remark}

\paragraph{Hermitian $K$-theory} If we start from unitary groups $\U_n(R)$, defined analogously to Section \ref{sec:ex_buildings_of_type_bn_cn} over a ring with involution $J$ and with an element $\varepsilon\in\{\pm 1\}$, we can also define a stable group $\U(R)$ via the direct limit. Then, if we copy Quillen's definition, we obtain \emph{hermitian $K$-theory}
\[
{}_\varepsilon L_i(R) \coloneq \pi_i( \BU(R)^+ \times L_0(R)),
\]
where $L_0(R)$ is the Grothendieck group over isomorphism classes of hyperbolic modules, in analogy to $K_0(R)$.

Both algebraic $K$-theory and hermitian $K$-theory are important tools in topology, taking the $K$- and $L$-groups over the group ring of the fundamental group of a space. Many topological obstructions are elements in $K$- or $L$-groups.

\chapter{Spectral sequences}\label{ch:spectral}

The main computational tools for this thesis will be spectral sequences. We will give a definition of spectral sequences and describe all the spectral sequences which are relevant for this thesis.

All of these can be deduced from the spectral sequence associated to a filtered complex and from the spectral sequence associated to a double complex. We will present the constructions of these two spectral sequences, but since this is a classical topic which is covered in various textbooks, we will not give convergence proofs and simply state most of the results.

Using these two spectral sequences, we will deduce all further required spectral sequences.

The standard guide to spectral sequences is of course \cite{McC:UGS:01}. The book by Weibel contains a chapter on spectral sequences, see \cite[Chapter 5]{Wei:IHA:94}. The definition of spectral sequences may seem strange when seen for the first time. A very nice introduction motivating this structure can be found in \cite{Cho:Ych:06}. Our examples should also show that spectral sequences arise quite naturally.

\section{The definition of spectral sequences}

The definition of spectral sequences starts with the definition of bicomplexes.

\begin{Definition}
	A \emph{bicomplex $(C_{p,q},\partial_{p,q})$ of bidegree $(a,b)$} is a bigraded module with differential
	\[
		\partial_{p,q} : C_{p,q} \rightarrow C_{p+a,q+b}
	\]
    such that $\partial^2=0$. The bicomplex is said to be \emph{first-quadrant} if we have $C_{p,q}=0$ for $p<0$ or $q<0$.
\end{Definition}

\noindent A bicomplex thus is made of infinitely many chain complexes aligned in a special fashion. We can calculate homology and obtain a new bigraded module.

Now, if we had a new differential on the homology bigraded module, we could continue forming homology modules of homology modules. The following definition makes this idea precise.

\begin{Definition}
    A \emph{spectral sequence $E$} is a sequence of bicomplexes $(E^r_{p,q},d^r_{p,q})_{r\geq r_0}$ with bidegrees $(-r,(r-1))$, such that
	\[
		E^{r+1} = H(E^r).
	\]
    The spectral sequence is said to be \emph{first-quadrant} if the bicomplex $E^{r_0}$ is (and then all others are) first-quadrant. Usually, we have $r_0\in\{0,1,2\}$.
\end{Definition}

\begin{figure}[hbt]
\[
\begin{array}{c|c}
\begin{xy}
	\xymatrix@=1.2em{
	&&&&&E^0_{\star,\star}\\
	\ar[d]\\
	E^0_{0,3}\ar[d]&\ar[d]\\
	E^0_{0,2}\ar[d]&E^0_{1,2}\ar[d]&\ar[d]\\
	E^0_{0,1}\ar[d]&E^0_{1,1}\ar[d]&E^0_{2,1}\ar[d]&\ar[d]\\
	E^0_{0,0}\ar@.[uuuuu]\ar@.[rrrrr]&E^0_{1,0}&E^0_{2,0}&E^0_{3,0}&&
}
\end{xy}&
\begin{xy}
\xymatrix@=1.2em{
	&&&&&E^1_{\star,\star}\\
	\\
	E^1_{0,3}&\ar[l]\\
	E^1_{0,2}&E^1_{1,2}\ar[l]&\ar[l]\\
	E^1_{0,1}&E^1_{1,1}\ar[l]&E^1_{2,1}\ar[l]&\ar[l]\\
	E^1_{0,0}\ar@.[uuuuu]\ar@.[rrrrr]&E^1_{1,0}\ar[l]&E^1_{2,0}\ar[l]&E^1_{3,0}\ar[l]&\ar[l]&
}
\end{xy}\\\hline
\begin{xy}
\xymatrix@=1.2em{
	&&&&&E^2_{\star,\star}\\
	\\
	E^2_{0,3}\\
	E^2_{0,2}&E^2_{1,2}&\ar[llu]\\
	E^2_{0,1}&E^2_{1,1}&E^2_{2,1}\ar[llu]&\ar[llu]\\
	E^2_{0,0}\ar@.[uuuuu]\ar@.[rrrrr]&E^2_{1,0}&E^2_{2,0}\ar[llu]&E^2_{3,0}\ar[llu]&\ar[llu]&
}
\end{xy}&

\end{array}
\]	\caption{The first three pages of a first-quadrant spectral sequence} \label{f:firstthreepages}
\end{figure}
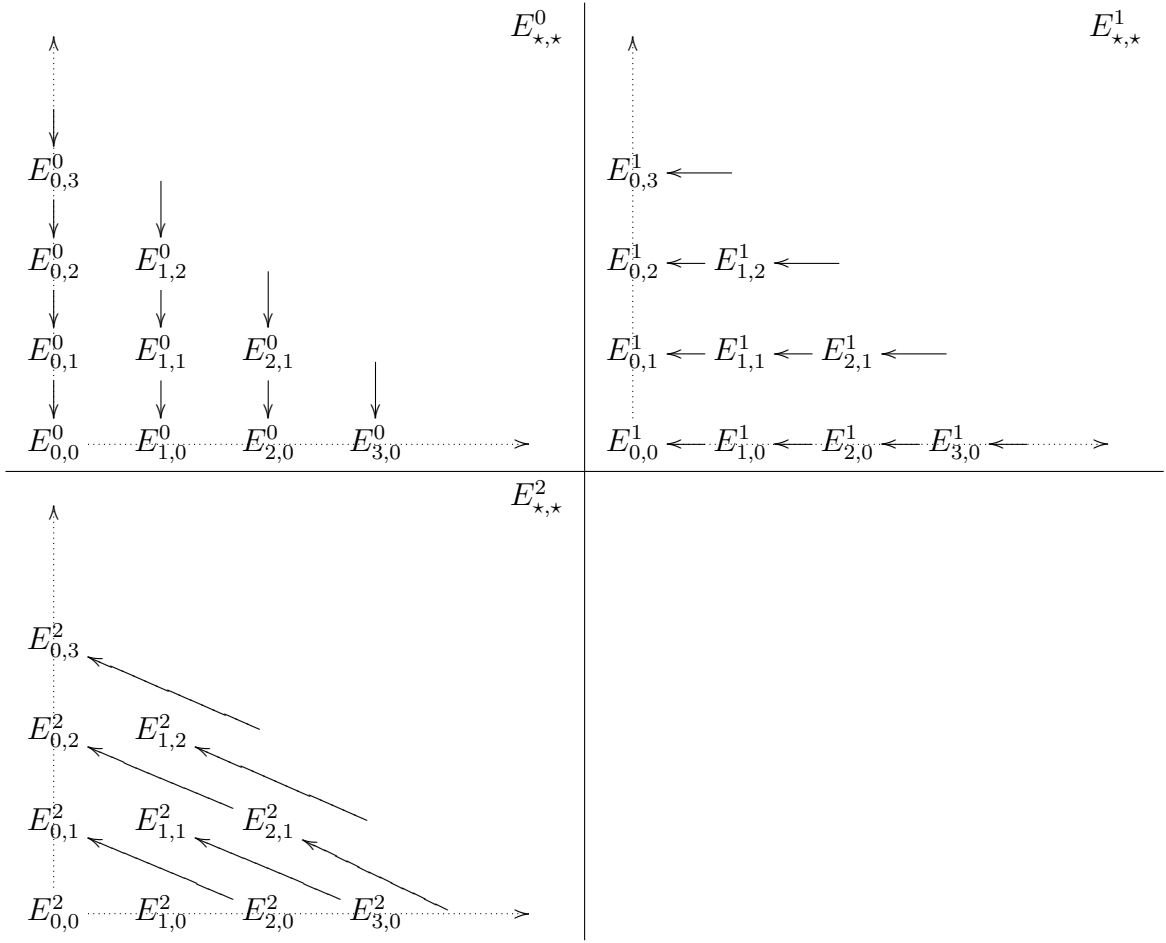

\paragraph{Visualisation} We picture this as a book whose pages are indexed by $r$. On each page, the graded module $E^r$ can be visualised as the modules $E^r_{p,q}$ sitting on the point $(p,q)$ in the real plane. Turning a page corresponds to calculating homology groups with respect to the differentials on this page. The differential arrows get longer and more slanted as the pages turn. For a sketch of the first three pages of a first-quadrant spectral sequence starting from $r_0=0$, see Figure \ref{f:firstthreepages}.

\paragraph{Convergence} We assume that for each $k\in\Z$ there are only finitely many non-zero modules $E^{r_0}_{p,q}\neq 0$ for $p+q=k$. In particular, this is true for first-quadrant spectral sequences.

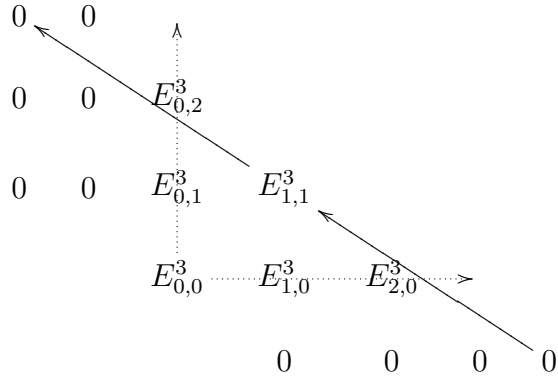
\begin{figure}[htb]
\[
	\xymatrix@=1.2em{
	0&0&\\
	0&0&E^3_{0,2}\\
	0&0&E^3_{0,1}& \ar[llluu]E^3_{1,1}\\
	&&E^3_{0,0}\ar@.[uuu]\ar@.[rrr] & E^3_{1,0} & E^3_{2,0}&\\
	&&&0&0&0&0\ar[llluu]
	}
\] \caption{Convergence of a first-quadrant spectral sequence} \label{f:convergence}
\end{figure}

In other words, for every pair $(p,q)$, there is hence a page number $R$ such that for $r\geq R$ the incoming and outgoing arrows of $E^{r}_{p,q}$ begin or end in trivial groups, see Figure \ref{f:convergence}. In particular, the module $E^R_{p,q}$ does not change any more, we call this stable module $E^\infty_{p,q}$. So for any pair $(p,q)$ there is an index $R(p,q)$, such that
\[
E^{R(p,q)}_{p,q}=E^{R(p,q)+1}_{p,q}=\cdots=:E^\infty_{p,q}.
\]
We call the resulting bigraded module $E^\infty$ the \emph{abutment} of the spectral sequence and we say that the spectral sequence $E$ \emph{converges} to $E^\infty$.
\paragraph{Collapse} If, for any reason, all differentials are zero from a page $E^R$ on, then obviously
\[
	E^R=E^{R+1}=\cdots=E^\infty,
\]
and we say that the spectral sequence \emph{collapses} on the $R$-th page.

\begin{figure}[hbt]
\[
	\xymatrix@=1.2em{
	0&&\\
	&0&\ar[l]\ar[llu]E_{0,2}\\
	&&E_{0,1}&E_{1,1}\\
	&&E_{0,0}\ar@.[uuu]\ar@.[rrr] & E_{1,0} & E_{2,0}&\\
	&&&&&0\ar[llu]\\
	&&&&&&0\ar[llluu]
	}
\]
	\caption{Edge homomorphisms: Incoming and outgoing differentials from the first and second page on.} \label{f:edgehom}
\end{figure}
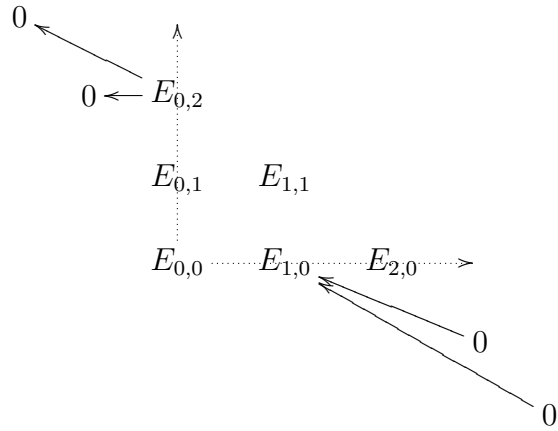

\paragraph{Edge homomorphisms} For a first-quadrant spectral sequence, let us consider the modules $E^r_{p,0}$ for $r\geq 2$. Since all `incoming' differentials from the second page on are zero (see Figure \ref{f:edgehom}), each $E^{r+k}_{p,0}$ is a submodule of $E^r_{p,0}$ for any $k\geq 0$. This of course is also true for $E^\infty_{p,0}$.

So for any $r\geq 2$ and $p\geq 0$ we have a canonical inclusion
\[
	\iota_p: E^\infty_{p,0} \hookrightarrow E^r_{p,0}.
\]

\noindent On the other hand, every `outgoing' differential of the groups $E^r_{0,q}$ is zero from the first page on. Then every $E^{r+k}_{0,q}$ is a quotient of $E^r_{0,q}$ for any $k\geq 0$. So we have a canonical projection
\[
	\pi_q: E^r_{0,q} \twoheadrightarrow E^\infty_{0,q}
\]
for any $r\geq 1$ and any $q\geq 0$. These maps are called \emph{edge homomorphisms}. In the `book picture' they can be imagined as arrows pointing up, respectively down at the `edges of the book'.

\section{The spectral sequence associated to a filtered chain complex}

Spectral sequences arise very naturally in the context of filtered chain complexes. In this section, we take all graded modules (chain complexes) $(M_q)_{q\in\Z}$ to be non-negative, that means $M_q=0$ for $q<0$.

\begin{Definition}\label{def:filtered_graded_module}
    A \emph{filtered graded module (chain complex)} is a sequence of graded submodules (subcomplexes) $(\cF^p M)_{p\in\Z}$ of a graded module (chain complex) $(M_q)_{q\in\Z}$ such that $\cF^p M\subseteq \cF^{p+1} M$.

    A filtration is called \emph{bounded} if for every $q$, there are integers $p$ and $P$ such that $\cF^p M_q = 0$ and $\cF^P M_q=M_q$.

    If $p=-1$ and $P=q$, that is $\cF^{-1} M = 0$ and $\cF^q M_q = M_q$ for all $q\in\Z$, the filtration is called \emph{canonically bounded}.
\end{Definition}

\paragraph{Remark} A filtration on a chain complex $C$ induces a natural filtration on the graded module $H(C)$ as follows:
\[
	\cF^p H(C) := \im\bigl( H(\cF^p C) \rightarrow H(C)\bigr).
\]

\begin{Definition}
    To any filtered graded module $(M_q)_{q\in\Z}$, especially to a filtered chain complex, we associate a bigraded module
\[
	\Gr_{p,q} M := \cF^p M_{p+q} / \cF^{p-1} M_{p+q}.
\]
Note that for a canonically bounded filtration $\Gr_{p,q}M=0$ for $p<0$ or $q<0$. So for a graded module $(M_q)_{q\in\Z}$ with a canonically bounded filtration, $\Gr_{p,q}M$ is a first-quadrant bigraded module.
\end{Definition}
\paragraph{Construction}For a filtered chain complex $(C,\partial^C)$, we have
\[
	E^0_{p,q} := \Gr_{p,q} C = \cF^p C_{p+q} / \cF^{p-1} C_{p+q}
\]
and the differential $\partial^C$ induces differentials of bidegree $(0,-1)$
\[
	d^0_{p,q}:\Gr_{p,q} C \rightarrow \Gr_{p,q-1}C.
\]

\noindent This means that $\Gr_{p,q}$ splits into `vertical' chain complexes, for which we can calculate homology groups. We obtain a bigraded module called
\[
	E^1_{p,q} := H_q(\Gr_{p,\stern}C) = \ker d^0_{p,q} /\im d^0_{p,q+1}.
\]
If we inspect the involved groups more closely, we can see that the differential $\partial^C$ induces new, `horizontal' differentials (of bidegree $(-1,0)$) on $E^1$. We can again calculate homology groups and obtain a graded module $E^2$. Here the differential $\partial^C$ induces new differentials, this time of bidegree $(-2,1)$. We again calculate homology and continue this process.

At every step, the differential $\partial^C$ induces new differentials with the right bidegree. A rather tedious calculation shows that this is a spectral sequence and that it converges to the graded module associated to the homology of $C$,
\[
	E^\infty_{p,q} = \Gr_{p,q} H(C) = \cF^p H_{p+q}(C) / \cF^{p-1} H_{p+q}(C).
\]
This results in the following theorem:

\begin{Theorem}[Section 5.4 in \cite{Wei:IHA:94}]\label{th:exspecseq}
	Associated to any chain complex $C$ with a bounded filtration, there is a spectral sequence with first terms
	\[
	E^0_{p,q}=\Gr_{p,q} C,\qquad E^1_{p,q} = H_{p+q}(\cF^p C / \cF^{p-1} C)
	\]
	converging to
	\[
		E^\infty_{p,q} = \Gr_{p,q} H(C).
	\]
	For a canonically bounded filtration, this is a first-quadrant spectral sequence. The standard notation for convergence is
	\[
		E^1_{p,q}= H_{p+q}(\cF^p C / \cF^{p-1} C) \quad\Rightarrow\quad H_{p+q}(C).
	\]
\end{Theorem}

\paragraph{Convergence} Note that, although the notation suggests that the abutment of this spectral sequence is $H_\stern(C)$, we actually only get something weaker, namely just the bigraded module $\Gr_{\stern,\stern}H(C)$. In general, it may not be possible to reconstruct $H_\stern(C)$ from $\Gr_{\stern,\stern}H(C)$, this is an extension problem. However, we will only deal with cases where this is possible.

\section{The spectral sequence associated to a filtered topological space}

The construction of the previous section can easily be applied to filtered topological spaces. We will only require the following simple version.

Let $X$ be a topological space filtered by a finite ascending chain of closed subspaces $\emptyset=X_{-1}\subseteq X_0\subseteq X_1\subseteq X_2\subseteq\cdots \subseteq X_n=X$. This implies that the induced filtration $C_*(X_{i})_{i\in\Z}$ on the simplicial chains of $X$ is bounded in the sense of Definition \ref{def:filtered_graded_module}. We apply Theorem \ref{th:exspecseq} to obtain

\begin{Theorem}\label{th:specseq_filtered_top_space}
	Associated to a filtered topological space $X$ as above, there is a spectral sequence
	\[
	E^1_{p,q} = H_{p+q}(X_p,X_{p-1}) \Rightarrow H_{p+q}(X).
	\]
\end{Theorem}

\begin{Proof}
	We apply Theorem \ref{th:exspecseq} to the filtered chain complex $C_*(X_p)$ and obtain a spectral sequence
	\[
	E^1_{p,q} = H_{p+q}(C_*(X_p)/C_*(X_{p-1})) = H_{p+q}(X_p,X_{p-1}) \Rightarrow H_{p+q}(X).
	\]
\end{Proof}

\begin{Remark} Note that this spectral sequence can be used to define cellular homology. If we take the filtration of a CW complex $X$ to be the filtration by skeletons, the spectral sequence is concentrated on the $x$-axis, that is $E^1_{p,q}=0$ for $q\neq 0$. The groups $E^1_{p,0}$ are precisely the groups of cellular chains and the differential on the first page is exactly the cellular boundary map. The convergence to $H_{p+q}(X)$ then proves that cellular and singular homology are isomorphic.
\end{Remark}

\section{The spectral sequences associated to a double complex}

The basis for all the following spectral sequences will be the spectral sequence associated to a double complex.

\begin{Definition}
	A \emph{double complex $D_{p,q}$} is a bigraded module with horizontal and vertical differentials $\partial^h:D_{p,q}\rightarrow D_{p-1,q}$ and $\partial^v:D_{p,q}\rightarrow D_{p,q-1}$ such that
\[
\partial^v\circ \partial^h - \partial^h\circ \partial^v = (\partial^v)^2 = (\partial^h)^2= 0.
\]
Associated to a double complex is the \emph{total complex}, defined to be
\[
\Tot(D_{p,q})_k \coloneq \bigoplus_{p+q=k} D_{p,q}
\]
with differential induced by
\[
\partial^{\Tot}(d_{p,q}) = \partial^h(d_{p,q}) + (-1)^p \partial^v(d_{p,q}).
\]
\end{Definition}

\noindent We assume that $D_{p,q}$ is \emph{bounded}, that is, for every $k\in\Z$, there are only finitely many non-vanishing modules $D_{p,q}$ with $p+q=k$. Then the homology of the total complex can be calculated via two spectral sequences, which are easily deduced from the spectral sequence of Theorem \ref{th:exspecseq} by using two different filtrations on the total complex $\Tot(D)$.

\begin{Proposition}[\cite{Bro:CoG:82}, VII.3]\label{prop:two_spectral_sequences}
 There are two spectral sequences both converging to the homology of the total complex
 \[
\begin{array}{l}
  E^1_{p,q} = H_q((D_{*,p},\partial^h)) \\
  L^1_{p,q} = H_q((D_{p,*},\partial^v))
\end{array} \Rightarrow H_{p+q}(\Tot(D)_*).
 \]
 The differentials on the first page $E^1$ are induced by $\pm\partial^v$, the differentials on $L^1$ are induced by $\partial^h$.
\end{Proposition}

\noindent For us, double complexes arise in the following situation: Fix a group $G$ and let $(F_*,\partial^F)$ and $(C_*,\partial^C)$ be chain complexes of $G$-modules. Assume that both $F$ and $C$ have only finitely many modules in negative dimensions. Then $D_{p,q}=F_p\otimes_G C_q$ with the differentials $(\partial^F \otimes_G\id)$ and $(\id\otimes_G \partial^C)$ is a bounded double complex.

The total complex of this double complex coincides with the \emph{tensor product of chain complexes $(F\otimes_G C)_*$}:
\[
\Tot(D_{p,q})_k = (F\otimes_G C)_k \coloneq \bigoplus_{p+q=k} F_p \otimes_G C_q
\]
with differential induced by
\[
\partial^{F\otimes_G C}(f_p \otimes_G c_q) = \partial^F(f_p) \otimes_G c_q + (-1)^p f_p \otimes_G \partial^C(c_q).
\]

\noindent Proposition \ref{prop:two_spectral_sequences} specialises to the following corollary.

\begin{Corollary}\label{cor:two_spectral_sequences_for_tensor_product}
 There are two spectral sequences both converging to the homology of the tensor product complex
\[
\begin{array}{r}
    E^1_{p,q} = H_q(F_* \otimes_G C_p) \\
    L^1_{p,q} = H_q(F_p \otimes_G C_*)
\end{array} \Rightarrow H_{p+q}(F\otimes_G C).
\]
The differentials on the first page are induced by $\pm(\id\otimes_G\partial^C)$ and $(\partial^F\otimes_G\id)$, respectively.
\end{Corollary}
Note that, if $F_*$ is a chain complex of free $G$-modules, each module $F_p$ is $\Z G$-free and hence $\Z G$-flat, and we obtain
\[
L^2_{p,q} = H_p(F_* \otimes_G H_q(C)) \Rightarrow H_{p+q}(F\otimes_G C).
\]
In particular, if $F_*=F_*(G)$ is the standard resolution of $\Z$ over $\Z G$ and if the chain complex $C_*$ is an exact complex of $G$-modules, we obtain

\begin{Corollary}\label{cor:ss_exact_coefficients}
        Let $G$ be a group and let $C_*$ be an exact chain complex of $G$-modules. Then there is a spectral sequence
        \[
        E^1_{p,q} = H_q(G;C_p) \Rightarrow 0
        \]
        which converges to zero.
\end{Corollary}

\noindent We will later also use transposed double complexes as follows: Let $F_*$ and $C_*$ be chain complexes of $G$-modules. We denote the \emph{transposed tensor complex} by $F \otimes^T_G C$:
\[
(F\otimes^T_G C)_k \coloneq \bigoplus_{p+q=k} F_p \otimes_G C_q
\]
with differential induced by
\[
\partial^{F\otimes^T_G C}(f_p \otimes_G c_q) = (-1)^q \partial^F(f_p) \otimes_G c_q + f_p \otimes_G \partial^C(c_q).
\]
This is obviously the total complex of the double complex $(F_q \otimes_G C_p)_{p,q}$.

\begin{Lemma}\label{lem:transposed_double_complex}
    The chain map $F\otimes_G C\rightarrow F\otimes^T_G C$ given on the basis by
    \[
    f_p \otimes_G c_q \mapsto (-1)^{pq} f_p \otimes_G c_q
    \]
    induces isomorphisms on homology $H_*(F\otimes_G C) \cong H_*(F\otimes^T_G C)$.
\end{Lemma}

\begin{Proof}
    It is easy to see that this map is a chain map. It is obviously bijective and induces hence an isomorphism on homology.
\end{Proof}

\section{Equivariant homology}\label{sec:equivariant_homology}

We will now specialise Corollary \ref{cor:two_spectral_sequences_for_tensor_product} to the case where $C_*=C_*(X)$ is the complex of cellular chains of a CW complex. The following construction can also be found in \cite[VII.7 and VII.8]{Bro:CoG:82}.

\begin{Definition}
	Let $X$ be a CW complex and let $G$ be a group acting on $X$ by homeomorphisms. The group action is called \emph{cellular} if the image of any cell is again a cell. We say that the group $G$ acts \emph{without inversions} if the setwise stabiliser of each cell in $G$ also fixes the cell pointwise.
\end{Definition}

\noindent A simple example for a cellular action is the action induced by a simplicial action on the geometric realisation of a simplicial complex. If the simplicial complex is a colourable chamber complex and the group action is type-preserving, then the induced cellular action is without inversions.

\begin{Theorem}\label{th:equiv_hom}
    If a group $G$ acts without inversions on a contractible CW complex $X$ and if $\Sigma_p$ is a system of representatives for the $p$-cells of $X$, there is a spectral sequence
    \[
    E^1_{p,q}= \bigoplus_{\sigma\in\Sigma_p} H_q(G_{\sigma}) \Rightarrow H_{p+q}(G).
    \]
    In particular, we have $E^1_{p,0}\cong \bigoplus_{\sigma\in\Sigma_p} \Z$ and the differential $d^1_{p,0}$ is induced by the differential $\partial_{p}$ on cellular chains. This results in
    \[
    E^2_{p,0}\cong H_p(G\backslash X).
    \]
\end{Theorem}

\begin{Proof}
    Let $F_*$ be a free resolution of $\Z$ over $\Z G$. We apply Corollary \ref{cor:two_spectral_sequences_for_tensor_product} to obtain two spectral sequences
    \[
    \begin{array}{l}
        E^1_{p,q} = H_q(G,C_p(X)) \\
        L^1_{p,q} = H_q(F_p \otimes_G C_*(X))
    \end{array} \Rightarrow H_{p+q}(F\otimes_G C(X)).
    \]
    As before, since $F_p$ is a free $G$-module, we obtain $L^1_{p,q} \cong F_p \otimes_G H_q(X)$, which vanishes except for $q=0$ since $X$ is contractible. The spectral sequence $L$ hence collapses on the second page and we obtain isomorphisms $H_*(G) \cong H_*(F \otimes_G C(X))$.

    On the other hand, notice that $C_p(X)$ decomposes as a direct sum
    \[
    C_p(X) = \bigoplus_{\sigma\in\Sigma_p} G \otimes_{G_\sigma} \Z,
    \]
    since the group acts without inversions (otherwise, we would have a module $\Z_\sigma$ on the right, which is twisted by elements of $G$ if the orientation of the cell $\sigma$ changes, see \cite[VII.7]{Bro:CoG:82}). By Shapiro's Lemma \ref{l:shapiro}, we obtain
    \[
        E^1_{p,q}\cong \bigoplus_{\sigma\in\Sigma_p} H_q(G_\sigma).
    \]
    Since in particular we have $E^1_{p,0} \cong \bigoplus_{\sigma\in\Sigma_p} \Z$ and the differential on the first page is induced by the cellular differential on $C_*(X)$, we obtain the desired result on $d^1$.
\end{Proof}

\begin{Corollary}\label{cor:rational_equiv_homology}
    If the stabilisers $G_\sigma$ are finite for all cells $\sigma\in X$, we have \[H_*(G;\Q)\cong H_*(G\backslash X;\Q).\]
\end{Corollary}

\begin{Proof}
    This follows easily from Theorem \ref{th:equiv_hom} and the fact that rational homology of finite groups vanishes except for $H_0$ by Proposition \ref{prop:rational_homology_of_finite_groups}.
\end{Proof}

\section{The Lyndon\slash Hochschild-Serre spectral sequence}

The Lyndon\slash Hochschild-Serre spectral sequence is a well-known tool in group homology. Nevertheless, we give a short summary of its construction in \cite[VII.6]{Bro:CoG:82} since we will need this in Chapter \ref{ch:stability}. Consider an exact sequence of groups
\[
    1\rightarrow H \rightarrow G \rightarrow Q \rightarrow 1.
\]
Write $F_*(G)$ for the standard resolution of $\Z$ over $\Z G$, this is a free resolution of $\Z$ over $\Z H$ as well. Let $M$ be a $G$-module. It is not difficult to see that
\[
    F_*(G) \otimes_G M \cong ( F_*(G) \otimes_H M)\otimes_Q \Z.
\]
Now set $C_* = (F_*(G) \otimes_H M)$ and consider the standard resolution $F_*(Q)$ of $\Z$ over $\Z Q$. By Corollary \ref{cor:two_spectral_sequences_for_tensor_product} applied to the tensor product $F(Q) \otimes_Q C$, there are two spectral sequences
\[
\begin{array}{l}
    E^1_{p,q} = H_q(Q; C_p) \\
    L^1_{p,q} = H_q(F_p(Q) \otimes_Q C_* )
\end{array} \Rightarrow H_{p+q}(F(Q)\otimes_Q C).
\]
Since $F_p(Q)$ is a free $Q$-module, the functor $(F_p(Q) \otimes_Q -)$ is exact and we obtain
\[
    L^2_{p,q} = H_p(Q;H_q(C)) = H_p(Q;H_q(H;M)).
\]
Now note that $\Z G \otimes_H M \cong \Z (G/H) \otimes_\Z M = \Z Q \otimes_\Z M$. Hence $H_q(Q; \Z G \otimes_H M)=0$ for $q>0$ by Corollary \ref{cor:induced_modules_are_acyclic} and hence also $H_q(Q,C_p)=0$ for all $p$ and all $q>0$, since $F_p(G)$ is a free $\Z G$-module for all $p$. We obtain
\[
E^1_{p,q}=\begin{cases}
    H_0(Q;C_p)\cong(F_p(G)\otimes_H M) \otimes_Q \Z \cong F_p(G)\otimes_G M\quad & \text{for }q=0 \\ 0 & \text{otherwise.}
\end{cases}
\]
The spectral sequence collapses hence on the second page and we have $E^2_{p,0}=H_p(G;M)$. This proves

\begin{BreakTheorem}[Lyndon\slash Hochschild-Serre]\label{th:lhs}
    Given a short exact sequence of groups
    \[
    1 \rightarrow H \rightarrow G \rightarrow Q \rightarrow 1,
    \]
    there is a convergent spectral sequence
    \[
    L^2_{p,q} = H_p(Q;H_q(H;M)) \Rightarrow H_{p+q}(G;M).
	\]
\end{BreakTheorem}

\chapter{Complexes of groups}\label{ch:complexes_of_groups}

We will define complexes of groups and their underlying geometric objects: small categories without loops. In addition, we consider metric properties of complexes of groups which ensure the existence of universal covers. The standard reference for this topic is \cite[Chapter III.$\cC$]{BH:NPC:99}.

\section{Small categories without loops}

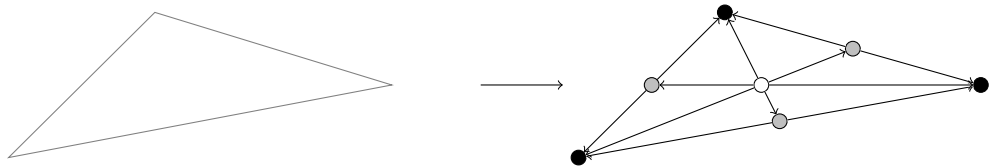
\begin{figure}[hbt]
\centering
    \begin{tikzpicture}[scale=1.25]
                \draw[color=gray] (-6,0,0) -- (-6,0,4) -- (-2.736,0,2) -- cycle;
                \draw[->] (-1.8,0,2) -- (-0.936,0,2);
                \tikzstyle{every node}=[draw,circle,fill=black,scale=0.5]
                \node (v1) at (0,0,0) {};
                \node (v2) at (0,0,4) {};
                \node (v3) at (3.464,0,2) {};
                % Draw middle triangle later
                \tikzstyle{every edge}=[draw]
                \tikzstyle{every node}=[draw,circle,fill=gray!50,scale=0.5]
                \node (l1) at (0,0,2) {} edge[->] (v1) edge[->] (v2);
                \node (l2) at (1.732,0,1) {} edge[color=white,very thick] (v1) edge[color=white,very thick] (v3) edge[->] (v1) edge [->] (v3);
                \node (l3) at (1.732,0,3) {} edge[color=white,very thick] (v2) edge[color=white,very thick] (v3)  edge[->] (v2) edge [->] (v3);
                \tikzstyle{every node}=[draw,circle,fill=white,scale=0.5]
                \node (f1) at (1.154,0,2) {} edge[color=white,very thick] (v1) edge[color=white,very thick] (v2) edge[color=white,very thick] (v3)
                                            edge[color=white,very thick] (l1) edge[color=white,very thick] (l2) edge[color=white,very thick] (l3)
                                            edge[->] (v1) edge[->] (v2) edge[->] (v3) edge[->] (l1) edge[->] (l2) edge[->] (l3);
    \end{tikzpicture}\caption{The scwol associated to a simplex}\label{fig:scwol_assoc_to_simplex}
\end{figure}

The quotient of a simplicial complex by a simplicial (even type-preserving) group action is, in general, no longer a simplicial complex. It is hence more natural to consider polyhedral complexes when studying quotients. It turns out that the right combinatorial object we consider to describe a polyhedral complex is a \emph{small category without loops}, in short \emph{scwol}. Intuitively, it can be thought of as the barycentric subdivision of the polyhedral complex which `remembers' the inclusion relation. Figure \ref{fig:scwol_assoc_to_simplex} can be used to visualise the following definitions.

We will introduce small categories without loops following \cite[III.$\cC$]{BH:NPC:99}, but changing the notation slightly by denoting vertices by $v$, $w$ instead of $\sigma$ and $\tau$.

\begin{Definition}
    A \emph{small category without loops (scwol)} is a set $\cX$ which is the disjoint union of a vertex set $V(\cX)$ and an edge set $E(\cX)$ together with initial and terminal vertex maps $i,t:E(\cX)\rightarrow V(\cX)$. We denote by $E^{(2)}(\cX)$ the set of pairs $(a,b)\in E(\cX)\times E(\cX)$ such that $i(a)=t(b)$ and require the existence of a composition map $E^{(2)}(\cX)\rightarrow E(\cX)$, $(a,b)\mapsto ab$, satisfying:
    \begin{itemize}
        \item For all $(a,b)\in E^{(2)}(\cX)$, we have $i(ab)=i(b)$ and $t(ab)=t(a)$.
        \item For all $a,b,c \in E(\cX)$, if $i(a)=t(b)$ and $i(b)=t(c)$, then $(ab)c=a(bc)$, which can hence be denoted by $abc$. (Associativity)
        \item For each $a\in E(\cX)$, we have $i(a)\neq t(a)$. (No loops condition)
    \end{itemize}
    We write $E^{(k)}(\cX)$ for the $k$-chains of composable edges $(a_1,\ldots,a_k)$ such that $(a_ia_{i+1})\in E^{(2)}(\cX)$. By convention $E^{(0)}(\cX)=V(\cX)$.

    The \emph{dimension of $\cX$} is the supremum of integers $k$ such that $E^{(k)}(\cX)\neq\emptyset$.
\end{Definition}

\begin{Example}
    We can associate a scwol to any partially ordered set $(P,\leq)$, in particular to any simplicial complex, as follows:

    The set of vertices is $P$, and the edges are ordered pairs $(p,q)\in P^2$ such that $p\leq q$, where $i(p,q)=q$, $t(p,q)=p$ and composition is given by $(p,q)(q,r)=(p,r)$.
\end{Example}

\begin{Construction}
    We can also associate a scwol $\cX$ to any $M_\kappa$-polyhedral complex $K$ by forming the partially ordered set of cells of $K$ ordered by inclusion and taking the associated scwol.
\end{Construction}

\noindent Conversely, associated to any scwol $\cX$ there is the \emph{geometric realisation $|\cX|$}, a polyhedral complex whose vertices are in bijection to $V(\cX)$ and whose $k$-simplices correspond to $k$-chains of composable edges in $E(\cX)$. We will now construct this polyhedral complex and we will endow it with a standard Euclidean metric.

\begin{Construction}[III.$\cC$.1.3 in \cite{BH:NPC:99}]
    Let $\Delta^k$ be the standard $k$-simplex, that is, the convex hull of the standard basis vectors in $\R^{k+1}$. For $k>1$ and $i=0,\ldots,k$, let $\partial_i:E^{(k)}(\cX)\rightarrow E^{(k-1)}(\cX)$ be the maps defined by
\begin{align*}
    \partial_0(a_1,\ldots,a_k) &= (a_2,\ldots,a_k) \\
    \partial_i(a_1,\ldots,a_k) &= (a_1,\ldots,a_ia_{i+1},\ldots,a_k) \qquad\text{for }0<i<k\\
    \partial_k(a_1,\ldots,a_k) &= (a_1,\ldots,a_{k-1})
\end{align*}
For $k=1$, we define $\partial_0(a)=i(a)$ and $\partial_1(a)=t(a)$. By the associativity relation, we have $\partial_i\partial_j=\partial_{j-1}\partial_i$ for $i<j$.

Correspondingly, we define maps $d_i:\Delta^{k-1}\rightarrow \Delta^k$ for $i=0,\ldots,k$ via
\[
d_i(t_0,\ldots,t_{k-1}) = (t_0,\ldots,t_{i-1},0,t_i,\ldots,t_{k-1}).
\]
In particular, $d_j(\Delta_{k-1})$ is the face of $\Delta^k$ consisting of points $(t_0,\ldots,t_k)$ such that $t_j=0$. More generally, the face of codimension $r$ given by $t_{i_1}=\cdots=t_{i_r}=0$ for $i_1>\cdots>i_r$ is given by $d_{i_1}\cdots d_{i_r}(\Delta^{k-r})$.

Obviously $d_id_j = d_{j+1}d_i$ for all $i\leq j$.
\end{Construction}

\begin{Definition}\label{def:geometric_realisation_scwol}
    The \emph{geometric realisation of a scwol $\cX$} is given by
    \[
    \lvert\cX\rvert \coloneq \Bigl(\bigcup_{\substack{k\geq 0\\ A\in E^{(k)}(\cX)}}\Delta^k \times \{A\}\Bigr)/\sim,
    \]
    where the equivalence relation $\sim$ is generated by
    \[
            (d_i(x),A) \sim (x,\partial_i(A))
    \]
    for $(x,A)\in \Delta^{k-1}\times E^{(k)}(\cX)$. The geometric realisation is a piecewise Euclidean complex whose $k$-cells are all isometric to $\Delta^k$.
\end{Definition}

\begin{Remark}[III.$\cC$.1.3.3 in \cite{BH:NPC:99}]
    As described above, we have endowed the geometric realisation $\lvert \cX\rvert$ with the intrinsic metric coming from the standard simplices $\Delta^k$. We will later identify each simplex $\Delta^k \times \{A\}$ with a simplex in some $M^n_0$ such that the inclusions $(x,\partial_i A)\mapsto (d_i(x),A)$ are isometries. In this fashion, we obtain a metric on $\lvert\cX\rvert$ as a different piecewise Euclidean complex. If the set of isometry classes of simplices we use is finite, then the topology associated to this new metric is the same as the topology defined above.
\end{Remark}

\noindent It can be shown that the geometric realisation of the scwol associated to a piecewise Euclidean complex is its barycentric subdivision, a concept we have not defined for piecewise Euclidean complexes but of which the intuition should be clear.

\begin{Definition}
    An \emph{automorphism of a scwol $\cX$} is a map $\cX\rightarrow \cX$ that maps $V(\cX)$ to $V(\cX)$ and $E(\cX)$ to $E(\cX)$ and that commutes with $i$, $t$ and the composition.

    An \emph{action} of a group $G$ on a scwol $\cX$ is a homomorphism from $G$ to the group of automorphisms of $\cX$ such that
    \begin{itemize}
        \item for all $a\in E(\cX)$ and $g\in G$, we have $g(i(a)) \neq t(a)$ and
        \item for all $g\in G$ and $a\in E(\cX)$, if $g(i(a))=i(a)$ then we have $g(a)=a$.
    \end{itemize}
    These conditions ensure that the \emph{quotient scwol $G\backslash \cX=\cY$} is well-defined, where $V(\cY)=G\backslash V(\cX)$ and $E(\cY)=G\backslash E(\cX)$ and where composition and the initial and terminal vertex maps factorise naturally.
\end{Definition}

\noindent We also need the definition of the lower and upper links.

\begin{Definition}
	Let $v\in V(\cX)$ be a vertex in a scwol $\cX$. The \emph{lower link} $\Lk^v(\cX)$ is given by
\begin{align*}
	V(\Lk^v(\cX)) &= \{ a\in E(\cX) : i(a)=v \}, \\
	E(\Lk^v(\cX)) &= \{ (a,b) \in E^{(2)}(\cX) : i(b)=v \},
\end{align*}
	with $i(a,b)\coloneq b$, $t(a,b)\coloneq ab$. Since all scwols in this thesis are at most two-dimensional, the lower link is at most one-dimensional and we do not define the composition of edges.
\end{Definition}

\begin{Definition}
    Let $v\in V(\cX)$ be a vertex in a scwol $\cX$. The \emph{upper link} $\Lk_v(\cX)$ is given by
\begin{align*}
	V(\Lk_v(\cX)) &= \{ a\in E(\cX) : t(a)=v \}, \\
	E(\Lk_v(\cX)) &= \{ (a,b) \in E^{(2)}(\cX) : t(a)=v \},
\end{align*}
    with $i(a,b)\coloneq ab$, $t(a,b)\coloneq a$. Again, we do not define composition.
\end{Definition}

\noindent To construct local developments, we need to define the join of scwols, which is in complete analogy to the join of simplicial complexes.

\begin{Definition}
    Let $\cX$ and $\cX'$ be two scwols. The \emph{join} $\cX * \cX'$ of $\cX$ and $\cX'$ is the scwol with vertex set
    \[
        V(\cX * \cX') = V(\cX) \sqcup V(\cX')
    \]
    and with edge set
    \[
        E(\cX * \cX') = E(\cX) \sqcup E(\cX') \sqcup \{ v * v' : v\in V(\cX), v'\in V'(\cX)\}.
    \]
    The maps $i,t$ are given by the respective maps on $E(\cX)$ and $E(\cX')$ and by $i(v * v')= v'$, $t(v * v') = v$. The new compositions are given by
    \[
        a ( i(a) * v') = (t(a) * v'),\qquad (v * t(a') ) a' = (v * i(a')).
    \]
    Note that the join is asymmetric: in general, we have $\cX * \cX'\not\cong \cX' * \cX$.
\end{Definition}

\begin{Definition}\label{def:scwol_closed_star}
    For any scwol $\cX$ and for any vertex $v\in V(\cX)$, we define
    \[
        \cX(v) \coloneq \Lk^v(\cX) * \{v\} * \Lk_v(\cX).
    \]
\end{Definition}

\noindent This scwol can be thought of as the closed star around $v$. However, in general, the canonical map $\cX(v)\rightarrow \cX$ is not injective. It is injective if and only if two distinct edges with terminal or initial vertex $v$ have different initial or terminal vertices, respectively. In any case, there is an affine isomorphism of the open star of $v$ in $|\cX(v)|$ onto the open star of $v\in|\cX|$. See \cite[III.$\cC$.1.17]{BH:NPC:99} for a more detailed discussion.

\section{Complexes of groups}\label{sec:complexes_of_groups}

Suppose a group $G$ acts on an $M_\kappa$-polyhedral complex $X$ by isometries preserving the polyhedral structure. We would like to endow the quotient polyhedral complex with additional structure which allows us to reconstruct the group $G$ and the complex $X$. If we pass to the associated scwol, the additional structure is exactly a complex of groups.

\begin{Definition}
	Let $\cY$ be a scwol. A \emph{complex of groups} $G(\cY)=(G_v,\psi_a,g_{a,b})$ over $\cY$ is given by the following data:
	\begin{enumerate}
		\item A `local group' $G_v$ for each vertex $v\in V(\cY)$.
		\item A monomorphism $\psi_a:G_{i(a)}\rightarrow G_{t(a)}$ for each edge $a\in E(\cY)$.
		\item A `twisting element' $g_{a,b}\in G_{t(a)}$ for each pair of composable edges $(a,b)\in E^{(2)}(\cY)$.
	\end{enumerate}
	We impose the following compatibility conditions:
\begin{itemize}
    \item $\Ad(g_{a,b})(\psi_{ab}) = \psi_a\psi_b$, where $\Ad(g_{a,b})$ is the conjugation by $g_{a,b}$ in $G_{t(a)}$.
    \item For each triple $(a,b,c)\in E^{(3)}(\cY)$, we have the \emph{cocycle condition}
        \[
        \psi(g_{b,c})g_{a,bc} = g_{a,b}g_{ab,c}.
        \]
\end{itemize}
    A complex of groups is called \emph{simple} if all twisting elements are trivial.
\end{Definition}

\begin{Remark}
    For all the complexes of groups we construct explicitly in this thesis, the underlying scwol is at most two-dimensional. In particular, the second compatibility condition will always be fulfilled. The first condition will also be trivially fulfilled, since all vertex groups whose vertices are initial vertices of a pair of composable edges $(a,b)$ will be trivial.
\end{Remark}

\noindent If a group $G$ acts on a scwol $\cX$, there is a natural way to construct a quotient complex of groups $G\backslash\!\backslash \cX$ over the quotient scwol $G\backslash\cX$, see also \cite[III.$\cC$.2.9]{BH:NPC:99}.

\begin{ConstructionN}\label{con:quotient_scwol}
    Let $G$ be a group acting on a scwol $\cX$. Let $\cY = G\backslash \cX$ be the quotient scwol and let $p:\cX\rightarrow \cY$ be the associated projection. We will construct a complex of groups with underlying scwol $\cY$, the \emph{quotient complex of groups $G\backslash\!\backslash \cX$}.

    For the local groups, we do the following: For each vertex $v\in V(\cY)$, we choose a preimage vertex $\bar v\in V(\cX)$ such that $p(\bar v)=v$. We set $G_v\coloneq G_{\bar v}$ to be the stabiliser of the vertex $\bar v$.

    Now by the definition of a group action on a scwol, for each edge $a\in E(\cY)$ with $i(a)=v$ there is a unique edge $\bar a\in E(\cX)$ such that $p(\bar a)=a$ and $i(\bar a) = \bar v$. But, if $u=t(a)$, then in general $\bar u \neq t(\bar a)$ and we choose $h_a\in G$ such that $h_a(t(\bar a))=\bar u$. See Figure \ref{fig:quotient_complex_of_groups} for a diagram illustrating these choices.

    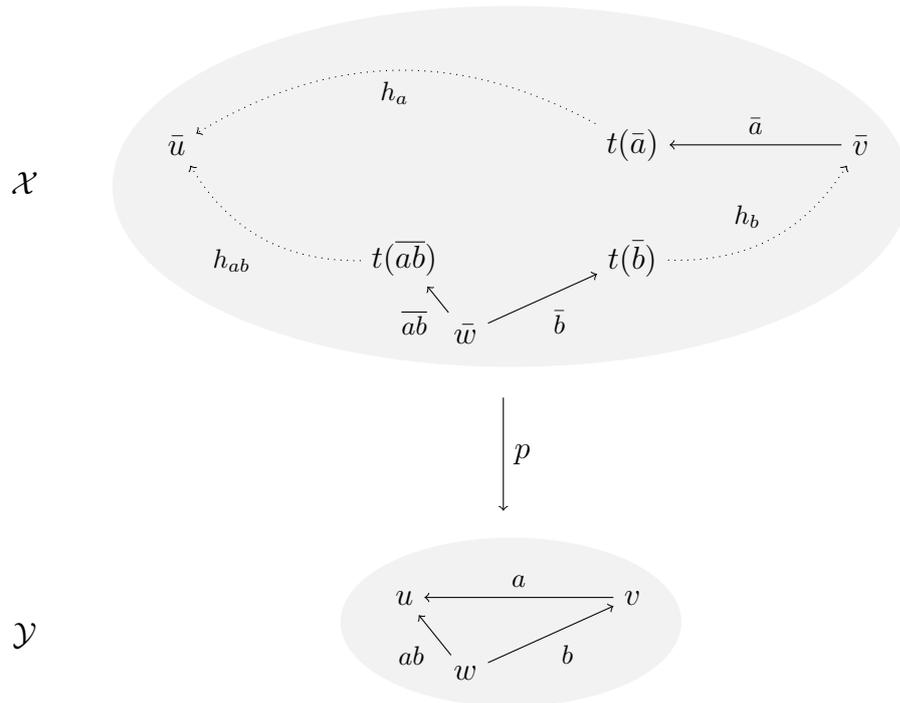
\begin{figure}[htb]
        \centering
            \begin{tikzpicture}[scale=1.5]
                \draw[draw=white,fill=lightgray!20] (1.6,0.45,1.732) ellipse (1.5 and .75);
                \draw[draw=white,fill=lightgray!20] (1.60,4.3,1.732) ellipse (3.5 and 1.60);
                \node (Y) at (-3,0,.866) {$\cY$};
                \node (X) at (-3,4,.866) {$\cX$};

                \node (u) at (0,0,0) {$u$};
                \node (v) at (2,0,0) {$v$};
                \node (w) at (1.2,0,1.732) {$w$};

                \draw[->] (w) -- node[auto,swap,font=\footnotesize] {$b$}  (v);
                \draw[->] (w) -- node[auto,font=\footnotesize] {$ab$} (u);
                \draw[->] (v) -- node[auto,swap,font=\footnotesize] {$a$}  (u);

                \draw[->] (1.2,2.1,0.866) -- node[right] {$p$} ++(0,-1,0);

                \node (wb) at (1.2,3,1.732) {$\bar w$};
                \node (tb) at (2,3,0) {$t(\bar b)$};
                \node (tab) at (0,3,0) {$t(\overline{ab})$};

                \draw[->] (wb) -- node[auto,swap,font=\footnotesize] {$\bar b$} (tb);
                \draw[->] (wb) -- node[auto,font=\footnotesize] {$\overline{ab}$} (tab);

                \node (vb) at (4,4,0) {$\bar v$};
                \node (ta) at (2,4,0) {$t(\bar a)$};

                \draw[->] (vb) -- node[auto,swap,font=\footnotesize] {$\bar a$} (ta);

                \draw[dotted,->] (tb) edge[bend right=30] node[auto,font=\footnotesize] {$h_b$} (vb);

                \node (wu) at(-2,4,0) {$\bar u$};
                \draw[dotted,->] (tab) edge[bend left=30] node[auto,font=\footnotesize] {$h_{ab}$} (wu);
                \draw[dotted,->] (ta) edge[bend right=30] node[auto,font=\footnotesize] {$h_a$} (wu);
            \end{tikzpicture}\caption{The construction of a quotient complex of groups}
        \label{fig:quotient_complex_of_groups}
    \end{figure}

    \noindent For the edge $a\in E(\cY)$ we then define the monomorphism $\psi_a:G_{i(a)}\rightarrow G_{t(a)}$ by
    \[
    \psi_a(g) = h_a g h_a^{-1}.
    \]
    For composable edges $(a,b)\in E^{(2)}(\cY)$, define
    \[
    g_{a,b}= h_ah_bh_{ab}^{-1}.
    \]
    It is easy to verify that this is a complex of groups.
\end{ConstructionN}

\begin{Definition}
    A complex of groups arising in this fashion from a group action on a scwol is called \emph{developable}.
\end{Definition}

\begin{Definition}
    A \emph{morphism $\varphi$ from a complex of groups $G(\cY)$ to a group $G$} consists of a homomorphism $\varphi_v : G_v\rightarrow G$ for each $v\in V(\cY)$ and an element $\varphi(a)\in G$ for each $a\in E(\cY)$ such that
    \[
    \varphi_{t(a)}\psi_a = \Ad(\varphi(a))\varphi_{i(a)}\qquad\text{and}\qquad \varphi_{t(a)}(g_{a,b})\varphi(ab)=\varphi(a)\varphi(b).
    \]
    The morphism $\varphi$ is \emph{injective on the local groups} if each homomorphism $\varphi_v$ is injective.
\end{Definition}

\begin{Example}
    Given a group $G$ acting on a scwol $\cX$, we construct the quotient complex of groups $G\backslash\!\backslash\cX$ as above. The natural inclusions $G_v\rightarrow G$ and the map $\varphi(a)=h_a$ induces a morphism from the complex of groups $G\backslash\!\backslash \cX$ to the group $G$ which is injective on the local groups.
\end{Example}

\noindent The notation in the following definition differs from the one in \cite{BH:NPC:99}. In addition, we have simplified the presentation slightly. It is very easy to see that this definition is equivalent to the one in the aforementioned book.

\begin{Definition}\label{def:fundamental_group}
	Let $G(\cY)$ be a complex of groups. We consider the associated graph with vertex set $V(\cY)$ and edge set $E(\cY)$. Choose a maximal tree $T$ in this graph. The \emph{fundamental group $\pi_1(G(\cY),T)$ of $G(\cY)$ with respect to $T$} is generated by the set
	\[
	\coprod_{v\in V(\cY)}G_v \sqcup \{ k_a : a\in E(\cY) \}
	\]
	subject to the relations in the groups $G_v$ and to
	\begin{align*}
		k_a k_b &= g_{a,b} k_{ab}& &\forall (a,b)\in E^{(2)}(\cY)\\
		\psi_a(g) &= k_agk_a^{-1}& &\forall a\in E(\cY),\,\forall g\in G_{i(a)}\\
		k_a &= 1& &\forall a\in T.
	\end{align*}
\end{Definition}

\noindent It can be shown that fundamental groups associated to different maximal trees are isomorphic. Consequently, we will usually omit $T$ and we will simply write $\pi_1(G(\cY))$. The fundamental group is particularly simple to describe in the following situation.

\begin{Remark}[III.$\cC$.3.11(1) in \cite{BH:NPC:99}]
    If $G(\cY)$ is a simple complex of groups and if $\cY$ is simply connected, then all elements $k_a$ are necessarily trivial and the fundamental group $\pi_1(G(\cY))$ is the direct limit (or colimit in the category of groups) of the system of groups consisting of the groups $G_v$ and the morphisms $\psi_a$. For a discussion of direct limits, see \cite[II.12.12]{BH:NPC:99} or \cite[\S 1]{Ser:Tre:80}.
\end{Remark}

\noindent In the case of a group acting on a scwol, we obtain a natural morphism.

\begin{Remark}
    Let a group $G$ act on a scwol $\cX$. If the quotient complex of groups $G\backslash\!\backslash\cX$ is constructed such that there is a maximal tree in $G\backslash \cX$ such that $h_a=1$ for all $a\in T$, the natural morphism $\varphi:G\backslash\!\backslash\cX\rightarrow G$ induces a homomorphism
    \[
        \pi_1(\varphi): \pi_1(G\backslash\!\backslash\cX) \rightarrow G
    \]
    where $\pi_1(\varphi)\vert_{G_v}$ is the inclusion into $G$ and where $\pi_1(\varphi)(k_a)=h_a$. For a detailed discussion of this construction see \cite[III.$\cC$.3.1]{BH:NPC:99}.
\end{Remark}

\begin{BreakTheorem}[Theorems III.$\cC$.2.13 and III.$\cC$.3.13 in \cite{BH:NPC:99}]
    Associated to each developable complex of groups $G(\cY)$, there is a simply connected scwol $\cX$, called the \emph{universal cover} or \emph{development} of $\cY$ on which $\pi_1(G(\cY))$ acts such that $G(\cY)$ is isomorphic to the complex of groups $\pi_1(G(\cY))\backslash\!\backslash \cX$.
\end{BreakTheorem}

\begin{ConstructionN}\label{con:universal_cover} We will require the precise construction of the universal cover $\cX$ later on. Abbreviate $\Gamma\coloneq \pi_1(G(\cY))$. Then write $\varphi_v:G_v\rightarrow \Gamma$ for the canonical homomorphisms, which are injective if and only if $\cY$ is developable. The universal cover $\cX$ is given by
\begin{align*}
	V(\cX)&\coloneq \Bigl\{ \bigl(g\varphi_v(G_v), v\bigr) : g\in \Gamma, v\in V(\cY)\Bigr\}, \\
	E(\cX)&\coloneq \Bigl\{ \bigl(g\varphi_{i(a)}(G_{i(a)}),a\bigr): g\in \Gamma, a\in E(\cY)\Bigr\},
\end{align*}
with initial and terminal vertex maps
\begin{align*}
	i\bigl( \bigl(g\varphi_{i(a)}(G_{i(a)}),a\bigr)\bigr) &\coloneq \bigl(g\varphi_{i(a)}(G_{i(a)}),i(a)\bigr)\\
	t\bigl( \bigl(g\varphi_{i(a)}(G_{i(a)}),a\bigr)\bigr) &\coloneq \bigl(gk_a^{-1}\varphi_{t(a)}(G_{t(a)}),t(a)\bigr).
\end{align*}
The composition is defined as follows
\[
(g\varphi_{i(a)}(G_{i(a)}),a)(h\varphi_{i(b)}(G_{i(b)}),b) = (h\varphi_{i(b)}(G_{i(b)}),ab).
\]
These edges are composable if $a$ and $b$ are composable and if
\[
    g\varphi_{i(a)}(G_{i(a)})=hk_b^{-1}\varphi_{i(a)}(G_{i(a)}).
\]
\end{ConstructionN}

\begin{Remark} If $\cY$ is finite and all vertex groups are finite, then the universal cover $\cX$ is locally finite. In particular, the geometric realisation $\lvert\cX\rvert$ is a locally compact piecewise Euclidean complex.

    The fundamental group then acts cocompactly on the geometric realisation with finite vertex stabilisers, see Section \ref{sec:lattices_in_autom_groups}.
\end{Remark}

\section{Developability and non-positive curvature}

\noindent However, not all complexes of groups are developable. If the polyhedral complex $|\cY|$ is endowed with a locally Euclidean metric then there is a criterion for developability assuming non-positive curvature on the \emph{local development}.

The local development is constructed in analogy to the complex $\cY(v)$, but using a modified upper link.

\begin{Definition}\label{def:upper_link}
	For a vertex $v\in V(\cY)$ in a complex of groups $G(\cY)$ we define the scwol $\Lk_{\tilde v}(\cY)$ as follows:
	\begin{align*}
		V(\Lk_{\tilde v}(\cY)) &= \bigl\{ \bigl(g\psi_a(G_{i(a)}),a\bigr) : a\in E(\cY), t(a)=v, g\psi_a(G_{i(a)})\in G_v/\psi_a(G_{i(a)}) \bigr\}, \\
		E(\Lk_{\tilde v}(\cY)) &= \bigl\{ \bigl(g\psi_{ab}(G_{i(b)}),a,b\bigr) : (a,b) \in E^{(2)}(\cY), t(a)=v, \\
        &\hspace{7.1cm} g\psi_{ab}(G_{i(b)}) \in G_v/\psi_{ab}(G_{i(b)}) \bigr\}.
	\end{align*}
	The maps $i,t : E(\Lk_{\tilde v}(\cY))\rightarrow V(\Lk_{\tilde v}(\cY))$ are defined by
	\begin{align*}
		i\bigl((g\psi_{ab}G_{i(b)},a,b)\bigr) &= \bigl(g\psi_{ab}(G_{i(b)}),ab\bigr),\\
		t\bigl((g\psi_{ab}G_{i(b)},a,b)\bigr) &= \bigl(gg^{-1}_{a,b}\psi_{ab}(G_{i(b)}),a\bigr).
	\end{align*}
    Again, we omit the definition of the composition of edges.
\end{Definition}

\begin{Definition}\label{def:local_development}
	For any vertex $v\in V(\cY)$ of a complex of groups $G(\cY)$, we define the \emph{local development of the complex of groups $G(\cY)$ at the vertex $v$} to be the scwol $\cY(\tilde v)\coloneq \Lk^v(\cY) * \{\tilde v\} * \Lk_{\tilde v}(\cY)$.

    Notation: We write $\St(\tilde v)=|\cY(\tilde v)|$ for the geometric realisation of $\cY(\tilde v)$ and $\st(\tilde v)$ for the open star of $\tilde v$ in $\St(\tilde v)$.
\end{Definition}

There is a natural projection of $\St(\tilde v)$ onto the closed star $\St(v)$ of the vertex $v$. Assume that the geometric realisation $|\cY|$ of $\cY$ is endowed with a locally Euclidean metric as described in the remark after Definition \ref{def:geometric_realisation_scwol}. Then $\St(v)$ is endowed with the induced metric, and we can endow $\St(\tilde v)$ and its subspace $\st(\tilde v)$ with a locally Euclidean metric such that the aforementioned projection is isometric on every simplex.

\begin{Proposition}[III.$\cC$.4.11 in \cite{BH:NPC:99}]\label{prop:local_isometry}
	Assume the complex of groups $G(\cY)$ is developable and let $\cX$ be the universal cover. Let $v\in V(\cY)$ and choose a preimage $\bar v\in V(\cX)$ of $v$. Then there is a $G_v$-equivariant isometry $\st(\tilde v) \rightarrow \st(\bar v)$.
\end{Proposition}

\noindent Hence we have a precise understanding of the local structure of the universal cover, even though the global structure may be complicated.

\begin{Definition}
	If $\lvert\cY\rvert$ can be endowed with a locally Euclidean metric such that for every vertex $v$ of $G(\cY)$ the metric space $\st(\tilde v)$ is CAT(0), we say that the complex of groups is \emph{non-positively curved}.
\end{Definition}

\begin{BreakTheorem}[Bridson-Haefliger, III.$\cC$.4.17 in \cite{BH:NPC:99}]\label{th:developability_of_nonpositively_curved_complexes}
	A non-positively curved complex of groups is developable.
\end{BreakTheorem}

\noindent We will use the following criterion later on.

\begin{Proposition}[III.$\cC$.4.18 in \cite{BH:NPC:99}]\label{prop:geometric_link_nonpos_curved}
	If for every vertex for every vertex $v\in V(\cY)$ the geometric link $\Lk(\tilde v,\st(\tilde v))$ is CAT(1), then the complex of groups is non-positively curved and hence developable.
\end{Proposition}

\section{Presentations by group actions}

If a group $G$ acts on a scwol $\cX$, there is an explicit presentation of the fundamental group $\pi_1(\lvert\cX\rvert)$ involving stabilisers and the fundamental group $\pi_1(G\backslash\!\backslash\cX)$ of the quotient complex of groups.

\begin{Theorem}\label{th:scwol_presentation_by_action}
    Let $\cX$ be a scwol and let $G$ act on $\cX$. Assume that the quotient scwol $G\backslash\cX$ is connected, and construct the quotient complex of groups $G\backslash\!\backslash\cX$. Let $\varphi: G\backslash\!\backslash\cX\rightarrow G$ be the natural morphism, and let
    \[
        \pi_1(\varphi) : \pi_1(G\backslash\!\backslash\cX)\rightarrow G
    \]
    be the induced morphism on the fundamental group. Then
    \begin{align*}
        \pi_0(\lvert\cX\rvert) &\cong G/\im(\pi_1(\varphi)).
        \intertext{If $\cX$ is connected, we have}
        \pi_1(\lvert\cX\rvert) &\cong \ker(\pi_1(\varphi)).
    \end{align*}
\end{Theorem}

\begin{Proof}
    This is the content of \cite[III.$\cC$.3.14]{BH:NPC:99}, using that the development $D(G\backslash\cX,\varphi)$ which we have not defined in this thesis is isomorphic to the scwol $\cX$ by \cite[III.$\cC$.2.13(2)]{BH:NPC:99}.
\end{Proof}

\noindent We can apply this theorem to a simplicial action of a group $G$ on a simplicial complex $X$.

\begin{Theorem}\label{th:simplicial_presentation_by_action}
    Let $X$ be a simplicial complex and let $G$ act on $X$ simplicially. Assume that there is a subcomplex $Y\subset X$ which is a strict fundamental domain, that is, $Y$ meets every orbit of simplices in exactly one simplex. Denote the subgroup of $G$ generated by all stabilisers of simplices in $Y$ by $G_0$. Then we have the following:
    \begin{itemize}
        \item If $Y$ is connected, we have $\pi_0(\lvert X\rvert) \cong G/G_0$.
        \item If $Y$ is also simply connected, we have $\pi_1(\lvert X\rvert) \cong \ker(\varinjlim \{G_s:s\in Y\} \twoheadrightarrow G_0)$, where the direct limit is taken over the system of stabilisers of simplices in $Y$ with the natural inclusions.
    \end{itemize}
\end{Theorem}

\begin{Proof}
    Consider the scwol $\cX$ associated to the simplicial complex $X$. Since $Y$ is a strict fundamental domain, the associated sub-scwol $\cY$ is isomorphic to the quotient scwol $G\backslash \cX$. In addition, the quotient complex of groups $G\backslash\!\backslash\cX$ is a simple complex of groups. In particular, the image of $\pi_1(\varphi):G\backslash\!\backslash\cX\rightarrow G$ is the subgroup $G_0$. Then apply Theorem \ref{th:scwol_presentation_by_action} to obtain the first claim.

    If $Y$ is simply connected, then we can apply the same theorem to a connected component $\cX_0$ of $\cX$ with the $G_0$-action to obtain the description of $\pi_1(\lvert X\rvert)$ as the kernel of 
    \[
    \pi_1(\varphi_0) : G_0\backslash\!\backslash \cX_0 \rightarrow G_0.
    \]
    By the remark after Definition \ref{def:fundamental_group}, the fundamental group $\pi_1(G_0\backslash\!\backslash\cX_0)$ is exactly the direct limit $\varinjlim\{G_s:s\in Y\}$.
\end{Proof}

\chapter{Wagoner complexes}\label{ch:wagoner}

In this chapter, we construct Wagoner complexes $\cW(G)$ associated to groups $G$ of Kac-Moody type. The main features of these complexes are their homotopy groups, which are related to the integral group homology of $G$ if the associated root datum is 2-spherical. The original motivation for the definition of Wagoner complexes was their strong connection to algebraic $K$-theory, which is also detailed below. Finally, we also introduce affine Wagoner complexes associated to groups having a root datum with valuation.

Our results in this chapter are also contained in \cite{Ess:OWC:09}. The reader might want to re-read Section \ref{sec:wagoner_results} for an overview of the results of this chapter.

\section{The definition of Wagoner complexes}\label{sec:definition}
Throughout this chapter, let $(W,S)$ be a Coxeter system of rank at least 2 such that the associated Coxeter diagram has no isolated nodes. Let $\Sigma=\Sigma(W,S)$ be the associated Coxeter complex. Denote by $\Phi(W,S)$ the set of roots of $\Sigma$ and let $G$ be a group with a twin root datum $(G,(U_\alpha)_{\alpha\in\Phi(W,S)},H)$ of type $(W,S)$ as in Definition \ref{def:root_datum}.

We will say that a simplex $s\in\Sigma$ is \emph{co-spherical} if it is not the empty simplex and if the corresponding residue $R_\Sigma(s)$ is of spherical type. We denote the set of co-spherical simplices by $\Sigma^{\mathrm c}$. The basic `building blocks' for Wagoner complexes will be the following groups.

\begin{Definition}
	For any simplex $s\in\Sigma^{\mathrm c}$ we define the group
    \[
		U_s \coloneq \langle U_\alpha : R_\Sigma(s) \subseteq \alpha \in \Phi(W,S)\rangle,
	\]
    where $R_\Sigma(s)$ is the residue of the simplex $s$ in $\Sigma$.
\end{Definition}

\noindent Remember that the parabolic subgroup $P_s$ is the stabiliser subgroup of $s$ in $G$. We have $U_s\leq P_s$.

By Proposition \ref{prop:root_datum_yields_bn_pair}, there is a thick building $X$ on which $G$ acts strongly transitively. We also have a standard apartment $\Sigma_0$ in the building which we identify with $\Sigma$, and a standard chamber $c_0\in\Sigma_0$. If the root datum is spherical, the building is strictly Moufang by Proposition \ref{prop:root_datum_building_is_moufang}. If the root datum is non-spherical, we still have a regular action of root groups at the boundary of roots by Proposition \ref{prop:general_moufang}, which is all we will use in the following. A simple consequence of this regular action is:

\begin{Remark}
        We have that $s\leq s'$ if and only if $R_\Sigma(s)\supseteq R_\Sigma(s')$ which is equivalent to $U_s\leq U_{s'}$.
\end{Remark}

\begin{Remark}
    By Theorem \ref{th:bn_pair_equiv_strongly_transitive}, the building $X$ is isomorphic to the set of standard cosets $\{ gP_s : s\subseteq c_0, g\in G\}$ partially ordered by reverse inclusion. In this chapter, we will identify the building with this partially ordered set.
\end{Remark}

\noindent Generalising the construction by Wagoner in \cite{Wag:BSK:73}, we define:

\begin{Definition}\label{def:wagonercomplex}
    The \emph{Wagoner complex $\cW(G)$} is the flag complex over the set of cosets $\{ gU_s : s\in \Sigma^{\mathrm c}, g\in G\}$ partially ordered by inclusion, which we write as follows:
	\[
    \cW(G)\coloneq \Flag(\{ gU_s : s\in\Sigma^{\mathrm c}, g\in G\}).
	\]
    We call a Wagoner complex \emph{spherical} if the root datum is of spherical type.
\end{Definition}

\section{Simple properties of Wagoner complexes}\label{sec:simple_properties}

Although we make no direct use of this, let us examine some interesting properties of spherical Wagoner complexes for further reference. Most of these results were already observed in \cite{Pet:HDW:04} or \cite{Wag:BSK:73}.

We restrict ourselves to spherical Wagoner complexes in this section, because the corresponding results are more interesting. We will comment on the general situation at the end of this section. For the rest of this section, assume that $\Sigma=\Sigma(W,S)$ is of spherical type. Then we have $\Sigma^{\mathrm c}=\Sigma\setminus\{\emptyset\}$.

\begin{Definition}
	We define the \emph{standard apartment of $\cW(G)$} to be
	\[
    \cA \coloneq \Flag(\{ U_s : s \in \Sigma^{\mathrm c} \}).
	\]
	Left translates of the standard apartment will be called \emph{apartments}.
\end{Definition}

\noindent A spherical Wagoner complex admits a natural projection onto the barycentric subdivision of the building, which is the flag complex over the building.

\paragraph{Projection} Consider the following map of partially ordered sets
\begin{align*}
    \tilde p_X:\{ gU_s: s\in\Sigma^{\mathrm c}, g\in G\} &\rightarrow X \\
    gU_s &\mapsto gs.
\end{align*}
This map is well-defined since $U_s\leq P_s$. In addition, it is easily seen to be order-preserving. In particular, it induces a surjective simplicial map $p_X:\cW(G)\rightarrow X'$ onto the barycentric subdivision of the building $X$.

The projection $p_X$ maps the standard apartment $\cA$ bijectively onto the barycentric subdivision of the standard apartment $(\Sigma_0)'$ which justifies the terminology. Furthermore, the map $p_X$ induces a surjective map on singular homology: by the Solomon-Tits Theorem \ref{th:solomon_tits}, the singular homology of the building $X$ is concentrated in top-dimension and it is generated by the apartments which are all images of apartments in $\cW(G)$.

\paragraph{Right $N$-action} There is a right action of the group $N$ on the vertices of the Wagoner complex $\{ gU_s: g\in G, s\in\Sigma^{\mathrm c}\}$ by right multiplication:
\[
	gU_sn = gnU_{n^{-1}s}.
\]
Since $U_{c}\cap N=\{1\}$ for any chamber $c$, this action is free. It induces a free right $N$-action on $\cW(G)$ and we have $\cA n = n\cA$. As already observed in \cite{Pet:HDW:04}, the projection $|\cW(G)|\rightarrow |\cW(G)/N|$ is a fibration with discrete fibre $N$.
The projection $p_X$ factors through $\cW(G)/N$, since
\[
	\tilde p_X( gnU_{n^{-1}s} ) = gn(n^{-1}s) = gs = \tilde p_X(gU_s).
\]
If we consider a complex analogous to the Wagoner complex but replacing the groups $U_s$ by the proper parabolic subgroups, we obtain another simplicial complex with a right $N$-action. The quotient complex is isomorphic to the barycentric subdivision of the building. The following diagram summarises all these constructions.
\[
\begin{xy}
	\xymatrix{
    \cW(G) = \Flag(\{ gU_s : s\in\Sigma^{\mathrm c}, g\in G\})\ar[d]^{\mod N}\ar[r]\ar[dr]^{p_X} & \Flag(\{ gP_s : s\in\Sigma^{\mathrm c}, g\in G\})\ar[d]^{\mod N} \\
		\cW(G)/N \cong \Flag(\{ gU_s : \emptyset\neq s\subseteq c_0, g\in G\}) \ar[r]& X' \cong \Flag(\{ gP_s : \emptyset\neq s\subseteq c_0, g\in G\})
	}
\end{xy}
\]

\paragraph{The general case} In the general case for non-spherical root data, the concept of apartments and the above projection still make sense. The projection is then onto the flag complex over all co-spherical simplices of $X$, which is the simplicial complex underlying the \emph{Davis realisation} of the building $X$ as in \cite[18.2]{Dav:GTC:08}. The Davis realisation is contractible, so the induced maps on homotopy and homology groups are necessarily trivial.

\section{The Steinberg group}\label{sec:steinberg}

From now on, we will often require the following technical condition. Denote by $\Pi$ the \emph{set of simple roots} of the root system $\Phi(W,S)$ which is given by all roots that contain the standard chamber $c_0$ but not all adjacent chambers. For roots $\alpha,\beta\in\Pi$, we let
\[
X_\alpha \coloneq \langle U_\alpha \cup U_{-\alpha}\rangle,\qquad X_{\alpha,\beta} \coloneq \langle X_\alpha \cup X_\beta \rangle.
\]
Following Caprace in \cite{Cap:O2S:07}, we consider the following condition, which excludes the possibility that $X_{\alpha,\beta}$ modulo its centre is isomorphic to some small finite groups.\let\storedtheequation=\theequation\renewcommand{\theequation}{Co${}^*$}
\begin{equation}\label{condition}
    X_{\alpha,\beta}/Z(X_{\alpha,\beta}) \not\cong B_2(2), G_2(2), G_2(3),{}^2F_4(2)\quad\text{ for all pairs }\{\alpha,\beta\}\subseteq \Pi.
\end{equation}
This condition will be required for most of the results later on. In addition, will will restrict ourselves to Coxeter complexes of $2$-spherical type.\renewcommand{\theequation}{\storedtheequation}

\begin{Definition}
    A Coxeter system is said to be \emph{$2$-spherical} if its diagram does not have any labels $\infty$, that is if all rank two residues are spherical.
\end{Definition}

\noindent The definition of the Steinberg group requires the concept of prenilpotent pairs of roots, root intervals and root group intervals defined in Sections \ref{subsec:roots} and \ref{subsec:root_data}. The following result describes the `shadow' of a root group interval in a group $U_s$.

\begin{Lemma}\label{lem:stabiliser_intersect_interval}
    Let $\{\alpha,\beta\}$ be a prenilpotent, non-nested pair of roots. Let $s\in\alpha\cap\beta$ be a co-spherical simplex. Then one of the following holds.
    \begin{itemize}
        \item There is a prenilpotent, non-nested pair of roots $\{\gamma,\delta\}$ such that $R_\Sigma(s)\subset \gamma\cap\delta$ and such that $U_s\cap U_{[\alpha,\beta]}\subseteq U_{[\gamma,\delta]}$.
        \item There is a root $\gamma$ containing $R_\Sigma(s)$ such that $U_s\cap U_{[\alpha,\beta]}\subseteq U_\gamma$.
    \end{itemize}
\end{Lemma}

\begin{Proof}
    The set of chambers in $R_\Sigma(s)$ is finite, and we denote the number of chambers in $R_\Sigma(s)\setminus (\alpha\cap\beta)$ by $n$. If $n=0$, then $\{\alpha,\beta\}$ is already a pair of roots we are looking for. Otherwise, we will replace $U_{[\alpha,\beta]}$ iteratively by smaller groups as follows, always decreasing $n$ strictly.

    If $n>0$, then $R_\Sigma(s)\not\subset \alpha$ or $R_\Sigma(s)\not\subset\beta$. We fix a chamber $c\in R_\Sigma(s)\cap\alpha\cap\beta$, and a chamber $d\in -\alpha\cap -\beta$. Denote by $e$ the projection of $d$ onto $R_\Sigma(s)$ as in Proposition \ref{prop:projections}. Since the roots $-\alpha$ and $-\beta$ are convex and $d\in-\alpha\cap-\beta$, we have $e\in-\alpha\cap R_\Sigma(s)$ or $e\in-\beta\cap R_\Sigma(s)$. By Lemma \ref{lem:root_intervals_order}, the minimal gallery from $c$ to $d$ through $e$ crosses one of the walls $\partial\alpha$ or $\partial\beta$ first among the walls of roots in $[\alpha,\beta]$. Assume that the gallery crosses the wall $\partial\alpha$ first. Then there is a chamber $c'\in-\alpha\cap R_\Sigma(s)$ such that $c'\in\gamma$ for all $\gamma\in(\alpha,\beta]$.

    By Lemma \ref{lem:root_group_interval_as_cyclic_product}, we have $U_{[\alpha,\beta]}=U_\alpha U_{(\alpha,\beta]}$. Then
    \[
    U_s \cap U_{[\alpha,\beta]} \subseteq U_s \cap U_{(\alpha,\beta]}
    \]
    since every non-trivial element of $U_\alpha$ cannot fix $c'$ by Proposition \ref{prop:general_moufang}, while all elements of $U_{(\alpha,\beta]}$ and of $U_s$ fix $c'$. Now, by Lemma \ref{l:root_intervals}, either there is a root $\alpha'$ such that $U_{(\alpha,\beta]}=U_{[\alpha',\beta]}$ or $(\alpha,\beta] = \{\beta\}$.

    In the first case, the number of chambers $n'$ in $R_\Sigma(s)\setminus(\alpha'\cap\beta)$ is strictly smaller than the number of chambers in $R_\Sigma(s)\setminus(\alpha\cap\beta)$, since the first one is a subset of the second and does not contain $c'$. If $n'=0$, we are finished, otherwise we repeat the procedure.

    In the second case, we have $U_{(\alpha,\beta]}=U_{\beta}$. Then either $R_\Sigma(s)\in\beta$ and we are done, or $U_{[\alpha,\beta]}\cap U_s \subseteq U_s \cap U_{\beta}=\{1\}$, in which case we can choose an arbitrary root $\gamma$ containing $R_\Sigma(s)$.
\end{Proof}

\noindent We will now introduce the Steinberg group following Tits in \cite{Tit:Uni:87}. The definition uses the concept of the direct limit of a system of groups as in \cite[\S 1]{Ser:Tre:80}. Note that this is nothing but the colimit in the category of groups, or alternatively, the fundamental group of an appropriate simple complex of groups as in Definition \ref{def:fundamental_group}.

\begin{Definition}
    The \emph{Steinberg group $\hat G$} associated to the group $G$ is the direct limit of the system of groups $U_\alpha$ and $U_{[\alpha,\beta]}$ with their natural inclusions for all roots $\alpha$ and all prenilpotent pairs of roots $\{\alpha,\beta\}$.
\end{Definition}

\noindent There is an improved characterisation of the Steinberg group due to Caprace, where we can restrict ourselves to non-nested pairs of roots.

\begin{BreakTheorem}[Proposition 3.9 in \cite{Cap:O2S:07}]
    Assume that the group $G$ is of $2$-spherical type and satisfies condition \eqref{condition}. Then the Steinberg group $\hat G$ is isomorphic to the direct limit of the system of groups $U_\alpha$ and $U_{[\alpha,\beta]}$ with their natural inclusions for all roots $\alpha$ and all prenilpotent, non-nested pairs of roots $\{\alpha,\beta\}$.
\end{BreakTheorem}

\noindent There is a similar characterisation for the groups $U_c$.

\begin{BreakTheorem}[Theorem 3.6 in \cite{Cap:O2S:07}]\label{th:chamber_case}
    Assume that the group $G$ is of $2$-spherical type and satisfies condition \eqref{condition}. For any chamber $c\in\Sigma$, the group $U_c$ is the direct limit of the system of groups $U_\alpha$ and $U_{[\alpha,\beta]}$ for all roots $\alpha$ containing $c$ and all prenilpotent, non-nested pairs of roots $\{\alpha,\beta\}$ containing $c$, respectively.
\end{BreakTheorem}

\noindent We need this statement also for the groups $U_s$.

\begin{Lemma}\label{l:g_sinjection}
    Assume again that $G$ is of $2$-spherical type and satisfies condition \eqref{condition}. For any simplex $s\in\Sigma^{\mathrm c}$, the group $U_s$ is the direct limit of the system of groups $U_\alpha$ and $U_{[\alpha,\beta]}$ for all roots $\alpha$ where $R_\Sigma(s)\subset\alpha$ and all prenilpotent, non-nested pairs of roots $\{\alpha,\beta\}$ with $R_\Sigma(s)\subset(\alpha\cap\beta)$, respectively.
\end{Lemma}

\begin{Proof}
    If $s=c$ is a chamber, this is Theorem \ref{th:chamber_case}. Now let $s\in\Sigma^{\mathrm c}$ be a simplex and fix a chamber $c$ that contains $s$. Consider the following systems of groups
    \begin{align*}
        \cU_c &\coloneq \{ U_\alpha : c\in\alpha\} \cup \{ U_{[\alpha,\beta]} : c \in \alpha\cap\beta,\, (\alpha,\beta) \text{ non-nested}\} \\
        \cU_s &\coloneq \{ U_\alpha : R_\Sigma(s)\subset\alpha\} \cup \{ U_{[\alpha,\beta]} : R_\Sigma(s) \subset \alpha\cap\beta,\, (\alpha,\beta) \text{ non-nested}\}
    \end{align*}
    with the canonical inclusions. Since $\cU_s\subseteq \cU_c$, by the universal property of $\varinjlim\cU_s$, we obtain a homomorphism
    \[
        \varphi:\varinjlim\cU_s \rightarrow \varinjlim\cU_c
    \]
    which is the identity on all groups $U\in \cU_s$. We have the following diagram:
    \begin{center}
    \begin{tikzpicture}
        \matrix (m) [matrix of math nodes, row sep=3em, column sep=2.5em, text height=1.5ex, text depth=0.25ex]
        { \varinjlim \cU_c & U_c \\
        \varinjlim \cU_s  & U_s \\
        };
        \path[->,shorten >=0.5ex] (m-2-1) edge node[auto,font=\scriptsize] {$\varphi$} (m-1-1);
        \path[->>] (m-2-1) edge (m-2-2);
        \path[->] (m-1-1) edge node[auto,font=\scriptsize,inner sep=1pt] {$\cong$} (m-1-2);
        \path[->] (m-2-1) edge (m-1-2);
        \path[right hook->] (m-2-2) edge (m-1-2);
    \end{tikzpicture}
    \end{center}
    where the bottom row is surjective since $U_s=\langle U: U\in\cU_s\rangle$. It remains to show that $\varphi$ is injective. Then the diagonal map is an isomorphism onto the image of $U_s$ in $U_c$.

    It is well known that the direct limit of a system of groups can be constructed explicitly as the free product of the involved groups modulo the relations coming from the inclusions in the system of groups, see \cite[\S 1]{Ser:Tre:80}. Assume that $\varphi$ is not injective, and fix $1\neq g\in\varinjlim \cU_s$ with $\varphi(g)=1$.

    By the above description, the element $g$ can be written as a word $g=g_1\cdots g_n$ with $g_i\in U_i\in\cU_s$ and with $n>0$ minimal.

    If $n=1$, then $g\in U_1\in \cU_s\subset \cU_c$. But then $\varphi(g)=1$ implies $g=1$, since $\varphi$ is the identity on $U_1$.

    If $n>1$, then $\varphi(g) = \varphi(g_1)\cdots \varphi(g_n)=1$. By the description of the direct limit $\varinjlim \cU_s$, there is $1\leq i <n$ such that $\varphi(g_i)$ and $\varphi(g_{i+1})$ are contained in a common group $V \in \cU_c$, since the word $\varphi(g_1)\cdots\varphi(g_n)$ can be reduced to $1$.

   The group $V$ is either a root group or a root group interval. In both cases, there is a prenilpotent, non-nested pair of roots $\{\alpha,\beta\}$ with $c\in\alpha\cap\beta$ such that $\varphi(g_i)$ and $\varphi(g_{i+1})$ are contained in $U_{[\alpha,\beta]}$.

    But observe that $\varphi(g_i)$ and $\varphi(g_{i+1})$ are contained in $U_s$, so by Lemma \ref{lem:stabiliser_intersect_interval}, there is a group $U\in \cU_s$ such that
    \[
    \{\varphi(g_i), \varphi(g_{i+1})\} \subset U_{[\alpha,\beta]}\cap U_s \subseteq U
    \]
    Hence the product $g_ig_{i+1}$ can already be formed in the group $\varinjlim\cU_s$ and
    \[
    g_1\cdots g_{i-1}(g_ig_{i+1})g_{i+2}\cdots g_n
    \]
    is a representation of $g$ as a word of shorter length in $\varinjlim \cU_s$, which contradicts the minimality of $n$. So $\varphi$ is injective.
\end{Proof}

\noindent All of these characterisations allow us to give a different description of the Steinberg group.

\begin{Proposition}\label{prop:tilde_g_steinberg}
	Assume that $G$ is of $2$-spherical type and satisfies condition \eqref{condition}. The direct limit $\tilde G$ of the system of groups $U_s$ for simplices $s\in\Sigma^{\mathrm c}$ with respect to their natural inclusions is isomorphic to the Steinberg group $\hat G$.
\end{Proposition}

\begin{Proof}
	The proof will be divided into two steps. In the first step, we construct a map $\hat G\rightarrow \tilde G$, the inverse map is constructed in the second step.
	\subparagraph{Step 1:} Of course, $U_\alpha\leq U_s$ for any co-spherical simplex $s$ such that $R_\Sigma(s)\subset\alpha$. So we have an inclusion $U_\alpha \hookrightarrow U_s$ for any such simplex $s$. Roots are connected by Lemma \ref{l:roots_are_convex} and chambers and panels are always co-spherical. So for any two co-spherical simplices $s,s'$ with $R_\Sigma(s),R_\Sigma(s') \subset \alpha$, there is a sequence of co-spherical simplices $s_1,\ldots,s_n$ such that
	\[
	U_s \alter{\subseteq}{\supseteq} U_{s_1} \alter{\subseteq}{\supseteq} \cdots \alter{\subseteq}{\supseteq} U_{s_n} \alter{\subseteq}{\supseteq} U_{s'},
	\]
	and the group $U_\alpha$ is contained in all of these. Hence the inclusion $U_\alpha\hookrightarrow \tilde G$, induced by any inclusion $U_\alpha\hookrightarrow U_s$ for any co-spherical simplex $s$ with $R_\Sigma(s)\subset \alpha$, does not depend on the choice of $s$. A similar argument works for the inclusions $U_{[\alpha,\beta]}\hookrightarrow U_s$ for any co-spherical simplex $s$ such that $R_\Sigma(s)\subset (\alpha \cap \beta)$.
	By the universal property of the colimit, there is hence a canonical homomorphism $\hat G \rightarrow \tilde G$ induced by these inclusions.

	\subparagraph{Step 2:} By Lemma \ref{l:g_sinjection}, every group $U_s$ can be written as a direct limit of a subsystem of the direct system for $\hat G$, hence there is a canonical homomorphism $U_s\rightarrow \hat G$. Since these homomorphisms are all induced by the inclusions of root groups, they are compatible with the direct system for $\tilde G$. So there is a canonical epimorphism $\tilde G \rightarrow \hat G$.

	A closer inspection yields that the composition of these maps is by construction trivial on the root groups $U_\alpha$, hence they invert each other.
\end{Proof}

\noindent The importance of the Steinberg group for this thesis comes from the following theorem by Caprace, for which we require some additional notation. For every root $\alpha$ write $X_\alpha\coloneq \langle U_\alpha \cup U_{-\alpha}\rangle$ and $H_\alpha=N_{X_\alpha}(U_\alpha)\cap N_{X_\alpha}(U_{-\alpha})$. For every pair of roots $\{\alpha,\beta\}\subset\Pi$ let $X_{\alpha,\beta}\coloneq \langle X_\alpha \cup X_\beta\rangle$ and let $\Phi_{\alpha,\beta}$ be the rank two root system corresponding to the Coxeter group $W_{\alpha,\beta}=\langle s_\alpha,s_\beta\rangle$, where $s_\alpha$ and $s_\beta$ are the reflections at the walls $\partial\alpha$ and $\partial\beta$. Set $H_{\alpha,\beta}=\bigcap_{\gamma\in\Phi_{\alpha,\beta}} N_{X_{\alpha,\beta}}(U_\gamma)$. Then $(X_{\alpha,\beta},(U_\gamma)_{\gamma\in\Phi_{\alpha,\beta}},H_{\alpha,\beta})$ is a root datum of rank two.

\begin{Theorem}[Caprace, Theorem 3.11 in \cite{Cap:O2S:07}]\label{th:caprace}
    Assume that the group $G$ is of $2$-spherical type and that its diagram has no isolated nodes. Suppose also that
	\begin{enumerate}
		\item For each $\alpha\in\Pi$, we have $[H_\alpha,U_\alpha]=U_\alpha$.
		\item For each pair of roots $\alpha,\beta$ in $\Pi$ such that the corresponding Coxeter group elements $s_\alpha$ and $s_\beta$ do not commute, the Steinberg group $\hat X_{\alpha,\beta}$ is the universal central extension of $X_{\alpha,\beta}$.
	\end{enumerate}
	Then $\hat G$ is the universal central extension of the little projective group $G^\dagger$.
\end{Theorem}

\noindent Condition \eqref{condition} is implied by the first condition, see the proof of Theorem 3.11 in \cite{Cap:O2S:07}. Furthermore, note that the second condition is essentially a condition on Moufang polygons. There is a complete classification of Moufang polygons satisfying this condition, see \cite{dMT:CE2:07}, stating that it is almost always valid except in the case of small fields. Using \cite[33.10-33.17]{TW:MP:02}, one can calculate that the first condition is also always fulfilled except for small fields.

\section{Homotopy groups of Wagoner complexes}\label{sec:homotopy}

\noindent By using the action of the group $G$ on the Wagoner complex $\cW(G)$, we can calculate the number of connected components and the fundamental group of its geometric realisation using Theorem \ref{th:simplicial_presentation_by_action}.

\begin{Theorem}\label{th:wag_main_result} Assume that $G$ is a group with a root datum of $2$-spherical type satisfying condition \eqref{condition}. Then we have
	\begin{align*}
        \pi_0(|\cW(G)|) &\cong G/G^\dagger, \\
        \intertext{If $G$ is of non-spherical type or of rank at least 3, we have}
		\pi_1(|\cW(G)|) &\cong \ker(\hat G\twoheadrightarrow G^\dagger).
    \end{align*}
\end{Theorem}

\begin{Proof}
    Consider the action of $G$ on the Wagoner complex $\cW(G)$. The stabiliser of a vertex $(U_s)$ is obviously the group $U_s$ itself, the stabiliser of a flag
    \[
    (U_{s_1}\subset U_{s_2} \subset \cdots \subset U_{s_r})
    \]
    is given by the smallest group $U_{s_1}$. The group $G_0 = \langle U_s : s\in\Sigma^{\mathrm c}\rangle$ is equal to the little projective group $G^\dagger$. The direct limit of all stabilisers of simplices in $\cA$ with the natural inclusions is hence the group $\tilde G$.
    
    A strict fundamental domain for the action on $\cW(G)$ is given by the standard apartment $\cA$. If $G$ is of spherical type of rank $n$, then $\lvert\cA\rvert$ is homotopy equivalent to a $(n-1)$-sphere, and hence connected and simply connected if $n\geq 3$.

    If $G$ is of non-spherical type, then $\lvert\cA\rvert$ is homeomorphic to the Davis realisation of the Coxeter system associated to $G$ as defined in \cite[Chapter 7]{Dav:GTC:08}, which is contractible by Theorem 8.2.13 in \cite{Dav:GTC:08}.
    
    So we can apply Theorem \ref{th:simplicial_presentation_by_action} and use that we have $\tilde G \cong \hat G$ by Proposition \ref{prop:tilde_g_steinberg}.
\end{Proof}

\begin{Corollary}\label{cor:wag_main_result}
    In particular, if the prerequisites of Theorems \ref{th:caprace} and \ref{th:wag_main_result} are satisfied, we have
    \[
        \pi_1(|\cW(G)|) \cong H_2(G^\dagger;\Z).
    \]
\end{Corollary}

\begin{Proof}
	By Theorem \ref{th:schur_multiplier}, the kernel of the universal central extension $\ker(\hat G \rightarrow G^\dagger)$ is isomorphic to the homology group $H_2(G^\dagger;\Z)$.
\end{Proof}

\section{The connection to algebraic \texorpdfstring{$K$}{K}-theory}\label{sec:ktheory}

Wagoner complexes were originally introduced by Wagoner in \cite{Wag:BSK:73} to provide a definition of higher algebraic $K$-theory. He then proved the following theorem, which is a special case of our result only if $R$ is a division ring.

\begin{Theorem}[Wagoner, Proposition 2 in \cite{Wag:BSK:73}]
	Let $R$ be a ring. If $n\geq 2$ then
	\[
    \pi_0(|\cW(\Gl_{n+1}(R))|) \cong \Gl_{n+1}(R)/\E_{n+1}(R) \cong H_1(\Gl_{n+1}(R);\Z) \cong K_1(R).
	\]
	If $n\geq 3$, then
	\[
    \pi_1(|\cW(\Gl_{n+1}(R))|) \cong \ker(\St_{n+1}(R)\twoheadrightarrow \E_{n+1}(R)) \cong H_2(\E_{n+1}(R);\Z) \cong K_2(R).
	\]
\end{Theorem}
In subsequent papers \cite{AKW:HAK:73} and \cite{AKW:HAK:77}, Anderson, Karoubi and Wagoner and in a final paper \cite{Wag:EAK:77} Wagoner proved that Wagoner's definition of higher algebraic $K$-theory coincides with Quillen's definition of $K$-theory in the sense that
\[
K_i(R) \cong \pi_{i-1}(|\cW(\Gl_{n+1}(R))|)\qquad\text{for }n\gg i, i>0.
\]

\noindent This leads naturally to the question whether
\[
\pi_{i-1}(|\cW(G)|) \cong \pi_i(BG^+)
\]
for any group with a spherical root datum with $\rank(G)\gg i$. This question would be interesting for hermitian $K$-theory, for example, by choosing the group $G$ to be a symplectic or orthogonal group. Unfortunately, we do not know whether these groups are isomorphic, in general.

\section{Affine Wagoner complexes}\label{sec:affine}

We also include the construction of Wagoner complexes related to affine buildings, imitating the corresponding construction in \cite{Wag:HTp:75}.

Let $(W,S)$ be a spherical Coxeter system and let $\Phi=\Phi(W,S)$ be its set of roots. Let $G$ be a group with a root datum with discrete valuation $(G,(U_\alpha)_{\alpha\in\Phi},(\varphi_\alpha)_{\alpha\in\Phi},H)$ as in Definition \ref{def:root_datum_valuation}. Let $X$ be the associated affine building and let $\Sigma$ be an apartment of $X$ such that the half-apartments of $\Sigma$ are parametrised as $H_{\alpha,k}$ with walls $\partial H_{\alpha,k}$ for $\alpha\in\Phi$, $k\in\Z$, compare Theorem \ref{th:root_data_with_valuations_have_affine_buildings}. Write $\Sigma^{\mathrm c} = \Sigma\setminus\{\emptyset\}$.

\begin{Definition}
    Let $s$ be a simplex in $\Sigma^{\mathrm c}$ and let $\geq 1$. The \emph{$n$-cell around $s$} is defined to be the set
	\[
    C_n(s) \coloneq \bigcap_{\substack{\alpha\in\Phi,\, k\in\Z \\R_\Sigma(s)\subset H_{\alpha,nk}}} H_{\alpha,nk}.
	\]
\end{Definition}

\begin{Definition}
    For any $n\geq 1$ and any simplex $s\in\Sigma^{\mathrm c}$, we write
	\[
	U^n_s \coloneq \langle U_{\alpha,k} : C_n(s)\subset H_{\alpha,k}, \alpha\in\Phi, k\in\Z \rangle.
	\]
	Analogously to Definition \ref{def:wagonercomplex}, for any $n\geq 1$, we define the \emph{affine Wagoner complex $\cW^n(G)$} to be 
    \[
        \cW^n(G)\coloneq\Flag(\{ gU^n_s : s\in \Sigma^{\mathrm c}, g\in G\}).
    \]
\end{Definition}

\noindent For every fixed $n$, we have by construction $s \subseteq s'$ implying that $U^n_s \leq U^n_{s'}$. There are also simplicial epimorphisms from $\cW^n(G)$ onto the barycentric subdivisions of certain `coarsenings' of the affine building $X$. For buildings of type $\tilde A_n$, this can be found in \cite[\S 5]{Wag:HTp:75}.

\begin{Remark}
    As above, we define the standard apartment of the Wagoner complex by
    \[
    \cA^n = \Flag(\{ gU^n_s : g\in G, s\in \Sigma^{\mathrm c} \}).
    \]
    Obviously, the apartment $\cA^n$ is isomorphic to the flag complex over the set of $n$-cells $\{C_n(s):s\in\Sigma^{\mathrm c}\}$. It is not hard to see that this complex is in turn isomorphic to the complex associated to the affine reflection group generated by the reflections at the walls $\partial H_{\alpha,nk}$ for $\alpha\in\Phi$ and $k\in\Z$. In particular, by rescaling by $\tfrac{1}{n}$, this complex is isomorphic to the Coxeter complex $\Sigma$, so $\cA^n$ is contractible.
\end{Remark}

\noindent In analogy to the construction above, we denote the direct limit of the system of groups $\{U^n_s : s\in\Sigma\}$ with their canonical inclusions by $\tilde G^n$. By the same argumentation as in the proof of Theorem \ref{th:wag_main_result}, using that $\cA^n$ is contractible, we obtain

\begin{Proposition} For any $n\geq 1$, we have
	\[
		\pi_1(|\cW^n(G)|) \cong \ker( \tilde G^n \twoheadrightarrow G^\dagger).
	\]
\end{Proposition}

\begin{Construction} If $n$ divides $m$, then $C_n(s)\subseteq C_m(s)$ and hence $U^n_s \supseteq U^m_s$ for any $s\in\Sigma$. Consequently, the map $p^m_n: \cW^m(G) \rightarrow \cW^n(G)$ induced by $gU^m_s\mapsto gU^n_s$ is well-defined.

We obtain an inverse system of the simplicial complexes $\cW^n(G)$ with the maps $p^m_n$ for $n$ dividing $m$. This induces an inverse system on the fundamental groups $\pi_1(|\cW^n(G)|)$.
\end{Construction}

\begin{Definition}
	The inverse limit of this system is denoted by
	\[
	\pi_1^{\text{aff}}(G) \coloneq \varprojlim_n \pi_1(|\cW^n(G)|)
	\]
	and is called the \emph{affine fundamental group of $G$}.
\end{Definition}

\noindent The following theorem was already proved by Wagoner in \cite[\S 2]{Wag:HTp:75} for $G=\Sl_k(\F)$ over a local field $\F$ in a very different way.

\begin{Theorem}
	Let $\cW(G)$ be the Wagoner complex associated to the spherical building at infinity. Assume that the rank of the spherical building is at least three. Then, for every $n\geq 1$, there is a homomorphism
	\[
		\pi_1(|\cW(G)|)\rightarrow \pi_1(|\cW^n(G)|)
	\]
	that is compatible with the inverse system above. In particular, we obtain a homomorphism
	\[
	\pi_1(|\cW(G)|)\rightarrow \pi_1^{\text{aff}}(G).
	\]
\end{Theorem}

\begin{Proof}
	We will first construct a homomorphism $\hat G\rightarrow \tilde G^n$. By the universal property of the direct limit, it is enough to construct homomorphisms $U_\alpha$, $U_{[\alpha,\beta]}\rightarrow \tilde G^n$ which are compatible with the natural inclusions.

	Note that by the definition of root data with valuations, we have $U_\alpha= \bigcup_{k\in\Z}U_{\alpha,k}$. This is a filtration and hence a direct limit. Again, it is enough to construct a homomorphism $U_{\alpha,k}\rightarrow \tilde G^n$ for all $k\in \Z$. The group $U_{\alpha,k}$ is contained in all groups $U^n_s$ such that $C_n(s)\subseteq H_{\alpha,k}$. There are obviously such simplices, and we consider the inclusions $U_{\alpha,k}\rightarrow U^n_s$. By an argument similar to the one in the proof of Proposition \ref{prop:tilde_g_steinberg}, using that half-apartments are connected and that all homomorphisms in the system of groups are inclusions, these maps are compatible with the system of groups for $\tilde G^n$, and we obtain homomorphisms $U_{\alpha,k}\hookrightarrow \tilde G^n$.

	For the groups $U_{[\alpha,\beta]}$, consider the filtration $U_{[\alpha,\beta],k}=\langle U_{\gamma,k} : \gamma\in[\alpha,\beta]\rangle$. Then the group $U_{[\alpha,\beta],k}$ is contained in all groups $U^n_s$ with $C_n(s)\subseteq \bigcap_{\gamma\in[\alpha,\beta]} H_{\gamma,k}$. This intersection contains a sector and hence such simplices exist. It is also connected and by the same argument as before, the inclusion is compatible with the system of groups and we obtain a homomorphism $U_{[\alpha,\beta]}\rightarrow \tilde G^n$.

	Since all of these maps are induced by inclusions, they are automatically compatible with the system of groups for $\hat G$. Hence we obtain a homomorphism $\hat G\rightarrow \tilde G^n$ which makes the following diagram commute:
	\[
\begin{xy}\xymatrix{
	1 \ar[r] & \pi_1(|\cW(G)|) \ar[r]\ar@.[d] & \hat G \ar[r]\ar[d] & G^\dagger \ar[r]\ar@{=}[d] & 1 \\
	1 \ar[r] & \pi_1(|\cW^n(G)|) \ar[r] & \tilde G^n \ar[r] & G^\dagger \ar[r] & 1.}
\end{xy}
	\]
	For standard reasons, the dotted map completing the above diagram exists. It is the desired homomorphism.
\end{Proof}

\noindent In the case of $G=\Sl_k(\F)$ over a local field $\F$, Wagoner proves in \cite{Wag:HTp:75} that $\pi_1^{\text{aff}}(G)$ is isomorphic to the kernel of the universal central topological extension of $G=G^\dagger=\Sl_k(\F)$, sometimes denoted by $H_2^{\text{top}}(\Sl_k(\F))$. We do not know whether this is true in the general case. Wagoner's proof relies heavily on the fact that $H_2(\Sl_k(\F)) \cong H_2^{\text{top}}(\Sl_k(\F)) \oplus D$, where $D$ is infinitely divisible, and that $H_2^{\text{top}}(\Sl_k(\F))$ is isomorphic to the group of roots of unity of $\F$.

\chapter{Homological stability}\label{ch:stability}

In this chapter, we will prove several homological stability results for classical groups. This will be done by investigating a spectral sequence obtained from the action of the classical group on the associated \emph{opposition complex}. The opposition complex is constructed out of the spherical building using the opposition relation.

We start with a general procedure to obtain a relative spectral sequence associated to a map of two tensor product complexes. As a first application, this can be used to construct a Lyndon\slash Hochschild-Serre spectral sequence for relative homology, which will be useful later on.

In the second section, we define the opposition complex associated to a weak spherical building with a strongly transitive group action of a group $G$. This complex has two properties which make it important for homological stability: First, it is homotopy equivalent to a bouquet of spheres by a result by von Heydebreck and secondly, all stabilisers are Levi subgroups. Levi subgroups can usually be split as direct or semidirect products of smaller groups, which makes them very useful for homological stability proofs.  Using the opposition complex, we construct a relative spectral sequence from the inclusion of a co-rank one subgroup into $G$ involving relative homology of Levi subgroups.

In the last section, we use the explicit structure of special linear and unitary groups and their Levi subgroups to prove homological stability results using the spectral sequence we have constructed.

The results of this chapter are also contained in \cite{Ess:SCG:09}. An overview of this chapter can be found in Section \ref{sec:stability_results}.

\section{Relative spectral sequences}

In this section, we will discuss a procedure to obtain relative spectral sequences from both the spectral sequence associated to a double complex as well as from the Lyndon\slash Hochschild-Serre spectral sequence.

\subsection{Relative spectral sequences associated to tensor complexes}

Fix a group $G$ and a subgroup $G'$. Let $F'$, $C'$ and $F$, $C$ be chain complexes of $G'$- and $G$-modules, respectively. Consider the two tensor product double complexes $(F_p \otimes_G C_q)_{p,q}$ and $(F'_p \otimes_{G'} C'_q)_{p,q}$. Assume that there is map of double complexes
\[
i: F'\otimes_{G'} C' \rightarrow F\otimes_G C.
\]
We denote the induced maps on the vertical and horizontal chain complexes by
    \begin{align*}
        i_{p,\bullet}: F'_p \otimes_{G'} C'_\bullet & \rightarrow F_p \otimes_G C_\bullet \\
        i_{\bullet,q}: F'_\bullet \otimes_{G'} C'_q & \rightarrow F_\bullet \otimes_G C_q.
    \end{align*}
    For every $q$, we then consider the two mapping cone chain complexes $\Cone_*(i_{\bullet,q})$ and $\Cone_*(i_{q,\bullet})$. It is a simple calculation to see that
    \[
    D_{p,q}=\Cone_p(i_{\bullet,q})\qquad\text{and}\qquad D^T_{p,q}=\Cone_p(i_{q,\bullet})
    \]
    are actually double complexes with respect to the cone differentials and the differentials
	\begin{align*}
        \partial: (F'_{p-1} \otimes_{G'} C'_q) \oplus (F_p \otimes_G C_q) &\longrightarrow (F'_{p-1} \otimes_{G'} C'_{q-1}) \oplus (F_p \otimes_G C_{q-1}) \\
        f' \otimes_{G'} c' + f \otimes_G c &\longmapsto f' \otimes_{G'} \partial^{C'}(c') + f\otimes_G\partial^C(c) \\
        \partial^T: (F'_q \otimes_{G'} C'_{p-1}) \oplus (F_q \otimes_G C_p) &\longrightarrow (F'_{q-1} \otimes_{G'} C'_{p-1}) \oplus (F_{q-1} \otimes_G C_p) \\
        f' \otimes_{G'} c' + f \otimes_G c &\longmapsto \partial^{F'}(f') \otimes_{G'} c' + \partial^F (f)\otimes_G c
	\end{align*}
	induced by the differentials $\partial^{C'}$ and $\partial^C$, respectively $\partial^{F'}$ and $\partial^F$.

\begin{Theorem}\label{th:relative_spectral_sequence}
  There are two spectral sequences corresponding to $D$ and $D^T$ satisfying
	\[
    \begin{array}{c} E^1_{p,q} = H_q\bigl(\Cone_*(i_{\bullet,p})\bigr) \\
        L^1_{p,q} = H_q\bigl(\Cone_*(i_{p,\bullet})\bigr) \end{array} \Rightarrow H_{p+q}\bigl(\Cone_*(i)\bigr).
	\]
    The differentials on the first pages are induced by $\pm\partial$ and $\pm\partial^T$, respectively.
\end{Theorem}
\begin{Proof}
    The first spectral sequence is the first one from Proposition \ref{prop:two_spectral_sequences} applied to the double complex $D$. We know that
    \[
    E^1_{p,q} = H_q\bigl(\Cone_*(i_{\bullet,p})\bigr) \Rightarrow H_{p+q}(\Tot(D)).
    \]
    But the formation of the mapping cone and of the total complex commute in the following sense: On the one hand
    \[
    \Tot(D)_k = \bigoplus_{p+q=k} (F'_{p-1} \otimes_{G'} C'_q) \oplus (F_p \otimes_G C_q).
    \]
	On the other hand
	\[
    \Cone_k(i:(F'_*\otimes_{G'} C'_*)\rightarrow (F_* \otimes_G C_*)) = \bigoplus_{p+q=k-1}(F'_p \otimes_{G'} C'_q) \oplus \bigoplus_{p+q=k}(F_p \otimes_G C_q).
	\]
	Hence the modules of the two chain complexes coincide. For the boundary maps, we have
    \begin{multline*}
        \partial^{\Tot(D)}\Bigl(\sum_{p+q=k} (f'_{p-1}\otimes_{G'} c'_q + f_p\otimes_G c_q)\Bigr) =
        \sum_{p+q=k} \bigl( -\partial f'_{p-1}\otimes_{G'} c'_q + i(f'_{p-1}\otimes_{G'} c'_q) + \partial f_p\otimes_G c_q + \\
        + (-1)^p f'_{p-1}\otimes_{G'} \partial c'_q + (-1)^p f_p\otimes_G \partial c_q\bigr)
	\end{multline*}
		and
	\begin{align*}
        \partial^{\Cone(i)}\Bigl(\sum_{p+q=k-1}f'_p\otimes_{G'} c'_q + \sum_{p+q=k}f_p\otimes_G c_q \Bigr) =
        &- \sum_{p+q=k-1} \bigl(\partial f'_p \otimes_{G'} c'_q + (-1)^p f'_p\otimes_{G'}\partial c'_q \bigr) +\\
        &+ i \bigl( \sum_{p+q=k-1} f'_p \otimes_{G'} c'_q \bigr) +\\
		&+ \sum_{p+q=k} \bigl(\partial f_p\otimes_G c_q + (-1)^p f_p\otimes_G\partial c_q\bigl)
    \end{align*}
	which are also easily seen to be identical.

    For the second spectral sequence, we take the first spectral sequence from Proposition \ref{prop:two_spectral_sequences} associated to the double complex $D^T$. We denote it by $L$, however, to make the notation consistent in the following. Note that we can also write $D^T_{p,q}\cong\Cone_p(i_{\bullet,q}^T)$, where $i^T: (F'_q\otimes_{G'}^T C'_p)_{p,q} \rightarrow (F_q\otimes_G^T C_p)_{p,q}$ is the map induced on the transposed double complexes. By applying the above calculation to $i^T$, we obtain
	\[
	L^1_{p,q} = H_q\bigl(\Cone_*(i_{p,\bullet})\bigr) \Rightarrow H_{p+q}(\Tot(D^T)) \cong H_{p+q}(\Cone_*(i^T)).
	\]
    We have seen in Lemma \ref{lem:transposed_double_complex} that there is a chain map inducing isomorphisms
	\[
	H_*(F \otimes_G C) \cong H_*(F \otimes_G^T C)
	\]
    on homology. Using these maps, one can easily construct a chain map $\Cone(i)\rightarrow \Cone(i^T)$ inducing an isomorphism on homology by the long exact sequence associated to the mapping cone from Proposition \ref{prop:les_mapping_cone} and the 5-Lemma \ref{l:5lemma}.
\end{Proof}

\begin{Remark}
  If $i$ is injective, we have
  \[
  H_{p+q}(\Cone(i))\cong H_{p+q}( F\otimes_G C, i(F'\otimes_{G'} C'))
  \]
  by Lemma \ref{l:mapping_cone_relative_homology}.
\end{Remark}

\subsection{The relative Lyndon\slash Hochschild-Serre spectral sequence}

If we apply this procedure to the Lyndon\slash Hochschild-Serre spectral sequence from Theorem \ref{th:lhs}, we obtain

\begin{Theorem}[Relative Lyndon\slash Hochschild-Serre]\label{th:relative_lyndon_hochschild_serre}
	Given a semidirect product $G = H \rtimes Q$, a subgroup $G' = H' \rtimes Q$ of $G$ such that $H'=H\cap G'$ and a $G$-module $M$, there is a spectral sequence
	\[
	L^2_{p,q} = H_p\bigl(Q;H_q(H,H';M)\bigr) \Rightarrow H_{p+q}(G,G';M).
	\]
	If both $H$ and $H'$ act trivially on $M$ and $M$ is $\Z$-free, we obtain
	\[
	L^2_{p,q} = H_p\bigl(Q;H_q(H,H';\Z)\otimes_\Z M\bigr) \Rightarrow H_{p+q}(G,G';M).
	\]
\end{Theorem}

\begin{Proof}
	Let $F_*(G)$, $F_*(G')$ and $F_*(Q)$ be the standard resolutions of $\Z$ over $\Z G$, $\Z G'$ and $\Z Q$, respectively. As in the construction of the Lyndon\slash Hochschild-Serre spectral sequence, we consider the double complexes $F_*(Q) \otimes_Q (F_*(G') \otimes_{H'} M)$ and $F_*(Q) \otimes_Q (F_*(G) \otimes_H M)$. Let
    \[
        i : F_*(Q) \otimes_Q ( F_*(G') \otimes_{H'} M) \rightarrow F_*(Q) \otimes_Q (F_*(G) \otimes_H M)
    \]
    be induced by the inclusions. Then apply Theorem \ref{th:relative_spectral_sequence} to obtain a spectral sequence with first page term
	\[
    L^1_{p,q}=H_q(\Cone_*(i_{p,\bullet}))\cong H_q\biggl(\frac{F_p(Q) \otimes_Q (F_*(G) \otimes_{H} M)}{F_p(Q) \otimes_Q (F_*(G') \otimes_{H'} M)}\biggr).
	\]
	The module $F_p(Q)$ is $\Z Q$-free and hence $\Z Q$-flat. We obtain
	\[
	L^1_{p,q}\cong F_p(Q) \otimes_Q H_q\biggl(\frac{F_*(G)\otimes_{H} M}{F_*(G') \otimes_{H'} M}\biggr) \cong F_p(Q) \otimes_Q H_q\biggl(\frac{F_*(H)\otimes_{H} M}{F_*(H') \otimes_{H'} M}\biggr),
	\]
    where the last isomorphism is the content of Proposition \ref{prop:relative_homology_of_subgroups}. This yields
    \[
    L^2_{p,q}\cong H_p(Q,H_q(H,H';M)) \Rightarrow H_{p+q}(\Cone_*(i)).
    \]

    \noindent On the other hand, consider the spectral sequence $E$ from Theorem \ref{th:relative_spectral_sequence} and apply Lemma \ref{l:mapping_cone_relative_homology} to obtain the following description of the first page.
    \[
    E^1_{p,q} \cong H_q \biggl( \frac{F_*(Q) \otimes_Q ( F_p(G) \otimes_H M)}{F_*(Q) \otimes_Q (F_p(G') \otimes_{H'} M)} \biggr) \Rightarrow H_{p+q}(\Cone_*(i)).
    \]
    In the proof of the regular Lyndon\slash Hochschild-Serre spectral sequence (Theorem \ref{th:lhs}), we have seen that $H_q(F_*(Q) \otimes_Q ( F_p(G) \otimes_H M))=0$ for $q\neq 0$ and that 
    \[
    H_0(F_*(Q) \otimes_Q ( F_p(G) \otimes_H M)) \cong F_p(G)\otimes_G M.
    \]
    The same is true if we replace $G$ and $H$ by $G'$ and $H'$. The map $i$ induces the map
    \[
    H_0(F_*(Q) \otimes_Q ( F_p(G') \otimes_{H'} M)) \rightarrow H_0(F_*(Q) \otimes_Q ( F_p(G) \otimes_H M))
    \]
    which under the above isomorphisms is just the inclusion $F_p(G')\otimes_{G'}M \rightarrow F_p(G)\otimes_G M$. In particular, it is injective. By the long exact sequence, we obtain $E^1_{p,q}=0$ for $q\neq 0$ and
    \[
    E^1_{p,0} \cong \frac{F_p(G) \otimes_G M}{F_p(G')\otimes_{G'} M},
    \]
    so the spectral sequence $E$ collapses on the second page and converges to $H_{p+q}(G,G';M)$, which proves the first statement. The second statement follows from Lemma \ref{l:trivial_coefficient_module}.
\end{Proof}

\section{The opposition complex and a relative spectral sequence}

Homological stability proofs usually consider the action of a group on some simplicial complex and then exhibit smaller groups of the same series of groups as stabiliser subgroups. In this part, we will introduce the opposition complex associated to a group with a weak spherical Tits system and construct a filtration of this complex which leads to a relative spectral sequence involving the group and its Levi subgroups. This will be used in the third part of this chapter to prove homological stability.

\subsection{The opposition complex}

Let $n\geq 2$ and let $X$ be a weak spherical building of rank $n$. Enumerate the type set $I=\{i_1,\ldots,i_n\}$ of $X$ arbitrarily. In addition, we fix a group $G$ that acts strongly transitively on $X$.

The basic geometry on which we will study the action of $G$ is not the building $X$, but its associated opposition complex.

\begin{Definition}
	The \emph{opposition complex $O(X)$} is the simplicial complex consisting of pairs of opposite simplices
	\[
		O(X) \coloneq \{ \sigma= (\sigma^+,\sigma^-) \in X\times X : \sigma^+ \text{ opposite } \sigma^- \}
	\]
	with the induced inclusion relations. Set $\type((\sigma^+,\sigma^-)) \coloneq \type(\sigma^+)$.
\end{Definition}

\begin{Remark}
    The opposition complex $O(X)$ is a $G$-simplicial complex since the $G$-action preserves opposition. The $G$-action is transitive on vertices of the same type of $O(X)$, since it is transitive on pairs of opposite vertices of a fixed type in the building.
\end{Remark}

\begin{Example}
    The vertices of the opposition complex associated to the general linear group $\Gl_{n+1}(D)$ acting on the associated projective space over $D^{n+1}$ as in Section \ref{sec:an_example} are pairs of complementary subspaces of $D^{n+1}$.
\end{Example}

\noindent The significance of the opposition complex for this thesis lies in the following theorem.

\begin{Theorem}[von Heydebreck, Theorem 3.1 in \cite{vH:HPC:03}]\label{th:opposition_complex_is_spherical}
	The opposition complex of a weak spherical building $X$ of rank $n$ is homotopy equivalent to a bouquet of $(n-1)$-spheres.
\end{Theorem}

\noindent We fix a set of representative vertices for the $G$-action.

\begin{Definition}\label{def:situation}
	We fix a standard apartment $\Sigma_0$ in $X$ and a chamber $c_0\in\Sigma_0$. In the following, we write $v_p^+$ for the vertex of $c_0$ of type $\{i_p\}$ and we denote the corresponding opposite vertex in $\Sigma_0$ by $v_p^-$. We write $v_p=(v_p^+,v_p^-)\in O(X)$ for the vertex in $O(X)$. 
\end{Definition}

\noindent For the inductive arguments to come, we investigate the structure of stabilisers.

\begin{Definition}
	Denote the stabilisers as follows:
\[
	L_p \coloneq G_{v_p} = G_{v_p^+} \cap G_{v_p^-}.
\]
These are intersections of two opposite parabolic subgroups, hence Levi subgroups of $G$ as in Definition \ref{def:parabolic_levi}.
\end{Definition}

\begin{Example}
	For the general linear group $\Gl_{n+1}(D)$ acting on the associated building $X$ of type $A_n$ as in Section \ref{sec:an_example}, the stabiliser of the vertex
    \[
        v=(\langle e_1,\ldots,e_k\rangle,\langle e_{k+1},\ldots,e_{n+1}\rangle)
    \] in $O(X)$ is the subgroup
    \[
        \begin{pmatrix}
        \Gl_{k}(D) & 0 \\ 0 & \Gl_{n+1-k}(D)
    \end{pmatrix}.
    \]
    Note that the Levi subgroup splits as a direct product of smaller general linear groups.
\end{Example}

\noindent A link in the opposition complex is isomorphic to the opposition complex of the corresponding link in the building. Compare this to Proposition \ref{prop:links_are_buildings}.

\begin{BreakProposition}[von Heydebreck, Proposition 2.1 in \cite{vH:HPC:03}]\label{prop:links_opposition_complex}
	For a $k$-simplex $(\sigma^+,\sigma^-)\in O(X)$ we have
	\[
		\lk_{O(X)}((\sigma^+,\sigma^-)) \cong O(\lk_X(\sigma^+)).
	\]
    This isomorphism is $G_{\{\sigma^+,\sigma^-\}}$-equivariant. In particular, the link $\lk_{O(X)}((\sigma^+,\sigma^-))$ is homotopy equivalent to a bouquet of $(n-1-k)$-spheres.
\end{BreakProposition}

\subsection{A filtration}

We construct an exact chain complex of $G$-modules associated to $O(X)$. The construction is similar to the construction of cellular chains of a CW complex. The filtration by skeletons is replaced by a filtration by types.

\begin{Definition}
	For $1\leq p \leq n$ let $I_p=\{i_1,\ldots,i_p\}$. Write
	\[
		O(X)_p \coloneq \{ \sigma\in O(X) : \type (\sigma) \subseteq I_p\},
	\]
	this is a $G$-invariant filtration of $O(X)$. We set $O(X)_0\coloneq \emptyset$.
\end{Definition}

\noindent Observe that $O(X)_p$ is of rank $p$ and hence of dimension $p-1$.

\begin{Definition}
    Let $v$ be a vertex in $O(X)_p$. We define the \emph{filtrated residue}, \emph{link} and \emph{star} by:
	\begin{align*}
        R(v)_p &= R_{O(X)}(v) \cap O(X)_p,\\
        \lk(v)_p &= \lk_{O(X)}(v) \cap O(X)_p,\\
		\st(v)_p &= \lk(v)_p \sqcup R(v)_p.
	\end{align*}
\end{Definition}

\begin{Remark}
	From the definition it is obvious that
	\[
	\st(v)_p \cap O(X)_{p-1} = \lk(v)_p = \lk(v)_{p-1}
	\]
	if $\type(v)=\{i_p\}$.
\end{Remark}

\begin{Proposition}\label{prop:filtered_homology}
	For $2\leq p\leq n$, we have
	\[
	H_i(O(X)_p,O(X)_{p-1}) \cong \bigoplus_{\type(v)=i_p} \tilde H_{i-1}(\lk(v)_{p-1}).
	\]
\end{Proposition}

\begin{Proof}
	We have
	\begin{align*}
		O(X)_p \setminus O(X)_{p-1} &= \{\sigma\in O(X) : i_p \in \type(\sigma) \subseteq I_p \} \\
		&= \coprod_{\type(v)=i_p} \{\sigma\in O(X)_p : v\in \sigma \} \\
		&= \coprod_{\type(v)=i_p} R(v)_p.
	\end{align*}
	Since $\st(v)_p\cap O(X)_{p-1} = \lk(v)_{p-1}$, we obtain the following push-out diagram
  \[
    \xymatrix{
    \coprod_{\type(v)=i_p} \lk(v)_{p-1}\ar[r]\ar[d] & O(X)_{p-1}\ar[d] \\
    \coprod_{\type(v)=i_p} \st(v)_p\ar[r] & O(X)_p.
    }
  \]
	By excision, we obtain
	\begin{align*}
		H_i(O(X)_p,O(X)_{p-1}) &\cong \bigoplus_{\type(v)=i_p} H_i(\st(v)_p,\lk(v)_{p-1}) \\
		&\cong \bigoplus_{\type(v)=i_p} \tilde H_i(\st(v)_p / \lk(v)_{p-1}) \\
		&\cong \bigoplus_{\type(v)=i_p} \tilde H_{i-1}(\lk(v)_{p-1}).
	\end{align*}
    The last line follows from the fact that $\st(v)_p$ is the simplicial cone over $\lk(v)_{p-1}$.
\end{Proof}

\noindent This description of relative homology allows us to show that each filtration subcomplex of the opposition complex is also homotopy equivalent to a bouquet of spheres.

\begin{Proposition}\label{prop:filtrated_opposition_complexes_are_spherical}
	For any $1\leq p \leq n$, the homology groups $\tilde H_i(O(X)_p)$ are trivial except for $i=p-1$. The group $\tilde H_{p-1}(O(X)_p)$ is $\Z$-free.
\end{Proposition}

\begin{Proof}
    Since the complexes $O(X)_p$ are $(p-1)$-dimensional, their top-dimensional homology group $H_{p-1}(O(X)_p)$ is automatically $\Z$-free.

	It remains to show that all other reduced homology groups vanish. We prove this by induction on $n=\rank(X)$. In any case, the statement is true for $O(X)_n=O(X)$ by Theorem \ref{th:opposition_complex_is_spherical} and it is trivial for $O(X)_1$.

	Combining these facts, we obtain the statement for $n=2$. Now assume $n\geq 3$. We prove the statement for $O(X)_p$ for all $2\leq p\leq n$ by reverse induction. The case $p=n$ is already known by Theorem \ref{th:opposition_complex_is_spherical}.

	Hence assume we have the statement for $O(X)_p$ and prove it for $O(X)_{p-1}$. First of all note that by Proposition \ref{prop:filtered_homology}, it follows that
	\[
	H_i(O(X)_p,O(X)_{p-1}) \cong \bigoplus_{\type(v)=i_p} \tilde H_{i-1}(\lk(v)_{p-1})
	\]
	and we have the induction hypothesis for $\lk(v)_{p-1}$, since $\rank(\lk_X(v^+))=\rank(X)-1$.
	We see in particular that $H_i(O(X)_p,O(X)_{p-1})$ vanishes for $i\neq p-1$.

    By the long exact sequence for the pair $(O(X)_p,O(X)_{p-1})$, we obtain that $\tilde H_i(O(X)_{p-1})$ vanishes for $i\neq p-2$.
\end{Proof}

\noindent The following modules $M_p$ will be the coefficient modules in the spectral sequence.

\begin{Definition}
	We define a sequence of $L_p$-modules as follows:
	\[
		M_p \coloneq\begin{cases}
			\Z & p=1 \\
			\tilde H_{p-2}(\lk(v_p)_{p-1}) & 2\leq p \leq n.
		\end{cases}
	\]
    These are $L_p$-modules, since $L_p$ fixes $\lk(v_p)$ and is type-preserving. It hence also fixes the subcomplex $\lk(v_p)_{p-1}$. Note additionally that $M_p$ is $\Z$-free by the previous proposition.
\end{Definition}

\noindent With this definition, we obtain a simple description of the relative homology modules.

\begin{Proposition}\label{prop:structure_filtered_homology}
	For $1\leq p\leq n$, we have
	\[
	H_i(O(X)_p,O(X)_{p-1}) \cong\begin{cases}
		0 & i \neq p-1\\ \Z G \otimes_{L_p} M_p & i=p-1.
	\end{cases}
	\]
\end{Proposition}

\begin{Proof}
	If $p=1$, then $O(X)_1 = \coprod_{\type(v)=i_1} \{v\}$. Then obviously
	\[
	H_0(O(X)_1) \cong \Z G \otimes_{L_1} \Z
	\]
	since $G$ acts transitively on pairs of opposite vertices of the same type.

	For $p>1$, by Proposition \ref{prop:filtered_homology}, we have
	\begin{align*}
        H_i(O(X)_p,O(X)_{p-1}) &\cong \bigoplus_{\type(v)=i_p} \tilde H_{i-1}(\lk(v)_{p-1}), \\
        &\cong \Z G \otimes_{L_p} \tilde H_{i-1}(\lk(v_p)_{p-1} ),
	\end{align*}
	again since $G$ acts transitively on pairs of opposite vertices of the same type. The claim now follows from Proposition \ref{prop:filtrated_opposition_complexes_are_spherical}.
\end{Proof}

\noindent The filtration allows us to obtain an exact complex of $G$-modules, which will be used to construct a relative spectral sequence.

\begin{Definition}\label{def:exact_chain_complexes}
	Consider the following sequence of $G$-modules:
	\[
		C_p\coloneq\begin{cases}
			H_{n-1}(O(X)) & p=n+1 \\
			H_{p-1}(O(X)_p,O(X)_{p-1}) & 1 \leq p \leq n \\
			\Z & p=0 \\
			0 & \text{otherwise.}
		\end{cases}
	\]
	We then have
	\[
	C_p \cong \Z G \otimes_{L_p}M_p
	\]
    for $1\leq p \leq n$ by Proposition \ref{prop:structure_filtered_homology}.
\end{Definition}

\noindent As for cellular chains, there is a chain complex structure on $C_*$. Note that, in contrast to the situation of cellular chains, we have added the modules $C_0$ and $C_{n+1}$ to make the chain complex $C_*$ exact.

\begin{Lemma}\label{l:kdelta_exact_chain_complex}
	There is a boundary map $\partial^C_*$ on $C_*$, which makes $C_*$ into an exact chain complex of $G$-modules.
\end{Lemma}

\begin{Proof}
    The filtration $(O(X)_p)_{p\in\Z}$ induces a $G$-equivariant filtration on cellular chains of $O(X)$. By Theorem \ref{th:specseq_filtered_top_space}, there is hence a spectral sequence of $G$-modules
	\[
	E^1_{p,q} = H_{p+q}(O(X)_{p+1},O(X)_p) \quad\Rightarrow\quad H_{p+q}(O(X)).
	\]
	By Proposition \ref{prop:structure_filtered_homology} we obtain
	\[
	E^1_{p,q} = \begin{cases}
		H_p(O(X)_{p+1},O(X)_p) & q=0 \\
		0 & q \neq 0.
	\end{cases}
	\]
	The spectral sequence hence collapses on the second page, the differential maps on the first page and the edge homomorphisms form the long exact sequence.
\end{Proof}

\begin{Remark}
    A closer comparison to the situation of cellular chains shows that this boundary map is given by the composition
    \[
    H_p(O(X)_{p+1},O(X)_p) \stackrel{\delta}{\rightarrow} H_{p-1}(O(X)_p) \stackrel{H(\pi)}{\rightarrow} H_{p-1}(O(X)_p,O(X)_{p-1})
    \]
    where $\delta$ is the connecting homomorphism of the long exact sequence associated to the pair $(O(X)_{p+1},O(X)_p)$ and where $H(\pi)$ is induced by the projection as in the long exact sequence associated to the pair $(O(X)_p,O(X)_{p-1})$. This is in complete analogy to the situation of cellular homology, for instance in \cite[Chapter 2]{Geo:Top:08}.
\end{Remark}

\subsection{The relative spectral sequence}

For $n\geq 2$, let $G$ be a group with a weak Tits system of rank $n+1$ with associated building $X$ whose type set $I=\{i_1,\ldots,i_{n+1}\}$ is ordered arbitrarily. We adapt the notation from the previous section.

In this general setting, of course, a part of the problem of homological stability is its precise formulation. Which subgroups of a given group $G$ should be considered as the ones to yield stability? These subgroups cannot be expressed explicitly in this generality, they depend on the chosen series of groups. Consequently, we allow for a certain amount of flexibility in the choice of the subgroup $G'$.

\begin{Definition}
    Let $G'$ be any subgroup of $L_{n+1}$ that still acts strongly transitively on the link $X'\coloneq\lk_X(v_{n+1}^+)$, which is a building of type $I_n=\{i_1,\ldots,i_n\}$. In particular, $G'$ admits a weak Tits system of type $I_n$.
	We call the pair $(G,G')$ a \emph{stability pair}.
\end{Definition}

\begin{Remark}
    We will always choose $G'$ such that this assumption of strong transitivity is obviously fulfilled. In fact, if the group $G$ admits a spherical root datum, then every Levi subgroup does as well by Proposition \ref{prop:levi_factors_have_root_datum}. The associated little projective group $L_{n+1}^\dagger$ also admits a spherical root datum, and we will always choose $G'$ to contain $L_{n+1}^\dagger$.
\end{Remark}

\noindent The aim of this section is to associate a spectral sequence to every stability pair. It will be a relative version of the spectral sequence in Corollary \ref{cor:ss_exact_coefficients}.

\begin{Definition}
	By Proposition \ref{prop:links_opposition_complex}, we have
	\[
	O(X)' \coloneq \lk_{O(X)}(v_{n+1}) \cong O(\lk_X(v_{n+1}^+))=O(X'),
	\]
	and we denote the complex on the left by $O(X)'$. This isomorphism is $G'$-equivariant.
	Define a type filtration on $O(X')$ analogously to the one defined on $O(X)$. In addition, the filtration on $O(X)$ induces a filtration on $O(X)'$. The isomorphism is then also filtration-preserving.
\end{Definition}

\noindent Now we consider the exact chain complexes $C_*$ and $C'_*$ associated to the groups $G$ and $G'$ as in Definition \ref{def:exact_chain_complexes}. We obtain the following descriptions:
\begin{align*}
		C_p&\coloneq\begin{cases}
			H_n(O(X)) & p=n+2 \\
			H_{p-1}(O(X)_p,O(X)_{p-1}) & 1 \leq p \leq n+1 \\
			\Z & p=0 \\
			0 & \text{otherwise}
		\end{cases}\\
	\intertext{and}
		C_p'&\coloneq\begin{cases}
			H_{n-1}(O(X)') & p=n+1 \\
			H_{p-1}(O(X)'_p,O(X)'_{p-1}) & 1 \leq p \leq n \\
			\Z & p=0 \\
			0 & \text{otherwise.}
		\end{cases}
\end{align*}
	We see that
	\[
  C_p' \cong \Z G' \otimes_{L_p'} M_p\qquad\text{ and }\qquad C_p\cong \Z G \otimes_{L_p} M_p
	\]
  for $1\leq p\leq n$ and $C'_{n+1} = H_{n-1}(\lk_{O(X)}(v_{n+1})) = M_{n+1}$, where
  \[
  M_p = \tilde H_{p-2}(\lk_{O(X)}(v_p)_{p-1}) \cong \tilde H_{p-2}(\lk_{O(X)'}(v_p)_{p-1})\qquad\text{for } 2\leq p\leq n
  \]
  are the same modules for both groups $G$ and $G'$.

\begin{Lemma}
	The inclusion $O(X)'\hookrightarrow O(X)$ induces a $G'$-equivariant chain map
  \[
  \iota: C_*' \rightarrow C_*.
  \]
\end{Lemma}

\begin{Proof}
    The inclusion $O(X)'\hookrightarrow O(X)$ is $G'$-equivariant and filtration-preserving. The inclusions of pairs then induce maps
	\[
	\iota_p: \underbrace{H_{p-1}(O(X)'_p,O(X)'_{p-1})}_{C_p'} \rightarrow \underbrace{H_{p-1}(O(X)_p,O(X)_{p-1})}_{C_p}
	\]
	for $1\leq p \leq n$ which are compatible with the boundary maps $\partial^C$ and $\partial^{C'}$. Of course $\iota_p=0$ for $p\leq -1$ and $p\geq n+2$. Consider
	\[
	\xymatrix{
		0 & \Z\ar[l] & H_0(O(X)_1) \ar[l]_-{\partial^C_1} & \ar[l]\cdots \\
		0 & \Z\ar[l]\ar@.[u]^{\iota_0} & H_0(O(X)'_1)\ar[l]_-{\partial^{C'}_1}\ar[u]^{\iota_1} & \ar[l]\cdots
	}
	\]
  Here $\iota_0\coloneq \partial^C_1 \circ \iota_1 \circ (\partial^{C'}_1)^{-1}$ is well-defined by a diagram chase. For $p=n+1$, consider
	\[
		\xymatrix{
        \cdots & H_{n-1}(O(X)_n,O(X)_{n-1})\ar[l] & H_n(O(X)_{n+1},O(X)_n)\ar[l]_-{\partial^C_{n+1}} & H_{n}(O(X)_{n+1})\ar[l] \\
        \cdots & H_{n-1}(O(X)'_n,O(X)'_{n-1})\ar[l]\ar[u]^{\iota_n} & \tilde H_{n-1}(O(X)'_n) \ar[l]_-{\partial^{C'}_{n+1}}\ar@.[u]_{\iota_{n+1}} & 0.\ar[l]
		}
	\]
    Note that $O(X)'_n=\lk_{O(X)}(v_{n+1}) = \lk(v_{n+1})_n$, so
    \[
    \tilde H_{n-1}(O(X)'_n) \cong H_n\bigl(\st(v_{n+1})_{n+1},\lk(v_{n+1})_n\bigr)
    \]
    as in the proof of Proposition \ref{prop:filtered_homology}. The inclusion of pairs
    \[
        \bigl(\st(v_{n+1})_{n+1},\lk(v_{n+1})_n\bigr) \hookrightarrow \bigl(O(X)_{n+1},O(X)_n\bigr)
    \]
    then induces the map $\iota_{n+1}$ on homology. A closer inspection using the explicit description of the boundary maps above shows that $\iota$ is a chain map.
\end{Proof}

\begin{Remark}
    It is not difficult to see that, under the above isomorphisms, the map $\iota_p$ is induced by the canonical inclusions
    \[
    \Z G'\otimes_{L'_p} M_p \hookrightarrow \Z G \otimes_{L_p} M_p
    \]
    for $1\leq p \leq n$ and that $\iota_{n+1}: M_{n+1} \rightarrow \Z G \otimes_{L_{n+1}} M_{n+1}$ is given by $m\mapsto 1\otimes_{L_{n+1}} m$.
\end{Remark}

\noindent With these chain complexes and the chain map, we are able to construct the following spectral sequence.

\begin{BreakTheorem}[Stability pair spectral sequence]\label{th:stability_pair_spectral_sequence}
    For each stability pair $(G,G')$, there is a first-quadrant spectral sequence
	\[
	E^1_{p,q}=\begin{cases}
		H_q(G,G';\Z) & p=0 \\
		H_q(L_p,L_p';M_p) & 1\leq p\leq n \\
    H_q(L_{n+1},G';M_{n+1}) & p=n+1
	\end{cases}
    \]
    which converges to zero.
\end{BreakTheorem}

\begin{Proof}
    Let $F_*(G)$ and $F_*(G')$ be the standard resolutions of $\Z$ over $\Z G$ and $\Z G'$, respectively. Consider the two double complexes $F_*(G')\otimes_{G'} C'_*$ and $F_*(G)\otimes_G C_*$ and the map of double complexes
    \begin{align*}
    i: F(G')\otimes_{G'} C' &\rightarrow F(G)\otimes_G C \\
    f \otimes_{G'} c &\mapsto f \otimes_G \iota(c)
\end{align*}
    induced by the inclusion $G'\hookrightarrow G$ and by $\iota$.
    We can apply Theorem \ref{th:relative_spectral_sequence} to obtain the relative spectral sequence
    \[
    E^1_{p,q} = H_q(\Cone_*(i_{\bullet,p})) \Rightarrow 0,
    \]
    converging to zero by Corollary \ref{cor:ss_exact_coefficients} and the long exact sequence associated to the mapping cone complex. All that remains is the description of the first page terms. For $1\leq p\leq n$, we apply Lemma \ref{l:mapping_cone_relative_homology} to obtain
    \[
    E^1_{p,q} \cong H_q\bigl( (F_*(G) \otimes_G C_p ) / (F_*(G') \otimes_{G'} C'_p)\bigr) \cong H_q\bigl( (F_*(G) \otimes_{L_p} M_p ) / (F_*(G') \otimes_{L'_p} M_p ) \bigr),
    \]
    which is isomorphic to $H_q(L_p,L'_p;M_p)$ by Proposition \ref{prop:relative_homology_of_subgroups}. For $p=n+1$, note that
    \begin{align*}
    E^1_{n+1,q} &= H_q\bigl( (F_*(G) \otimes_G C_{n+1} ) / (F_*(G') \otimes_{G'} C'_{n+1})\bigr) \\
    &\cong H_q\bigl( (F_*(G) \otimes_{L_{n+1}} M_{n+1} ) / (F_*(G') \otimes_{G'} M_{n+1} ) \bigr),
\end{align*}
    which is isomorphic to $H_q(L_{n+1},G';M_{n+1})$.
\end{Proof}

\begin{Remark}\label{rem:e1_vanishes_at_n+1}
  Note that by Lemma \ref{l:relative_h0}, we have $E^1_{p,0}=0$ for all $0\leq p\leq n+1$.
\end{Remark}

\section{Homological stability}

In the third part, we will apply the results of the first two parts to specific series of groups. We will only sketch the original application to general linear groups by Charney in \cite{Cha:HSD:80}. Instead, we will use the strong stability theorem for general linear groups by Sah (Theorem \ref{th:sah}) for the following proofs. In the following sections, we then prove homological stability for special linear groups, for unitary groups and for special orthogonal groups.

For the definitions of these groups and of their associated buildings, we will refer to Chapter \ref{ch:building_examples}.

\subsection{General linear groups}\label{sec:general_linear_groups}

Let $D$ be any division ring. We consider the group $G=\Gl_{n+2}(D)$ for $n\geq 2$. The associated building $X$ is the flag complex over $(n+1)$-dimensional projective space and hence of type $A_{n+1}$, see Section \ref{sec:an_example}. The vertices of the opposition complex $O(X)$ are then pairs of complementary vector subspaces of $D^{n+2}$ by Lemma \ref{lem:an_opposition}.

\begin{figure}[hbt]
  \centering
    \begin{tikzpicture}[font=\small]
      \node (G) at (-1,2.4) {$G$};
      \node (G) at (-1,1.6) {$L_p$};
      \node (G) at (-1,.8) {$G'$};
      \node (L) at (-1,0) {$L'_p$};
      \foreach \y in {0,.8,1.6,2.4} {
      \foreach \x in {0,1,2,3,4,6,7,8,9} { \fill (\x,\y) circle (.7mm);}
      \draw (1,\y) -- (2,\y);
      \draw[dotted] (2,\y) -- (3,\y);
      \draw (3,\y) -- (4,\y);
      \draw (6,\y) -- (7,\y);
      \draw[dotted] (7,\y) -- (8,\y);
      \draw (8,\y) -- (9,\y);
      \draw (0,\y) -- (1,\y);
      }
      \draw (9,2.4) -- (10,2.4);
      \draw (9,1.6) -- (10,1.6);
      \draw (4,.8) -- (6,.8);
      \draw (4,2.4) -- (6,2.4);
      \fill (5,.8) circle (.7mm);
      \fill (5,2.4) circle (.7mm);
      \fill (10,2.4) circle (.7mm);
      \fill (10,1.6) circle (.7mm);
      \node (0) at (0,-.6) {$1$}; \node (1) at (1,-.6) {$2$}; \node (3) at (4,-.6) {$p-1$}; \node (4) at (5,-.63) {$p$}; \node (5) at (6,-.6) {$p+1$}; \node (7) at (9,-.62) {$n$}; \node (8) at (10,-.6) {$n+1$};
    \end{tikzpicture}\caption{The Coxeter diagrams of the groups $G$, $L_p$, $G'$, $L'_p$ for type $A_n$.} \label{fig:an_coxeter_diagrams}
\end{figure}
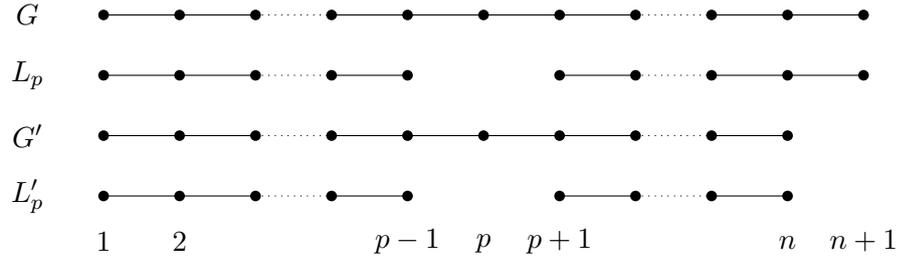

As in Section \ref{sec:an_example}, we choose a basis $e_1,e_2,\ldots,e_{n+2}$ of $D^{n+2}$ and construct the associated frame $\cF_0=\{\langle e_1\rangle,\ldots,\langle e_{n+2}\rangle\}$ which leads to the standard apartment $\Sigma_0=\Sigma(\cF_0)$. We choose the standard chamber $c_0$ in $\Sigma_0$ to be:
\[
    c_0 = ( \langle e_1\rangle \subset \langle e_1,e_2\rangle\subset\cdots\subset\langle e_1,\ldots,e_{n+1}\rangle).
\]
We enumerate the type set $I=\{i_1,\ldots,i_{n+1}\}$ such that a vertex of type $i_k$ is a subspace of dimension $k$. Consequently, we have $v_p^+ = (\langle e_1,\ldots,e_p \rangle)$ for the vertices of $c_0$, and hence $v_p^- = (\langle e_{p+1},\ldots,e_{n+2}\rangle)$. We obtain
\[
	L_p = \begin{pmatrix}
		\Gl_p(D) & 0 \\ 0 & \Gl_{n+2-p}(D)
	\end{pmatrix}
\]
for $1\leq p \leq n+1$ as in Section \ref{sec:an_example}. In particular
\[
	L_{n+1}\cong \begin{pmatrix} \Gl_{n+1}(D) & 0 \\ 0 & D^\times\end{pmatrix}.
\]
We choose $G' = \begin{pmatrix}
	\Gl_{n+1}(D) & 0 \\ 0 & 1
\end{pmatrix}$ which implies that
\[
	L'_p = \begin{pmatrix}
		\Gl_p(D) & 0 & 0 \\
		0 & \Gl_{n+1-p}(D) & 0 \\
	  0 & 0 & 1 \\
	\end{pmatrix}.
\]
Obviously, we obtain $L_p\cong \Gl_p(D) \times \Gl_{n+2-p}(D)$ and $L'_p \cong \Gl_p(D)\times \Gl_{n+1-p}(D)$. The Coxeter diagrams of $G$, $L_p$, $G'$ and $L'_p$ can be found in Figure \ref{fig:an_coxeter_diagrams}.

By Theorem \ref{th:stability_pair_spectral_sequence}, we obtain the following spectral sequence, which has originally been constructed by Charney in \cite{Cha:HSD:80}. We obtain a slightly better description of the first page, but her result is for arbitrary principal ideal domains.

\begin{Theorem}[Charney, Theorem 2.2 in \cite{Cha:HSD:80}]
	For any $n\geq 2$, there is a first-quadrant spectral sequence
	\[
	E^1_{p,q} =\begin{cases}
		H_q\bigl(\Gl_{n+2}(D),\Gl_{n+1}(D);\Z\bigr) & p=0 \\
		H_q\bigl(\Gl_{p}(D)\times \Gl_{n+2-p}(D),\Gl_p(D)\times \Gl_{n+1-p}(D); M_p\bigr) & 1\leq p \leq n\\
        H_q\bigl(\Gl_{n+1}(D)\times D^\times, \Gl_{n+1}(D);M_{n+1}\bigr) & p=n+1
	\end{cases}
	\]
	which converges to zero.
\end{Theorem}

\begin{Remark}
    It is not difficult to see that the maps
	\[
  E^1_{1,q} = H_q\bigl(D^\times \times \Gl_{n+1}(D), D^\times \times \Gl_n(D);\Z\bigr) \rightarrow H_q\bigl(\Gl_{n+2}(D),\Gl_{n+1}(D);\Z\bigr) = E^1_{0,q}
	\]
	are induced by the inclusion of pairs, which is an important ingredient in the homological stability proof by Charney described below.
\end{Remark}

\noindent Note that the vertices of $\lk(v_p)_{p-1}$ are given by
\[
\lk(v_p)_{p-1}^0 = \bigl\{ (V,W) : V\subsetneq \langle e_1,\ldots,e_p\rangle, W \supsetneq \langle e_{p+1},\ldots,e_{n+2}\rangle, V\oplus W=D^{n+2}\bigr\}.
\]
In particular, these factors of the Levi subgroups
\[ \begin{pmatrix}
	\one_p & 0 \\ 0 & \Gl_{n+2-p}(D)
\end{pmatrix}\qquad\text{and}\qquad\begin{pmatrix}
	\one_p & 0 & 0\\ 0 & \Gl_{n+1-p}(D) & 0 \\ 0 & 0 & 1
\end{pmatrix}
\]
act trivially on $\lk(v_p)_{p-1}$ and hence on $M_p$. We can hence apply the relative Lyndon\slash Hochschild-Serre spectral sequence (Theorem \ref{th:relative_lyndon_hochschild_serre}) to obtain a description of the first page of the aforementioned spectral sequence in terms of integral relative group homology of smaller general linear groups. This can be used to apply an ingenious ``bootstrap procedure'' to prove homological stability inductively for general linear groups as in \cite{Cha:HSD:80}. The result, originally proved for principal ideal domains, is

\begin{Theorem}[Charney, Theorem 3.2 in \cite{Cha:HSD:80}]
	For $n\geq 3k$, we have
	\[
	H_k(\Gl_{n+1}(D),\Gl_{n}(D);\Z) = 0.
	\]
\end{Theorem}

\noindent This theorem is far from optimal, to the author's knowledge, the results by Sah and van der Kallen in Section \ref{sec:stability} have the best known stability ranges of $n\geq k$ for division rings with infinite centre and of $n\geq 2k$ for other division rings.

\subsection{Special linear groups}

In this section, we consider special linear groups $\Sl_{n+2}(D)$ over infinite fields $D$. By the arguments of Section \ref{sec:an_example}, special linear groups act strongly transitively on the same building $X$ we have considered in the previous section. The opposition complex $O(X)$ and hence also the modules $M_p$ are the same as for general linear groups. In the case $G=\Sl_{n+2}(D)$ for $n\geq 2$, the Levi subgroups admit the following structure:
\[
L_p = \biggl\{\begin{pmatrix}
	A & 0 \\ 0 & B
\end{pmatrix} : A \in \Gl_{p}(D), B \in \Gl_{n+2-p}(D), \det(A)\det(B)=1 \biggr\},
\]
which splits as a semidirect product
\begin{align*}
L_p &= \biggl\{\begin{pmatrix}
  A & 0 & 0 \\ 0 & \det(A)^{-1}  & 0 \\ 0 & 0 & \one_{n+1-p}
\end{pmatrix} : A \in \Gl_p(D) \biggr\}\ltimes\begin{pmatrix}
	\one_{p} & 0 \\ 0 & \Sl_{n+2-p}(D)
      \end{pmatrix}\\
		& \cong \Gl_{p}(D)  \ltimes \Sl_{n+2-p}(D),
\end{align*}
and the group $\Sl_{n+2-p}(D)$ acts trivially on $M_p$ by the same argument as in Section \ref{sec:general_linear_groups}. In particular
\[
L_{n+1} = \biggl\{\begin{pmatrix}
  A & 0 \\ 0 & \det(A)^{-1}
\end{pmatrix} : A\in\Gl_{n+1}(D) \biggr\}\text{, and we choose } G'=\begin{pmatrix}
  \Sl_{n+1}(D) & 0 \\ 0 & 1
\end{pmatrix}.
\]
We obtain
\begin{align*}
L_{p}' &=\biggl\{\begin{pmatrix}
  A & 0 & 0 \\ 0 & \det(A)^{-1} & 0 \\ 0 & 0 & \one_{n+1-p}
    \end{pmatrix} : A \in \Gl_p(D) \biggr\}\ltimes\begin{pmatrix}
	\one_{p} & 0 &0 \\ 0 & \Sl_{n+1-p}(D) &0 \\ 0 & 0 & 1
  	\end{pmatrix}\\
		&\cong \Gl_{p}(D) \ltimes \Sl_{n+1-p}(D),
	\end{align*}
again with $\Sl_{n+1-p}(D)$ acting trivially on $M_p$. In particular, we have
\begin{align*}
  L_1 &= \biggl\{\begin{pmatrix}
    \det(A)^{-1} & 0 \\ 0 & A
  \end{pmatrix} : A\in\Gl_{n+1}(D) \biggr\}\cong \Gl_{n+1}(D) \\
  L'_1&=\biggl\{\begin{pmatrix}
    \det(A)^{-1} & 0 & 0 \\ 0 & A & 0 \\ 0 & 0 & 1
  \end{pmatrix} : A \in \Gl_{n}(D) \biggr\}\cong \Gl_{n}(D).
\end{align*}

\noindent Figure \ref{fig:an_coxeter_diagrams} also shows the Coxeter diagrams of the groups which are involved here. As in \cite{Cha:HSD:80}, we apply Theorem \ref{th:stability_pair_spectral_sequence} to obtain a spectral sequence converging to zero. Again, we obtain a slightly better description of the first page, while restricting to fields.

\begin{Theorem}[Charney, Theorem 2.3 in \cite{Cha:HSD:80}]\label{th:sln_spectral_sequence}
	For any $n\geq 2$, there is a first-quadrant spectral sequence
	\[
	E^1_{p,q} =\begin{cases}
		H_q\bigl(\Sl_{n+2}(D),\Sl_{n+1}(D);\Z\bigr) & p=0 \\
		H_q\bigl(\Gl_{p}(D)\ltimes \Sl_{n+2-p}(D),\Gl_p(D)\ltimes \Sl_{n+1-p}(D); M_p\bigr) & 1\leq p \leq n\\
    H_q\bigl(\Gl_{n+1}(D),\Sl_{n+1}(D);M_{n+1}\bigr) & p=n+1
	\end{cases}
	\]
	which converges to zero.
\end{Theorem}

\noindent Charney then uses this spectral sequence and a version of the relative Lyndon\slash Hochschild-Serre spectral sequence (Theorem \ref{th:relative_lyndon_hochschild_serre}) to prove homological stability for $n\geq 3k$. One can improve this result using Theorem \ref{th:sah}, however.

\begin{Theorem}\label{th:sln_result}
  If $D$ is an infinite field, then $n\geq 2k-1$ implies
  \[
	  H_k(\Sl_{n+1}(D),\Sl_{n}(D);\Z)=0.
  \]
\end{Theorem}

\begin{Proof}
  In spite of the formulation of the theorem we will prove equivalently that $n\geq 2k-2$ implies
  \[
  H_k(\Sl_{n+2}(D),\Sl_{n+1}(D);\Z)=0.
  \]
  Since relative $H_0$ vanishes always by Lemma \ref{l:relative_h0} and since $H_1(\Sl_{n+2}(D);\Z)=0$ for all $n\geq 0$ by \cite[2.2.3]{HoM:CGK:89}, we can start an induction over $k$.

  Let $k\geq 2$. As induction hypothesis, assume that
  \[
  H_l(\Sl_{n+2}(D),\Sl_{n+1}(D);\Z)=0\quad\text{ for all $l<k$ and all $n\geq 2l-2$.}
  \]
  We will show that $H_k(\Sl_{n+2}(D),\Sl_{n+1}(D);\Z)=0$ for all $n\geq 2k-2$. For any $n\geq 2k-2\geq 2$, by Theorem \ref{th:sln_spectral_sequence}, we have the spectral sequence $E^1_{p,q}$ converging to zero. Note first of all that
  \[
  E^1_{0,q} = H_q(\Sl_{n+2}(D),\Sl_{n+1}(D);\Z).
  \]
  In particular, we have to prove that $E^1_{0,k}=0$. Since the spectral sequence converges to zero, it is enough to prove $E^1_{p,q}=0$ for all $p+q\leq k+1$ and $p\geq 1$, which will then imply $E^1_{0,q}=0$ for all $0\leq q\leq k$. The situation is illustrated in Figure \ref{fig:sln}. We will prove the vanishing of the modules separately for the regions I, II and III.

\begin{figure}
\centering
    \begin{tikzpicture}[scale=0.5,font=\small]
    \draw[->] (0,0) -- (0,9);
    \draw[->] (0,0) -- (10,0);
    \draw (0,8) -- (2,8) -- (2,7) -- (3,7) -- (3,6) -- (4,6) -- (4,5) -- (5,5) -- (5,4) -- (6,4) -- (6,3) -- (7,3) -- (7,2) -- (8,2) -- (8,1) -- (9,1) -- (9,0);
    \node (h1) at (-.5,.5) {$0$};
    \node (hk) at (-.5,7.5) {$k$};
    \node (v1) at (.5,-.5) {$0$};
    \node (vk) at (7.5,-.5) {$k$};
    \node (p) at (10.5,0) {$p$};
    \node (q) at (0,9.5) {$q$};
    \node (leq) at (5,7.5) {$2\leq p \leq n$};
    \draw[dotted] (1,8.5) -- (1,1);
    \draw[dotted] (2,8.5) -- (2,1);
    \draw[dotted] (8,8.5) -- (8,1);
    \draw[dotted] (-.5,1) -- (9,1) -- (9,-.5);
    \node (1) at (.5,7.5) {$*$};
    \node (2) at (1.5,4) {I};
    \node (3) at (3.5,2) {II};
    \node (4) at (4.5,.5) {III};
  \end{tikzpicture}\caption{$E^1_{p,q}$ in the proof of Theorem \ref{th:sln_result}}\label{fig:sln}
\end{figure}
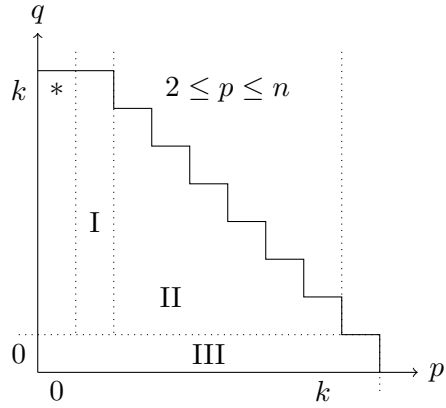
  For region I, that is for $p=1$ and $1\leq q\leq k$, note that
  \[
  E^1_{1,q} \cong H_q(\Gl_{n+1}(D),\Gl_{n}(D);\Z).
  \]
  Since $n\geq 2k-2\geq k\geq q$, we know that $E^1_{1,q}=0$ by Theorem \ref{th:sah}.

  For region II, that is $p+q\leq k+1$, $q\geq 1$ and $2\leq p\leq n$, we have in particular $q\leq k-1$. Additionally
  \[
    2q-2+p = (p+q) + q -2 \leq 2k-2 \leq n,
  \]
  so $n-p\geq 2q-2$. By the induction hypothesis, we hence have
  \begin{equation}
      H_q(\Sl_{n-p+2}(D), \Sl_{n-p+1}(D);\Z) = 0 \label{eq:vanishing_terms_of_relative_lhs}
  \end{equation}
  for all $p,q$ with $p+q\leq k+1$ and $p\geq 2$.

  For the terms $E^1_{p,q}=H_q(L_p,L'_p;M_p)$ in region II, consider the relative Lyndon\slash Hochschild-Serre spectral sequence from Theorem \ref{th:relative_lyndon_hochschild_serre}:
  \[
  L^2_{i,j} = H_i(\Gl_p(D);H_j(\Sl_{n+2-p}(D),\Sl_{n+1-p}(D);\Z)\otimes_\Z M_p) \Rightarrow H_{i+j}(L_p,L'_p;M_p).
  \]
  Note that we proved in \eqref{eq:vanishing_terms_of_relative_lhs} that $L^2_{i,j}=0$ for $j\leq q\leq k+1-p$, hence $H_q(L_p,L'_p;M_p)=0$ for all $p$, $q$ in region II. Hence $E^1_{p,q}=0$ for region II.

  Finally, note that region III, that is, all $E^1_{p,0}$ with $0\leq p\leq k+1\leq n+1$ vanish by the remark after Theorem \ref{th:stability_pair_spectral_sequence}.

  Inspecting the spectral sequence $E^1_{p,q}$ yields $E^1_{0,k}=E^{\infty}_{0,k}=0$, which proves the theorem.
\end{Proof}

\begin{Remark}
    If we compare this to the results of Section \ref{sec:stability}, we see that the stability result by Hutchinson and Tao (Theorem \ref{th:hutchinson_tao}) for special linear groups improves this to $n\geq k$ in the case where $D$ is a field of characteristic zero. Our result above improves the known stability range for all other infinite fields by one, however. Up to now, the best result known to the author is due to van der Kallen (Theorem \ref{th:vdk}) proving stability for $n\geq 2k$ for division rings.
\end{Remark}

\subsection{Unitary groups}

In this section, let $D$ be a division ring. Let $J:D\rightarrow D$ be an involution and let $\varepsilon\in\{\pm 1\}$. Let $V$ be a finite-dimensional right $D$-vector space. Choose a (trace-valued) $(\varepsilon,J)$-hermitian form $h:V\times V\rightarrow D$ of Witt index $\ind(h)=n+1$ and consider the associated unitary group $\U(V)$ as in Section \ref{sec:ex_buildings_of_type_bn_cn}.

From this section, we obtain that $\U(V)$ admits a Tits system of type $C_{n+1}$ which is not thick if and only if $J=\id$, $\dim(V)=2\ind(h)$ and $\varepsilon\neq -1$, which is the case of orthogonal groups. We also allow this case explicitly. The corresponding building $X$ is isomorphic to the flag complex over all totally isotropic subspaces of $V$. The opposition complex $O(X)$ is then isomorphic to the flag complex over pairs of opposite totally isotropic subspaces of $V$, where $A\leq V$ is \emph{opposite} $B\leq V$ if $A\oplus B^\perp = V$.

By Proposition \ref{prop:bases_for_spaces_with_hermitian_form}, there is a decomposition $V= \caH_{n+1}\oplus E$, where $E$ is anisotropic, and there is a basis
\[
  e_{-(n+1)},\ldots,e_{-1},e_1,\ldots,e_{n+1}\]
  of the hyperbolic module $\caH_{n+1}$ such that
\[
  h(e_i,e_j) =\begin{cases}
    \varepsilon\qquad & \text{if } i+j=0 \text{ and } i>0 \\
    1 & \text{if } i+j=0\text{ and } i<0 \\
    0 & \text{otherwise.}
  \end{cases}
\]
The associated frame determines an apartment $\Sigma_0$ of the building $X$. We choose the chamber $c_0\in\Sigma_0$ and the type enumeration such that
\[
v_p = ( \langle e_{-(n+1)},\ldots,e_{-p}\rangle, \langle e_{p},\ldots,e_{n+1}\rangle )
\]
and we write $\caH_p=\langle e_{-p},\ldots,e_{-1},e_1,\ldots,e_p\rangle$ and $\caH_0=\{0\}$. In Section \ref{sec:ex_buildings_of_type_bn_cn}, we have already seen that the Levi subgroups $L_p$ then admit the following structure
\begin{align*}
L_p &=\biggl\{\begin{pmatrix}
    S & 0 & 0 & 0 \\
    0 & A & 0 & B \\
    0 & 0 & S^{-J} & 0 \\
    0 & C & 0 & D
  \end{pmatrix} : S\in \Gl_{n+2-p}(D), \begin{pmatrix}
    A & B \\
    C & D
  \end{pmatrix}\in \U(\caH_{p-1}\perp E)\biggr\}\\
  &\cong \Gl_{n+2-p}(D)\times \U(\caH_{p-1}\perp E).
\end{align*}
In particular, we have
\[
L_{n+1} = \biggl\{\begin{pmatrix}
    s & 0 & 0 & 0 \\
    0 & A & 0 & B \\
    0 & 0 & s^{-J} & 0 \\
    0 & C & 0 & D
\end{pmatrix} : s\in D^\times, \begin{pmatrix}
    A & B \\
    C & D
  \end{pmatrix}\in \U(\caH_n\perp E)\biggr\}.
\]
We choose $G'$ to be
\[
  G' = \biggl\{\begin{pmatrix}
    1 & 0 & 0 & 0 \\
    0 & A & 0 & B \\
    0 & 0 & 1 & 0 \\
    0 & C & 0 & D
\end{pmatrix} : \begin{pmatrix}
    A & B \\
    C & D
  \end{pmatrix}\in \U(\caH_n\perp E)\biggr\}\cong \U(\caH_n\perp E),
\]
hence
\begin{align*}
L'_p &= \biggl\{\begin{pmatrix}
  1 & 0 & 0 & 0 & 0 & 0 \\
  0 & S & 0 & 0 & 0 & 0 \\
  0 & 0 & A & 0 & 0 & B \\
  0 & 0 & 0 & S^{-J} & 0 & 0 \\
  0 & 0 & 0 & 0 & 1 & 0 \\
  0 & 0 & C & 0 & 0 & D
\end{pmatrix} : S \in\Gl_{n+1-p}(D), \begin{pmatrix}
    A & B \\
    C & D
  \end{pmatrix}\in \U(\caH_{p-1}\perp E)\biggr\}\\
&\cong \Gl_{n+1-p}(D) \times \U(\caH_{p-1}\perp E).
\end{align*}

\noindent The diagrams of these groups can be seen in Figure \ref{fig:cn_coxeter_diagrams}.

\begin{figure}[hbt]
  \centering
    \begin{tikzpicture}[font=\small]
      \node (G) at (-1,2.4) {$G$};
      \node (G) at (-1,1.6) {$L_p$};
      \node (G) at (-1,.8) {$G'$};
      \node (L) at (-1,0) {$L'_p$};
      \foreach \y in {0,.8,1.6,2.4} {
      \foreach \x in {0,1,2,3,4,6,7,8,9} { \fill (\x,\y) circle (.7mm);}
      \draw (1,\y) -- (2,\y);
      \draw[dotted] (2,\y) -- (3,\y);
      \draw (3,\y) -- (4,\y);
      \draw (6,\y) -- (7,\y);
      \draw[dotted] (7,\y) -- (8,\y);
      \draw (8,\y) -- (9,\y);
      \draw (0,\y + 0.05) -- ++(1,0);
      \draw (0,\y - 0.05) -- ++(1,0);
      }
      \draw (9,2.4) -- (10,2.4);
      \draw (9,1.6) -- (10,1.6);
      \draw (4,.8) -- (6,.8);
      \draw (4,2.4) -- (6,2.4);
      \fill (5,.8) circle (.7mm);
      \fill (5,2.4) circle (.7mm);
      \fill (10,2.4) circle (.7mm);
      \fill (10,1.6) circle (.7mm);
      \node (0) at (0,-.6) {$1$}; \node (1) at (1,-.6) {$2$}; \node (3) at (4,-.6) {$p-1$}; \node (4) at (5,-.63) {$p$}; \node (5) at (6,-.6) {$p+1$}; \node (7) at (9,-.62) {$n$}; \node (8) at (10,-.6) {$n+1$};
    \end{tikzpicture}\caption{The Coxeter diagrams of the groups $G$, $L_p$, $G'$, $L'_p$ for type $C_n$.} \label{fig:cn_coxeter_diagrams}
\end{figure}

For $n\geq 2$, we apply Theorem \ref{th:stability_pair_spectral_sequence} to obtain a spectral sequence
\[
E^1_{p,q} \cong\begin{cases}
  H_q\bigl(\U(\caH_{n+1}\perp E),\U(\caH_{n}\perp E);\Z\bigr) & p=0 \\
  H_q\bigl(\Gl_{n+2-p}(D)\times \U(\caH_{p-1}\perp E),\\ \hspace{4cm} \Gl_{n+1-p}(D) \times \U(\caH_{p-1}\perp E); M_p\bigr) \hspace{1cm} & 1\leq p \leq n\\
  H_q\bigl(D^\times\times \U(\caH_n\perp E), \U(\caH_n\perp E);M_{n+1}\bigr) & p=n+1
\end{cases}
\]
which converges to zero. Using this spectral sequence, we obtain a generalisation of the theorems by Vogtmann for orthogonal groups in \cite{Vog:HSO:79} and \cite{Vog:SPH:81}, which were later generalised to Dedekind rings by Charney in \cite{Cha:gtV:87}.

\begin{Theorem}\label{th:un_result}
  For a division ring $D$ with infinite centre, the relative homology modules
  \[
    H_k\bigl(\U(\caH_{n+1}\perp E), \U(\caH_n\perp E);\Z\bigr)
  \]
  vanish for $n\geq 2$ if $k=1$ and for $n\geq k\geq 2$. If the centre of $D$ is finite, relative homology vanishes for $n\geq 2k$.
\end{Theorem}

\begin{Proof}
    Note first that all assumptions imply $n\geq 2$, so we can construct the spectral sequence $E^1_{p,q}$ as above. For $2\leq p\leq n$, the vertices of $\lk(v_p)_{p-1}$ are given by
  \[
  \lk(v_p)_{p-1}^0 = \bigl\{ (A,B) : \langle e_{-(n+1)},\ldots,e_{-p}\rangle \subsetneq A, B \supsetneq \langle e_p,\ldots,e_{n+1}\rangle, A\oplus B^{\perp}=V\bigr\}.
  \]
  In particular, the general linear group factors of $L_p$ and $L'_p$ act trivially on this filtrated link and hence on $M_p$. For $p\geq 1$, we can hence apply the relative Lyndon\slash Hochschild-Serre spectral sequence (Theorem \ref{th:relative_lyndon_hochschild_serre}) to obtain
  \[
  L^2_{i,j} = H_i( \U(\caH_{p-1}\perp E); H_j(\Gl_{n+2-p}(D),\Gl_{n+1-p}(D);\Z)\otimes_\Z M_p) \Rightarrow H_{i+j}(L_p,L'_{p};M_p).
  \]

\begin{figure}
\centering
    \begin{tikzpicture}[scale=0.5,font=\small]
    \draw[->] (0,0) -- (0,9);
    \draw[->] (0,0) -- (10,0);
    \draw (0,8) -- (2,8) -- (2,7) -- (3,7) -- (3,6) -- (4,6) -- (4,5) -- (5,5) -- (5,4) -- (6,4) -- (6,3) -- (7,3) -- (7,2) -- (8,2) -- (8,1) -- (9,1) -- (9,0);
    \node (h1) at (-.5,.5) {$0$};
    \node (hk) at (-.5,7.5) {$k$};
    \node (v1) at (.5,-.5) {$0$};
    \node (vk) at (7.5,-.5) {$k$};
    \node (p) at (10.5,0) {$p$};
    \node (q) at (0,9.5) {$q$};
    \node (leq) at (5,7.5) {$1\leq p \leq k\leq n$};
    \draw[dotted] (1,8.5) -- (1,1);
    \draw[dotted] (8,8.5) -- (8,1);
    \draw[dotted] (-.5,1) -- (9,1) -- (9,-.5);
    \node (1) at (.5,7.5) {$*$};
    \node (3) at (3,3) {I};
    \node (4) at (4.5,.5) {II};
  \end{tikzpicture} \caption{$E^1_{p,q}$ in the proof of Theorem \ref{th:un_result}}\label{fig:un}
\end{figure}
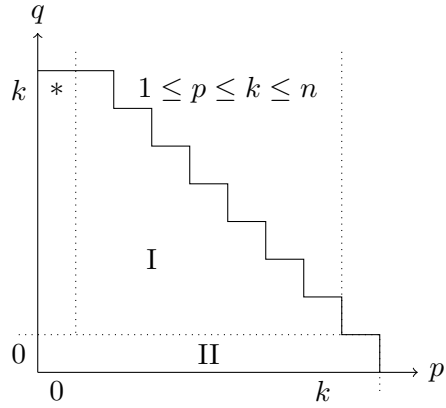

Again, we want to show that $E^1_{0,k}=0$. We do this by showing that $E^1_{p,q}=0$ for $p+q\leq k+1$ and $p\geq 1$. The situation is illustrated in Figure \ref{fig:un}. We will prove the vanishing of these modules separately for the regions I and II.

First of all, consider region II, which consists of the modules $E^1_{p,0}$ with
\[
1\leq p\leq k+1\leq n+1.
\]
These modules vanish by the remark after Theorem \ref{th:stability_pair_spectral_sequence}.

Region I is given by $p+q\leq k+1$, $q\geq 1$ and $1\leq p\leq k\leq n$. We distinguish cases.

If the centre of $D$ is infinite, by Sah's result (Theorem \ref{th:sah}), we have
\[
    H_j(\Gl_{n+2-p}(D),\Gl_{n+1-p}(D);\Z)=0
\]
for $j\leq n+1-p$. So $E^1_{p,q}$ vanishes for $p+q\leq n+1$ and $1\leq p\leq n$ by the relative Lyndon\slash Hochschild-Serre spectral sequence. Since $n\geq k$, the modules $E^1_{p,q}$ vanish in region I.

  If the centre of $D$ is finite, van der Kallen's result (Theorem \ref{th:vdk}) implies $E^1_{p,q}=0$ for $p+2q\leq n+1$. By hypothesis, $n\geq 2k$, so $p+q\leq k+1$ and $q\leq k$ imply $p+2q\leq 2k+1\leq n+1$, which shows again that $E^1_{p,q}=0$ in region I. This proves the theorem.
\end{Proof}

\begin{Remark}
    As already mentioned in Section \ref{sec:stability}, Mirzaii and van der Kallen have proved homological stability for unitary groups over local rings with infinite residue fields in \cite{MaB:HSU:02} and \cite{Mir:HSU:05} with a slightly weaker stability range of $n\geq k+1$. We restrict ourselves to division rings with infinite centre, but on the other hand we allow an anisotropic kernel.
\end{Remark}

\noindent Note that, if $E=\{0\}$ and $J=\id$, which forces $D$ to be a field, we obtain stability results for the following two special cases.

\begin{Theorem}
  For the \emph{symplectic} and \emph{orthogonal groups over an infinite field $D$}, we obtain
  \begin{align*}
    H_k\bigl(\Sp_{2n+2}(D),\Sp_{2n}(D);\Z\bigr) &=0 \\
    H_k\bigl(\Orth_{n+1,n+1}(D),\Orth_{n,n}(D);\Z\bigr) &=0
  \end{align*}
  if $k=1$ and $n\geq 2$ or $n\geq k\geq 2$. If $D$ is a finite field, then the relative homology groups vanish for $k\geq 2n$.
\end{Theorem}

\noindent By a combination of the methods for Theorem \ref{th:sln_result} and Theorem \ref{th:un_result}, results on special unitary groups can also be obtained. The following result on special orthogonal groups is particularly interesting:

\begin{Theorem}\label{th:son_result}
  For an infinite field $D$, we have
  \[
    H_k\bigl(\SO_{n+1,n+1}(D),\SO_{n,n}(D);\Z\bigr) =0
  \]
  for $k=1$ and $n\geq 2$ or $n \geq k\geq 2$. If $D$ is a finite field, then the relative homology groups vanish for $n\geq 2k$.
\end{Theorem}

\begin{Proof}
  In this case, note that $J=\id$. The Levi subgroups hence have the following simple structure:
  \[
  L_p =\biggl\{\begin{pmatrix}
    S & 0 & 0 & 0 \\
    0 & A & 0 & B \\
    0 & 0 & S^{-1} & 0 \\
    0 & C & 0 & D
  \end{pmatrix} : S\in \Gl_{n+2-p}(D), \begin{pmatrix}
    A & B \\
    C & D
  \end{pmatrix}\in \SO_{p-1,p-1}(D)\biggr\}.
  \]
  Since by assumption $n\geq 2$, we can apply Theorem \ref{th:stability_pair_spectral_sequence} to obtain a relative spectral sequence converging to zero. We have to prove $E^1_{0,k}=0$.
  As in the proof of Theorem \ref{th:un_result}, the general linear factors of $L_p$, $L'_p$ act trivially on $M_p$ and we consider the relative Lyndon\slash Hochschild-Serre spectral sequence to obtain
  \[
  L^2_{i,j} = H_i\bigl( \SO_{p-1,p-1}(D); H_j(\Gl_{n+2-p}(D),\Gl_{n+1-p}(D);\Z)\otimes_\Z M_p\bigr) \Rightarrow H_{i+j}(L_p,L'_{p};M_p).
  \]
  Now, for infinite fields, note that by Theorem \ref{th:sah}, for $j\leq n+1-p$, we have $L^2_{i,j}=0$. In particular, $p+q\leq k+1\leq n+1$ and $p\geq 1$ implies $E^1_{p,q}=0$, except for $p=n+1=k+1$, where we apply the remark after Theorem \ref{th:stability_pair_spectral_sequence} again. We obtain the theorem by inspection of the spectral sequence.

  For finite fields, by Theorem \ref{th:vdk}, for $2j\leq n+1-p$, we have $L^2_{i,j}=0$. So $p+q\leq k+1$ and $q\leq k$ imply $p+2q\leq 2k+1\leq n+1$, so $E^1_{p,q}=0$, and we obtain the result.
\end{Proof}

\chapter{Lattices in buildings}\label{ch:lattices}

This chapter deals with the construction of lattices in groups acting on affine buildings. More specifically, we will construct uniform lattices in the automorphism groups of two-dimensional affine buildings. This is done by constructing a complex of groups with the right local developments whose universal cover is a building and whose fundamental group acts cocompactly with finite stabilisers. The fundamental group is hence a uniform lattice in the full automorphism group of the building.

In the first section, we recall some basic facts about lattices in the automorphism groups of polyhedral complexes.

In order to construct the required complexes of groups, we study groups acting regularly on points or lines of finite generalised polygons in Section \ref{sec:singer_polygons}. Using polygons with such an automorphism group, so-called \emph{Singer polygons}, we then give a construction of complexes of groups whose local developments are generalised polygons.

In the next two sections, we construct complexes of groups that give rise to uniform lattices on buildings of type $\tilde A_2$ and $\tilde C_2$, respectively, and investigate some of their properties.

In the last section, we calculate group homology of these lattices.

The results of this chapter are also contained in \cite{Ess:PRL:09}. The reader might want to consult Section \ref{sec:lattices_results} for an overview of this chapter.

\section{Lattices in automorphism groups of polyhedral complexes}\label{sec:lattices_in_autom_groups}

Let us first introduce the general concept of a lattice. We will only require uniform lattices in this thesis. Hence we will not give detailed definitions of the Haar measure, for example.

Let $G$ be a locally compact topological group. Then there is a left invariant measure on $G$ which is unique up to a multiplicative constant. It is called the \emph{Haar measure}. This measure induces an invariant measure on every quotient space of $G$.

\begin{Definition}
    A \emph{lattice} $\Gamma \leq G$ is a discrete subgroup of $G$ such that the measure space $\Gamma \backslash G$ has finite measure.

    A \emph{uniform lattice $\Gamma$} is a discrete and cocompact subgroup of $G$, that means $\Gamma\backslash G$ is compact.
\end{Definition}

\begin{Examples}
    The additive group $\Z\leq\R$ is a uniform lattice. The group $\Sl_n(\Z)$ is a non-uniform lattice in $\Sl_n(\R)$.
\end{Examples}

\noindent We will only be interested in uniform lattices in this thesis. Locally compact groups and uniform lattices arise naturally in the following setting.

\paragraph{Situation}
Let $X$ be a locally finite piecewise Euclidean complex. Its full automorphism group $G=\Aut(X)$ is a topological group, a neighbourhood basis for the identity is given by the stabilisers of compact sets. The vertex stabilisers $G_x$ are then open and compact, since
\[
G_x = \varprojlim_{r\rightarrow\infty} G_x\vert_{B_r(x)},
\]
where $B_r(x)$ is the closed ball around $x$. Obviously, since $X$ is locally finite, every restriction $G_x\vert_{B_r(x)}$ is a finite group, so $G_x$ is profinite and hence compact.
In particular, $G$ is locally compact. We assume additionally that $G$ acts cocompactly on $X$, that is $G\backslash X$ is compact.

See \cite[Chapter 3]{BL:TL:01} for a description of this situation in case $X$ is a locally finite tree.

\begin{Lemma}
    A subgroup $\Gamma\leq G$ is discrete if and only if all its stabiliser subgroups of vertices are finite. It is cocompact in $G$ if and only if it acts cocompactly on $X$.
\end{Lemma}

\begin{Proof}
    If $\Gamma$ is discrete, then $\Gamma\cap G_x$ is compact and discrete and hence finite for all $x\in X$. If $\Gamma\cap G_x$ is finite, then $\{1\}$ is open in $\Gamma\cap G_x$ and hence in $\Gamma$, since $G_x$ is open in $G$. So $\Gamma$ is discrete.

    For cocompactness, observe that, since $G$ acts cocompactly on $X$, there is a compact set $C\subseteq X$ such that $G\cdot C = X$. Hence $(\Gamma\backslash G)\cdot C = \Gamma\backslash X$ and so $\Gamma\backslash G$ is compact if and only if $\Gamma\backslash X$ is compact.
\end{Proof}

\noindent In particular, $\Gamma\leq G$ is a uniform lattice if and only if it acts cocompactly with finite stabilisers on $X$.

\begin{Construction}
    Let $G(\cY)$ be a developable complex of groups with finite vertex groups over a finite scwol. Then the geometric realisation of the universal cover $|\cX|$ is a locally finite, piecewise Euclidean complex.

    Its automorphism group $\Aut(|\cX|)$ is then a locally compact group with compact open stabilisers. The fundamental group $\pi_1(G(\cY))$ acts on $|\cX|$, the stabilisers are the vertex groups and the quotient is the geometric realisation of the scwol $\cY$. So $\pi_1(G(\cY))$ is a uniform lattice in $\Aut(\lvert\cX\rvert)$.
\end{Construction}

\noindent Hence, to construct lattices in the automorphism groups of locally finite buildings, all we have to do is the following:

\paragraph{Goal} Construct a finite scwol with finite vertex groups whose universal cover is a Euclidean building.

\section{Complexes of groups with generalised polygons as local developments}\label{sec:singer_polygons}

To achieve that goal of constructing a suitable scwol, in particular the local developments have to have the right structure. Remember that, by the list of diagrams of affine Coxeter groups (Figure \ref{fig:affine_coxeter_diagrams}), a two-dimensional affine building must be of type $\tilde A_2$, $\tilde C_2$ or $\tilde G_2$. In particular, a vertex link in a two-dimensional affine building must be a generalised polygon as in Section \ref{sec:generalised_polygons}.

In order to obtain a complex of groups whose universal cover is a two-dimensional building, we have to ensure that the local developments at vertices are finite generalised polygons. The main ingredients for our constructions are generalised polygons admitting \emph{Singer groups}.

\begin{Definition}
    A \emph{Singer polygon} is a generalised polygon which admits an automorphism group acting regularly, that is transitively and freely, on its set of points. Such a group of automorphisms is called a \emph{Singer group}.
\end{Definition}

\noindent We will investigate examples for Singer polygons in the following sections. The easiest case is a generalised 2-gon, which is nothing but a complete bipartite graph.

\begin{Lemma}
	Any group of order $k$ acts point- and line-regularly on a bipartite graph of order $(k,k$).
\end{Lemma}

\begin{Proof}
	Identify the set of points and the set of lines of the bipartite graph each with the elements of the group and consider the left regular representation on both sets.
\end{Proof}

\subsection{Projective planes}\label{subsec:proj_planes}

A generalised triangle is a projective plane. A projective plane is \emph{classical} if it is the usual projective plane over a division ring as in Section \ref{sec:an_example}. The following classical result gives rise to the name of Singer groups.

\begin{Theorem}[Singer, \cite{Si:FPG:38}]
    Every finite classical projective plane admits a cyclic Singer group. A generator of such a group is called a \emph{Singer cycle}.
\end{Theorem}

\noindent It is conjectured that, as a stronger converse statement, every finite projective plane admitting a group acting transitively on points is already classical. This has been an open problem for a long time, see \cite{Gil:TPP:07} for a recent contribution. It is not even known whether a projective plane with an abelian Singer group must be classical. On the other hand, all Singer groups on finite classical projective planes are classified by Ellers and Karzel:

\begin{Theorem}[Ellers-Karzel, 1.4.17 in \cite{De:FG:68}]
    Fix a finite classical projective plane of order $p^e$ with a Singer group $S$.

    Then there is a divisor $n$ of $3e$ such that $n(p^e-1)$ divides $p^{3e}-1$ and $S$ can be presented in the following way:
    \[
    S = \langle \gamma, \delta \,|\, \gamma^s=1, \delta^n=\gamma^t, \delta\gamma\delta^{-1}= \gamma^{p^{3e/n}} \rangle,
    \]
    where $s = \frac{p^{3e}-1}{n(p^e-1)}$ and $t = \frac{p^{3e}-1}{n(p^{3e/n}-1)}$.

    Conversely, every group satisfying these relations acts point-regularly on a classical projective plane. In addition, the group $S$ is contained in $\PGl_3(\F_{p^e})$ if and only if $n\in\{1,3\}$.
\end{Theorem}

\noindent Finally, we require the following simple fact about arbitrary finite projective planes.

\begin{Proposition}[4.2.7 in \cite{De:FG:68}]\label{prop:regular_on_lines}
	A Singer group of a finite projective plane also acts regularly on lines.
\end{Proposition}

\subsection{Generalised quadrangles}\label{subsec:GQs}

There is no classification of all Singer quadrangles known to the author. However, there is a list of all Singer quadrangles among all known finite generalised quadrangles, see the book \cite{STW:SQ:09}. All known Singer quadrangles except one arise by Payne derivation from translation generalised quadrangles and symplectic quadrangles.

Since we require rather explicit descriptions of the quadrangles, we consider only one large class of examples, the so-called \emph{slanted symplectic quadrangles}. For a detailed investigation of these quadrangles, see \cite{GJS:SSQ:94} and \cite{Str:PSQ:03}.

We first define symplectic quadrangles, which are a special case of the polar spaces we have constructed in Section \ref{sec:ex_buildings_of_type_bn_cn}. Nevertheless, we give the concrete definition here to fix the notation.

\begin{Definition}
     Fix a prime power $q>2$. Let $V=\langle e_{-2},e_{-1},e_1,e_2\rangle$ be a four-dimensional vector space over the finite field $\F_q$ endowed with the symplectic form $h$ satisfying
    \[
        h(e_i,e_j) =\begin{cases}
            -1 & i+j=0, i>j \\
            1 & i+j=0, i<j \\
            0 & \text{else.}
        \end{cases}
    \]
    The polar space consisting of all non-trivial totally isotropic subspaces of $V$, that is, subspaces $U$ satisfying $U\subseteq U^\perp$, forms a generalised quadrangle, called the \emph{symplectic quadrangle $W(q)$}. More specifically, set
    \begin{align*}
        \cP(W(q)) &= \{ U\subseteq V : \dim(U)=1, U\subseteq U^\perp \}, \\
        \caL(W(q)) &= \{ W\subseteq V : \dim(W)=2, W\subseteq W^\perp \}, \\
        \cF(W(q)) &= \{ (U,W) \in \cP(W(q)) \times \caL(W(q)) : U\subseteq W \}.
    \end{align*}
\end{Definition}

We denote the full symplectic group by $\Sp_4(q)$, which is the group of all linear maps preserving $h$, and we fix the following notation
\[
    x(a) =\begin{pmatrix}
        1 & a & 0 & 0 \\
        0 & 1 & 0 & 0 \\
        0 & 0 & 1 & -a \\
        0 & 0 & 0 & 1
    \end{pmatrix},\qquad y(b) =\begin{pmatrix}
        1 & 0 & b & 0 \\
        0 & 1 & 0 & b \\
        0 & 0 & 1 & 0 \\
        0 & 0 & 0 & 1
    \end{pmatrix},\qquad z(c) =\begin{pmatrix}
        1 & 0 & 0 & c \\
        0 & 1 & 0 & 0 \\
        0 & 0 & 1 & 0 \\
        0 & 0 & 0 & 1
    \end{pmatrix}
\]
for specific elements in $\Sp_4(q)$, where $a,b,c\in\F_q$. We abbreviate $x=x(1)$, $y=y(1)$ and $z=z(2)$. It is not hard to see that $[x,y]=z$.

\begin{Lemma}\label{lemma:heisenberg_group}
    The subgroup $E=\langle x(a),y(b),z(c) : a,b,c\in\F_q \rangle\leq \Sp_4(q)$ acts regularly on all points not collinear with the point $p_0=\F_q (1,0,0,0)^T$.
    \begin{itemize}
	    \item If $q$ is odd, then $E$ is a three-dimensional Heisenberg group over $\F_q$. If $q$ is prime, we have the following simple presentation:
    \[
        E = \langle x,y \,|\, z=xyx^{-1}y^{-1}, x^q=y^q=z^q=1, xz=zx, yz=zy \rangle.
    \]
        \item If $q$ is even, then $E\cong \F_q^3$, which is hence an elementary abelian 2-group.
    \end{itemize}
\end{Lemma}

\begin{Proof}
    It is an easy calculation to show that all points not collinear with $p_0$ have the form $\F_q(x,y,z,1)^T$. Simple matrix calculations then imply that $E$ acts regularly on these points.

    If $q$ is odd, then $E$ is a three-dimensional Heisenberg group by \cite[4.3]{STW:SQ:09}. It is a simple calculation to verify the presentation in case that $q$ is prime. By calculating three commutators, one can see that $E$ is isomorphic to $\F_q^3$ if $q$ is even.
\end{Proof}

\paragraph{Payne derivation} We will give a definition of \emph{Payne derivation} as described in \cite{STW:SQ:09}. The same procedure is called \emph{slanting} in \cite{GJS:SSQ:94}. We write $p\sim r$ for collinear points $p$ and $r$. Then we define
\begin{align*}
    p^\perp & = \{ r\in \cP(W(q)) : p \sim r \} \\
    \{p,p'\}^{\perp\perp} &= \{ r\in\cP(W(q)) : r \in s^\perp \text{ for all } s\in p^\perp \cap p'^\perp \}.
\end{align*}

\begin{Definition}
    The \emph{slanted symplectic quadrangle $W(q)^\Diamond$} is given as follows:
\begin{itemize}
    \item The points of $W(q)^\Diamond$ are all points of $W(q)$ not collinear with $p_0$.
    \item The lines of $W(q)^\Diamond$ are all lines of $W(q)$ not meeting $p_0$ as well as the sets
	    \[\{p_0,r\}^{\perp\perp} \setminus \{p_0\}\] for all points $r\in W(q)^\Diamond$.
\end{itemize}
\end{Definition}

\noindent The following characterisation will be used later on to construct lattices in buildings of type $\tilde C_2$.

\begin{Theorem}\label{th:slanted_sympl_quadrangle}
	The slanted symplectic quadrangle $W(q)^\Diamond$ is a Singer quadrangle of order $(q-1,q+1)$ with Singer group $E$. A set of representatives $L$ for the $E$-action on lines is given by all lines through the point $p_1=(0,0,0,1)^T$:
	\[
    L = \bigl\{ l_{[a:b]} = \bigl(\F_q(0,0,0,1)^T + \F_q (0,b,a,0)^T\bigr) \,\big|\, [a:b]\in \bP\F_q^2 \bigr\} \sqcup \bigl\{l_0 = \{p_0,p_1\}^{\perp\perp}\setminus\{p_0\}\bigr\}.
	\]
	The line stabilisers have the form
	\[
	E_{l_{[a:b]}}=\biggl\{\begin{pmatrix}
		1 & fa & fb & 0 \\
		0 & 1 & 0 & fb \\
		0 & 0 & 1 & -fa \\
		0 & 0 & 0 & 1
    \end{pmatrix}: f\in \F_q \biggr\}\quad\text{and}\quad E_{l_0} = \{ z(f) : f\in\F_q \}.
	\]
	In particular, all line stabilisers are isomorphic to the additive group of $\F_q$. If $q$ is an odd prime, the line stabilisers have the form
	\[
    E_{l_{[a:b]}} = \langle x^ay^bz^{-\tfrac{1}{2}ab} \rangle\quad\text{and}\quad E_{l_0} = \langle z \rangle,
	\]
    where of course $2$ is invertible in $\F_q^\times$. If $q$ is even, we have
    \[
    E_{l_{[a:b]}} = \langle x(a)y(b) \rangle\quad\text{and}\quad E_{l_0} = \langle z(1) \rangle
    \]
    as $\F_q$-subspaces in $\F_q^3$.
\end{Theorem}

\begin{Proof}
	The first claim follows from Lemma \ref{lemma:heisenberg_group}. For the set of line representatives and the structure of the stabilisers, see \cite[Lemma 4.1]{Str:PSQ:03}. The last two statements are simple calculations.
\end{Proof}

\subsection{A local complex of groups construction}\label{subsec:local_cplx_of_groups_construction}

Using Singer polygons, we will now construct a complex of groups whose local development at one vertex is a generalised polygon. With this local piece of data, we will later construct complexes of groups whose local developments are buildings.

See Chapter \ref{ch:complexes_of_groups} for the definitions of scwols and of complexes of groups.

\begin{Definition}
	Associated to any generalised $n$-gon $\cI=(\cP,\caL,\cF)$, there is a scwol $\cZ(\cI)$ with vertex set
	\begin{align*}
		V(\cZ(\cI)) &= \cP\sqcup\caL\sqcup\cF \\
		\intertext{and with edges, suggestively written as follows,}
		E(\cZ(\cI)) &= \{ p \leftarrow f : p\in\cP, f\in \cF, p\in f \} \sqcup \{ l\leftarrow f: l\in\caL, f\in \cF,l\in f\},
\end{align*}
where we mean $i(x\leftarrow y)=y$ and $t(x\leftarrow y)=x$.
\end{Definition}

\noindent The geometric realisation of this scwol is isomorphic as a simplicial complex to the barycentric subdivision of the geometric realisation $|\cI|$ of the generalised polygon.

\paragraph{Notation} Let $\cI=(\cP,\caL,\cF)$ be a Singer polygon with Singer group $S$. Fix a point $p\in \cP$ and a set $L$ of line representatives for the $S$-action on $\caL$. For all lines $l\in L$, we consider the sets
\[
	D_l = \{ d\in S: (p,dl)\in\cF \}.
\]
The sets $D_l$ are called \emph{difference sets}. Furthermore, let $F_l\coloneq \{ (p,dl) : d\in D_l\}$ for all $l\in L$. This is a partition of the set of all flags containing $p$, which we denote by
\[
F \coloneq \{f\in\cF : p\in f\} = \coprod_{l\in L}F_l = \coprod_{l\in L} \{ (p,dl):d\in D_l\}.
\]
In particular, the set $F$ is a set of representatives of flags for the $S$-action on $\cF$.

\paragraph{Construction} Let a complex of groups $G(\cY)$ be given. Let $v\in V(\cY)$ be a minimal element (that is: no edges start in $v$), and assume that the scwol $\cY(v)$ as in Definition \ref{def:scwol_closed_star} has the following structure:
\begin{align*}
	V(\cY(v)) \cong&\,\{v \} \sqcup \{p\} \sqcup L \sqcup F. \\
	\intertext{Again, we will write the edges in the following suggestive fashion:}
	E(\cY(v))\cong&\,\{ v\leftarrow p \} \sqcup \{ v\leftarrow l : l\in L \} \sqcup \{ v\leftarrow f : f  \in F \} \\
	& \quad\sqcup \{ p\leftarrow f : f\in F \} \sqcup \coprod_{l\in L}\Bigl\{ l\leftarrow f: f\in F_l\Bigr\}.
\end{align*}

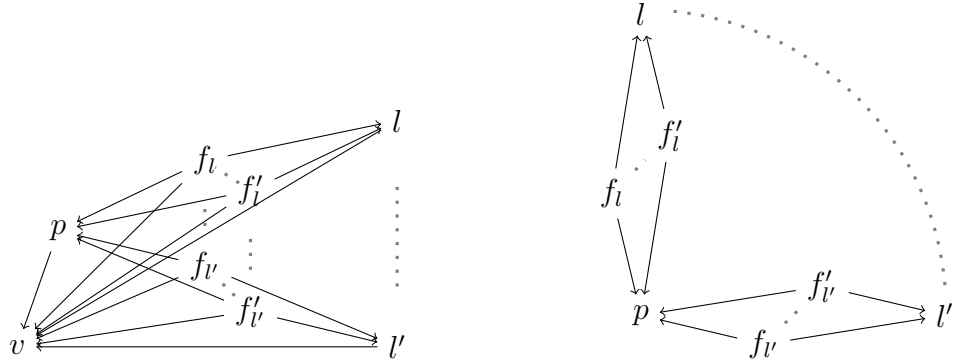
\begin{figure}[hbt]
\centering
    \begin{tikzpicture}
	\node (v) at (0,0,0) {$v$};
	\node (p) at (-1,0,-4) {$p$} edge[->] (v);
	\node (l1) at (5,3,0) {$l$} edge[->] (v);
	\node (l2) at (5,0,0) {$l'$} edge[->] (v) edge[color=gray, loosely dotted,very thick, shorten <=.5cm, shorten >=.5cm] (l1);
	\node (f1a) at (1.7,1.7,-2) {$f_{l}$} edge[->] (v) edge[->] (p) edge[->] (l1);
	\node (f1b) at (2.3,1.3,-2) {$f'_{l}$} edge[->] (v) edge[->] (p) edge[->] (l1) edge[color=gray, loosely dotted,very thick, shorten <=.2cm, shorten >=.2cm] (f1a);
	\node (f2a) at (1.7,0.3,-2) {$f_{l'}$} edge[->] (v) edge[->] (p) edge[->] (l2) edge[color=gray, loosely dotted, very thick, shorten <=.2cm, shorten >=.2cm] (f1a);
	\node (f2b) at (2.3,-0.3,-2) {$f'_{l'}$} edge[->] (v) edge[->] (p) edge[->] (l2) edge[color=gray, loosely dotted,very thick, shorten <=.2cm, shorten >=.2cm] (f2a) edge[color=gray, loosely dotted, very thick, shorten <=.2cm, shorten >=.2cm] (f1b);
\end{tikzpicture} \hspace{2cm}
\begin{tikzpicture}[scale=2]
	\node (p) at (0,0,0) {$p$};
	\node (l1) at (2,0,0) {$l'$};
	\node (l2) at (0,2,0) {$l$};
	\node (f1a) at (1,0,.5) {$f_{l'}$} edge[->] (p) edge[->] (l1);
    \node (f1b) at (1,0,-.5) {$f'_{l'}$} edge[->] (p) edge[->] (l1) edge[color=gray,loosely dotted,very thick, shorten <=.2cm, shorten >=.2cm] (f1a);
	\node (f2a) at (0,1,.5) {$f_l$} edge[->] (p) edge[->] (l2);
    \node (f2b) at (0,1,-.5) {$f'_{l}$} edge[->] (p) edge[->] (l2) edge[color=gray,loosely dotted,very thick, shorten <=.2cm, shorten >=.2cm] (f2a);
	\draw[color=gray,loosely dotted,very thick] (2,.2,0) arc (5:85:2);
\end{tikzpicture} \caption{The scwol $\cY(v)$ on the left, the link $\Lk_v(\cY)$ on the right.}\label{fig:local_development}
\end{figure}

\noindent The scwol $\cY(v)$ can be visualised as in Figure \ref{fig:local_development}. The left image shows the full scwol $\cY(v)$, while the right picture illustrates the upper link of the vertex $v$.

\noindent Assume now that the vertex groups in $G(\cY)$ are as follows
\[
G_v \cong S,\qquad G_p=G_{f}=\{1\},\qquad G_l \cong S_l \qquad \forall l\in L, \forall f\in F,
\]
that the monomorphisms are the obvious inclusions and that the twist elements are
\[
g_{v\leftarrow p\leftarrow (p,dl)} = 1, \qquad g_{v\leftarrow l \leftarrow (p,dl) } = d^{-1}\qquad\forall l\in L, d\in D_l.
\]
The significance of this construction lies in the following

\begin{Proposition}\label{prop:local_development_v}
	The local development $\cY(\tilde v)$ isomorphic to the cone over the scwol $\cZ(\cI)$ associated to the generalised polygon $\cI$.
\end{Proposition}

\begin{Proof}
	Since $v$ is a minimal element in $\cY$, by Definition \ref{def:local_development}, we have $\cY(\tilde v) = \{\tilde v\} * \Lk_{\tilde v}(\cY)$. By the construction of the scwol $\Lk_{\tilde v}(\cY)$ in Definition \ref{def:upper_link}, we have
	\begin{align*}
		V(\Lk_{\tilde v}(\cY)) =& \,\bigl\{ (g\psi_a(G_{i(a)}),a) : a\in E(\cY), t(a)=v, g\psi_a(G_{i(a)})\in G_v/\psi_a(G_{i(a)}) \bigr\} \\
		=&\, \{ ( \{g\}, v \leftarrow p ) : g\in S\} \sqcup \{ (g S_l, v\leftarrow l) : g\in S, l\in L \}\\
		&\,\quad \sqcup \{ (\{g\}, v\leftarrow f) : g\in S, f\in F \},
	\end{align*}
	which is in bijection to $V(\cZ(\cI))=\cP \sqcup \caL \sqcup \cF$ via
    \[
        (\{g\},v\leftarrow p)\mapsto gp,\, (gS_l,v\leftarrow l)\mapsto gl\,\text{ and }\,(\{g\}, v\leftarrow f)\mapsto gf.
    \]
    For the edge set, we obtain
	\begin{align*}
		E(\Lk_{\tilde v}(\cY)) =& \,\bigl\{ (g\psi_{ab}(G_{i(b)}),a,b) : (a,b) \in E^{(2)}(\cY), t(a)=v, g\psi_{ab}(G_{i(b)}) \in G_v/\psi_{ab}(G_{i(b)}) \bigr\} \\
		=&\, \{ ( \{g\}, v\leftarrow p \leftarrow f) : g\in G, f\in F \} \\
		&\,\quad \sqcup \{ ( \{g\}, v\leftarrow l \leftarrow f) : g\in G, l\in L, f\in F_l \}
	\end{align*}
	which is in bijection to the set of edges $E(\cZ(\cI))$ via
	\begin{align*}
        ( \{g\}, v\leftarrow p \leftarrow f ) &\mapsto (gp \leftarrow gf ) \\
        ( \{g\}, v\leftarrow l \leftarrow (p,dl) ) &\mapsto (gdl \leftarrow (gp,gdl) ).
	\end{align*}
	The twist elements precisely guarantee that these bijections commute with the maps $i$ and $t$, respectively. We have
	\[\begin{array}{rclcrcl}
		i( (\{g\},v\leftarrow p \leftarrow f)) &=& (\{g\}, v\leftarrow f) &\,\mapsto\,& gf &=& i(gp\leftarrow gf) \\
		t( (\{g\},v\leftarrow p \leftarrow f)) &=& (\{g\}, v\leftarrow p) &\,\mapsto\,& gp &=& t(gp\leftarrow gf) \\
		i( (\{g\},v\leftarrow l \leftarrow (p,dl))) &=& (\{g\}, v\leftarrow (p,dl)) &\,\mapsto\,& (gp,gdl) &=& i(gdl\leftarrow (gp,gdl)) \\
		t( (\{g\},v\leftarrow l \leftarrow (p,dl))) &=& (gdS_l, v\leftarrow l) &\,\mapsto\,& gdl &=& t(gdl\leftarrow (gp,gdl)),
	\end{array}\]
    where we use Definition \ref{def:upper_link} for the maps $i$ and $t$. So $\Lk_{\tilde v}(\cY)\cong\cZ(\cI)$.
\end{Proof}

\noindent Note that the special situation is always simpler.
    \begin{description}
        \item[Bipartite graphs] If a group $S$ acts regularly on points and lines of a complete bipartite graph of order $(k,k)$, there is only one line orbit, hence $L=\{l\}$ and $S_l=\{1\}$. Since the bipartite graph is complete, we have $|F|=|F_l|=k$ and $D_l=S$.
        \item[Projective planes] If a group $S$ is a Singer group on a projective plane, there is only one line orbit by Proposition \ref{prop:regular_on_lines}, hence $L=\{l\}$ and $S_l=\{1\}$. In particular, there is only one difference set $D\coloneq D_l$.
        \item[Generalised quadrangles] For slanted symplectic quadrangles of order $(q-1,q+1)$, there are $q$ line orbits, hence $|L|=|F|=q$. For simplicity, we choose the set $L$ to be the set of all lines incident to $p$. Then we have $D_l=\{1\}$ for all $l\in L$.
    \end{description}

\noindent We endow the geometric realisation $|\cY|$ with a piecewise Euclidean metric such that the angles between $v\leftarrow p$ and $v\leftarrow f$ as well as the angles between $v\leftarrow f$ and $v\leftarrow l$ are $\pi/2n$ for all lines $l\in L$ and all flags $f\in F$.

\begin{Proposition}\label{prop:key}
    The geometric link $\Lk(\tilde v,\st(\tilde v))$ is isometric to the barycentric subdivision of the generalised $n$-gon $|\cI|$ with its standard metric. In particular, it is a connected CAT(1) space with diameter $\pi$.
\end{Proposition}

\begin{Proof}
    By Proposition \ref{prop:local_development_v}, the local development $\cY(\tilde v)$ is isomorphic to the cone over the scwol $\cZ(\cI)$. By Definition \ref{def:geometric_link}, the polyhedral complex structure on $\Lk(\tilde v,\cY(\tilde v))$ is then the barycentric subdivision of $\cI$.
    By construction, the angles at the vertex $\tilde v$ are $\pi/2n$. Since the edge length in the geometric link $\Lk(\tilde v,\st(\tilde v))$ is given by the angles, the geometric link is hence isometric to the barycentric subdivision of $|\cI|$ with the natural metric of a generalised polygon. The CAT(1) condition follows from Theorem \ref{th:metrics_on_buildings}.
\end{Proof}

\section{Lattices in buildings of type \texorpdfstring{$\tilde A_2$}{\textasciitilde A2}}\label{sec:a2}

In this section, we will give a construction of small complexes of groups whose universal covers are affine buildings of type $\tilde A_2$. We will first state a general construction. Later, we specialise to cyclic Singer groups to obtain lattices with very simple presentations. For these special lattices, we investigate the structure of spheres of radius two in the corresponding buildings.

\subsection{A general construction of panel-regular lattices in \texorpdfstring{$\tilde A_2$}{\textasciitilde A2}-buildings}\label{subsec:general_a2_construction}

For the whole construction, we fix three finite Singer projective planes $\cI_1$, $\cI_2$ and $\cI_3$ of order $q$ with three (possibly isomorphic) Singer groups $S_1$, $S_2$ and $S_3$. Note that the projective planes need not be classical, even though this is likely in view of what we have said in Section \ref{subsec:proj_planes}. By Proposition \ref{prop:regular_on_lines}, these Singer groups act regularly on lines as well. As in Section \ref{subsec:local_cplx_of_groups_construction}, we obtain a difference set for each of these groups by choosing a point and a line in the corresponding projective plane. We denote these three difference sets by $D_1$, $D_2$ and $D_3$.

Write $J=\{0,1,\dots,q\}$ and choose bijections $d_\alpha : J\rightarrow D_\alpha$ for $\alpha\in \{1,2,3\}$, which we call \emph{ordered difference sets}.

\paragraph{Construction} Associated to this piece of data, we consider the complex of groups $G(\cY)$ over the scwol $\cY$ with vertices
\begin{align*}
       V(\cY) &\coloneq \{ v_1, v_2, v_3\} \sqcup \{e_1,e_2,e_3\} \sqcup \{ f_j : j \in J \} \\
       \intertext{and edges}
       E(\cY) &\coloneq \{ v_\alpha \leftarrow e_\beta: \alpha\neq\beta\} \sqcup \{ v_\alpha \leftarrow f_j : j\in J\} \sqcup \{e_\beta \leftarrow f_j : j\in J\}.
\end{align*}

\noindent Figure \ref{fig:cYfor2} illustrates this scwol for $q=2$.

\begin{figure}[hbt]
\centering
\begin{tikzpicture}[scale=4,font=\small]
       \node (p1) at (0,0,0) {$v_1$};
       \node (p2) at (2,0,0) {$v_2$};
       \node (p3) at (1,0,.866) {$v_3$};
       \node (l3) at (1,0,0) {$e_3$} edge[->] (p1) edge[->] (p2);
       \node (l1) at (1.5,0,.433) {$e_1$} edge[->] (p2) edge[->] (p3);
       \node (l2) at (.5,0,.433) {$e_2$} edge[->] (p1) edge[->] (p3);
       \node (f1) at (1,.5,.433) {$f_1$} edge[->] (p1) edge[->] (p2) edge[->] (p3) edge[->] (l1) edge[->] (l2) edge[->] (l3);
       \node (f2) at (1,0,.433) {$f_2$} edge[->] (p1) edge[->] (p2) edge[->] (p3) edge[->] (l1) edge[->] (l2) edge[->] (l3);
       \node (f3) at (1,-.5,.433) {$f_3$} edge[->] (p1) edge[->] (p2) edge[->] (p3) edge[->] (l1) edge[->] (l2) edge[->] (l3);
\end{tikzpicture} \caption{The scwol $\cY$ for $q=2$}\label{fig:cYfor2}
\end{figure}

Choose the three vertex groups $G_{v_\alpha}$ to be isomorphic to $S_\alpha$ for $\alpha\in\{1,2,3\}$. All other vertex groups are trivial. The twist elements are defined as follows:
\[
g_{v_\alpha\leftarrow e_\beta \leftarrow f_j} \coloneq\begin{cases}
       1 & \beta - \alpha \equiv_3 1 \ \\
       d_\alpha(j)^{-1}& \beta - \alpha \equiv_3 2.
\end{cases}
\]

\noindent We endow $|\cY|$ with a locally Euclidean metric as follows:

Let $\Delta$ be the geometric realisation of one triangle in the affine Coxeter complex of type $\tilde A_2$. For each $j\in J$, we map each subcomplex spanned by $\{v_1,v_2,v_3,e_1,e_2,e_3,f_j\}$ onto the barycentric subdivision of $\Delta$ in the obvious way and pull back the metric. We obtain a locally Euclidean metric on $|\cY|$. In particular, the angles at every vertex $v_\alpha$ are $\pi/3$.

\begin{Proposition}\label{prop:a2_is_developable}
       The complex of groups $G(\cY)$ is developable.
\end{Proposition}

\begin{Proof} We want to apply Proposition \ref{prop:geometric_link_nonpos_curved}, so we have to check that the geometric link of every vertex in its local development is CAT(1).

    For any vertex $v_\alpha$, the subcomplex $\cY(v_\alpha)$ obviously has the appropriate structure for the construction in Section \ref{subsec:local_cplx_of_groups_construction}. By construction, the induced complex of groups has the right form to apply Proposition \ref{prop:local_development_v}, and we obtain that $\cY(\tilde v_\alpha)$ is isomorphic to the cone over the scwol $\cZ(\cI_\alpha)$. Then, by Proposition \ref{prop:key}, the geometric link $\Lk(\tilde v_\alpha, \st(\tilde v_\alpha))$ is CAT(1).

       For any vertex $e_\beta$, by Definition \ref{def:local_development} the local development has the form
       \[
       \St(\tilde e_\beta) = \lvert\Lk^{e_\beta}(\cY) \ast\, \{\tilde e_\beta\} \ast\, \Lk_{\tilde e_\beta}(\cY)\rvert.
       \]
       The lower link $\Lk^{e_\beta}(\cY)$ is isomorphic to the scwol consisting of the two adjacent vertices $\{v_\alpha: \alpha\neq \beta\}$, the upper link of the local development $\Lk_{\tilde e_\beta}(\cY)$ is isomorphic to a scwol consisting of the vertices $\{f_j: j\in J\}$. The geometric realisation of the local development is hence isometric to $q+1$ triangles joined along one edge. In particular, the geometric link is CAT(1).

       The upper link of the local development $\Lk_{\tilde f_j}(\cY)$ is trivial. The lower link $\Lk^{f_j}(\cY)$ consists of the sub-scwol of $\cY$ spanned by the vertices $\{v_1,v_2,v_3,e_1,e_2,e_3\}$. The local development $\cY(\tilde f_j)$ is hence just one triangle, so the geometric link is also CAT(1). By Proposition \ref{prop:geometric_link_nonpos_curved}, the complex $G(\cY)$ is then non-positively curved and hence developable by Theorem \ref{th:developability_of_nonpositively_curved_complexes}.
\end{Proof}

\begin{Proposition}\label{prop:fundamental_group_a2}
       The fundamental group $\Gamma$ of the complex $G(\cY)$ admits the following presentation:
    \[
    \Gamma = \Bigl\langle S_1,S_2,S_3 \,\Big| \begin{array}{c}\text{ all relations in the groups }S_1,S_2,S_3,\\ d_1(j)d_2(j)d_3(j) = d_1(j')d_2(j')d_3(j') \quad\forall j,j'\in J\end{array}\Bigr\rangle.
    \]
\end{Proposition}

\begin{Remark}
    We show in Theorem \ref{th:a2_classification} that $D_\alpha$ and $d_\alpha$ can always be chosen such that $d_\alpha(0)=1$, which simplifies this presentation even further.
\end{Remark}

\begin{Proof}
       Consider the following maximal spanning subtree $T$:
       \begin{align*}
               E(T)= \{ v_1\leftarrow e_2, v_1\leftarrow e_3, v_2\leftarrow e_3, v_2\leftarrow e_1, v_3\leftarrow e_1 \} \sqcup \{e_3\leftarrow f_j : j \in J \}.
       \end{align*}
    By taking a look at Definition \ref{def:fundamental_group}, we obtain a presentation of $\Gamma$ with generating sets $S_1$, $S_2$ and $S_3$ and generators for all edges of the scwol $\cY$ not contained in $T$. The relations imposed are first of all the relations from the groups $S_\alpha$. For the relations of the form $k_a k_b = g_{a,b}k_{ab}$, consider Figure \ref{fig:triangle_for_presentation} showing the $j$-th triangle in the complex of groups $G(\cY)$. There, black arrows indicate edges contained in $T$ and all other edges are drawn dotted. Group elements $g$ written on dotted edges $a$ indicate that $k_a=g$.

\begin{figure}[hbt]
   \centering
   \begin{tikzpicture}[scale=7.5,font=\small]
           \node (v1) at (0,0) {$v_1$};
           \node (v3) at (1,0) {$v_3$};
           \node (v2) at (.5,.866) {$v_2$};
           \node (e2) at (.5,0) {$e_2$} edge[->] (v1) edge[->,dotted] node[below,color=gray]{$k_{v_3\leftarrow e_2}$} (v3);
           \node (e3) at (.25,.466) {$e_3$} edge[->] (v1) edge[->] (v2);
           \node (e1) at (.75,.466) {$e_1$} edge[->] (v2) edge[->] (v3);
           \node (fj) at (.5,.289) {$f_j$} edge[->] (e3) edge[->,dotted] node[color=gray,pos=0.7]{$d_1(j)$} (v1) edge[->,dotted] node[color=gray,above]{$1$} (v2) edge[->,dotted] node[pos=0.2,color=gray]{$d_2(j)^{-1}$} (v3) edge[->,dotted] node[color=gray]{$d_2(j)^{-1}$} (e1) edge[->,dotted] node[color=gray]{$d_1(j)$} (e2);
           \node (m12) at (.593,.540) {$d_2(j)^{-1}$};
           \node (m32) at (.417,.540) {$1$};
           \node (m31) at (.25,.252) {$d_1(j)^{-1}$};
           \node (m23) at (.75,.252) {$1$};
           \node (m21) at (.333,.096) {$1$};
           \node (m21) at (.666,.096) {$d_3(j)^{-1}$};
    \end{tikzpicture}\caption{One triangle in $G(\cY)$}\label{fig:triangle_for_presentation}
\end{figure}

All edge elements but $k_{v_3\leftarrow e_2}$ have already been replaced by elements in the groups $S_\alpha$. For $k_{v_3\leftarrow e_2}$, we obtain the relation $k_{v_3\leftarrow e_2}=d_1(j)d_2(j)d_3(j)$ for all $j\in J$. Since this edge is contained in all triangles for all $j\in J$, we obtain the additional relations for all pairs of triangles. All other relations from Definition \ref{def:fundamental_group} do not show up here, since almost all vertex groups are trivial.
\end{Proof}

\begin{Theorem}\label{th:a2_is_building}
The universal cover $\cX$ is a building of type $\tilde A_2$, where the vertex links are isomorphic to $\cI_1$, $\cI_2$ and $\cI_3$. The fundamental group $\Gamma=\pi_1(G(\cY))$ with presentation
    \[
    \Gamma = \Bigl\langle S_1,S_2,S_3 \,\Big| \begin{array}{c}\text{ all relations in the groups }S_1,S_2,S_3,\\ d_1(j)d_2(j)d_3(j) = d_1(j')d_2(j')d_3(j') \quad\forall j,j'\in J\end{array}\Bigr\rangle
    \]
is hence a uniform lattice in the full automorphism group of the building. The set of all chambers containing a fixed panel is a fundamental domain for the action.
\end{Theorem}

\begin{Proof}
    The universal cover $|\cX|$ is a simply connected space which is locally CAT(0). So by Theorem \ref{th:cartan_hadamard}, the space $\lvert\cX\rvert$ is CAT(0). Consider the space $|\cX|$ as a polyhedral complex consisting only of simplices where each cell is the preimage of subscwols spanned by $\{v_1,v_2,v_3,e_1,e_2,e_3,f_j\}_{j\in J}$. The vertices of this polyhedral complex are the preimages of the vertices $v_{\alpha}$ of $\cX$. This is, a priori, not a simplicial complex since the intersection of two simplices might not be a simplex.

However, in view of Theorem \ref{th:recognition} it remains to see that $|\cX|$ is thick and that the geometric links of all vertices are connected and have diameter $\pi$. Thickness of $|\cX|$ is clear by construction. Again by Proposition \ref{prop:key}, the geometric link of every vertex is connected and of diameter $\pi$.

By Theorem \ref{th:recognition}, the geometric realisation $|\cX|$ is then either a two-dimensional affine building or a product of two trees. Of course, the universal cover cannot be a product of trees and it has to be of type $\tilde A_2$, since all vertex links are projective planes.
\end{Proof}

\begin{Remark}
    Unfortunately, there are many examples of buildings of type $\tilde A_2$. It is clear that $\cX$ has an automorphism group which is transitive on vertices of the same type, but there are known examples of non-classical buildings with vertex-transitive automorphism groups, see \cite{HvM:NTB:90}. However, some information can be obtained by considering spheres of radius two in $\cX$ and comparing these to corresponding spheres in the classical buildings of type $\tilde A_2$ associated to $\Q_p$ and $\F_p(\!(t)\!)$, respectively. We will explicitly calculate these spheres of radius two in Section \ref{subsec:spheres}.
\end{Remark}

\paragraph{Construction} Since we know that $|\cX|$ is a building, it is a flag complex. By using the construction of the universal cover in Construction \ref{con:universal_cover}, we obtain a very explicit description of the building as the flag complex whose 1-skeleton is the following graph:
\begin{align*}
        V(X) &\coloneq \Gamma/S_1 \sqcup \Gamma/S_2 \sqcup \Gamma/S_3 \\
        E(X) &\coloneq \bigl\{ (gS_1,gS_2) : g\in\Gamma\bigr\} \sqcup
        \bigl\{ (gS_2,gS_3) : g\in\Gamma\bigr\}
		\\&\qquad\quad\sqcup \bigl\{ (gS_1,gk^{-1}S_3) : g\in\Gamma\bigr\},
\end{align*}
where $k^{-1}=d_1(0)d_2(0)d_3(0)$. Note again that we can choose $D_\alpha$ and $d_\alpha$ such that $k=1$ by Theorem \ref{th:a2_classification} to obtain a symmetric description.

\begin{Remark}
    Finally, it is not clear in which way the building $\cX$ and the lattice $\Gamma$ depend on the given data, namely on the Singer groups $S_\alpha$, on the difference sets $D_\alpha$ and on the ordered difference sets $d_\alpha$. In Section \ref{subsec:cyclic_lattices}, we give examples where different orderings $d_\alpha$ lead to non-isomorphic lattices. We do not know whether the lattices arising from different orderings are all commensurable. In addition, it is not even clear whether this construction can lead to different buildings for fixed projective planes.
\end{Remark}

\noindent Several strong properties of these lattices can be obtained without knowing the exact building they act on.

\begin{Remark}
    By unpublished work of Shalom and Steger, lattices in arbitrary buildings of type $\tilde A_2$ follow a version of Margulis' normal subgroup theorem. In particular, the lattices $\Gamma$ are almost perfect, their first homology group is finite. For the special case of lattices generated by cyclic Singer groups, one can calculate group homology quite explicitly, as we show in Section \ref{sec:group_homology}.
\end{Remark}

\noindent In addition, by using the so-called spectral criterion, one can see that these lattices have property (T).

\begin{Proposition}[Cartwright-M{\l}otkowski-Steger, \.Zuk]
	Cocompact lattices in buildings of type $\tilde A_2$ have property (T).
\end{Proposition}

\begin{Proof}
	See \cite{CMS:TA2:94} or \cite{Zuk:TGP:96}. An detailed description of the latter result can be found in \cite[Theorem 5.7.7]{BHV:KPT:08}.
\end{Proof}

\subsection{The classification of panel-regular lattices}\label{subsec:a2_classification}

In this section, we will classify all panel-regular lattices on locally finite buildings of type $\tilde A_2$. In this process, we will show that the examples from the previous section cover all possible lattices, and that the presentations can even be simplified to a more symmetric form. We start with a simple observation.

\begin{Lemma}
    Let $X$ be a locally finite building of type $\tilde A_2$. Let $\Gamma\leq \Aut(X)$ be a lattice acting regularly on one type of panel. Then $\Gamma$ already acts regularly on all types of panels.
\end{Lemma}

\begin{Proof}
    Let $v,w,u$ be the vertices of a chamber in the building. Assume that $\Gamma$ acts regularly on the panels of the same type as the panel $(v,w)$. In particular, the vertex stabiliser $\Gamma_v$ is a group acting regularly on points (or lines) of the finite projective plane $\lk_X(v)$. By Proposition \ref{prop:regular_on_lines}, $\Gamma_v$ then also acts regularly on lines (or points) of $\lk_X(v)$, hence $\Gamma$ also acts regularly on the panels of the same type as $(v,u)$. By repeating the argument, we see that $\Gamma$ acts regularly on all types of panels.
\end{Proof}

\noindent In particular, the lattice $\Gamma$ acts transitively on vertices of the same type, which implies that vertex links of vertices of the same type are isomorphic. There are hence three finite projective planes of the same order, denoted by $\cI_1$, $\cI_2$ and $\cI_3$ associated to $X$. Pick a chamber in $X$ with vertices $\bar v_1, \bar v_2$ and $\bar v_3$. By the same argumentation as before, the three vertex stabilisers $\Gamma_{\bar v_1}$, $\Gamma_{\bar v_2}$ and $\Gamma_{\bar v_3}$ are then Singer groups on the vertex links and we denote these by $S_1$, $S_2$ and $S_3$. In the following, we write $J=\{0,1,\ldots,q\}$.

\begin{Theorem}\label{th:a2_classification}
    Let $X$ be a locally finite building of type $\tilde A_2$. We assume that there is a lattice $\Gamma\leq \Aut(X)$ acting regularly on one and hence on all types of panels. Then there are three Singer groups $S_1$, $S_2$ and $S_3$ associated to the three types of possible vertex links. Furthermore, there are three ordered difference sets $d_1$, $d_2$ and $d_3$ corresponding to these Singer groups satisfying $d_\alpha(0)=1\in S_\alpha$. The lattice $\Gamma$ admits the following presentation:
    \[
        \Gamma \cong \langle S_1,S_2,S_3 \,|\, \text{ all relations in the groups }S_1,S_2,S_3,\quad d_1(j)d_2(j)d_3(j) = 1 \quad\forall j \in J\rangle.
    \]
\end{Theorem}

\begin{Proof}
    Consider the canonical scwol $\cX$ associated to the building. Choose a vertex $\bar f_0\in V(\cX)$ corresponding to a chamber in $X$. Denote the vertices corresponding to vertices of the chamber by $\bar v_1$, $\bar v_2$ and $\bar v_3$ and the vertices corresponding to panels of the chamber by $\bar e_1$, $\bar e_2$ and $\bar e_3$. Choose $q$ chambers representing the other $\Gamma$-orbits of chambers and denote the corresponding vertices in $\cX$ by $\bar f_1,\ldots,\bar f_q$.

    Now we construct the associated quotient complex of groups as in Construction \ref{con:quotient_scwol}. For this, note that, since the $\Gamma$-action is regular on panels, the vertices we just named form a system of representatives for all orbits of vertices of $\cX$. It is not hard to see that the quotient scwol then has the following structure
    \begin{align*}
       V(\Gamma\backslash\!\backslash \cX) &\coloneq \{ v_1, v_2, v_3\} \sqcup \{e_1,e_2,e_3\} \sqcup \{ f_j : j \in J \} \\
       E(\Gamma\backslash\!\backslash \cX) &\coloneq \{ v_\alpha \leftarrow e_\beta: \alpha\neq\beta\} \sqcup \{ v_\alpha \leftarrow f_j : j\in J\} \sqcup \{e_\beta \leftarrow f_j : j\in J\},
    \end{align*}
where we denote the orbits of vertices by removing the bar. The quotient scwol is also shown in Figure \ref{fig:quotient_scwol} for $q=2$.

    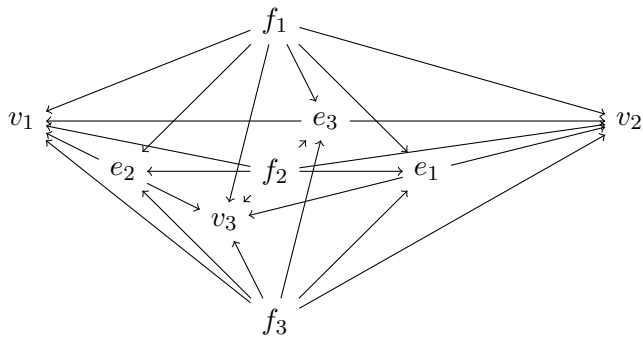
\begin{figure}[hbt]
    \centering
\begin{tikzpicture}[scale=4,font=\small]
	\node (p1) at (0,0,0) {$v_1$};
	\node (p2) at (2,0,0) {$v_2$};
	\node (p3) at (1,0,.866) {$v_3$};
	\node (l3) at (1,0,0) {$e_3$} edge[->] (p1) edge[->] (p2);
	\node (l1) at (1.5,0,.433) {$e_1$} edge[->] (p2) edge[->] (p3);
	\node (l2) at (.5,0,.433) {$e_2$} edge[->] (p1) edge[->] (p3);
	\node (f1) at (1,.5,.433) {$f_1$} edge[->] (p1) edge[->] (p2) edge[->] (p3) edge[->] (l1) edge[->] (l2) edge[->] (l3);
	\node (f2) at (1,0,.433) {$f_2$} edge[->] (p1) edge[->] (p2) edge[->] (p3) edge[->] (l1) edge[->] (l2) edge[->] (l3);
	\node (f3) at (1,-.5,.433) {$f_3$} edge[->] (p1) edge[->] (p2) edge[->] (p3) edge[->] (l1) edge[->] (l2) edge[->] (l3);
\end{tikzpicture}\caption{The quotient scwol $\Gamma\backslash\!\backslash\cX$ for $q=2$}\label{fig:quotient_scwol}
\end{figure}

The vertex groups are just the stabilisers. All of these are hence trivial except the stabilisers $\Gamma_{\bar v_\alpha}$, which are the Singer groups $S_\alpha$. All monomorphisms are obviously trivial.

It remains to determine the twist elements as in Construction \ref{con:quotient_scwol}, which is done by choosing the elements $h_a$ for $a\in E(\Gamma\backslash\!\backslash\cX)$ as follows: For every edge $a\in E(\Gamma\backslash\!\backslash \cX)$, there is precisely one preimage $\bar a\in E(\cX)$ satisfying $i(\bar a) = \overline{i(a)}$, but in general $t(\bar a)\neq \overline{t(a)}$. We choose $h_a$ such that $h_a(t(\bar a)) = \overline{t(a)}$.

We first consider the edges of the form $v_\alpha\leftarrow e_\beta$. For these, we obviously have $h_{v_\alpha\leftarrow e_\beta}=1$. Now consider the edges of the form $e_\beta\leftarrow f_j$. There is exactly one preimage $\bar e'_\beta\leftarrow \bar f_j$ where $\bar e_\beta'$ is in the orbit $e_\beta$. Since the $\Gamma$-action is regular on panels, there is exactly one element $h_{e_\beta\leftarrow f_j}\in\Gamma$ satisfying
\[
h_{e_\beta\leftarrow f_j}(t(\bar e_\beta'\leftarrow \bar f_j) ) = \bar e_\beta.
\]

For preimage edges of the type $\bar v_\alpha' \leftarrow \bar f_j$, where $\bar v_\alpha'$ is in the orbit $v_\alpha$, set $h_{v_\alpha\leftarrow f_j}=h_{e_\beta\leftarrow f_j}$ where $\beta \equiv_3 \alpha + 1$. Since by construction $g_{a,b}= h_ah_b h_{ab}^{-1}$, we have
\[
g_{v_\alpha\leftarrow e_\beta\leftarrow f_j} = \underbrace{h_{v_\alpha\leftarrow e_\beta}}_{=1} h_{e_\beta\leftarrow f_j}h_{v_\alpha\leftarrow f_j}^{-1}.
\]
This implies
\[
g_{v_\alpha\leftarrow e_\beta \leftarrow f_j}= 1\quad\text{for all } \beta-\alpha\equiv_3 1,\, j\in J.
\]
The other twist elements then necessarily form three ordered difference sets $d_\alpha$ as follows
\[
g_{v_\alpha\leftarrow e_\beta \leftarrow f_j}= d_\alpha(j) \quad\text{for } \beta-\alpha\equiv_3 2,\, j\in J,
\]
since the element $h_{e_\beta\leftarrow f_j}h^{-1}_{e_{\beta-1}\leftarrow f_j}$ transports the unique vertex in the orbit $f_j$ adjacent to $\bar e_{\beta-1}$ to the unique vertex in the orbit $f_j$ adjacent to $\bar e_\beta$. Since the vertices $\{\bar v_1,\bar v_2,\bar v_3,\bar e_1,\bar e_2,\bar e_3,\bar f_0\}$ span a subcomplex, we have
\[
d_1(0)=d_2(0)=d_3(0)=1.
\]
By the same calculation as in the proof of Proposition \ref{prop:fundamental_group_a2} we obtain the required presentation.
\end{Proof}

\begin{Remark}
    In particular, this means that any lattice constructed as in Section \ref{subsec:general_a2_construction} admits a presentation as in Theorem \ref{th:a2_classification}. This will simplify our considerations in the following sections.
\end{Remark}

\begin{Corollary}\label{cor:a2_building_description}
    In this situation, the building $X$ is isomorphic to the flag complex over the graph with vertices
    \[
        V(X) = \Gamma/S_1 \sqcup \Gamma/S_2 \sqcup \Gamma/S_3
    \]
    and edges
    \[
        E(X) = \{ (gS_1,gS_2), (gS_2,gS_3), (gS_3,gS_1) : g\in \Gamma\}.
    \]
\end{Corollary}

\subsection{Lattices generated by cyclic Singer groups}\label{subsec:cyclic_lattices}

In this section, we start with cyclic Singer groups to obtain very simple lattices. Let $\cI$ be the classical projective plane of order $q$. Take three cyclic Singer groups $S_\alpha$ of $\cI$ of order $q^2+q+1$ with respective generators $\sigma_\alpha$. We choose the points and lines for the construction of the difference sets $D_\alpha$ to be incident, such that $1\in D_\alpha$. We write
\[
	\Delta_\alpha \coloneq \{ \delta\in\Z/(q^2+q+1) : \sigma_\alpha^\delta \in D_\alpha \}.
\]
The sets of numbers $\Delta_\alpha$ are then \emph{difference sets} in the classical sense. In addition, we choose bijections $\delta_\alpha:J = \{0,1,\ldots,q\}\rightarrow \Delta_\alpha$ satisfying $\delta_\alpha(0)=0$. If we apply the construction from Section \ref{subsec:general_a2_construction}, we obtain the following very simple presentation.

\begin{BreakTheorem}[Lattices generated by cyclic Singer groups]\label{th:cyclic_lattices}
	For any prime power $q$ and any three classical difference sets $\Delta_1$, $\Delta_2$ and $\Delta_3$ containing 0, and for any bijections $\delta_\alpha:\{0,1,\ldots,q\}\rightarrow \Delta_\alpha$ satisfying $\delta_\alpha(0)=0$, the group $\Gamma$ with presentation
	\[
		\Gamma = \bigl\langle \sigma_1, \sigma_2, \sigma_3 \,\big|\, \sigma_1^{q^2+q+1}=\sigma_2^{q^2+q+1}=\sigma_3^{q^2+q+1}=1,\, \sigma_1^{\delta_1(j)}\sigma_2^{\delta_2(j)}\sigma_3^{\delta_3(j)} = 1 \quad\forall j \in J\bigr\rangle
\]
	is a uniform lattice in a building of type $\tilde A_2$.
\end{BreakTheorem}

\begin{Examples}
	We give explicit presentations for three lattices in the two smallest cases. More difference sets $\Delta$ can be obtained from \cite{LaJolla}.
	\begin{enumerate}
		\item For $q=2$, we choose all difference sets to be $\Delta=\{0,1,3\}$. We obtain the lattice
			\[
				\Gamma_2 \coloneq \langle \sigma_1,\sigma_2,\sigma_3 \,|\, \sigma_1^7=\sigma_2^7=\sigma_3^7 = \sigma_1\sigma_2\sigma_3 = \sigma_1^3\sigma_2^3\sigma_3^3 = 1\rangle.
			\]
			It can easily be seen that $H_1(\Gamma_2)\cong (\Z/7)^2$.
		\item By changing the order of the first difference set, we obtain a new lattice which is not isomorphic to the first one.
			\[
				\Gamma_2' \coloneq \langle \sigma_1,\sigma_2,\sigma_3 \,|\, \sigma_1^7=\sigma_2^7=\sigma_3^7 = \sigma_1^3\sigma_2\sigma_3 = \sigma_1\sigma_2^3\sigma_3^3 = 1\rangle.
			\]
			For the abelianisation, we obtain $H_1(\Gamma_2') \cong \Z/7$.
		\item For $q=3$, we choose the difference set $\Delta=\{0,1,3,9\}$. We obtain
			\[
			\Gamma_3 \coloneq \langle \sigma_1,\sigma_2,\sigma_3 \,|\, \sigma_1^{13}=\sigma_2^{13}=\sigma_3^{13} = \sigma_1\sigma_2\sigma_3 = \sigma_1^3\sigma_2^3\sigma_3^3 =\sigma_1^9\sigma_2^9\sigma_3^9 = 1\rangle.
			\]
			Here, we have $H_1(\Gamma_3)\cong (\Z/13)^2$.
	\end{enumerate}
    For $q=5$, it is possible to construct perfect lattices with this method.
    % With the matrix
    %
    % [ 0  4 10 23 24 26]
    % [ 4 10 23 24 26  0]
    % [10 23 24 26  0  4]
    % [23 24 26  0  4 10]
    % [24 26  0  4 10 23]
    % [26  0  4 10 23 24]
\end{Examples}

\begin{Remark}
	In \cite{KMW:A2l:84}, Köhler, Meixner and Wester showed that the following group is a  chamber-regular lattice in the building associated to $\Sl_3(\F_2(\!(t)\!))$:
    \[
    \Gamma = \langle a,b,c\,|\, a^3=b^3=c^3=1,\, (ab)^2=ba, (ac)^2=ca, (c^{-1}b)^2=bc^{-1} \rangle.
    \]
    It is not hard to verify that the subgroup generated by $a^{-1}b^{-1}$, $bc^{-1}$ and $ca$ is isomorphic to the lattice $\Gamma_2$ from the example above.

    Computer searches done by the author using \cite{SAGE} and \cite{GAP4} have not yielded any embeddings of our lattices in $\Sl_3(\F_q(\!(t)\!))$ for $q>2$.
\end{Remark}

\subsection{Spheres of radius two}\label{subsec:spheres}

For any vertex of an affine building of type $\tilde A_2$, the combinatorial sphere of radius one, usually called the link, forms a projective plane. The spheres of larger radii also admit interesting incidence structures, the so-called \emph{Hjelmslev planes}, see \cite{HvM:PHP:89}. It is a well-known fact that the two families of classical buildings of type $\tilde A_2$ associated to the groups $\Sl_3(\Q_p)$ and $\Sl_3(\F_p(\!(t)\!))$ can already be distinguished by inspecting spheres of radius two, respectively Hjelmslev planes of level two.

We will use a criterion by Cartwright, Mantero, Steger and Zappa, see \cite[Section 8]{CMSZ2}, to distinguish between these two affine buildings. Consider the following construction.

Denote the type set of a building $X$ of type $\tilde A_2$ by \{1,2,3\}. Fix a chamber $c_0$ and denote its vertices of type $1$, $2$ and $3$ by $v_1$, $v_2$ and $v_3$, respectively. Denote by $\cP$ and $\caL$ the sets of vertices adjacent to $v_1$ of types $3$ and $2$, respectively. Then of course $\cP$ and $\caL$ with the adjacency relation in the building form a projective plane.

Now consider the vertices of type $2$ and $3$ at distance $2$ from the vertex $v_1$ and denote these sets by $\cP^2$ and $\caL^2$, respectively.

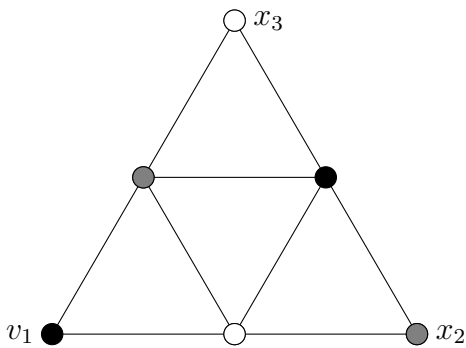
\begin{figure}[hbt]
    \centering
    \begin{tikzpicture}[scale=0.6]
        \tikzstyle{every node}=[circle, draw, inner sep=0pt, minimum width=8pt]
        \draw (0,0) node[fill=black] {} -- (4,0) node[fill=white] {} -- (8,0) node[fill=gray] {};
        \draw (0,0) -- (2,3.464) node[fill=gray] {} -- (4,6.928) node[fill=white] {};
        \draw (2,3.464) node[fill=gray] {} -- (4,0) node[fill=white] {} -- (6,3.464) -- cycle;
        \draw (4,6.928) node[fill=white] {} -- (6,3.464) node[fill=black] {} -- (8,0) node[fill=gray] {};
        \tikzstyle{every node}=[]
        \node (v1) at (-.7,0) {$v_1$};
        \node (v2) at (8.75,0) {$x_2$};
        \node (v3) at (4.75,6.928) {$x_3$};
    \end{tikzpicture}\caption{Adjacency in the Hjelmslev plane of level 2}\label{fig:adjacency}
    \end{figure}

\begin{Definition}
    Two vertices $x_2 \in \cP^2$, $x_3\in \caL^2$ are \emph{adjacent}, $x_2\sim x_3$, if the configuration in Figure \ref{fig:adjacency},
    where the colours of the vertices indicate the type, exists in the building. The incidence structure $(\cP^2,\caL^2,\sim)$ is called a \emph{Hjelmslev plane} of level two.
\end{Definition}

\noindent We will give an explicit construction of the Hjelmslev plane in the case of a building $X$ of type $\tilde A_2$ admitting a panel-regular lattice $\Gamma$.

Denote the vertex stabilisers of the vertices $v_\alpha$ by $S_\alpha$, these are Singer groups in the corresponding vertex links. For every vertex $v_\alpha$, the link $\lk_X(v_\alpha)$ is a projective plane of order $q$, where call the vertices of type $\alpha+1 \pmod{3}$ lines and the vertices of type $\alpha+2 \pmod{3}$ points.

For each of these projective planes, we construct a difference set $D_\alpha$ with respect to the points and lines given by the other vertices of $c_0$, respectively. We take the description of the building from Corollary \ref{cor:a2_building_description} and prove a simple lemma.

\begin{Lemma}\label{la:triangles} We have
    \begin{itemize}
        \item If $h,f\in S_1$ and $(v_1,hv_2,fv_3)$ is a triangle in $X$, then $f^{-1}h\in D_1$.
        \item If $g,f\in S_2$ and $(gv_1,v_2,fv_3)$ is a triangle in $X$, then $g^{-1}f\in D_2$.
        \item If $g,h\in S_3$ and $(gv_1,hv_2,v_3)$ is a triangle in $X$, then $h^{-1}g\in D_3$.
    \end{itemize}
\end{Lemma}

\begin{Proof}
    We will only prove the first claim, the other two are analogous. If $(v_1, hv_2,fv_3)$ is a triangle, then so is $(v_1,f^{-1}hv_2,v_3)$ which we obtain by multiplying with $f^{-1}$. Since $S_1$-translates of $v_2$ correspond to lines in the projective plane $\lk_X(v_1)$ by the definition of the difference set $D_1$, the line $f^{-1}h v_2$ is incident to the point $v_3$ if and only if $f^{-1}h\in D_1$.
\end{Proof}

\begin{Lemma}\label{lem:points_and_lines}
    We have
    \[
    \cP = \{ s_1 v_3 : s_1 \in S_1 \},\qquad \caL = \{ t_1 v_2 : t_1 \in S_1 \},
    \]
    as well as
    \[
    \cP^2 = \{ s_1s_3 v_2: s_1\in S_1, s_3\in S_3\backslash D_3^{-1} \},\qquad \caL^2 = \{ t_1t_2 v_3 : t_1\in S_1, t_2\in S_2\backslash D_2 \}.
    \]
\end{Lemma}

\begin{Proof}
    The first claim is obvious, since $S_1$ acts regularly on points and lines of the projective plane $\lk(v_1)$. The second claim is a simple consequence of Lemma \ref{la:triangles} --- the vertex $t_2v_3$ is not adjacent to $v_1$ if and only if $t_2\not\in D_2$, and similarly for $s_3$.
\end{Proof}

\noindent Adjacency in the Hjelmslev plane is relatively complex to describe but becomes manageable for cyclic Singer groups. The following observation will be used frequently in this section:

From the presentation of $\Gamma$ in Theorem \ref{th:a2_classification}, we see that for $d_1\in D_1$, there are always elements $d_2\in D_2$, $d_3\in D_3$ such that $d_1d_2d_3=1$.

\begin{Lemma}\label{lem:hjelmslev_adjacency}
    Fix a point $s_1s_3v_2\in\cP^2$ and a line $t_1t_2v_3\in\caL^2$, where $s_1,t_1\in D_1$, $t_2\in D_2$ and $s_3\in D_3$. In the Hjelmslev plane of level 2 around $v_1$, we have $s_1s_3v_2 \sim t_1t_2v_3$ if and only if the following conditions hold.
    \begin{description}
        \item[(C1)] We have $s_1^{-1}t_1 \in D_1$.
    \end{description}
     There are hence elements $d_2\in D_2$, $d_3\in D_3$ such that $(s_1^{-1}t_1)d_2d_3=1$.
     \begin{description}
        \item[(C2)] There is an element $n_2\in t_2D_2^{-1} \cap d_2D_2^{-1}$ such that
            \[
            s_3^{-1}d_3^{-1}e_3 \in D_3,
            \]
            where $e_1(n_2^{-1}d_2)e_3 =1$ for some elements $e_1\in D_1$ and $e_3\in D_3$.
    \end{description}
    Note that the elements $e_1$ and $e_3$ exist since $n_2\in d_2D_2^{-1}$.

    If the lattice arises from cyclic Singer groups $S_\alpha=\langle \sigma_\alpha\rangle$, $\alpha\in\{1,2,3\}$, as in Section \ref{subsec:cyclic_lattices} and if the associated unordered difference sets $\Delta_\alpha=\Delta$ are all equal, we obtain:
    \[
    \sigma_1^{j_1}\sigma_3^{j_3}v_2\sim \sigma_1^{k_1}\sigma_2^{k_2}v_3 \quad\Leftrightarrow\quad\begin{cases}
        k_1 - j_1 \in\Delta \text{ and}\\
        \exists n\in (k_2 - \Delta) \cap (-j_3-\Delta) \cap ( k_1-j_1-\Delta).
    \end{cases}
    \]
\end{Lemma}

\begin{Proof}
	Assume that $s_1s_3v_2\sim t_1t_2v_3$. Then there must be an element $n_2\in S_2$ such that we have the configuration of Figure \ref{fig:proof_adjacency} in $X$. We will now investigate the triangles in the order given by Roman numerals and apply Lemma \ref{la:triangles} repeatedly to obtain the required relations.
\begin{figure}[hbt]
    \centering
    \begin{tikzpicture}[scale=0.6]
        \tikzstyle{every node}=[circle, draw, inner sep=0pt, minimum width=8pt]
        \draw (0,0) node[fill=black] {} -- (4,0) node[fill=white] {} -- (8,0) node[fill=gray] {};
        \draw (0,0) -- (2,3.464) node[fill=gray] {} -- (4,6.928) node[fill=white] {};
        \draw (2,3.464) node[fill=gray] {} -- (4,0) node[fill=white] {} -- (6,3.464) -- cycle;
        \draw (4,6.928) node[fill=white] {} -- (6,3.464) node[fill=black] {} -- (8,0) node[fill=gray] {};
        \tikzstyle{every node}=[]
        \node (v1) at (-.7,0) {$v_1$};
        \node (v2) at (9.3,0) {$s_1s_3v_2$};
        \node (v3) at (5.3,6.928) {$t_1t_2v_3$};
        \node (v4) at (4,-.7) {$s_1v_3$};
        \node (v5) at (1.0,3.464) {$t_1v_2$};
        \node (v6) at (7.4,3.464) {$t_1n_2v_1$};
        \node (i) at (2,1.154) {I};
        \node (ii) at (4,4.618) {II};
        \node (iii) at (4,2.309) {III};
        \node (iv) at (6,1.154) {IV};
    \end{tikzpicture}\caption{The configuration in the proof of Lemma \ref{lem:hjelmslev_adjacency}}\label{fig:proof_adjacency}
    \end{figure}
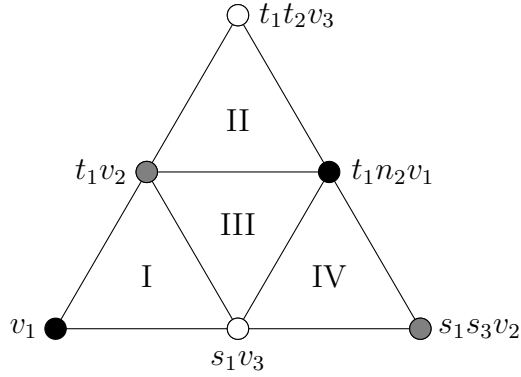

    \noindent Since $(v_1,t_1v_2,s_1v_3)$ form triangle I, we obtain
    \begin{equation}\label{eq:I}
        s_1^{-1}t_1\in D_1.
    \end{equation}

    \noindent There are hence elements $d_2\in D_2$ and $d_3\in D_3$ such that $(s_1^{-1}t_1)d_2d_3=1$. By looking at triangle II translated by $t_1^{-1}$, we obtain
    \begin{equation}\label{eq:II}
        n_2^{-1}t_2\in D_2.
    \end{equation}

    \noindent From triangle III, by observing that $s_1v_3=s_1d_3^{-1}v_3$ and by translating by $t_1^{-1}$, we obtain the triangle \[
    (n_2v_1,v_2,t_1^{-1}s_1d_3^{-1}v_3)=(n_2v_1,v_2,d_2v_3),
    \]
    which again by Lemma \ref{la:triangles} yields
    \begin{equation}
        n_2^{-1}d_2 \in D_2,\label{eq:III}.
    \end{equation}

    \noindent There are hence elements $e_1\in D_1$ and $e_3\in D_3$ such that $e_1(n_2^{-1}d_2)e_3=1$. For the last triangle IV, we first write $t_1n_2v_1=t_1n_2e_1^{-1}v_1$. By translating by $s_1^{-1}$, we obtain the triangle
    \begin{align}
        (s_1^{-1}t_1n_2e_1^{-1}v_1,s_3v_2,v_3).\nonumber \\
        \intertext{Substituting $s_1^{-1}t_1$ by $d_3^{-1}d_2^{-1}$ yields the triangle}
        (d_3^{-1}d_2^{-1}n_2e_1^{-1}v_1,s_3v_2,v_3).\nonumber \\
        \intertext{Finally, we substitute $d_2^{-1}n_2=e_3e_1$ and obtain the triangle}
        (d_3^{-1}e_3v_1,s_3v_2,v_3).\nonumber \\
        \intertext{We can now apply Lemma \ref{la:triangles} again to obtain}
        s_3^{-1}d_3^{-1}e_3 \in D_3. \label{eq:IV}
    \end{align}

    \noindent By assembling equations \eqref{eq:I},\eqref{eq:II}, \eqref{eq:III} and \eqref{eq:IV}, we obtain conditions (C1) and (C2). For the reverse direction, observe that (C1) and (C2) imply directly that the configuration of Figure \ref{fig:proof_adjacency} exists in the building.

    The second claim is then a simple calculation.
\end{Proof}

\noindent There is an obvious projection $\psi:(\cP^2,\caL^2,\sim) \rightarrow (\cP,\caL,\cF)$ given by
\begin{align*}
    \psi: \cP^2 &\rightarrow \cP: s_1s_3v_2 \mapsto s_1v_3\\
         \caL^2 &\rightarrow \caL: t_1t_2v_3 \mapsto t_1v_2.
\end{align*}

\begin{Lemma}[Lemma 8.1 in \cite{CMSZ2}]\label{l:cmsz}
    For $p,p'\in \cP^2$, consider the images $\psi(p)$ and $\psi(p')$. If $\psi(p)\neq \psi(p')$, there is a unique $l\in \caL^2$ such that $p\sim l\sim p'$. If $\psi(p)=\psi(p')$, there are $q$ distinct such lines. The same is true if the roles of points and lines are reversed.
\end{Lemma}

\begin{Proposition}\label{prop:hjelmslev_splitting}
    If the lattice $\Gamma$ and the building $X$ arise from cyclic Singer groups as in Section \ref{subsec:cyclic_lattices} and if in addition all difference sets $\Delta_\alpha=\Delta$ are equal, there is a splitting map $\iota: (\cP,\caL,\cF) \rightarrow (\cP^2,\caL^2,\sim)$ satisfying $\psi\circ\iota=\id_{(\cP,\caL,\cF)}$.
\end{Proposition}

\begin{Proof}
    With the same notation as in Lemma \ref{lem:hjelmslev_adjacency}, we choose an element
    \[
    m\in(\Z/(q^2+q+1))\setminus(\Delta\cup(-\Delta)).
    \]
    This is always possible for cardinality reasons, since $q^2+q+1 > 2q+2$ for $q\geq 2$. The map $\iota$ is given by
    \begin{align*}
        \cP & \rightarrow \cP^2&&\text{and}             & \caL & \rightarrow \caL^2 \\
        s_1v_3 & \mapsto s_1\sigma_3^{-m}v_2&&&t_1v_2 & \mapsto t_1\sigma_2^{m}v_3.
    \end{align*}
    This map preserves incidence by Lemma \ref{lem:hjelmslev_adjacency}, since for any $j_1, k_1$, the set
    \[
        (m - \Delta) \cap ( k_1 - j_1 - \Delta)
    \]
    has exactly one element, which corresponds to the intersection of lines in the projective plane $\lk_X(v_1)$.
\end{Proof}

\begin{Definition}[Section 8 in \cite{CMSZ2}]
    Consider $\cP'\subset \cP^2$ and $\caL'\subset \caL^2$. We say that $(\cP',\caL')$ is a \emph{substructure} of $(\cP^2,\caL^2)$ if for all $p,p'\in\cP'$ satisfying $\psi(p)\neq\psi(p')$, the unique line $l\in\caL^2$ incident with both $p$ and $p'$ is already contained in $\caL'$ and if the same condition holds with points and lines exchanged. The substructure \emph{generated by a set $\cP''\subset \cP^2$} is the smallest substructure containing $\cP''$.
\end{Definition}

\noindent We will use the following characterisation to distinguish spheres of radius $2$ in buildings of type $\tilde A_2$ by Cartwright-Mantero-Steger-Zappa in \cite{CMSZ2}.

\begin{Proposition}[Proposition 8.6 in \cite{CMSZ2}]\label{prop:cmsz}
   Assume that $X$ is the building associated to the special linear group $\Sl_3(K)$, where $K$ is a local, non-Ar\-chi\-me\-de\-an field with residue field of prime order $p$. Choose a vertex $v_1\in X$ and construct the plane $(\cP,\caL,\cF)$ and the Hjelmslev plane of level 2 denoted by $(\cP^2,\caL^2,\sim)$ as above. Consider four points $p_1,\ldots,p_4\in\cP^2$ such that no three points of $\psi(p_1),\ldots,\psi(p_4)$ are collinear.
   \begin{itemize}
       \item If $K=\Q_p$, the substructure of the Hjelmslev plane generated by the set of points $\{p_1,\ldots,p_4\}$ is all of $(\cP^2,\caL^2,\sim)$.
       \item Otherwise, the substructure generated by the set of points $\{p_1,\ldots,p_4\}$ is a projective plane of order $p$.
   \end{itemize}
\end{Proposition}

\begin{Corollary}
    The panel-regular lattices constructed out of cyclic Singer groups with identical difference sets cannot be contained in the building associated to $\Sl_3(\Q_p)$.
\end{Corollary}

\begin{Proof}
    Take any four points $p_1,\ldots,p_4$ in $(\cP,\caL,\cF)$ such that no three of these points are collinear. Then $\iota(p_1),\ldots, \iota(p_4)$ satisfy the conditions of Proposition \ref{prop:cmsz}. But by Lemma \ref{l:cmsz} and by Proposition \ref{prop:hjelmslev_splitting}, the substructure generated by these four points is exactly $\im(\iota)$, which is a projective plane of order $p$. By Proposition \ref{prop:cmsz}, if the building $X$ is associated to $\Sl_3(K)$, then necessarily $K\neq \Q_p$.
\end{Proof}

\begin{Remark}
    Conjecturally, these very simple lattices should be contained in $\Sl_3(\F_p(\!(t)\!))$. It might be possible to prove this by investigating the structure of all Hjelmslev planes of all levels and prove a splitting lemma for every step. This seems tedious and difficult. Surprisingly, we have not found any realisation of such a lattice in $\Sl_3(\F_p(\!(t)\!))$ for $p\neq 2$.
\end{Remark}

\section{Lattices in buildings of type \texorpdfstring{$\tilde C_2$}{\textasciitilde C2}}\label{sec:c2}

In this section, we will construct lattices in buildings of type $\tilde C_2$. This will be done using the slanted symplectic quadrangles discussed in Section \ref{subsec:GQs}. We will give two different constructions of panel-regular lattices in buildings of type $\tilde C_2$.

\paragraph{Situation} Fix a prime power $q>2$, and set $J=\{0,1,\ldots, q+1\}$. Let $\cI$ and $\cI'$ be two copies of the slanted symplectic quadrangle of order $(q-1,q+1)$. For the quadrangle $\cI$, we fix a point $p$ and denote the set of lines through $p$ by $L$. The set of flags  $F$ is given by $F=\{ (p,l):l\in L\}$. The associated Singer group will be denoted by $S$, the line stabilisers by $S_l$ for lines $l\in L$. We fix the same objects for $\cI'$ and add a prime $'$ to the notation.

\subsection{Lattices acting regularly on two types of panels}\label{subsec:c2_two_panels}

Since $q$ determines the quadrangle and the Singer group uniquely up to isomorphism, there is a bijection between the sets of line representatives $L$ and $L'$ and the corresponding line stabilisers $S_l$ and $S'_{l'}$ are pairwise isomorphic.

\paragraph{Construction} We choose two bijections $\lambda:J\rightarrow L$, $\lambda':J\rightarrow L'$ and an abstract group $S_j\cong \Z/q$ with isomorphisms $\psi_j: S_j\rightarrow S_{\lambda(j)}$ and $\psi'_j: S_j\rightarrow S'_{\lambda'(j)}$ for every $j\in J$. Consider the complex of groups $G(\cY)$ over the scwol $\cY$ with vertices
\begin{align*}
	V(\cY) \coloneq&\, \{v,v',w\} \sqcup \{e,e'\} \sqcup \{e_j : j\in J\} \sqcup \{f_j : j\in J \} \\
	\intertext{and edges}
	E(\cY) \coloneq&\, \{w \leftarrow e, w\leftarrow e', v\leftarrow e, v'\leftarrow e'\} \sqcup \{ v\leftarrow e_j, v'\leftarrow e_j : j\in J\} \\
	& \quad\sqcup \{ v\leftarrow f_j, v'\leftarrow f_j, w\leftarrow f_j : j\in J \} \sqcup \{e\leftarrow f_j, e'\leftarrow f_j, e_j\leftarrow f_j : j \in J\}.
\end{align*}
Figure \ref{fig:cYfor3} illustrates the scwol $\cY$ in the case $q=3$.
\begin{figure}[hbt]
\centering
\begin{tikzpicture}[scale=4,font=\small]
	\node (p1) at (0,0,0) {$w$};
	\node (p2) at (2,0,0) {$v'$};
	\node (p3) at (1,0,.866) {$v$};
	\node (l3) at (1,0,0) {$e'$} edge[->] (p1) edge[->] (p2);
	\node (l2) at (.5,0,.433) {$e$} edge[->] (p1) edge[->] (p3);
	\node[color=gray] (la) at (1.5,1,.433) {$e_0$};
	\node[color=gray] (lb) at (1.5,.5,.433) {$e_1$};
	\node (lc) at (1.5,0,.433) {$e_2$} edge[->] (p2) edge[->] (p3);
	\node[color=gray] (ld) at (1.5,-.5,.433) {$e_3$};
	\node[color=gray] (le) at (1.5,-1,.433) {$e_4$};
	\node[color=gray] (fa) at (1,.66,.433) {$f_0$} edge[->,color=lightgray] (p1) edge[->,color=lightgray] (p2) edge[->,color=lightgray] (p3) edge[->,color=lightgray] (la) edge[->,color=lightgray] (l2) edge[->,color=lightgray] (l3);
	\node[color=gray] (fb) at (1,.33,.433) {$f_1$} edge[->,color=gray] (p1) edge[->,color=gray] (p2) edge[->,color=gray] (p3) edge[->,color=gray] (lb) edge[->,color=gray] (l2) edge[->,color=gray] (l3);
	\node (fc) at (1,0,.433) {$f_2$} edge[->] (p1) edge[->] (p2) edge[->] (p3) edge[->] (lc) edge[->] (l2) edge[->] (l3);
	\node[color=gray] (fd) at (1,-.33,.433) {$f_3$} edge[->,color=gray] (p1) edge[->,color=gray] (p2) edge[->,color=gray] (p3) edge[->,color=gray] (ld) edge[->,color=gray] (l2) edge[->,color=gray] (l3);
	\node[color=gray] (fe) at (1,-0.66,.433) {$f_4$} edge[->,color=lightgray] (p1) edge[->,color=lightgray] (p2) edge[->,color=lightgray] (p3) edge[->,color=lightgray] (le) edge[->,color=lightgray] (l2) edge[->,color=lightgray] (l3);
	\draw[->,color=lightgray] (la) .. controls (1.75,.7,.216) .. (p2);
	\draw[->,color=lightgray] (la) .. controls (1.25,.7,.649) .. (p3);
	\draw[->,color=lightgray] (le) .. controls (1.75,-.7,.216) .. (p2);
	\draw[->,color=lightgray] (le) .. controls (1.25,-.7,.649) .. (p3);
	\draw[->,color=gray] (lb) .. controls (1.75,.35,.216) .. (p2);
	\draw[->,color=gray] (lb) .. controls (1.25,.35,.649) .. (p3);
	\draw[->,color=gray] (ld) .. controls (1.75,-.35,.216) .. (p2);
	\draw[->,color=gray] (ld) .. controls (1.25,-.35,.649) .. (p3);
\end{tikzpicture} \caption{The scwol $\cY$ for $q=3$.}\label{fig:cYfor3}
\end{figure}

Now choose the vertex groups to be $G_{v}= S$, $G_{v'}= S'$ and $G_{w}= \langle c\,|\, c^{q+2}=1\rangle$. In addition, set $G_{e_j}= S_j$ for all $j\in J$. All other vertex groups are trivial. The only non-trivial monomorphisms are chosen to be the maps $\psi_j$ and $\psi'_j$ for $j\in J$. All twist elements are trivial except
\[
g_{w\leftarrow e\leftarrow f_j} \coloneq c^j.
\]

\noindent We endow $|\cY|$ with a locally Euclidean metric as follows:

Let $\Delta$ be the geometric realisation of one triangle in the affine Coxeter complex of type $\tilde C_2$. For each $j\in J$, we map the subcomplex spanned by $\{v,v',w,e,e',e_j,f_j\}$ onto the barycentric subdivision of $\Delta$ in the obvious way and pull back the metric. We obtain a locally Euclidean metric on $|\cY|$. In particular, the angles at the vertices $v$ and $v'$ are $\pi/4$, the angle at $w$ is $\pi/2$.

\begin{Proposition}\label{prop:c2_developable}
	The complex of groups $G(\cY)$ is developable.
\end{Proposition}

\begin{Proof}
    The proof is analogous to the one of Proposition \ref{prop:a2_is_developable}. By construction of the complex, the geometric links in the local developments at $v$, $v'$ and $w$ are $\cZ(\cI)$, $\cZ(\cI')$ and the scwol associated to a complete bipartite graph of order $(q+2,q+2)$. In particular, they are CAT(1) by Proposition \ref{prop:key}.

	The local developments of the edges $e$ and $e'$ are the same as the closed stars in $\cY$: $(q+2)$ 2-simplices glued along one edge. For the edges $e_j$, the local development is isometric to $q$ triangles glued along one edge. The local developments of the vertices $f_j$ are just flat triangles.

    By Proposition \ref{prop:geometric_link_nonpos_curved}, the complex of groups $G(\cY)$ is then developable.
\end{Proof}

\begin{Proposition}\label{prop:c2_presentation}
	The fundamental group $\Gamma=\pi_1(G(\cY))$ admits the following presentation:
	\[
        \Gamma = \bigl\langle S,S',c \,\big|\, \text{ all relations in } S,S',\, c^{q+2}=1,\, c^j\psi_j(s)c^{-j}= \psi'_j(s)\quad \forall j \in J, s\in S_j \bigr\rangle.
	\]
\end{Proposition}

\begin{Proof}
	Consider the following maximal spanning subtree $T$:
	\[
		E(T)=\{w\leftarrow e, w\leftarrow e', v\leftarrow e, v'\leftarrow e'\} \sqcup \{w\leftarrow f_j, e_j\leftarrow f_j: j\in J\}.
	\]
	As in the proof of Proposition \ref{prop:fundamental_group_a2}, the group is generated by $S$, $S'$, $c$ and group elements for each edge not contained in $T$. Consider, however Figure \ref{fig:c2_triangle} showing one triangle in $\cY$.
	\begin{figure}
	\centering
	\begin{tikzpicture}[scale=6,font=\small]
		\node (v3) at (0,0) {$w$};
		\node (v2) at (0,1) {$v$};
		\node (v1) at (1,0) {$v'$};
		\node (e1) at (0,.5) {$e$} edge[->] (v3) edge [->] (v2);
		\node (e2) at (.5,0) {$e'$} edge[->] (v3) edge [->] (v1);
        \node (el) at (.5,.5) {$e_j$} edge[->,dotted] node[color=gray] {$1$} (v1) edge[->,dotted] node[color=gray]{$c^j$} (v2);
		\node (fl) at (.33,.33) {$f_j$} edge[->] (v3) edge[->,dotted] node[color=gray]{$1$} (e2) edge[->,dotted] node[color=gray]{$1$} (v1) edge[->,dotted] node[color=gray]{$c^j$} (e1) edge[->,dotted] node[color=gray]{$c^j$} (v2) edge[->] (el);
		\node (m12) at (.11,.61) {$1$};
		\node (m21) at (.61,.11) {$1$};
		\node (m13) at (.11,.29) {$c^j$};
		\node (m32) at (.29,.11) {$1$};
		\node (ml2) at (.27,.61) {$1$};
		\node (ml1) at (.61,.27) {$1$};
	\end{tikzpicture}\caption{One triangle in $\cY$}\label{fig:c2_triangle}
\end{figure}
	As before, edges in $T$ are drawn black, the other edges are drawn dotted. Group elements $g$ written on an edge $a$ indicate that $k_a=g$. From the presentation in Definition \ref{def:fundamental_group}, one can see that all additional relations are of the form stated in the result.
\end{Proof}

\begin{Theorem}\label{th:c2_building}
        The universal cover $\cX$ is a building of type $\tilde C_2$. The fundamental group $\Gamma=\pi_1(G(\cY))$ with presentation
	\[
        \Gamma = \bigl\langle S,S',c \,\big|\, \text{ all relations in } S,S',\, c^{q+2}=1,\, c^j\psi_j(s)c^{-j}= \psi'_j(s)\quad \forall j \in J, s\in S_j \bigr\rangle
	\]
        is hence a uniform lattice in the automorphism group of $\cX$, which is panel-regular on two types of panels.
\end{Theorem}

\begin{Proof}
    The proof is analogous to the one of Theorem \ref{th:a2_is_building}, using Proposition \ref{prop:key} and Theorem \ref{th:recognition}.
\end{Proof}

\begin{Remark}
    Except for the slanted symplectic quadrangle $W(3)^\Diamond$, all Singer quadrangles are necessarily exceptional. Hence, except for possibly this case, all buildings constructed in this way are exceptional buildings of type $\tilde C_2$.

    It is not hard to see that $\langle c\vert c^{q+2}=1\rangle$ can be replaced by any group $C$ of order $q+2$. In the above presentation, the term $c^j$ has to be replaced by an arbitrary bijection $J\rightarrow C$.
\end{Remark}

\noindent Finally, we also have a very explicit description of this building as the flag complex whose 1-skeleton is the following graph
\begin{align*}
        V(X) &\coloneq \Gamma/S \sqcup \Gamma/S' \sqcup \Gamma/\langle c\rangle \\
        E(X) &\coloneq \{ (gS,gS'), (gS,g\langle c\rangle), (gS',g\langle c \rangle) : g\in \Gamma\},
\end{align*}
	which could be used to study these exotic affine buildings.

	Instead, we will focus on the two cases where we obtain relatively simple presentations of these lattices.

\subsubsection{The prime case}

If $q=p$ is prime, we have very simple presentations of the Heisenberg groups $S$ and $S'$. Denote their generators by $x, y$ and $x', y'$, respectively and set $z=[x,y]$ and $z'=[x',y']$. The sets of line representatives $L$ and $L'$ can be parametrised by $\bP\F_p^2\sqcup\{0\}$ as in Theorem \ref{th:slanted_sympl_quadrangle}. We pick the following generators for the line stabilisers
\[
c_{[a:b]}= x^ay^bz^{-\tfrac{1}{2}ab},\quad c_0 = z,\quad c'_{[a:b]}= x'^ay'^bz'^{-\tfrac{1}{2}ab},\quad c'_0 = z'.
\]

\begin{Theorem}\label{th:c2_prime_result1}
    For any odd prime $p$, for $J=\{0,1,\ldots,p+1\}$ and any two bijections $\lambda, \lambda' : J \rightarrow \bP\F_p^2\sqcup\{0\}$, consider the group $\Gamma$ presented by
    \[
    \Gamma = \Biggl\langle x,y,x',y',c \,\Bigg|\, \begin{array}{c} z= xyx^{-1}y^{-1}, z'=x'y'x'^{-1}y'^{-1}, \\
            x^p=y^p=x'^p=y'^p=c^p=z^p=z'^p=1, \\ xz=zx, yz=zy, x'z'=z'x', y'z' = z'y', \\
            c^j c_{\lambda(j)} c^{-j} = c'_{\lambda'(j)} \qquad\forall j\in J \end{array}\Biggr\rangle.
    \]
    Then $\Gamma$ is a uniform lattice in a building of type $\tilde C_2$, where all vertex links which are quadrangles are isomorphic to the slanted symplectic quadrangle of order $(p-1,p+1)$.
\end{Theorem}

\begin{Proof}
    This is a direct application of Theorem \ref{th:c2_building}.
\end{Proof}

\subsubsection{The abelian case}

If $q$ is even, then the associated Singer groups $S$ and $S'$ are isomorphic to $\F_q^3$. We write $J=\{0,1,\ldots,q+1\}$ and fix two bijections
\[
\lambda, \lambda': J \rightarrow \bP\F_q^2\sqcup\{0\}.
\]
We write $S_{[a:b]}=\F_q(a,b,0)^T$ and $S_{0}= \F_q(0,0,1)^T$ for the line stabilisers in $S$ as in Theorem \ref{th:slanted_sympl_quadrangle} and add a prime $'$ for the respective groups in $S'$. Finally, we fix isomorphisms of abstract groups $\psi_j,\psi_j': \Z/q \rightarrow S_{\lambda(j)},S'_{\lambda'(j)}$ for any $j\in J$.

\begin{Theorem}
    Consider the group $\Gamma$ given by
    \[
    \Gamma= \bigl(S * S' * \langle c \rangle \bigr) / \langle c^{q+2}=1, c^j (\psi_j(x)) c^{-j} = \psi'_j(x) : x\in\Z/q, j\in J\rangle.
    \]
    Then $\Gamma$ is a uniform lattice in an exotic building of type $\tilde C_2$, where all vertex links which are quadrangles are isomorphic to the slanted symplectic quadrangle of order $(q-1,q+1)$.
\end{Theorem}

\begin{Proof}
	Again, this is a simple application of Theorem \ref{th:c2_building}.
\end{Proof}

\subsection{Lattices acting regularly on one type of panel}\label{subsec:c2_one_panel}

In this section, we will consider lattices acting regularly on only one type of panel in a building of type $\tilde C_2$. We will concentrate on the simpler case where the type of this panel corresponds to the two extremal vertices in the $\tilde C_2$-diagram.

As before, let $\cI$ and $\cI'$ be two slanted symplectic quadrangles of the same order $(q-1,q+1)$ and let $S$ and $S'$ be the associated Singer groups. Pick points $p$ and $p'$ and denote the sets of lines incident to $p$ and $p'$ by $L$ and $L'$, respectively. Then write $F=\{(p,l):l\in L\}$ and $F'=\{(p',l') : l'\in L'\}$.

\paragraph{Construction} We set $J=\{0,1,\ldots,q+1\}$ and choose two bijections $\lambda:J\rightarrow L$, $\lambda':J\rightarrow L'$. We consider the scwol $\cY$ given by
\begin{align*}
    V(\cY) &\coloneq \{ v,v'\} \sqcup \{v_j : j\in J\} \sqcup \{e\} \sqcup \{e_j, e'_j : j\in J\} \sqcup \{f_j : j\in J\} \\
    E(\cY) &\coloneq \{ v\leftarrow e, v'\leftarrow e\} \sqcup \{v \leftarrow e_j, v_j \leftarrow e_j: j\in J\} \sqcup \{v'\leftarrow e'_j, v_j\leftarrow e'_j : j\in J\} \\
    &\qquad\quad\sqcup \{v\leftarrow f_j, v'\leftarrow f_j, v_j\leftarrow f_j:j\in J\}\sqcup \{ e\leftarrow f_j, e_j\leftarrow f_j, e'_j\leftarrow f_j: j\in J\}.
\end{align*}
Figure \ref{fig:cYfor32} illustrates this scwol for $q=3$.
\begin{figure}
\centering
\begin{tikzpicture}[scale=4,x=30,y=30,z=10]
    \makeatletter
    \define@key{cylindricalkeys}{angle}{\def\myangle{#1}}
    \define@key{cylindricalkeys}{radius}{\def\myradius{#1}}
    \define@key{cylindricalkeys}{z}{\def\myz{#1}}
    \tikzdeclarecoordinatesystem{cylindrical}
    {
       \setkeys{cylindricalkeys}{#1}
       \pgfpointadd{\pgfpointxyz{0}{0}{\myz}}{\pgfpointpolarxy{\myangle}{\myradius}}
    }

    \tikzstyle{every edge}=[draw,->]
    \node (vv) at (0,0,0) {$v'$};
    \node (v) at (0,0,2) {$v$};
    \node (e) at (0,0,1) {$e$} edge (v) edge (vv);

    \foreach \a in {4,3,2,1,0}
    {
        \ifnum \a > 1
            \tikzstyle{every node}=[color=gray];
            \tikzstyle{every edge}=[draw,->,color=lightgray];
        \fi

        \node (v1) at (cylindrical cs:angle=\a*72+20,radius=1,z=1) {$v_\a$};
        \node (e1) at (cylindrical cs:angle=\a*72+20,radius=.5,z=1.5) {$e_\a$} edge (v) edge (v1);
        \node (ee1) at (cylindrical cs:angle=\a*72+20,radius=.5,z=.5) {$e'_\a$} edge (vv) edge (v1);

        \ifnum \a > 0 % Fix tikz quirk -- short arrows are drawn reversed!
        \node (f1) at (cylindrical cs:angle=\a*72+20,radius=.333,z=1) {$f_\a$} edge (v) edge (vv) edge (e) edge (e1) edge (ee1) edge (v1);
        \else
        \node (f1) at (cylindrical cs:angle=\a*72+20,radius=.333,z=1) {$f_\a$} edge (v) edge (vv) edge (e) edge (e1) edge (v1);
        \draw[->] (ee1) -- (f1);
        \fi
    }
\end{tikzpicture}\caption{The scwol $\cY$ for $q=3$.}\label{fig:cYfor32}
\end{figure}

We construct a complex of groups $G(\cY)$ over $\cY$ by setting the vertex groups to be
\[
G_v = S,\qquad G_{v'}=S',\qquad G_{v_j}= S_{\lambda(j)} \times S'_{\lambda'(j)},\qquad G_{e_j}=S_{\lambda(j)},\qquad G_{e'_j}=S'_{\lambda'(j)}.
\]
All other vertex groups are trivial. The inclusions are the obvious ones. All twist elements are trivial, so $G(\cY)$ is a \emph{simple} complex of groups.

We endow $|\cY|$ with a locally Euclidean metric as follows: Let $\Delta$ be the geometric realisation of one triangle in the affine Coxeter complex of type $\tilde C_2$. For each $j\in J$, we map each subcomplex spanned by $\{v,v',v_j,e,e_j,e'_j,f_j\}$ onto the barycentric subdivision of $\Delta$ in the obvious way and pull back the metric. We obtain a locally Euclidean metric on $|\cY|$ with angle $\pi/4$ at the vertices $v,v'$ and with angle $\pi/2$ at the vertex $v_j$.

\begin{Lemma}\label{la:c2_local_development}
	For any $j\in J$, the local development $\cY(\tilde v_j)$ is isomorphic to the cone over the barycentric subdivision of a complete bipartite graph of order $(q,q)$. The geometric link $\Lk(\tilde v_j,\st(\tilde v_j))$ is hence isometric to the barycentric subdivision of a generalised $2$-gon, in particular a connected, CAT(1) polyhedral complex of diameter $\pi$.
\end{Lemma}

\begin{Proof}
	We abbreviate $H=S_{\lambda(j)}$, $F=S'_{\lambda'(j)}$ and $G=H\times F$. We have $\cY(\tilde v) = \{\tilde v\} * \Lk_{\tilde v_j}(\cY)$ by the definition of the local development, where
	\begin{align*}
		V(\Lk_{\tilde v_j}(\cY)) &= \{ (gH, e_j) : g\in G \} \sqcup \{ (gF,e'_j) : g\in G\} \sqcup \{ (g,f_j) : g \in G \} \cong F \sqcup H \sqcup G\\
		E(\Lk_{\tilde v_j}(\cY)) &= \{ (g, e_j\leftarrow f_j) : g\in G\} \sqcup  \{ (g, e'_j\leftarrow f_j) : g\in G\},
	\end{align*}
    which is easily seen to be the barycentric subdivision of a complete bipartite graph on the vertex set $F\sqcup H$. The rest of the argument is analogous to the proof of Proposition \ref{prop:key}.
\end{Proof}

\begin{Theorem}\label{th:c2_result2}
    The complex $G(\cY)$ is developable. Its fundamental group $\Gamma = \pi_1(G(\cY))$ admits the presentation
    \[
    \Gamma = \langle S,S' \,|\, \text{ all relations in the groups } S,S',\; [S_{\lambda(j)},S'_{\lambda'(j)}]=1\quad \forall j\in J \rangle.
    \]
    The universal cover $\cX$ of this complex of groups is a building of type $\tilde C_2$, and $\Gamma$ is a uniform lattice acting regularly on one type of panels of $\cX$.
\end{Theorem}

\begin{Proof}
    This proof is analogous to the proofs of Proposition \ref{prop:c2_developable} and Theorem \ref{th:c2_building}, where we use Lemma \ref{la:c2_local_development} for the local developments at the vertices $v_j$.
\end{Proof}

\noindent In the next two small sections, we will make this explicit in the cases where we have simple presentations.

\subsubsection{The prime case}

If $q=p$ is prime and if we start with two slanted symplectic quadrangles $\cI$, $\cI'$ of order $(p-1,p+1)$ with corresponding Heisenberg groups $S$, $S'$ with respective generator pairs $x$, $y$ and $x'$, $y'$, the sets of line representatives $L$ and $L'$ can be parametrised by $\bP\F_p^2\sqcup\{0\}$. The stabilisers then have the form 
\[
S_{[a:b]}=\langle x^ay^bz^{-\tfrac{1}{2}ab}\rangle\qquad\text{and}\qquad S_0=\langle z\rangle,
\]
where $z=[x,y]$ and similarly for $S'$. Let $J=\{0,1,\ldots, p+1\}$.

\begin{Theorem}\label{th:c2_prime_result2}
	Let $\lambda, \lambda':J\rightarrow \bP\F_p^2\sqcup \{0\}$ be two bijections. The group
    \[
    \Gamma = \Biggl\langle x,y,x',y' \,\Bigg|\, \begin{array}{c} z= xyx^{-1}y^{-1}, z'=x'y'x'^{-1}y'^{-1}, \\
            x^p=y^p=x'^p=y'^p=z^p=z'^p=1, \\ xz=zx, yz=zy, x'z'=z'x', y'z' = z'y', \\\space
            [S_{\lambda(j)}, S'_{\lambda'(j)}] = 1\qquad \forall j\in J \end{array}\Biggr\rangle
    \]
    is a uniform lattice in a building of type $\tilde C_2$.
\end{Theorem}

\begin{Proof}
	This is a direct application of Theorem \ref{th:c2_result2}.
\end{Proof}

\subsubsection{The abelian case}

If $q$ is even, then the associated Singer groups $S$ and $S'$ are isomorphic to $\F_q^3$. In this case, we obtain a very simple presentation:

Let $J=\{0,1,\ldots,q+1\}$ and fix two bijections $\lambda,\lambda': J \rightarrow \bP\F_q^2\sqcup \{0\}$. We write $S_{[a:b]}=\F_q(a,b,0)^T$ and $S_{0}= \F_q(0,0,1)^T$ for the line stabilisers in $S$ and add a prime $'$ for the respective groups in $S'$.

\begin{Theorem}
	The group
	\[
	\Gamma = \bigl(S * S' \bigr) / \langle [ S_{\lambda(j)}, S'_{\lambda'(j)}] : j\in J\rangle
	\]
	is a uniform lattice in an exotic building of type $\tilde C_2$.
\end{Theorem}

\begin{Proof}
    Again, this is a direct application of Theorem \ref{th:c2_result2}.
\end{Proof}

\begin{Remark}
    Note that this is a finitely presented group generated by involutions such that all other relations make elements commute. It would be very interesting if one of these groups were a right-angled Coxeter group, but this seems very unlikely.
\end{Remark}

\subsection{Property (T)}

Even though most of these buildings are necessarily exotic, the lattices have automatically property (T).

\begin{Proposition}[\.Zuk]
	All cocompact lattices in buildings of type $\tilde C_2$ have property (T) if the order of the generalised quadrangles which appear as vertex links is not $(2,2)$. In particular, all lattices we construct here have property (T).
\end{Proposition}

\begin{Proof}
	See \cite{Zuk:TGP:96}. A more detailed exposition can be found in \cite[5.6]{BHV:KPT:08}.
\end{Proof}

\section{Group homology of lattices}\label{sec:group_homology}

We will calculate group homology of the lattices we have constructed using the action on the building. Since the building $X$ is contractible and the lattice acts without inversions, we can use the spectral sequence associated to equivariant homology we have obtained in Section \ref{sec:equivariant_homology}.

\subsection{Group homology of lattices in buildings of type \texorpdfstring{$\tilde A_2$}{\textasciitilde A2}}

We will calculate the full group homology of the lattices generated by cyclic Singer groups which we have constructed in Section \ref{subsec:cyclic_lattices}. The abelianisation can be calculated directly.

\begin{Proposition}
    Let $\Gamma_1$ be a lattice constructed as in Theorem \ref{th:cyclic_lattices} using a classical projective plane of order $q$ and three ordered classical difference sets $\delta_1$, $\delta_2$ and $\delta_3$. We have
    \[
    H_1(\Gamma_1) = \Bigl\langle s_1,s_2,s_3 \,\Big|\, \begin{array}{c} (q^2+q+1)s_1 = (q^2+q+1)s_2 = (q^2+q+1)s_3 = 0 \\ \delta_1(j)s_1+\delta_2(j)s_2+\delta_3(j)s_3=0\quad\forall j\in J \end{array}\Bigr\rangle.
    \]
    In particular, $\Gamma_1$ is perfect if and only if the $(q\times 3)$-matrix $\caD= (\delta_j(i))_{i,j}$ over the ring $\Z/(q^2+q+1)$ has full rank in the sense of \cite[Chapter 4]{Bro:MCR:93}.
\end{Proposition}

\begin{Proof}
    Since $H_1(\Gamma_1)$ is isomorphic to the abelianisation of $\Gamma_1$, we just have to abelianise the presentation from Theorem \ref{th:cyclic_lattices}. The resulting group is the kernel of the linear map $(\Z / (q^2+q+1))^3 \rightarrow (\Z / (q^2+q+1))^q$ described by $\caD$. In particular, by \cite[Theorem 5.3]{Bro:MCR:93}, the kernel is trivial if and only if $\caD$ has full rank.
\end{Proof}

\begin{Theorem}\label{th:a2_homology}
    Again, let $\Gamma_1$ be a lattice constructed as in Theorem \ref{th:cyclic_lattices}. Let $q$ be the order of the associated projective plane. Then
    \[
        H_j(\Gamma_1) \cong\begin{cases}
            \Z & j=0\\
            \ker(\caD) & j=1 \\
            \Z^q & j=2 \\
            (\Z/(q^2+q+1))^3 & j\geq 3 \text{ odd} \\
            0 & \text{ else.}
        \end{cases}
    \]
    In addition
    \[
    H_2(\Gamma_2;\Q)=\Q^q,\qquad H_j(\Gamma_2;\Q)=0 \text{ for }j\not\in\{ 0,2\}
    \]
    for any lattice $\Gamma_2$ constructed as in Theorems \ref{th:a2_is_building} or \ref{th:cyclic_lattices}.
\end{Theorem}

\begin{Proof}
    Consider the spectral sequence $E^1_{i,j}$ from Theorem \ref{th:equiv_hom} for the $\Gamma_1$-action on the associated building $X$. The quotient space $\Gamma_1\backslash X$ is the geometric realisation of the scwol $\cY$ from the construction in Section \ref{subsec:general_a2_construction}, which has hence the homotopy type of a bouquet of $2$-spheres. The structure of the spectral sequence is particularly simple since most stabiliser subgroups are trivial. We have
    \[
    E^1_{0,j} \cong \bigoplus_{k=1}^3 H_j( \Z/(q^2+q+1)), \qquad E^1_{1,0} \cong \Z^3, \qquad E^1_{2,0} \cong \Z^{q+1}.
    \]
    All other groups on the first page are trivial. The differential on the bottom row is induced by the cellular differential, so we have
    \[
    E^2_{0,0}\cong \Z,\qquad E^2_{1,0}=0, \qquad E^2_{2,0}\cong \Z^q.
    \]
    All other groups remain unchanged. Since the only remaining non-trivial differential is
    \[
    d^2_{2,0}: E^2_{2,0}\cong \Z^q \rightarrow (\Z/(q^2+q+1))^3 = E^2_{0,1},
    \]
    whose kernel is always isomorphic to $\Z^q$, we obtain the full description of $H_*(\Gamma_1)$ by inspecting the third page $E^3_{i,j}$ and using Proposition \ref{prop:hom_of_cycl_groups} for the homology of cyclic groups. The rational homology can easily be calculated using Corollary \ref{cor:rational_equiv_homology}.
\end{Proof}

\subsection{Group homology of lattices in buildings of type \texorpdfstring{$\tilde C_2$}{\textasciitilde C2}}

We will first calculate rational group homology for the lattices of type $\tilde C_2$ we have constructed in Section \ref{sec:c2}.

\begin{Theorem}\label{th:c2_rational_homology}
    For a lattice $\Gamma$ constructed as in Sections \ref{subsec:c2_two_panels} or \ref{subsec:c2_one_panel}, we have $H_j(\Gamma;\Q)=0$ for $j\neq 0$.
\end{Theorem}

\begin{Proof}
    The geometric realisations of the quotient scwols are contractible in both cases. We can hence apply Corollary \ref{cor:rational_equiv_homology} to obtain the result.
\end{Proof}

\noindent For the lattices acting regularly on two types of panels as in Section \ref{subsec:c2_two_panels}, the situation with integral coefficients is rather complicated and probably depends on the identification functions $\lambda$ and $\lambda'$ given in the construction.

Hence we will consider a lattice $\Gamma$ acting regularly on one type of panel as in Section \ref{subsec:c2_one_panel}. Assume that the associated slanted symplectic quadrangle has order $(q-1,q+1)$, denote the associated Singer group by $S$.

\begin{Theorem}\label{th:c2_homology}
    For the lattice $\Gamma$, we have
    \[
    H_1(\Gamma)\cong (\Z/q)^6\qquad\text{and}\qquad H_2(\Gamma) \cong H_2(S) \oplus H_2(S).
    \]
\end{Theorem}

\begin{Proof}
    The structure of the abelianisation can easily be seen using the presentation from Theorem \ref{th:c2_result2}. For $H_2$, we inspect the corresponding spectral sequence from Theorem \ref{th:equiv_hom}. As a set of representatives for the action, we use the subcomplex induced by the set of all chambers containing a common panel of the type the lattice acts regularly on. Then we have
    \begin{align*}
        E^1_{0,2} &\cong H_2(S)^2,& E^1_{1,2} &=0,\\
        E^1_{0,1} &\cong (\Z/q)^6 \oplus (\Z/q)^{2q+4},& E^1_{1,1} &\cong (\Z/q)^{2q+4},\\
        E^1_{0,0} &\cong \Z^{q+4},& E^1_{1,0}&\cong \Z^{2q+5},& E^1_{2,0}&\cong \Z^{q+2}.
    \end{align*}
    In this calculation, we use Proposition \ref{prop:hom_of_cycl_groups} as well as $H_1(S)=(\Z/q)^3$, which is easily verified. Since the quotient space $\Gamma\backslash X$ is contractible, we obtain
    \[
    E^2_{0,0} \cong \Z,\qquad E^2_{1,0}=0,\qquad E^2_{2,0}=0.
    \]
    Since we already know that $H_1(\Gamma)\cong (\Z/q)^6$, the map $d^1_{1,1}:E^1_{1,1}\rightarrow E^1_{0,1}$ must be injective and we obtain $E^2_{1,1}=0$ which proves the theorem.
\end{Proof}

\chapter{Open problems}\label{ch:problems}

We will end this thesis by stating some of the problems which could be considered next.

\minisec{Wagoner complexes}
For Wagoner complexes, we think that the following question is natural to ask: Given a Wagoner complex $\cW$ associated to a group with a spherical root datum $G$, is there an isomorphism $\pi_i(BG^+) \cong \pi_{i-1}(\cW)$? This would naturally extend Wagoner's results concerning $K$-theory and could lead to a new characterisation of hermitian $K$-theory.

In addition, one might hope to generalise Wagoner's result on the homotopy groups of affine Wagoner complexes to arbitrary groups with root data with valuations, but this will require a new technique.

\minisec{Homological stability}
We think that the construction of the spectral sequence in Theorem \ref{th:stability_pair_spectral_sequence} associated to any stability pair can be applied to many other situations. For instance, one could investigate homology groups of the (short) series of groups of types $D_5$, $E_6$, $E_7$, $E_8$.

Also, one could try to investigate relative homology of groups of types $C_n$ and $A_{n-1}$, which might lead to relations between hermitian $K$-theory and algebraic $K$-theory.

Additionally, one should be able to prove homological stability results for reductive algebraic groups and their Levi subgroups given only the types of the groups, albeit with a weak stability range, by using semidirect decompositions of the rational points of Levi subgroups.

Finally, there is an analogous version of the opposition complex for arbitrary twin buildings. Unfortunately, not much is known about its homotopy type. If this opposition complex turns out to be highly connected, one could try to generalise our method to groups acting on twin buildings.

\minisec{Lattice constructions}
For lattices, the largest open question is of course whether the $\tilde A_2$-buildings we have constructed are exotic. In addition, it might be very interesting to consider certain finite-index subgroups of the lattices and investigate their properties.

Finally, one could hope to extend the construction we have given to higher dimensions. This seems to be increasingly difficult, though, and would of course not yield any new buildings, but instead (new) presentations of arithmetic lattices.

\appendix

\chapter{Coxeter diagrams}

We include the lists of all spherical and affine Coxeter diagrams here for reference.

\begin{figure}[ht]
    \vspace{1em}
	\begin{center}
	\begin{tikzpicture}
		\node (an) at (-2,11) {$A_n$};
		\node (bn) at (-2,10) {$B_n=C_n$};
		\node (dn) at (-2,8.5) {$D_n$};
		\node (e6) at (-2,6.5) {$E_6$};
		\node (e7) at (-2,4.5) {$E_7$};
		\node (e8) at (-2,2.5) {$E_8$};
		\node (f4) at (-2,1.5) {$F_4$};
		\node (g2) at (-2,0.5) {$G_2$};
		\node (g2) at (-2,-0.5) {$H_3$};
		\node (g2) at (-2,-1.5) {$H_4$};
		\node (g2) at (-2,-2.5) {$I_2(m)$};

		% A, B, D

		\foreach \y in {11,10,8.5} {
			\tikzstyle{every node}=[draw,circle,fill,inner sep=0pt,minimum width=4pt];
            \draw (1,\y) node {} -- (2,\y) node {};
			\draw[dotted] (2,\y) -- (3,\y);
			\draw (3,\y) node {} -- (4,\y) node {};
			\draw (4,\y) -- (5,\y) node {};
		}
		\draw (1,11) -- (0,11);
		\draw (1,10.05) -- (0,10.05); \draw (1,9.95) -- (0,9.95);
		\tikzstyle{every node}=[draw,circle,fill,inner sep=0pt,minimum width=4pt];
		\node (x) at (0,11) {}; \node (x) at (0,10) {};
		\draw (1,8.5) -- (0,9) node {};
		\draw (1,8.5) -- (0,8) node {};

		% E6 - E8

		\foreach \y in {6.5,4.5,2.5} {
		\foreach \x in {0,1,2,3} {
			\draw (\x,\y) node {} -- ++(1,0);
		}
			\draw (2,\y) -- ++(0,1) node {};
			\node (x) at (4,\y) {};
		}
		\draw (4,4.5) -- ++(1,0) node {};
		\draw (4,2.5) -- ++(1,0) node {}; \draw (5,2.5) -- ++(1,0) node {};

		% Rest

		\foreach \y in {1.5,0.5,-0.5,-1.5,-2.5} {
		\draw (0,\y) node {} -- (1,\y) node {};
		}
		\draw (1,1.55) -- ++(1,0); \draw (1,1.45) -- ++(1,0);
		\draw (2,1.5) node {} -- (3,1.5) node {};

		\draw (1,-0.5) -- ++(1,0) node{};
		\draw (1,-1.5) -- ++(1,0) node{};
		\draw (2,-1.5) -- ++(1,0) node{};

		\tikzstyle{every node}=[];
		\node (g2) at (0.5,0.8) {$6$};
		\node (h3) at (1.5,-0.2) {$5$};
		\node (h4) at (2.5,-1.2) {$5$};
		\node (im) at (0.5,-2.2) {$m$};

		% for

		\node (a) at (8,11) {$n\geq 1$};
		\node (b) at (8,10) {$n\geq 2$};
		\node (d) at (8,8.5) {$n\geq 4$};

	\end{tikzpicture}
	\end{center}
	\caption{Spherical Coxeter diagrams}
	\label{fig:spherical_coxeter_diagrams}
\end{figure}
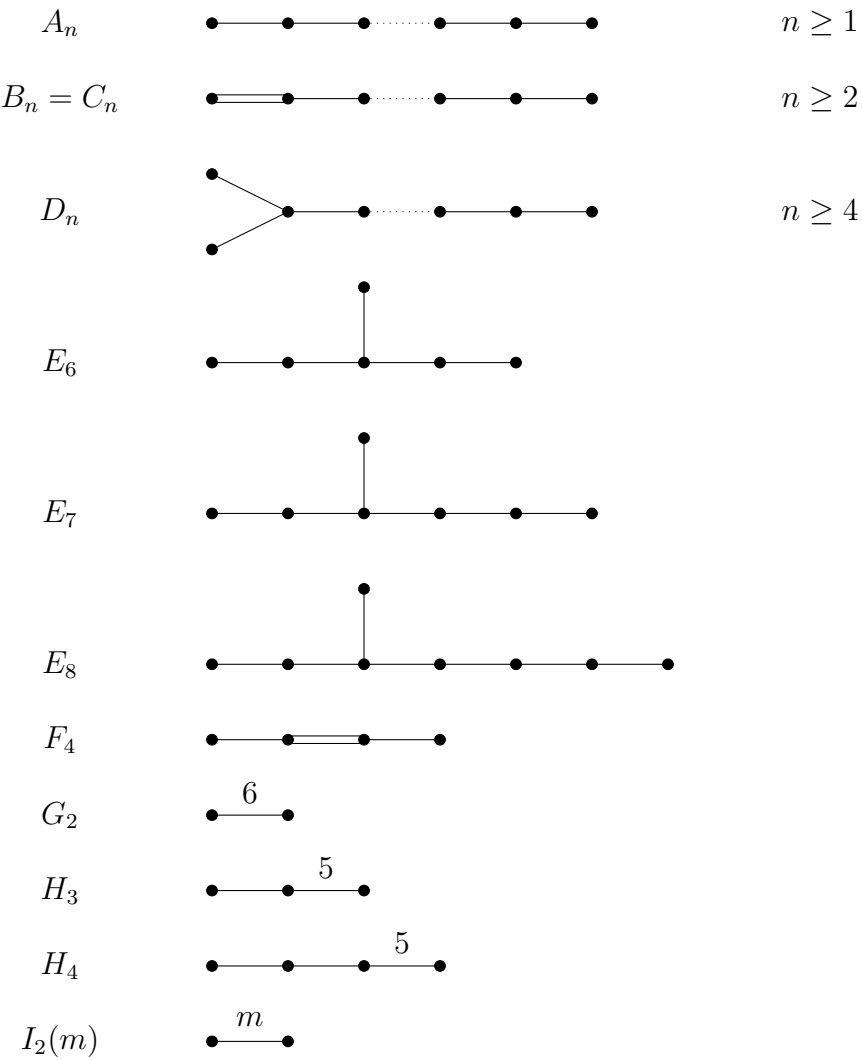

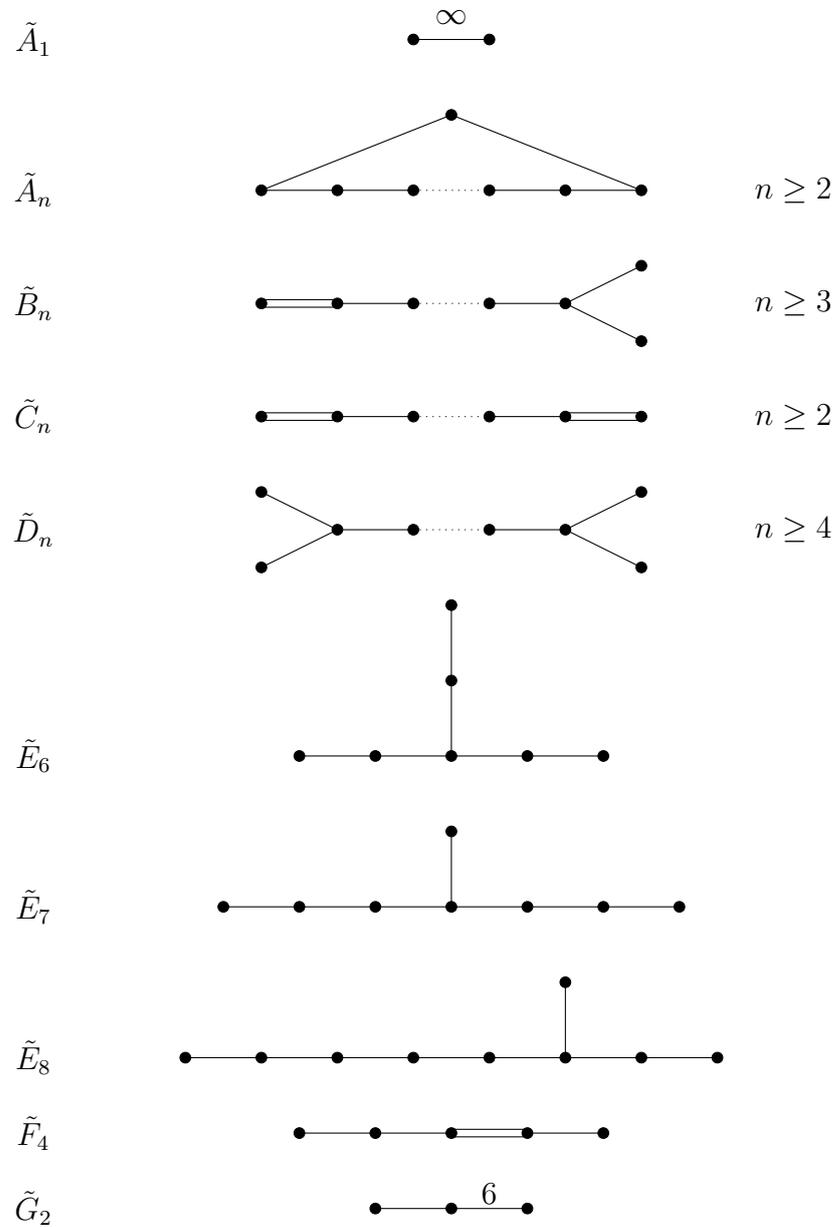
\begin{figure}[ht]
	\begin{center}
	\begin{tikzpicture}
		\node (a1) at (-2,15) {$\tilde A_1$};
		\node (an) at (-2,13) {$\tilde A_n$};
		\node (bn) at (-2,11.5) {$\tilde B_n$};
		\node (cn) at (-2,10) {$\tilde C_n$};
		\node (dn) at (-2,8.5) {$\tilde D_n$};
		\node (e6) at (-2,5.5) {$\tilde E_6$};
		\node (e7) at (-2,3.5) {$\tilde E_7$};
		\node (e8) at (-2,1.5) {$\tilde E_8$};
		\node (f4) at (-2,0.5) {$\tilde F_4$};
		\node (g2) at (-2,-0.5) {$\tilde G_2$};

        % A1

        \node (a1) at (3.5,15.3){$\infty$};
		\tikzstyle{every node}=[draw,circle,fill,inner sep=0pt,minimum width=4pt];
        \draw (3,15) node {} -- ++(1,0) node {};

		% A, B, D

		\foreach \y in {13,11.5,10,8.5} {
			\tikzstyle{every node}=[draw,circle,fill,inner sep=0pt,minimum width=4pt];
            \draw (2,\y) node {} -- ++(1,0) node {};
			\draw[dotted] (3,\y) -- ++(1,0);
			\draw (4,\y) node {} -- ++(1,0) node {};
		}
        % An
        \draw (1,13) node {} -- ++(1,0);
        \draw (5,13) -- ++(1,0) node {};
        \draw (6,13) -- (3.5,14) node{} -- (1,13);
        \tikzstyle{every node}=[]
		\draw (2,10.05) -- ++(-1,0) node{}; \draw (2,9.95) -- ++(-1,0) node{};
		\draw (5,10.05) -- ++(1,0) node{}; \draw (5,9.95) -- ++(1,0) node{};
		\draw (1,11.55) -- ++(1,0); \draw (1,11.45) -- ++(1,0);
		\tikzstyle{every node}=[draw,circle,fill,inner sep=0pt,minimum width=4pt];
		\node (x) at (1,11.5) {}; \node (x) at (6,10) {};
		\node (x) at (1,10) {};
        \draw (6,8) node {} -- ++(-1,0.5) -- ++(1,0.5) node {};
        \draw (1,8) node {} -- ++(1,0.5) -- ++(-1,0.5) node {};
        \draw (6,11) node {} -- ++(-1,0.5) -- ++(1,0.5) node {};

		% E6 - E8

		\foreach \y in {5.5,3.5} {
		\foreach \x in {1.5,2.5,3.5,4.5} {
			\draw (\x,\y) node {} -- ++(1,0);
		}
			\node (x) at (5.5,\y) {};
		}
        % E6
        \draw (3.5,5.5) -- ++(0,1) node{} -- ++(0,1) node{};
        % E7
        \draw (3.5,3.5) -- ++(0,1) node{};
        \draw (0.5,3.5) node {} -- ++(1,0); \draw (5.5,3.5) -- ++(1,0) node {};
		% E8
        \draw (5,1.5) -- ++(0,1) node{};
        \foreach \x in {0,1,2,3,4,5,6} {
        \draw (\x,1.5) node {} -- ++(1,0);
        }
        \node (e8) at (7,1.5) {};

		% Rest

		\draw (1.5,0.5) node {} -- ++(1,0) node {} -- ++(1,0) node {};
		\draw (3.5,0.55) -- ++(1,0); \draw (3.5,0.45) -- ++(1,0);
		\draw (4.5,0.5) node {} -- ++(1,0) node {};

        \draw (2.5,-0.5) node {} -- ++(1,0) node {} -- ++(1,0) node {};
		\tikzstyle{every node}=[];
		\node (g2) at (4,-0.3) {$6$};

		% for

		\node (a) at (8,13) {$n\geq 2$};
		\node (b) at (8,11.5) {$n\geq 3$};
		\node (c) at (8,10) {$n\geq 2$};
		\node (d) at (8,8.5) {$n\geq 4$};

	\end{tikzpicture}
	\end{center}
	\caption{Affine Coxeter diagrams}
	\label{fig:affine_coxeter_diagrams}
\end{figure}

\bibliographystyle{alpha}
\bibliography{../biblio}
\end{document}